\documentclass[aos,preprint]{imsart}

\usepackage[numbers]{natbib}
\usepackage{amsthm,amsmath,amssymb}
\usepackage{xcolor}
\usepackage{graphicx}
\usepackage[margin = 1.5 in]{geometry}
\startlocaldefs
\newtheorem{theorem}{Theorem}
\newtheorem{lemma}[theorem]{Lemma}
\newtheorem{corollary}[theorem]{Corollary}
\theoremstyle{remark}
\newtheorem{remark}[theorem]{Remark}
\endlocaldefs

\newcommand{\mR}{\mathbb{R}}

\DeclareMathOperator{\etr}{etr}
\DeclareMathOperator{\diag}{diag}
\DeclareMathOperator{\newre}{Re}
\DeclareMathOperator{\newim}{Im}

\setattribute{journal}{name}{}

\begin{document}

\begin{frontmatter}

% "Title of the paper"
\title{Testing in high-dimensional spiked models}
\runtitle{Testing in spiked models}

% indicate corresponding author with \corref{}
% \author{\fnms{John} \snm{Smith}\corref{}\ead[label=e1]{smith@foo.com}\thanksref{t1}}
% \thankstext{t1}{Thanks to somebody} 
% \address{line 1\\ line 2\\ printead{e1}}
% \affiliation{Some University}

\author{\fnms{Iain M.} \snm{Johnstone}\ead[label=e1]{imj@stanford.edu}\thanksref{t1}}
\address{\printead{e1}}
\thankstext{t1}{Supported in part by NSF DMS 0906812 and 1407813.} 
\and
\author{\fnms{Alexei} \snm{Onatski}\ead[label=e2]{ao319@cam.ac.uk}\thanksref{t2}}
\address{\printead{e2}}
\thankstext{t2}{Supported by the J.M. Keynes Fellowship Fund,
  University of Cambridge.} 
\affiliation{Stanford University and University of Cambridge}
%\affiliation{???}

\runauthor{I. M. Johnstone and A. Onatski}

\begin{abstract}
We consider the five classes of multivariate statistical problems
identified by James (1964), 
which together cover much of classical multivariate analysis,
plus a simpler limiting case,
symmetric matrix denoising. 
Each of James' problems  involves the 
eigenvalues of $E^{-1}H$ where $H$ and $E$ are proportional to
high dimensional Wishart matrices. Under the null hypothesis, both Wisharts
are central with identity covariance. Under the alternative, the
non-centrality or the covariance parameter of $H$ has a single eigenvalue, a
spike, that stands alone. When the spike is smaller than a case-specific phase
transition threshold, 
none of the sample eigenvalues separate from the bulk, making the
testing problem challenging. 
Using a unified strategy for the six cases,
we show that the log likelihood
ratio processes parameterized by the value of the sub-critical spike converge
to Gaussian processes with logarithmic correlation. We then
derive asymptotic power envelopes for tests for the presence of a
spike.

% We consider five different classes of multivariate statistical problems
% identified by James (1964). Each of these problems is related to the
% eigenvalues of $E^{-1}H$ where $H$ and $E$ are proportional to
% high-dimensional Wishart matrices. Under the null hypothesis, both Wisharts
% are central with identity covariance. Under the alternative, the
% non-centrality or the covariance parameter of $H$ has a single eigenvalue, a
% spike, that stands alone. When the spike is larger than a case-specific phase
% transition threshold, one of the eigenvalues of $E^{-1}H$ separates from the
% bulk. This makes the alternative easily detectable, so that reasonable
% statistical tests have asymptotic power one. In contrast, when the spike is
% sub-critical, that is lies below the threshold, none of the eigenvalues
% separates from the bulk, which makes the testing problem more interesting from
% the statistical perspective. In such cases, we show that the log likelihood
% ratio processes parameterized by the value of the sub-critical spike converge
% to Gaussian processes with logarithmic correlation. We use this result to
% derive the asymptotic power envelopes for tests for the presence of a spike in
% the data representing each of the five cases in James' classification.

\end{abstract}



\end{frontmatter}

%AOS,AOAS: If there are supplements please fill:
% \begin{supplement}[id=suppA]
%  \sname{Supplement A}
%  \stitle{Title}
%  \slink[doi]{10.1214/00-AOASXXXXSUPP}
%  \sdatatype{.pdf}" 
%  \sdescription{Some text}
% \end{supplement}

\addtocontents{toc}{\protect\setcounter{tocdepth}{0}}
\section{Introduction}

High-dimensional multivariate models and methods, such as regression,
principal components, and canonical correlation analysis, 
repay study 
% are the subject of much recent research 
% In contrast to the classical framework where the
% dimensionality is fixed, the current focus is on 
in frameworks where the
dimensionality diverges to infinity together with the sample size. 
``Spiked'' models that deviate from a reference model along a
small fixed number of unknown directions have proven to be a fruitful
abstraction and research
tool in this context.
 A basic statistical question that arises in the analysis of such
models is how to test for the presence of spikes in the data.

James (1964) arranges multivariate statistical problems in five different
groups with broadly similar features. His remarkable classification,
recalled in Table \ref{tab:HandE}, 
relies on the five most common hypergeometric functions
 $_{\mathsf{p}}F_{\mathsf{q}}$. In this paper, we describe 
rank-one spiked models that
represent each of James' classes in a high dimensional setting. 
We derive the asymptotic behavior of the
corresponding likelihood ratios
in a regime where the dimensionality $p$ of the data and the degrees of freedom
$n_1, n_2$ increase proportionally. 
Specifically, we study the ratios of the joint densities of
the relevant data under the alternative hypothesis, which assumes the presence
of a spike, to that under the null of no spike. 
%In each of the cases, 
The relevant data consist, in each case,
of the maximal invariant statistic represented by
eigenvalues of a large random matrix. 

We find that the joint distributions of the
eigenvalues under the alternative hypothesis and under the null are mutually
contiguous when the values of the spike is below a phase transition
threshold. The value of the threshold depends on the problem type.
Furthermore, we find that the log likelihood ratio processes parametrized by
the value of the spike are asymptotically Gaussian, with logarithmic mean
and autocovariance functions. These findings allow us to compute the
asymptotic power envelopes for the tests for the presence of spikes in five
multivariate models representing each of James' classes.

\begin{table}[t]
  \centering
%\caption{}
\caption{The five cases of James (1964) }
\begin{tabular}[h]{llll}
    & {Statistical method} & $n_1H$ & $n_2 E$ \\[8pt]
$_{0}F_{0}$ \ \ PCA & Principal components analysis & $W_{p}(n_1, \Sigma + \Phi) $ &
$n_2 \Sigma$ \\
       & \qquad [latent roots of covariance matrix] & & \\
$_{1}F_{0}$ \ \ SigD & Signal Detection  & $W_{p}(n_1, \Sigma + \Phi) $ &
$ W_{p}(n_2,\Sigma)$ \\ 
       & \qquad [equality of covariance matrices] & & \\
$_{0}F_{1}$ \ \ REG$_0$ & Multivariate regression, known error & $W_{p}(n_1,
\Sigma, n_{1}\Phi) $ &  $n_2 \Sigma $ \\ 
       & \qquad covariance \ [non-central means] & & \\
$_{1}F_{1}$ \ \ REG & Multivariate regression, unknown error  & $W_{p}(n_1,
\Sigma, n_{1}\Phi) $ &  $ W_{p}(n_2,\Sigma)$ \\ 
       & \qquad covariance \ [non-central latent roots] & & \\
$_{2}F_{1}$ \ \ CCA & Canonical correlation analysis & $W_{p}(n_1, \Sigma,
\Phi(Y)) $ & 
$ W_{p}(n_2,\Sigma)$ 
  \end{tabular}
\\[12pt]{\it \noindent James' names for the cases,
    when different from ours, are shown in brackets.
Final two columns interpret $H$ and $E$ of \eqref{basic} 
% below in each case, 
for Gaussian data, so that $W_{p}$ denotes a $p$-variate central
    or noncentral Wishart distribution, see Definitions. Matrix
    $\Phi$ has low rank, equal to one in this paper. 
For CCA, $\Phi(Y)$ is a random noncentrality matrix, see Supplementary Material (SM)
\ref{Identification of the parameters} for
definition. In
    cases 1 and 3, $E$ is deterministic, $\Sigma$ is
    known, and $n_2$ disappears. Otherwise $E$ is assumed independent of $H$.}
  \label{tab:HandE}
\end{table}

% \begin{table}[t]
%   \centering
% \caption{}
% %  \caption{The five cases of James (1964) }
% \begin{tabular}[h]{llll}
%     & {\it The five cases of James (1964)} & $n_1H$ & $n_2 E$ \\[8pt]
% 1. PCA & Principal components analysis & $W(n_1, \Sigma_0 + \Phi) $ &
% $n_2I$ \\
%        & \qquad [latent roots of covariance matrix] & & \\
% 2. SigD & Signal Detection  & $W(n_1, \Sigma + \Phi) $ &  $ W(n_2,\Sigma)$ \\
%        & \qquad [equality of covariance matrices] & & \\
% 3. REG$_0$ & Multivariate regression, known error & $W(n_1,
% \Sigma_0, \Phi) $ &  $n_2 \Sigma_0 $ \\ 
%        & \qquad covariance \ [non-central means] & & \\
% 4. REG & Multivariate regression, unknown error  & $W(n_1,
% \Sigma, \Phi) $ &  $ W(n_2,\Sigma)$ \\ 
%        & \qquad covariance \ [non-central latent roots] & & \\
% 5. CCA & Canonical correlation analysis & $W(n_1, \Sigma, \Phi(X)) $ &
% $ W(n_2,\Sigma)$ 
%   \end{tabular}
% \\[12pt]{\it \noindent James' terms,
%     when different from ours, are shown in brackets.
% Final two columns interpret $H$ and $E$ of \eqref{basic} below 
%  in each case, for Gaussian data, while $W$ denotes a $p$-variate central
%     or noncentral Wishart distribution, see Definitions. Matrix
%     $\Phi$ has low rank. In
%     cases 1 and 3, $E$ is deterministic and $\Sigma_0$
%     known. Otherwise $E$ is assumed independent of $H$.}
%   \label{tab:HandE}
% \end{table}

Our analysis is based on classical results that assume Gaussian data. All
the likelihood ratios that we study correspond to the joint densities of the
solutions to the basic equation of classical multivariate statistics,%
\begin{equation}
\det\left(  H- \lambda E\right)  =0, \label{basicJO}
\end{equation}
where the hypothesis $H$ and error sums of squares $E$ are proportional to Wishart matrices, as summarized
for the various cases in Table \ref{tab:HandE}.
% and the generalized eigenvalue $f$ is the multivariate extension of
% Fisher's $F$ ratio.
% The five different cases that we study are: 1) $E$ is a known deterministic
% matrix, and $H$ is a central Wishart matrix with covariance equal to a
% low-rank perturbation of $E$; 2) both $E$ and $H$ are central Wisharts with
% unknown covariance matrices that differ by a matrix of low rank; 3) $E$ is a
% known deterministic matrix, and $H$ is a non-central Wishart matrix with
% covariance equal to $E$ and with a low-rank non-centrality; 4) $E$ is a
% central Wishart matrix, while $H$ is a non-central one with the same unknown
% covariance matrix and with a low-rank non-centrality; 5) $E$ is a central
% Wishart, while $H$ is a non-central Wishart conditionally on a random low-rank
% non-centrality parameter. 
The five cases can be linked via sufficiency and
invariance arguments to the statistical problems listed in the table.
% a principal components problem, a signal detection
% problem, hypotheses testing in multivariate regression with known and with
% unknown error covariance, and a canonical correlation problem,
% respectively. 
We briefly discuss these links in the next section.

James' classification suggests common features that call for a
systematic approach.
Thus the main steps of our asymptotic analysis are the same for all the
five cases.
The likelihood ratios have explicit forms that involve hypergeometric
functions of two high-dimensional matrix arguments. 
However, one of the arguments has low rank under our spiked model
alternatives.  
Indeed, for tractability, we focus on the rank one setting.
We can then represent the hypergeometric function of two
high-dimensional matrix arguments in the form of a contour integral that
involves a \textit{scalar} hypergeometric function of the same
type, Lemma \ref{DJLemma}, using the recent result of 
%Dharmawansa and Johnstone (2014).
\cite{dhjo14}.
Then we deform the 
contour of integration so that the integral becomes amenable to Laplace
approximation analysis, 
%(see Olver (1997), chapter 4).
extending \cite[ch. 4]{olve97}.

Using the Laplace approximation technique, we show that the log likelihood
ratios are asymptotically equivalent to simpler random functions of the
spike parameters, Theorems \ref{Icase} and \ref{LR asymptotics}. The randomness 
% in the quadratic function 
enters via a linear
spectral statistic of a large random matrix of either sample covariance or
$F$-ratio type. 
Using central limit theorems for the two cases, due to
%Bai and Silverstein (2004) and Zheng (2012) respectively,
\cite{basi04} and \cite{zheng12} respectively,
% Using Central Limit Theorems (CLT) for the linear spectral statistics, established by
% Bai and Silverstein (2004) for the sample-covariance-type random matrices and
% by Zheng (2012) for the $F$-ratio-type random matrices, 
we derive the
asymptotic Gaussianity and obtain the mean and the autocovariance functions of
the log likelihood ratio processes, Theorem \ref{AsymptoticNormality}.

These asymptotics of the log likelihood processes show that the
corresponding statistical experiments do not converge to Gaussian shift
models. In other words, the experiments that consist of observing the
solutions to equation (\ref{basicJO}) parameterized by the value of the spike
under the alternative hypothesis are not of Locally Asymptotically Normal
(LAN) type. This implies that there are no ready-to-use\ optimality results
associated with LAN experiments that can be applied in our setting. However at
the fundamental level, the derived asymptotics of the log likelihood ratio
processes is all that is needed for the asymptotic analysis of the risk of the
corresponding statistical decisions.

In this paper, we use the derived asymptotics together with the Neyman-Pearson
lemma and Le Cam's third lemma to find simple
analytic expressions for the asymptotic power envelopes for the statistical
tests of the null hypothesis of no spike in the data,
Theorem \ref{th:PE}. The form of the
envelope depends only on whether both $H$ and $E$ in equation
(\ref{basicJO}) are Wisharts or only $H$ is Wishart 
whereas $E$ is deterministic.

For most of the cases, as the value of the spike under the alternative
increases, the envelope, at first, rises very slowly. Then, as the spike
approaches the phase transition, the rise quickly accelerates and the envelope
`hits' unity at the threshold. However, in cases of two Wisharts and when the
dimensionality is not much smaller than the degrees of freedom of $E,$ the
envelope rises 
%much faster. 
more rapidly.
In such cases, the information in all the
eigenvalues of $E^{-1}H$ might be useful for detecting population spikes which
lie far below the phase transition threshold.

A type of the analysis performed in this paper has been previously implemented
in the study of the principal components case by 
%Onatski et al (2013).
\cite{omh13}. 
Our work here identifies common features in James' classification of
multivariate statistical problems and uses them to extend the analysis
to the full system.
One of the hardest challenges in such an
extension is the rigorous implementation of the Laplace approximation step.
With this goal in mind, we have developed asymptotic approximations to the
hypergeometric functions $_{1}F_{1}$ and $_{2}F_{1}$ which are uniform in
certain domains of the complex plane, Lemma \ref{1F1approximation}.

The simple observation that the solutions to equation (\ref{basicJO}) can be
interpreted as the eigenvalues of random matrix $E^{-1}H$ relates our work to
the vast literature on the spectrum of large random matrices. 
% We refer the reader to Bai and Silverstein (2006) for a recent
% book-long treatment of the subject.
Three extensively studied classical ensembles of random matrices are
the Gaussian, Laguerre and Jacobi ensembles, e.g. \cite{meht04}.
%Mehta (2004). 
However, only
the Laguerre and Jacobi ensembles appear in high-dimensional analysis of 
James' five-fold classification. This prompts us to look for a
\textquotedblleft missing\textquotedblright\ class in James' system that could
be linked to the Gaussian ensemble.

Such a class is easy to obtain by taking the limit of $\sqrt{n_{1}}\left(
H-\Sigma\right)  $ with $\Sigma=I_{p}$ as $n_{1}\rightarrow\infty,$ for $p$ fixed.
% where $n_{1}$ and $p$ are $H$'s
% degrees of freedom and dimensionality, respectively. 
We call the corresponding
statistical problem \textquotedblleft symmetric matrix
denoising\textquotedblright (SMD). Under the null hypothesis, the
observations are 
given by a $p\times p$ matrix $Z/\sqrt{p}$ with $Z$ from the Gaussian
Orthogonal Ensemble. Under the alternative, the observations are given by
$Z/\sqrt{p}+\Phi,$ where $\Phi$ is a deterministic symmetric matrix of low
rank, again of rank one for this paper.
 We add this \textquotedblleft case zero\textquotedblright to James' 
classification and derive the asymptotics of the
corresponding log likelihood ratio and 
power envelope. 

% Many existing results in the random matrix literature do not require that the
% data are Gaussian. This suggests that some results about tests for the
% presence of the spikes in the data may remain valid without the
% Gaussian assumption.
% % One may for example consider $H$ and $E$ in (\ref{basic}) that, although have
% % the form of sample covariance matrices, do not come from the underlying
% % Gaussian distribution, and study the properties of the corresponding tests. We
% % leave this line of research to the future.
% % Since the explicit form of the joint distribution of the solutions to
% % (\ref{basic}) is only known in the Gaussian case, it seems unlikely that one
% % would be able to completely summarize the asymptotic behavior of the
% % corresponding non-Gaussian statistical experiments. 
% We hope that the results of this paper
% %, that provide such a summary under the Gaussianity, can 
% might provide 
% a benchmark for such future studies.
% % that would relax our assumptions.

To summarize, the contributions of this paper are as follows.
\begin{itemize}
\item We revisit James' classification, which covers a large part of
  classical multivariate analysis,
now in the setting of
high-dimensional data and show that the classification accommodates low
rank structures as departures from the classical null hypotheses.
\item We show that in such high dimensional settings with rank-one
  structure, random matrix theory 
  allows tractable approximations to the joint eigenvalue density functions, in
  place of slowly converging zonal polynomial series.
\item We show that the log likelihood ratio processes, when parametrized by
  spike magnitude, converge to Gaussian process limits in the
  sub-critical interval.
\item Hence, we show that informative tests are possible based on
  \textit{all} the eigenvalues whereas tests based on the largest
  eigenvalue alone are uninformative.
\item As a tool, we develop new uniform approximations to certain
  hypergeometric functions.
\item We identify symmetric matrix denoising as a limiting case of
  each of James' models. It is the simplest model displaying all the
  phenomena seen in the paper. We clarify the manner in which the
  simpler cases are limits of the more complex ones.

\end{itemize}

The rest of the paper fleshes out this program and its conclusions.
The proofs are largely deferred to the extensive Supplementary
Material (SM).
They reflect substantial effort to identify and exploit
common structure in the six cases.
Indeed some of this common structure appears remarkable and not yet
fully explained.

%All this detail is deferred to the extensive Supplementary Material.
%, which is arranged to parallel the main text.

% is organized as follows. 
% In the next section, we relate
% the five different cases of equation (\ref{basic}) to the classical
% multivariate statistical problems representing different cells of James'
% (1964) five-fold classification system. In Section 3, we obtain explicit
% expressions for the likelihood ratios. Section 4 represents the likelihood
% ratios in the form of contour integrals. Section 5 performs the Laplace
% approximation analysis. Section 6 derives the asymptotic power envelopes.
% Section 7 concludes. Technical proofs are given in the Appendix and in the
% Supplementary Material (SM).

\smallskip
\textit{Definitions and global assumptions.} \ Let $Z$ be an $n \times p$
data matrix with rows drawn i.i.d. from $N_p(0,\Sigma)$, a
$p$-dimensional normal distribution with mean $0$ and covariance
$\Sigma$. Suppose that $M$ is also $n \times p$, but deterministic. 
If $Y = M + Z$, then $H = Y'Y$ has a $p$ dimensional Wishart
distribution $W_p(n, \Sigma, \Psi)$ with $n$ degrees of freedom,
covariance matrix $\Sigma$ 
and non-centrality matrix $\Psi = \Sigma^{-1} M'M$. 
The central Wishart distribution, corresponding to $M = 0$, is denoted 
$W_p(n,\Sigma)$. 

Throughout the paper, we shall assume that
\[
p\leq\min\left\{  n_{1},n_{2}\right\}  ,
\]
where $p$ is the dimensionality of matrices $H$ and $E$, and $n_{1},n_{2}$ 
are the degrees of freedom of the corresponding Wishart distributions,
as summarized in Table \ref{tab:HandE}.
The assumption $p\leq n_{2}$ ensures almost sure invertibility of
matrix $E$ in (\ref{basicJO}), whereas the assumption $p\leq n_{1}$
while not essential, is made for brevity, as it reduces the number of
various situations which need to be considered.

\section{Links to statistical problems}

We briefly review examples of statistical problems, old and new, that
lead to each of James' five cases, plus symmetric matrix denoising,
and explain our choice of labels
for those cases.

\smallskip
\textit{PCA.} \ In the first case $n_{1}$
i.i.d. $N_{p}\left( 0,\Omega\right)  $  observations are used to test the null
hypothesis that the population covariance $\Omega$ equals a given matrix
$\Sigma$. The alternative of interest is
\[
\Omega=\Sigma+\Phi \quad \text{with} \quad \Phi=\theta \mathbb{\psi}\mathbb{\psi}^{\prime},
\]
where $\theta>0$ and $\mathbb{\psi}$ are unknown, and $\mathbb{\psi}$ is
normalized so that $ \| \Sigma^{-1/2}\mathbb{\psi} \| =1$.

Without loss of generality (wlog), we may assume that $\Sigma=I_{p}$. Then under the
null, the data are isotropic noise, whereas under the alternative, the first
principal component explains a larger portion of the variation than the other
principal components. 
% We therefore label Case 1 as the `principal components
% analysis' (\rm{PCA}) case.

The null and the alternative hypotheses can be formulated in terms of the
spectral `spike' parameter $\theta$ as
\begin{equation}
H_{0}:\theta_{0}=0\text{ and }H_{1}:\theta_{0}=\theta>0,
\label{Basic hypothesesJO}
\end{equation}
where $\theta_{0}$ is the true value of the `spike'. This testing problem
remains invariant under the multiplication of the $p\times n_{1}$ data matrix
from the left and from the right by orthogonal matrices, and under the
corresponding transformation in the parameter space. A maximal invariant
statistic consists of the solutions $\lambda_{1}\geq...\geq\lambda_{p}$ of
equation (\ref{basicJO}) with $n_1 H$ equal to the sample covariance matrix and
$E=\Sigma$. We restrict attention to the invariant tests. Therefore, the
relevant data are summarized by $\lambda_{1},...,\lambda_{p}$.
For convenience, details of the invariance and sufficiency arguments
for all cases are in SM \ref{Sufficiency and invariance considerations}.

\smallskip
\textit{SigD.} \ Consider  testing the equality of covariance
matrices, $\Omega$ and $\Sigma$, corresponding to two independent
$p$-dimensional zero-mean Gaussian samples of sizes $n_{1}$ and $n_{2}$.
The alternative hypothesis is the same as for case PCA.
Invariance considerations lead to tests based on the eigenvalues
of the $F$-ratio of the sample covariance matrices. 
% \footnote{For the reader's convenience, we 
% give details of the invariance and sufficiency arguments for all the cases in the
% SM.} 
Matrix $H$ from
(\ref{basicJO}) equals the sample covariance corresponding to the observations
that might contain a `signal' responsible for the covariance spike, whereas
matrix $E$ equals the other `noise' sample covariance matrix.
We again can assume that the population covariance of the `noise'
$\Sigma=I_{p}$, although this time it is unknown to the statistician
(SM \ref{Sufficiency and invariance considerations} explains
why such an assumption involves no loss of generality).
% We label Case 2 as the `signal detection' (SigD) case. 
Here, we find it more convenient to
work with the $p$ solutions to the equation%
\begin{equation}
\det\left(  H-\lambda\left(  E+\frac{n_{1}}{n_{2}}H\right)  \right)  =0,
\label{basic beta formJO}
\end{equation}
which we also denote $\lambda_{1}\geq...\geq\lambda_{p}$ to make the notations
as uniform across the different cases as possible. 
Note that as the second sample size $n_{2} \to \infty$, while $n_{1}$
and $p$ are held constant, equation (\ref{basic beta formJO}) reduces to
equation (\ref{basicJO}), $E$ converges to $\Sigma$, and SigD reduces to PCA.
% Note that as the number of
% observations in the second sample, $n_{2}$, diverges to infinity while $n_{1}$
% and $p$ are held constant, equation (\ref{basic beta formJO}) reduces to
% equation (\ref{basic}), $E$ converges to $\Sigma$, and SigD reduces to PCA.

\smallskip
\textit{REG$_0$, REG.} \ 
Now consider linear regression with multivariate response
\[
Y=X\beta+\varepsilon
\]
when the goal is to test linear restrictions on the matrix of coefficients
$\beta$. 
In case REG$_0$ the covariance matrix
$\Sigma$ of the i.i.d. Gaussian rows of the error matrix $\varepsilon$
is assumed known. 
%We label this case as `regression with known variance' (REG$_{0}$).
REG corresponds to unknown $\Sigma.$ 
%and we label it as `regression with unknown variance' (REG).

As explained in \cite[pp. 433--434]{muir82},
%Muirhead (1982), pp. 433--434, 
the problem of testing linear
restrictions on $\beta$ can be cast in the canonical form, where the matrix of
transformed response variables is split into three parts, $Y_{1},$
$Y_{2},$ and $Y_{3}$. Matrix $Y_{1}$ is $n_{1}\times p,$
where $p$ is the number of response variables and $n_{1}$ is the number of
linear restrictions (per each of the $p$ columns of matrix $\beta$). Under the null hypothesis, $\mathbb{E}Y_{1}=0,$ whereas
under the alternative,
\begin{equation}
\mathbb{E}Y_{1}=\sqrt{n_{1}\theta}\varphi\psi^{\prime},
\label{REG0mean}%
\end{equation}
where $\theta > 0,$ $\| \Sigma^{-1/2}\psi\| =1,$ and
$\left\Vert \varphi\right\Vert =1.$ Matrices $Y_{2}$ and $Y_{3}$
are $\left(  q-n_{1}\right)  \times p$ and $\left(  T-q\right)  \times p,$
respectively, where $q$ is the number of regressors and $T$ is the number of
observations. These matrices have, respectively, unrestricted and zero means
under both the null and the alternative.
SM \ref{Sufficiency and invariance considerations}
contains a
discussion of the relationship between alternative (\ref{REG0mean}) and a 
corresponding constraint on the coefficients of the 
untransformed regression model.

In the important example of comparison of $q$ group means, i.e. one-way
MANOVA, the null hypothesis imposes equality of all means, while a
rank one alternative would posit that the $q$ mean vectors lie along a
line, for example $\mu_k = \mu_1 + s_k \psi$ for scalar $s_k, k =2,
\ldots, q$ and $\psi \in \mR^p$. This will be a plausible reduction of
a global alternative hypothesis in some applications.

For REG$_{0}$, sufficiency and invariance arguments lead to tests based on the
solutions $\lambda_{1},...,\lambda_{p}$ of (\ref{basicJO}) with
\[
H=Y_{1}^{\prime}Y_{1}/n_{1}\text{ and }E=\Sigma.
\]
These solutions represent a multivariate analog of the difference between the
sum of squared residuals in the restricted and unrestricted
regressions. Under the null hypothesis, $n_{1}H$ is distributed as $W_{p}(n_{1},\Sigma)$
whereas under the alternative, it is distributed as $W_{p}(n_{1},\Sigma,n_{1}\Phi)$,
where $\Phi=\theta\Sigma^{-1}\psi\psi^{\prime}$.
Without loss of generality, we may assume that $\Sigma=I_{p}$.

The canonical form of REG$_0$ is essentially equivalent to the
recently studied setting of \textit{matrix denoising}
\begin{equation*}
  Y_1 = M + Z.
\end{equation*}
References which point to a variety of applications include
\cite{ccs10,shno13,nada14,doga14}.
Often $M$ is assumed to have low rank, and the matrix valued noise $Z$
to have i.i.d. Gaussian entries. Here we test $M = 0$ versus a rank
one alternative.

For REG, similar arguments lead to tests based on the $p$ solutions $\lambda
_{1},...,\lambda_{p}$ of (\ref{basic beta formJO}) with
\[
H=Y_{1}^{\prime}Y_{1}/n_{1}\text{ and }E=Y_{3}^{\prime}%
Y_{3}/n_{2},
\]
where the error d.f. $n_{2}=T-q$. These solutions represent a
multivariate analog of the $F$ ratio: the difference between the sum
of squared residuals in the restricted and unrestricted regressions to
the sum of squared residuals in the restricted regression. Again, 
we may assume wlog that $\Sigma$, although unknown to the statistician, 
equals $I_{p}$. 
Note that,
as $n_{2}\rightarrow\infty$ while $n_{1}$ and $p$ are held constant,
REG reduces to REG$_{0}$.

\smallskip
\textit{CCA.} \  Consider testing for 
independence between Gaussian vectors $x_{t}\in\mathbb{R}^{p}$ and $y_{t}%
\in\mathbb{R}^{n_{1}},$ given zero mean observations with $t=1,...,n_{1}%
+n_{2}.$ Partition the population and sample covariance matrices of the
observations $\left(  x_{t}^{\prime},y_{t}^{\prime}\right)  ^{\prime}$ into
\[
\left(
\begin{array}
[c]{cc}%
\Sigma_{xx} & \Sigma_{xy}\\
\Sigma_{yx} & \Sigma_{yy}%
\end{array}
\right)  \text{ and }\left(
\begin{array}
[c]{cc}%
S_{xx} & S_{xy}\\
S_{yx} & S_{yy}%
\end{array}
\right)  ,
\]
respectively. Under $H_0: \Sigma_{xy}=0,$ the alternative of
interest is
\begin{equation}
  \label{eq:ccapar}
\Sigma_{xy}=\sqrt{\frac{n_{1}\theta}{n_{1}\theta+n_{1}+n_{2}}}\psi
\varphi^{\prime},
\end{equation}
where the vectors of nuisance parameters $\psi\in\mathbb{R}^{p}$ and
$\varphi\in\mathbb{R}^{n_{1}}$ are normalized so that
\[
\lVert \Sigma_{xx}^{-1/2}\psi\rVert =\lVert \Sigma_{yy}%
^{-1/2}\varphi\rVert =1.
\]
The peculiar parameterizations of the alternative $\theta \neq 0$ in 
\eqref{REG0mean} and \eqref{eq:ccapar} are chosen to allow unified treatments of 
PCA and REG$_0$ and of SigD, REG and CCA in our main results,
Theorems \ref{LR asymptotics} and \ref{AsymptoticNormality} below.
% The peculiar form of the expression under the square root is chosen so as to
% simplify various expressions in the analysis that follow

The test can be based on the squared sample canonical correlations $\lambda
_{1},...,\lambda_{p}$, which are solutions to (\ref{basicJO}) with
\[
H=S_{xy}S_{yy}^{-1}S_{yx}\text{ and }E=S_{xx}.
\]
%We label Case 5 as the `canonical correlation analysis' (CCA) case.
Remarkably, the squared sample canonical correlations also solve
(\ref{basic beta formJO}) with different $H$ and $E,$ such that $E$ is a central
Wishart matrix and $H$ is a non-central Wishart matrix conditionally on a
random non-centrality parameter (see SM \ref{Identification of the parameters}).
%(for details, see Theorem 11.3.2 of Muirhead (1982)).

\smallskip
\textit{SMD.} \ 
% Finally, as discussed in the Introduction, we also consider Case 0, which we
% label as the `symmetric matrix denoising' (SMD) case. 
We observe a $p\times p$
matrix $X=\Phi+Z/\sqrt{p},$ where $Z$ is a noise matrix from the Gaussian
Orthogonal Ensemble (GOE), i.e. it is symmetric and 
%  Recall, that a symmetric
% matrix $Z$ belongs to GOE if
% its diagonal and sub-diagonal entries are independently distributed as%
\[
Z_{ii}\sim N\left(  0,2\right)  \text{ and }Z_{ij}\sim N\left(  0,1\right)
\text{ if }i>j\text{.}%
\]
We seek to make inference about a
symmetric rank-one \textquotedblleft signal\textquotedblright\ matrix
$\Phi=\theta \mathbb{\psi} \mathbb{\psi}^{\prime}$.
The null and the alternative hypotheses are given by (\ref{Basic hypothesesJO}).
The nuisance vector $\mathbb{\psi}\in\mathbb{R}^{p}$ is normalized so that
$\left\Vert \mathbb{\psi}\right\Vert =1$. The problem remains invariant under
the multiplication of $X$ from the left by an orthogonal matrix, and from the
right by its transpose. A maximal invariant statistic consists of the
solutions $\lambda_{1},...,\lambda_{p}$ to (\ref{basicJO}) with $H=X$ and
$E=I_{p}.$ We consider tests based on $\lambda_{1},...,\lambda_{p}$.

The SMD case can be viewed as a degenerate version of each of the above
cases.
For example, consider PCA with $p$ held fixed and $n_1 \to \infty$. 
Take $\Sigma = I_{p}$ for convenience and set 
$\Omega = I_{p} + \sqrt{p/n_1} \Phi$ 
%$\Omega = I + (p/n_1)^{1/2} \Phi$ 
with $\Phi = \theta \psi \psi'$, so
that the original value of the spike is rescaled to
be a local perturbation. 
Now  write $H$ in the form $\Omega^{1/2} \check H \Omega^{1/2}$ 
where $\check H \sim W_p(n_{1},I_{p})$.
A standard matrix central limit theorem for $p$ fixed, 
e.g. \cite[Th. 2.5.1]{fus10}, says that
\begin{equation*}
  \check H = I_p + Z/\sqrt{n_1} + o_{\mathrm{P}}(n_1^{-1/2}),
\end{equation*}
where $Z$ belongs to GOE. 
Writing $\Omega^{1/2} = I_{p} + \tfrac{1}{2} \sqrt{p/n_1} \Phi +
o(n_{1}^{-1/2})$, and 
introducing $\mu = \sqrt{n_1/p} \, (\lambda - 1)$, we can rewrite
\begin{equation*}
  \det (H - \lambda I_{p}) 
    = (p/n_1)^{p/2} \det [ \Phi +Z/\sqrt{p} - \mu I_{p} + o_{\rm P}(1)],
\end{equation*}
so that PCA degenerates to SMD.
Compare also \cite{bayi88}.

% For example, consider REG$_{0}$ with%
% \[
% \mathbb{E}Y_{1}^{}=\sqrt{\left(  p/n_{1}\right)  ^{1/2}n_{1}\theta}%
% \varphi\psi^{\prime},
% \]
% so that the original value of the spike $\theta$ (see equation (\ref{REG0mean}%
% )) is scaled by $\left(  p/n_{1}\right)  ^{1/2}$. Suppose now that $n_{1}$
% diverges to infinity while $p$ is held constant. Then, by a CLT,%
% \begin{equation}
% \Sigma^{-1/2}H\Sigma^{-1/2}-I_{p}=Z/\sqrt{n_{1}}+\sqrt{p/n_{1}}\eta\theta
% \eta^{\prime}+o_{\mathrm{P}}\left(  n_{1}^{-1/2}\right)  ,
% \label{equality link}%
% \end{equation}
% where $Z$ belongs to GOE and $\eta=\Sigma^{-1/2}\psi$. On the other hand,
% equation (\ref{basic}) is equivalent to%
% \begin{equation}
% \det\left(  \Sigma^{-1/2}H\Sigma^{-1/2}-\lambda I_{p}\right)  =0.
% \label{modified equation 1}%
% \end{equation}
% Multiplying it by $\sqrt{n_{1}/p}$ and using (\ref{equality link}), we see
% that equation (\ref{modified equation 1}) degenerates to%
% \[
% \det\left(  Z/\sqrt{p}+\eta\theta\eta^{\prime}-\mu I_{p}\right)  =0\text{ with
% }\mu=\sqrt{n_{1}/p}\left(  \lambda-1\right)  .
% \]
% Hence, REG$_{0}$ degenerates to SMD.

\begin{figure}[htbp]
\centering
\includegraphics[height=3.85in,width=3.85in]{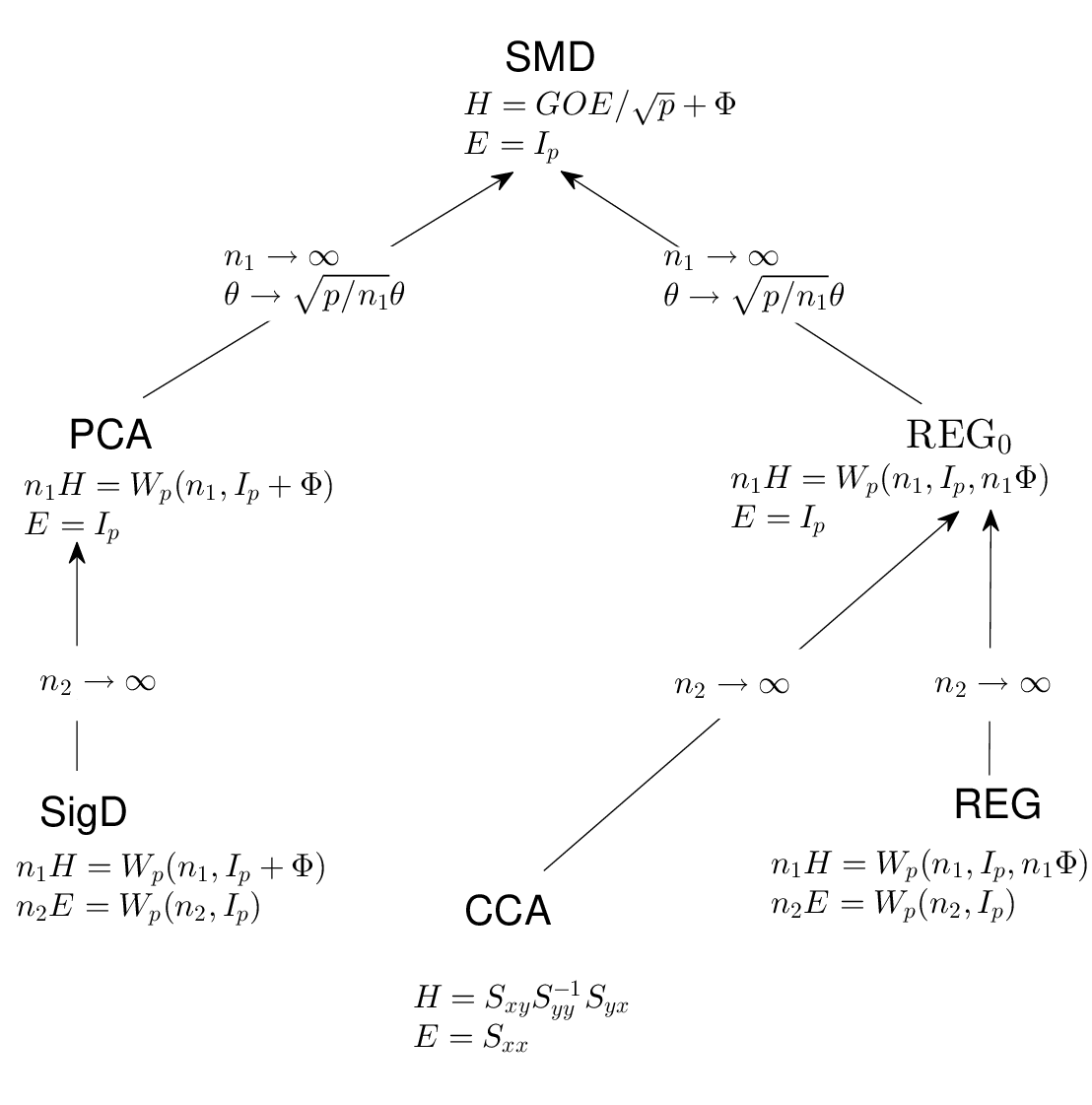}
\caption{Matrices $H$ and $E,$ and links between the different cases. 
Without loss of generality, matrix $E$ or, in SigD, REG, and CCA cases,
its population counterpart $\Sigma$ is assumed to be equal to $I_{p}$.
Matrix
$\Phi$ has the form $\theta \psi\psi^{\prime}$ with $\theta\geq0$ and
$\left\Vert \psi\right\Vert =1$.}
\label{diagram1}
\end{figure}

Indeed, each of the cases eventually degenerate to SMD via sequential
asymptotic links (SM \ref{Sequential asymptotic links between the cases} has details). 
For convenience, we summarize links between the different cases
and the definitions of the corresponding matrices $H$ and $E$ in Figure
\ref{diagram1}. 
We note that the SMD model has been studied recently,
e.g. \cite{bawa13,lemi16} and references therein,
though not with our techniques.
% We denote the $p$-dimensional Wishart distribution with $n$
% degrees of freedom, covariance parameter $\Sigma$, and non-centrality
% parameter $\Psi$ as $W_{p}\left(  n,\Sigma,\Psi\right)  .$ Recall that, if
% $A=B^{\prime}B,$ where the $n\times p$ matrix $B$ is $N\left(  M,I_{n}%
% \otimes\Sigma\right)  ,$ then $A\sim W_{p}\left(  n,\Sigma,\Psi\right)  $ with
% the non-centrality $\Psi=\Sigma^{-1}M^{\prime}M$. Notation $W_{p}\left(
% n,\Sigma\right)  $ is used for the central Wishart distribution. Without loss
% of generality, we assume that $\Sigma=I_{p}$.

Cases SMD, PCA, and REG$_{0}$, forming the upper half of the diagram,
correspond to random $H$ and deterministic $E.$ The cases in the lower half of
the diagram correspond to both $H$ and $E$ being random. Cases PCA and SigD
are \textquotedblleft parallel\textquotedblright\ to cases REG$_{0}$ and REG
in the sense that the alternative hypothesis is characterized by a rank one
perturbation of the covariance and of the non-centrality parameter of $H$ for
the former and for the latter two cases, respectively. Case CCA
\textquotedblleft stands alone\textquotedblright\ because of the different
structure of $H$ and $E.$ As discussed above, CCA can be reinterpreted in
terms of $H$ and $E$ such that $E$ is Wishart, but $H$ is a non-central
Wishart only after conditioning on a random non-centrality parameter.%

\section{The likelihood ratios}

Our goal is to study the asymptotic behavior of likelihood ratios
based on the observed eigenvalues
\[
\Lambda=\operatorname*{diag}\left\{  \lambda_{1},...,\lambda_{p}\right\}.
\]
Let $p(\Lambda; \theta)$ be the joint density of the eigenvalues under
the alternative and $p(\Lambda; 0)$ the corresponding density under
the null. 
James' formulas for these joint densities lead to our starting
point, which is a unified form for the likelihood ratio
% , which
% are defined as the ratios of the joint density of $\lambda_{1},...,\lambda
% _{p}$ under the alternative to that under the null hypothesis, where both
% densities are evaluated at the observed values of the $\lambda$'s. Let
% \[
% \Lambda=\operatorname*{diag}\left\{  \lambda_{1},...,\lambda_{p}\right\}  ,
% \]
% and let us denote the likelihood ratio corresponding to particular case `Case'
% = `SMD', `PCA', etc. as $L^{(\rm{Case})}\left(  \theta;\Lambda\right)$. Then
\begin{equation}
%L^{(\rm{Case})}\left(  \theta;\Lambda\right)  
L(\theta; \Lambda)
 = \frac{p(\Lambda;\theta)}{p(\Lambda; 0)}
=\alpha\left(  \theta\right) 
\left. _{\mathsf{p}}F_{\mathsf{q}} (a, b; \Psi, \Lambda)  \right. ,
\label{LRgeneral}
\end{equation}
where $\Psi=\Psi(\theta)$ is a $p$-dimensional matrix $\operatorname*{diag}\left\{
\Psi_{11},0,...,0\right\}  ,$ and the values of $\Psi_{11}$, $\alpha\left(
\theta\right)$, $\mathsf{p}$, $\mathsf{q}$, $a$, and $b$ are as given in Table
\ref{Table 2}.

\begin{table}[b] 
\centering
\caption{Parameters of the explicit expression (\ref{LRgeneral}) for the likelihood ratios. Here $n\equiv n_1+n_2$.}
\vspace{.1in}
\begin{tabular}
[c]{llllll}%
Case & $_{\mathsf{p}}F_{\mathsf{q}}$ & $\alpha\left(  \theta\right)  $ & $a$ & $b$ &
$\Psi_{11}$\\[8pt]
SMD & $_{0}F_{0}$ & $\exp\left(  -p\theta^{2}/4\right)  $ & \_ & \_ & $\theta p/2$\\[6pt]
PCA & $_{0}F_{0}$ & $\left(  1+\theta\right)  ^{-n_{1}/2}$ & \_ & \_ & $\theta
n_{1}/(2\left(  1+\theta\right)  )$\\[6pt]
SigD & $_{1}F_{0}$ & $\left(  1+\theta\right)  ^{-n_{1}/2}$ & $n/2$ & \_ &
$\theta n_{1}/\left(  n_{2}\left(  1+\theta\right)  \right)  $\\[6pt]
REG$_{0}$ & $_{0}F_{1}$ & $\exp\left(  -n_{1}\theta/2\right)  $ & \_ &
$n_{1}/2$ & $\theta n_{1}^{2}/4$\\[6pt]
REG & $_{1}F_{1}$ & $\exp\left(  -n_{1}\theta/2\right)  $ & $n/2$ & $n_{1}/2$
& $\theta n_{1}^{2}/\left(  2n_{2}\right)  $\\[6pt]
CCA & $_{2}F_{1}$ & $\left(  1+n_{1}\theta/n\right)  ^{-n/2}$ & $\left(
n/2,n/2\right)  $ & $n_{1}/2$ & $\theta n_{1}^{2}/\left(  n_{2}^{2}+n_{2}%
n_{1}\left(  1+\theta\right)  \right)  $\\
\end{tabular}
\label{Table 2}
\end{table}

For SMD, we prove that $L\left(  \theta;\Lambda\right)  $ is as in
(\ref{LRgeneral}) in SM \ref{SMD entry}. For PCA, the explicit form of the
likelihood ratio is derived in 
%Onatski et al (2013). 
\cite{omh13}.
For SigD, REG$_{0}$, and
REG, the expressions (\ref{LRgeneral}) 
%with the parameters given in Table \ref{Table 2} 
follow, respectively, from equations (65), (68), and (73) of
%James (1964). 
\cite{jame64}.
For CCA, the expression is a corollary of \cite[Th. 11.3.2]{muir82}.
%Theorem 11.3.2 of Muirhead (1982). 
Further details appear in SM \ref{Identification of the parameters}.

Recall that hypergeometric functions of two matrix arguments $\Psi$ and
$\Lambda$ are defined as%
\[
_{\mathsf{p}}F_{\mathsf{q}}\left(  a,b;\Psi,\Lambda\right)  =\sum_{k=0}^{\infty
}\frac{1}{k!}\sum_{\kappa\vdash k}\frac{\left(  a_{1}\right)  _{\kappa
}...\left(  a_{\mathsf{p}}\right)  _{\kappa}}{\left(  b_{1}\right)  _{\kappa
}...\left(  b_{\mathsf{q}}\right)  _{\kappa}}\frac{C_{\kappa}\left(  \Psi\right)
C_{\kappa}\left(  \Lambda\right)  }{C_{\kappa}\left(  I_{p}\right)  },
\]
where $a=\left(  a_{1},...,a_{\mathsf{p}}\right)  $ and $b=\left(  b_{1}%
,....,b_{\mathsf{q}}\right)  $ are parameters, $\kappa$ are partitions of the
integer $k$, $\left(  a_{j}\right)  _{\kappa}$ and $\left(  b_{i}\right)
_{\kappa}$ are the generalized Pochhammer symbols, and $C_{\kappa}$ are the
zonal polynomials, e.g. \cite[Def. 7.3.2.]{muir82}.
% As mentioned in the
% Introduction, James' (1964) classification of the multivariate statistical
% problems is based on the type of $_{\mathsf{p}}F_{\mathsf{q}}$ that occur in related
% probability distributions. 
%The function $_{0}F_{0}$ of \textit{exponential
%type} corresponds to the first class represented by PCA; the function
%$_{1}F_{0}$ of \textit{binomial type} corresponds to the second class
%represented by SigD; the function $_{0}F_{1}$ of \textit{Bessel type} is
%associated with the third class represented by REG$_{0}$; the
%\textit{confluent hypergeometric function} $_{1}F_{1}$ is associated with the
%fourth class represented by REG; and the \textit{Gaussian hypergeometric
%function} $_{2}F_{1}$ corresponds to the fifth class represented by CCA. 
Note that some links between the cases illustrated in Figure \ref{diagram1} can
also be established via asymptotic relations between 
the hypergeometric functions.
%in the different rows of Table \ref{Table  2}. 
For example, the confluence relations
%(see,for example, Olver et al (2010), equation 35.8.9)
\begin{align*}
\left.  _{0}F_{0}\left(  \Psi,\Lambda\right)  \right.   &  =\lim
_{a\rightarrow\infty}\left.  _{1}F_{0}\left(  a;a^{-1}\Psi,\Lambda\right)
\right. \qquad \text{and}\\ 
\left.  _{0}F_{1}\left(  b;\Psi,\Lambda\right)  \right.   &  =\lim
_{a\rightarrow\infty}\left.  _{1}F_{1}\left(  a,b;a^{-1}\Psi,\Lambda\right)
\right.  
\end{align*} e.g. \cite[eq. 35.8.9]{ollobocl10},
imply the links
SigD $\mapsto$ PCA and REG $\mapsto$ REG$_{0}$ 
as $n_{2}\rightarrow\infty$ for 
$p$ and $n_{1}$ held constant.

In the next section, we shall study the asymptotic behavior of the likelihood
ratios (\ref{LRgeneral}) as $n_{1},n_{2},$ and $p$ go to infinity so that%
\begin{equation}
c_{1}\equiv p/n_{1}\rightarrow\gamma_{1}\in\left(  0,1\right)  \text{
and }c_{2}\equiv p/n_{2}\rightarrow\gamma_{2}\in\left(  0,1\right]  .
\label{asymptotic regime}%
\end{equation}
We denote this asymptotic regime by $\mathbf{n},p\rightarrow_{\boldsymbol{\gamma}%
}\infty,$ where $\mathbf{n}=\left\{  n_{1},n_{2}\right\}  $ and
$\boldsymbol{\gamma}=\left\{  \gamma_{1},\gamma_{2}\right\}  .$ To make our
exposition as uniform as possible, we use this notation for all the cases,
even though the simpler ones, such as SMD, do not refer to $\mathbf{n}$. We briefly 
discuss possible extensions of our analysis to the
situations with $\gamma_{1}\geq1$ in Section \ref{sec-Concluding remarks}.

We are interested in the asymptotics of the likelihood ratios under the null
hypothesis, that is when the true value of the spike, $\theta_{0}$, equals
zero. 
First, some background on the eigenvalues. 
% Before turning to the next section, let us provide a relevant background
% on the asymptotics of $\Lambda$. 
Under the null, $\lambda_{1},...,\lambda_{p}$
are the eigenvalues of $GOE/\sqrt{p}$ in the SMD case; of $W_{p}\left(
n_{1},I_{p}\right)  /n_{1}$ for PCA and REG$_{0}$; and of a 
% a scaled (by a factor of $n_{2}/n_{1}$) 
$p$-dimensional multivariate beta matrix,
e.g. \cite[p. 110]{muir78},  with parameters $n_{1}/2$ and $n_{2}/2$ 
and here scaled by a factor of $n_2/n_1$, in the SigD, REG, and CCA cases.
% For a
% definition of the multivariate beta, see Muirhead (1982), p. 110.
%In each of these cases, 
The empirical distribution
 of  $\lambda_{1},...,\lambda_{p}$
\begin{equation*}
\hat F 
= \frac{1}{p} \sum_{j=1}^p I\{ \lambda_j \leq \lambda \}
\end{equation*}
is well known, \citep{bai99}, to converge weakly almost
surely (a.s.) in each case:
\begin{equation*}
\hat F 
%= \frac{1}{p} \sum_{j=1}^p I\{ \lambda_j \leq \lambda \}
\Rightarrow F_{\boldsymbol{\gamma}} =
\begin{cases}
  F^{\rm SC} & \text{for  SMD} \\
  F^{\rm MP} & \text{for  PCA, REG}_0 \\
  F^{\rm W} & \text{for  SigD, REG, CCA}, 
\end{cases}
\end{equation*}
the semi-circle, Marchenko-Pastur and (scaled) Wachter distributions
respectively. 
% Let%
% \[
% \hat{F}^{(\rm{Case})}\left(  \lambda\right)  =\frac{1}{p}%
% %TCIMACRO{\dsum _{j=1}^{p}}%
% %BeginExpansion
% {\displaystyle\sum_{j=1}^{p}}
% %EndExpansion
% 1\left\{  \lambda_{j}\leq\lambda\right\}
% \]
% be the empirical distribution of $\lambda_{1},...,\lambda_{p}$. As is well
% known (see Bai (1999)), as $\mathbf{n},p\rightarrow_{\boldsymbol{\gamma}}\infty,$
% $\hat{F}^{(\rm{Case})}$ almost surely (a.s.) weakly converges
% \[
% \hat{F}^{(\rm{Case})}\Rightarrow F_{\boldsymbol{\gamma}}^{\lim},
% \]
% where $F_{\boldsymbol{\gamma}}^{\lim}$ is the Semi-circle distribution $F^{\rm{SC}}$ in
% SMD case; the Marchenko-Pastur distribution $F_{\gamma_{1}}^{\rm{MP}}$ in PCA and
% REG$_{0}$ cases; and the (scaled) Wachter distribution $F_{\boldsymbol{\gamma}%
% }^{\rm{W}}$ in SigD, REG, and CCA cases. 
Table \ref{Table 3} recalls the explicit
forms of these limiting distributions. The cumulative distribution
functions $F_{\boldsymbol{\gamma}}^{\lim}\left(  \lambda\right)  $ are linked in
the sense that 
% \begin{alignat*}{2}
%   F_\gamma^{\rm W}(\lambda) & \to F_{\gamma_1}^{\rm MP}(\lambda) &
%   \qquad 
%       & \text{when } \gamma_2 \to 0 \\
%   F_{\gamma_1}^{\rm MP}(\sqrt{\gamma_1} \lambda + 1) 
%       &  \to F^{\rm SC} (\lambda) & & \text{when } \gamma_1 \to 0.
% \end{alignat*}
\begin{alignat*}{3}
  & F_{\boldsymbol{\gamma}}^{\rm W}(\lambda) & &\to F_{\gamma_1}^{\rm MP}(\lambda) 
  \qquad 
      & & \text{as } \gamma_2 \to 0, % \text{  and } 
\\
  & F_{\gamma_1}^{\rm MP}(\sqrt{\gamma_1} \lambda + 1) 
      & & \to F^{\rm SC} (\lambda) & & \text{as } \gamma_1 \to 0.
\end{alignat*}

\begin{table}[t]
  \centering
  \caption{Semi-circle, Marchenko-Pastur and scaled Wachter
    distributions} 
\vspace{.1in}
  \begin{tabular}[h]{p{.3in}cccc}
Case & \qquad $F_{\boldsymbol{\gamma}}^{\rm lim}$ & Density, $\lambda \in [\beta_-,\beta_+]$
& $\beta_\pm$ & Threshold $\bar \theta$ \\[12pt]
SMD  & \qquad SC & $\dfrac{R(\lambda)}{2 \pi} $ & $\pm 2$ & $1$  \\[12pt]
PCA  REG$_0$ & \qquad MP & $\dfrac{R(\lambda)}{2 \pi \gamma_1 \lambda} $ 
     & $(1 \pm \sqrt{\gamma_1})^2 $ & $\sqrt{\gamma_1}$  \\[12pt]
SigD REG CCA& \qquad W  & \qquad \quad $\dfrac{(\gamma_1 +
  \gamma_2)R(\lambda)}{2 \pi 
  \gamma_1 \lambda   (\gamma_1 - \gamma_2 \lambda)} $  \qquad \quad
     & $\gamma_1 \left( \dfrac{\rho \pm 1}{\rho \pm \gamma_2} \right)^2$ 
     & $\dfrac{\rho + \gamma_2}{1 - \gamma_2}$ \\[30pt]
& & $R(\lambda) = \sqrt{(\beta_+ - \lambda)(\lambda -
    \beta_-)}$ & 
\multicolumn{2}{c}{$\rho = \sqrt{\gamma_1 + \gamma_2 - \gamma_1 \gamma_2} $}  
% \multicolumn{5}{r}{$R(\lambda) = \sqrt{(\beta_+ - \lambda)(\lambda -
%     \beta_-)},  \qquad  
%     \rho = \sqrt{\gamma_1 + \gamma_2 - \gamma_1 \gamma_2} $}  
  \end{tabular}
  \label{Table 3}
\end{table}

%\newpage
%For what follows it will be important that 
If $\varphi$ is a `well-behaved' function, the centered \textit{linear
spectral statistic}
\begin{equation}%
%TCIMACRO{\dsum _{j=1}^{p}}%
%BeginExpansion
{\displaystyle\sum_{j=1}^{p}}
%EndExpansion
\varphi\left(  \lambda_{j}\right)  -p\int\varphi\left(  \lambda\right)
\mathrm{d}F_{\mathbf{c}}^{\lim}\left(  \lambda\right)  ,
\label{linear spectral statistics}%
\end{equation}
 converges in distribution to a
Gaussian random variable
in each of the semicircle \cite{bayo05},
Marchenko-Pastur \cite{basi04}
and Wachter \cite{zheng12} cases.
% The corresponding CLTs are established in Bai and
% Yao (2005), Bai and Silverstein (2004), and Zheng (2012) for the cases of the
% Semi-circle, Marchenko-Pastur, and Wachter limiting distributions,
% respectively. 
Note that the centering constant is defined in terms of
$F_{\mathbf{c}},$ where $\mathbf{c}=\left\{  c_{1},c_{2}\right\}  .$
That is, the \textquotedblleft correct centering\textquotedblright\ can be
computed using the densities from Table \ref{Table 3}, where $\gamma_{1}$ and
$\gamma_{2}$ are replaced by $c_{1}\equiv p/n_{1}$ and $c_{2}\equiv p/n_{2}$, respectively.

Finally, let us recall the behavior of the largest eigenvalue $\lambda_{1}$
under the alternative hypothesis. 
As long as $\theta \leq \bar \theta$,
the phase transition threshold 
% The value of the threshold is 
reported in Table \ref{Table 3},
the top eigenvalue
$\lambda_{1} \to \beta_+$,
the upper boundary of support of $F_{\boldsymbol{\gamma}}$, 
almost surely.
When $\theta>\bar{\theta}$, $\lambda_{1}$ separates from `the
bulk' of the other eigenvalues and a.s. converges to a point strictly
above $\beta_+$.
%the upper boundary of the support of $F_{\boldsymbol{\gamma}}.$ 
For details, we refer to
\cite{mai07,basi06,nadasil10,on07,dhjo14,bahupazh14}
%  to Ma\"{\i}da (2007), Baik and Silverstein
% (2006), Nadakuditi and Silverstein (2010), Onatski (2007), Dharmawansa et al
% (2014a), and Bao et al (2014) 
for the respective cases SMD, PCA, SigD, REG$_{0}$, REG, and
CCA.

The fact that $\lambda_{1}$ converges to different limits under the null and
under the alternative hypothesis sheds light on the behavior of the likelihood
ratio when $\theta$ is above the phase transition threshold $\bar
\theta$. 
In such \textit{super-critical} cases, the likelihood ratio
degenerates. The sequences of measures corresponding to the distributions of
$\Lambda$ under the null and under super-critical alternatives are
asymptotically mutually singular as $\mathbf{n},p\rightarrow_{\boldsymbol{\gamma}%
}\infty$, as shown in \cite{mai07} and \cite{omh13}
% Montanari et al (2014) and Onatski et al (2013)  
for SMD and PCA respectively. In contrast, as we show below, the
sequences of measures corresponding to the distributions of $\Lambda$ under
the null and under \textit{sub-critical} alternatives $\theta < \bar \theta$ 
%(is below the threshold)
are mutually contiguous, and the likelihood ratio converges to a
Gaussian process. 
In the super-critical setting, an analysis of the likelihood ratios
under local alternatives appears in \cite{dhjoon14}.

\section{Contour integral representation}

The asymptotic behavior of the likelihood ratios (\ref{LRgeneral})
depends on that 
of $_{\mathsf{p}}F_{\mathsf{q}}\left(  a,b;\Psi,\Lambda\right)  $. 
When the dimension of the matrix arguments remains fixed, 
there is a large
and well established literature on the asymptotics of $_{\mathsf{p}}F_{\mathsf{q}%
}\left(  a,b;\Psi,\Lambda\right)  $ for large parameters and  norm of the
matrix arguments, see \cite{muir78} for a review. 
%while the dimensionality of the latter remains fixed
%(see Muirhead (1978) for a review). 
In contrast, relatively little is known
about when the dimensionality of the matrix
arguments $\Psi,\Lambda$ diverge to infinity.
It is this regime we study in this paper, noting that
%We exploit the fact that, since we study
in single-spiked models, the matrix argument $\Psi$ has rank one. This allows us
to represent $_{\mathsf{p}}F_{\mathsf{q}}\left(  a,b;\Psi,\Lambda\right)  $ in the
form of a contour integral of a hypergeometric function with a single scalar
argument. Such a representation implies contour integral representations for
the corresponding likelihood ratios.
%  which we summarize in the following
% lemma. 
% The results of the lemma are used below to derive the asymptotics of
% the likelihood ratios via the Laplace approximation.

% In what follows, we omit the superscripts `(Case)' and `lim' for quantities
% such as $L^{(\rm{Case})}\left(  \theta;\Lambda\right)  ,$ $\hat{F}^{(\rm{Case})}\left(
% \lambda\right)  ,$ and $F_{\mathbf{c}}^{\lim}\left(  \lambda\right)  $ to
% simplify our notation. However, we shall use these superscripts to identify
% particular instances, when necessary.

\begin{lemma}
\label{DJLemma}Assume that $p\leq\min\left\{  n_{1},n_{2}\right\}  .$ Let
$\mathcal{K}$ be a contour in the complex plane $\mathbb{C}$ that starts at
$-\infty$, encircles $0$ and $\lambda_{1},...,\lambda_{p}$ counterclockwise,
and returns to $-\infty$. Then%
\begin{equation}
L\left(  \theta;\Lambda\right)  =\frac{\Gamma\left(  s+1\right)  \alpha\left(
\theta\right)  q_{s}}{\Psi_{11}^{s}2\pi\mathrm{i}}\int_{\mathcal{K}}\left.
_{\mathsf{p}}F_{\mathsf{q}}\right.  \left(  a-s,b-s;\Psi_{11}z\right)
%TCIMACRO{\dprod _{j=1}^{p}}%
%BeginExpansion
{\displaystyle\prod_{j=1}^{p}}
%EndExpansion
\left(  z-\lambda_{j}\right)  ^{-1/2}\mathrm{d}z, \label{integral}%
\end{equation}
where $s=p/2-1,$ the values of $\alpha\left(  \theta\right)  ,$ $\Psi_{11},$
$a,$ $b,$ $\mathsf{p},$ and $\mathsf{q}$ for the different cases are given in Table
\ref{Table 2}; $a-s$ and $b-s$ denote vectors with elements $a_{j}-s$ and
$b_{j}-s,$ respectively; 
% the hypergeometric function under the integral is the
% standard hypergeometric function of a scalar argument; 
and
\[
q_{s}=%
%TCIMACRO{\dprod _{j=1}^{\bar{p}}}%
%BeginExpansion
{\displaystyle\prod_{j=1}^{\mathsf{p}}}
%EndExpansion
\frac{\Gamma\left(  a_{j}-s\right)  }{\Gamma\left(  a_{j}\right)  }%
%TCIMACRO{\dprod _{i=1}^{\bar{q}}}%
%BeginExpansion
{\displaystyle\prod_{i=1}^{\mathsf{q}}}
%EndExpansion
\frac{\Gamma\left(  b_{i}\right)  }{\Gamma\left(  b_{i}-s\right)  }.
\]
In cases SigD and CCA, we require, in addition, that the contour $\mathcal{K}$
does not intersect $\left[  \Psi_{11}^{-1},\infty\right)  $, which ensures the
analyticity of the integrand in an open subset of $\mathbb{C}$ that includes
$\mathcal{K}$.
\end{lemma}

The statement of the lemma immediately follows from 
%Proposition 1 of Dharmawansa and Johnstone (2014) 
\cite[Prop. 1]{dhjo14}
and from equation (\ref{LRgeneral}). Our next
step is to apply the Laplace approximation to integrals (\ref{integral}). To
this end, we shall transform the right hand side of (\ref{integral}) so that
it has a \textquotedblleft Laplace form\textquotedblright%
\begin{equation}
L\left(  \theta;\Lambda\right)  =\sqrt{\pi p}\frac{1}{2\pi\mathrm{i}}%
\int_{\mathcal{K}}\exp\left\{  -(p/2)f(z;\theta)\right\}  g(z;\theta)\mathrm{d}z.
\label{Laplace form}%
\end{equation}
The dependence on $\theta$ will usually not be shown explicitly.
Leaving $\sqrt{\pi p}/\left(  2\pi\mathrm{i}\right)  $ separate from $g(z)$
allows us to choose $f(z)$ and $g(z)$ that are bounded in probability, and
makes some of the expressions below more compact. In order to apply the
Laplace approximation, we shall deform the contour of integration so that it
passes through a critical point $z_{0}$ of $f(z)$ and is such that
$\operatorname{Re}f(z)$ is strictly increasing as $z$ moves away from $z_{0}$
along the contour, at least in a vicinity of $z_{0}$.

\subsection{The Laplace form}

We shall transform (\ref{integral}) to (\ref{Laplace form}) in three steps. As
a result, functions $f$ and $g$ will have the forms of a sum and a product,%
\begin{align}
f\left(  z\right)   &  =f_{\mathrm{c}}+f_{\mathrm{e}}\left(  z\right)  +f_{\mathrm{h}}\left(
z\right)  \text{ and}  \label{eq:fdecomp}\\
g(z)  &  =g_{\mathrm{c}}\times g_{\mathrm{e}}\left(  z\right)  \times
g_{\mathrm{h}}(z), \notag
\end{align}
where $f_{\mathrm{c}}$ and $g_{\mathrm{c}}$ do not depend on $z$.
The subscripts (c,e,h) are mnemonic for `coefficient', `eigenvalues' and
`hypergeometric'. 

First, using the definitions of $\alpha\left(  \theta\right)  ,$ $q_{s},$
$\Psi_{11}$ and employing Stirling's approximation, we obtain a decomposition%
\begin{equation}
\frac{\Gamma\left(  s+1\right)  \alpha\left(  \theta\right)  q_{s}}{\sqrt{\pi
p}\Psi_{11}^{s}}=\exp\left\{  -(p/2)f_{\mathrm{c}}\right\}  g_{\mathrm{c}}, \label{decomposition 0}
\end{equation}
where $g_{\mathrm{c}}$ remains bounded as $\mathbf{n},p\rightarrow_{\boldsymbol{\gamma}%
}\infty$. The values of $f_{\mathrm{c}}$ and $g_{\mathrm{c}}$ are given in Table \ref{Table 2a}.
Details of the derivation are given in SM \ref{Derivations for Table 4}. 
% It should be noted that 
% $f_{\mathrm{c}}^{\rm{REG}},f_{\mathrm{c}}^{\rm{CCA}}\rightarrow
% f_{\mathrm{c}}^{\rm{REG_{0}}}$ and $f_{\mathrm{c}}^{\rm{SigD}}\rightarrow f_{\mathrm{c}}^{\rm{PCA}}$ as
% $c_{2}\rightarrow0$.%
% $f_{\mathrm{c}}^{(\rm{REG})},f_{\mathrm{c}}^{(\rm{CCA})}\rightarrow
% f_{\mathrm{c}}^{(\rm{REG_{0}})}$ and $f_{\mathrm{c}}^{(\rm{SigD})}\rightarrow f_{\mathrm{c}}^{(\rm{PCA})}$ as
% $c_{2}\rightarrow0$.%

\begin{table}[tbp] 
\centering
\caption{Values of $f_{\mathrm{c}}$ and $\check g_{\mathrm{c}} = g_{\mathrm{c}}/(1+o(1))$ for the different cases.
The terms $o(1)$ do not depend on $\theta$ and converge to zero as
$\mathbf{n},p\rightarrow _{\boldsymbol{\gamma }}\infty $. In the table, $l(\theta)=1+(1+\theta)c_2/c_1$ 
and $r^{2}=c_{1}+c_{2}-c_{1}c_{2}$.}
\vspace{.1in}
\begin{tabular}[h]{lll}
   
Case & $f_c$ & $\check g_c = g_c/(1+o(1))$ \\
\hline
   \\

SMD  & $1 + \theta^2/2 + \log \theta$ &  $\theta$ \\[8pt]

    PCA  & $1 + \dfrac{1-c_1}{c_1} \log(1+\theta) + \log 
    \dfrac{\theta}{c_1}$ & $\theta (1+\theta)^{-1} c_1^{-1}$ \\[8pt]

    SigD & $f_{\rm c}^{\rm PCA} + f_{10}$ & $\check g_c^{\rm PCA} \check g_{10}$ \\[8pt]

    REG$_0$ & $1 + \dfrac{\theta + c_1}{c_1} + \log \dfrac{\theta}{c_1}
    + \dfrac{1-c_1}{c_1} \log(1-c_1)$  & $\theta c_1^{-1}
      (1-c_1)^{-1/2}$  \\[5pt] 

   REG   & $f_{\rm c}^{{\rm REG}_0} + f_{10}$ & $\check g_c^{{\rm REG}_0}
   \check g_{10}$
   \\[8pt]
   
CCA   & $f_{\rm c}^{\rm REG} + f_{21}$ & $\check g_c^{\rm REG} \check
   g_{10}/l(\theta)$ \\[6pt] 
   \hline \\[-2pt]
& $f_{10} = -1  - \dfrac{r^2}{c_1 c_2} \log \dfrac{r^2}{c_1 + c_2}
 + \log \dfrac{c_1 + c_2}{c_1}$ \qquad
 & 
$\check g_{10} = c_1^{-1} r (c_1 + c_2)^{1/2}$  \\[5pt]
& $f_{21} = -1 - \dfrac{\theta}{c_1} - \dfrac{r^2}{c_1 c_2} \log
\dfrac{r^2}{c_1 l(\theta)}$ 
\\[8pt]
\hline
\end{tabular}
  
\label{Table 2a}

\end{table}

Second, we consider the decomposition%
\begin{equation}
{\displaystyle\prod_{j=1}^{p}}
\left(  z-\lambda_{j}\right)  ^{-1/2}=\exp\left\{  -(p/2)f_{\mathrm{e}}(z)\right\}
g_{\mathrm{e}}(z), \label{step1}
\end{equation}
where%
\begin{equation}
f_{\mathrm{e}}(z)=\int\ln\left(  z-\lambda\right)  \mathrm{d}F_{\mathbf{c}}(\lambda),
\label{f1Case}
\end{equation}
and%
\begin{equation}
g_{\mathrm{e}}(z)=\exp\left\{  -(p/2)\int\ln\left(  z-\lambda\right)
\mathrm{d}\left(  \hat{F}\left(  \lambda\right)  -F_{\mathbf{c}}\left(
\lambda\right)  \right)  \right\}  . \label{Deltap}%
\end{equation}
For $f_{\mathrm{e}}(z)$ and $g_{\mathrm{e}}(z)$ to be well-defined we need $z$ not to belong
to the support of $F_{\mathbf{c}},$ which we assume.
%\footnote{}
In addition, $z \notin \text{supp} (\hat F)$ since by definition 
contour $\mathcal{K}$ encircles it.
%$z\in\mathcal{K}$ does not belong to such support.] 
Note that $g_{\mathrm{e}}(z)$ is
the exponent of a linear spectral statistic, which converges to a Gaussian
random variable as $\mathbf{n},p\rightarrow_{\boldsymbol{\gamma}}\infty$ under the
null hypothesis. 
% Since $F_{\mathbf{c}}^{\rm{W}}(\lambda)\rightarrow F_{c_{1}}%
% ^{\rm{MP}}(\lambda)$ as $c_{2}\rightarrow0,$ we have $f_{\mathrm{e}}^{(\rm{SigD})}%
% (z)=f_{\mathrm{e}}^{(\rm{REG})}(z)=f_{\mathrm{e}}^{(\rm{CCA})}(z)$ converging to $f_{\mathrm{e}}^{(\rm{PCA})}%
% (z)=f_{\mathrm{e}}^{(\rm{REG_{0}})}(z)$.

Third and finally, we describe a decomposition%
\begin{equation}
\left.  _{\mathsf{p}}F_{\mathsf{q}}\right.  \left(  a-s,b-s;\Psi_{11}z\right)
=\exp\left\{  -(p/2)f_{\mathrm{h}}(z)\right\}  g_{\mathrm{h}}(z). 
\label{LaplaceForm}
\end{equation}
For the $\mathsf{q}=0$ cases, the corresponding $\left.  _{\mathsf{p}}F_{\mathsf{q}%
}\right.  $ can be expressed in terms of elementary functions.
Indeed, $_{0}F_{0} (z) = e^z$ and $_{1}F_{0} (a; z) = (1 - z)^{-a}$. 
We set%
\begin{equation}
f_{\mathrm{h}}(z)=\left\{
\begin{array}
[c]{ll}%
-z\theta & \text{for SMD}\\
-z\theta/\left(  c_{1}\left(  1+\theta\right)  \right)  \text{ } & \text{for
PCA}\\
\ln\left[  1-c_{2}z\theta/\left\{  c_{1}\left(  1+\theta\right)  \right\}
\right]  r^{2}/\left(  c_{1}c_{2}\right)  & \text{for SigD,}%
\end{array}
\right.   \label{elementary}%
\end{equation}
and%
\begin{equation}
g_{\mathrm{h}}(z)=\left\{
\begin{array}
[c]{ll}%
1 & \text{for SMD and PCA}\\
\left[  1-c_{2}z\theta/\left\{  c_{1}\left(  1+\theta\right)  \right\}
\right]  ^{-1} & \text{for SigD.}%
\end{array}
\right.   
\label{component2SMD}
\end{equation}
% As $c_{2}\rightarrow0$, $f_{\mathrm{h}}^{(\rm{SigD})}(z)$ converges to $f_{\mathrm{h}}%
% ^{(\rm{PCA})}(z).$ Since, as has been shown above, a similar convergence holds for
% $f_{\mathrm{c}}$ and $f_{\mathrm{e}},$ we have $f^{(\rm{SigD})}(z)\rightarrow f^{(\rm{PCA})}(z)$ as
% $c_{2}\rightarrow0$. Combining (\ref{f1Case}) and (\ref{elementary}) with the
% information supplied by Table \ref{Table 2a}, we also see that $f^{(\rm{PCA})}%
% (z)\rightarrow f^{(\rm{SMD})}(z)$ as $c_{1}\rightarrow0$ after the transformations
% $\theta\mapsto\sqrt{c_{1}}\theta$ and $z\mapsto\sqrt{c_{1}}z+1$.

Unfortunately, for the $\mathsf{q}=1$ cases, the corresponding $\left.
_{\mathsf{p}}F_{\mathsf{q}}\right.  $ do not admit exact representations in terms of
elementary functions. Therefore, we shall consider their asymptotic
approximations instead. Let%
\[
m=\left(  n_{1}-p\right)  /2\text{ and }\kappa=\left(  n-p\right)
/\left(  n_{1}-p\right)  .
\]
Further, let%
\begin{equation}
\eta_{j}=\left\{
\begin{array}
[c]{ll}%
z\theta/\left(  1-c_{1}\right)  ^{2} & \text{for }j=0\\
z\theta c_{2}/\left[  c_{1}\left(  1-c_{1}\right)  \right]  & \text{for }j=1\\
z\theta c_{2}^{2}/\left[  c_{1}^{2}l\left(  \theta\right)  \right]  &
\text{for }j=2
\end{array}
\right.  , \label{etas}%
\end{equation}
where%
\begin{equation}
l\left(  \theta\right)  =1+\left(  1+\theta\right)  c_{2}/c_{1}.
\label{loftheta}%
\end{equation}
With this notation, we have%
\begin{equation}
\left.  _{\mathsf{p}}F_{\mathsf{q}}\right.  =\left\{
\begin{array}
[c]{ll}%
\left.  _{0}F_{1}\right.  \left(  m+1;m^{2}\eta_{0}\right)  \equiv F_{0} &
\text{for REG}_{0}\\
\left.  _{1}F_{1}\right.  \left(  m\kappa+1;m+1;m\eta_{1}\right)  \equiv
F_{1} & \text{for REG}\\
\left.  _{2}F_{1}\right.  \left(  m\kappa+1,m\kappa+1;m+1;\eta
_{2}\right)  \equiv F_{2} & \text{for CCA.}%
\end{array}
\right.   \label{pFq}%
\end{equation}

The function $F_{0}(z)$ can be expressed in terms of the modified Bessel function
of the first kind $I_{m}\left(  \cdot\right)$, see \cite[eq. 9.6.47]{abst64}, as 
%(see Abramowitz and Stegun (1964), equation 9.6.47)%
\begin{equation}
F_{0}=\Gamma\left(  m+1\right)  \left(  m^{2}\eta_{0}\right)  ^{-m/2}%
I_{m} \big(  2m\eta_{0}^{1/2}\big)  . \label{bessel rep}%
\end{equation}
This representation allows us to use a known uniform asymptotic approximation
of the Bessel function \cite[eq. 9.7.7]{abst64} 
%(see Abramowitz and Stegun (1964), equation 9.7.7) to
to obtain Lemma \ref{uniform1}, proven in SM \ref{Approximation to 0F1}. To state it let%
\begin{equation}
\varphi_{0}\left(  t\right)  =\ln t-t-\eta_{0}/t+1\text{ and }t_{0}=\left(
1+\sqrt{1+4\eta_{0}}\right)  /2. \label{phi0}%
\end{equation}
Further, for any $\delta>0,$ let $\Omega_{0\delta}$ be the set of $\eta_{0}%
\in\mathbb{C}$ such that%
\[
\left\vert \arg\eta_{0}\right\vert \leq\pi-\delta,\text{ and }\eta_{0}\neq0.
\]

\begin{lemma}
\label{uniform1}As $m\rightarrow\infty,$ we have%
\begin{equation}
F_{0}=\left(  1+4\eta_{0}\right)  ^{-1/4}\exp\left\{  -m\varphi_{0}\left(
t_{0}\right)  \right\}  \left(  1+o(1)\right)  . \label{BesselApprox}%
\end{equation}
The convergence $o(1)\rightarrow0$ holds uniformly with respect to $\eta
_{0}\in\Omega_{0\delta}$ for any $\delta>0$.
\end{lemma}

%We would like to point out 
To foreshadow our results for $F_1(z)$ and $F_2(z)$, we note
that the right hand side of (\ref{BesselApprox})
can be formally linked, via (\ref{bessel rep}), to the saddle-point
approximation of the integral representation, 
%(see Watson (1944), p. 181)%
see \cite[p. 181]{watson4},
\[
I_{m}\left(  2m\eta_{0}^{1/2}\right)  =\frac{\eta_{0}^{m/2}e^{m}}%
{2\pi\mathrm{i}}\int_{-\infty}^{(0+)}\exp\left\{  -m\varphi_{0}\left(
t\right)  \right\}  t^{-1}\mathrm{d}t.
\]
Point $t_{0}$ can be interpreted as a saddle point of $\varphi_{0}\left(
t\right)  ,$ and the term $\left(  1+4\eta_{0}\right)  ^{-1/4}$ in
(\ref{BesselApprox}) can be interpreted as a factor of $\left(  \varphi
_{0}^{\prime\prime}\left(  t_{0}\right)  \right)  ^{-1/2}$.

Turning now to functions $F_1(z)$ and $F_2(z)$, 
to obtain uniform asymptotic approximations,
we use the contour integral representations, 
%(see Olver et al (2010), equations 13.4.9 and 15.6.2)
see \cite[eqs. 13.4.9 and 15.6.2]{ollobocl10}, 
\begin{equation}
F_{j}=\frac{C_{m}}{2\pi\mathrm{i}}\int_{0}^{(1+)}\exp\left\{  -m\varphi
_{j}\left(  t\right)  \right\}  \psi_{j}\left(  t\right)  \mathrm{d}t,
\label{integralrepF11}%
\end{equation}
where%
\begin{equation}
C_{m}=\frac{\Gamma\left(  m+1\right)  \Gamma\left(  m\left(  \kappa
-1\right)  +1\right)  }{\Gamma\left(  m\kappa+1\right)  }, \label{Calpha}%
\end{equation}%
\begin{equation}
\varphi_{j}(t)=\left\{
\begin{array}
[c]{ll}%
-\eta_{j}t-\kappa\ln t+\left(  \kappa-1\right)  \ln\left(
t-1\right)  & \text{for }j=1\\
-\kappa\ln\left(  t/\left(  1-\eta_{j}t\right)  \right)  +\left(
\kappa-1\right)  \ln\left(  t-1\right)  & \text{for }j=2
\end{array}
\right.  , \label{phij}%
\end{equation}
and%
\begin{equation}
\psi_{j}\left(  t\right)  =\left\{
\begin{array}
[c]{ll}%
\left(  t-1\right)  ^{-1} & \text{for }j=1\\
\left(  t-1\right)  ^{-1}\left(  1-\eta_{j}t\right)  ^{-1} & \text{for }j=2
\end{array}
\right.  . \label{gammaj}%
\end{equation}
For $j=2$, the contour does not encircle $1/\eta_{2},$ and the representation
is valid for $\eta_{2}$ such that $\left\vert \arg\left(  1-\eta_{2}\right)
\right\vert <\pi$. We derive
%obtain the following lemma by deriving 
a saddle-point
approximation to the integral in (\ref{integralrepF11})
to be summarized in Lemma \ref{1F1approximation} below. The relevant saddle
points are%
\begin{equation}
t_{j}=\left\{
\begin{array}
[c]{cl}%
\frac{1}{2\eta_{j}}\left\{  \eta_{j}-1+\sqrt{\left(  \eta_{j}-1\right)
^{2}+4\kappa\eta_{j}}\right\}  & \text{for }j=1\\
\frac{1}{2\eta_{j}\left(  \kappa-1\right)  }\left\{  -1+\sqrt
{1+4\kappa\left(  \kappa-1\right)  \eta_{j}}\right\}  & \text{for
}j=2
\end{array}
\right.  . \label{tj}%
\end{equation}

We shall need the following additional notation. Let
\begin{equation}
\omega_{j}=\arg\varphi_{j}^{\prime\prime}\left(  t_{j}\right)  +\pi
\quad \text{ and} \quad 
\omega_{0j}=\arg\left(  t_{j}-1\right)  , \label{omegaj}%
\end{equation}
where the branches of $\arg\left(  \cdot\right)  $ are chosen so that
$\left\vert \omega_{j}+2\omega_{0j}\right\vert \leq\pi/2.$ 
% Further, for any
% small $\delta>0$ let $\Omega_{1\delta}$ be the set of $\left(  \kappa
% ,\eta_{1}\right)  \in\mathbb{R}\times\mathbb{C}$ such that $\delta
% \leq\kappa-1\leq1/\delta,$ and%
% \[
% \operatorname{Re}\eta_{1}\geq-2\kappa+1,\text{ }\mathrm{dist}\left(
% \eta_{1},\mathbb{R}\backslash\left[  0,\infty\right)  \right)  \geq
% \delta,\text{ }\left\vert \eta_{1}\right\vert \leq1/\delta.
% \]
% Similarly, let $\Omega_{2\delta}$ be the set of $\left(  \kappa,\eta
% _{2}\right)  \in\mathbb{R}\times\mathbb{C}$ such that $\delta\leq
% \kappa-1\leq1/\delta,$ and%
% \[
% \mathrm{dist}\left(  \eta_{2},\mathbb{R}\backslash\left[  0,1\right]  \right)
% \geq\delta,\text{ }\left\vert \eta_{2}\right\vert \leq1/\delta.
% \]
% Here, for any $A\subseteq\mathbb{C}$ and $B\subseteq\mathbb{C}$,
% $\mathrm{dist}\left(  A,B\right)  =\inf_{a\in A,b\in B}\left\vert
% a-b\right\vert .$ Figure \ref{uniform} shows cross-sections of $\Omega
% _{1\delta}$ and $\Omega_{2\delta}$ for fixed $\kappa.$%

%\begin{figure}[h]
%\centering
%\includegraphics[height=2.19in,width=4.74in]{uniform.eps}
%\caption{Cross-sections of the sets $\Omega_{j\delta}$ for $\kappa=2$ and
%$\delta=0.1.$ The horizontal and vertical axes correspond to the real and
%purely imaginary numbers, respectively.}
%\label{uniform}
%\end{figure}

\begin{lemma}
\label{1F1approximation}As $m\rightarrow\infty,$ we have for $j=1,2$%
\begin{equation}
F_{j}=C_{m}\psi_{j}\left(  t_{j}\right)  e^{-\mathrm{i}\omega_{j}/2}\left\vert
2\pi m\varphi_{j}^{\prime\prime}\left(  t_{j}\right)  \right\vert ^{-1/2}%
\exp\left\{  -m\varphi_{j}\left(  t_{j}\right)  \right\}  \left(
1+o(1)\right)  . \label{ConfluentApprox}%
\end{equation}
The convergence $o(1)\rightarrow0$ holds uniformly with respect to $\left(
\kappa,\eta\right)  \in\Omega_{j\delta}$ for any $\delta>0$,
where $\Omega_{j\delta}$ are as defined in Table \ref{Table Def Omegajdelta}.
\end{lemma}

\begin{table}[b]
\centering
\caption{Definition of $\Omega_{j\delta}$ from Lemma \ref{1F1approximation}.}
\vspace{.1in}
\begin{tabular}
[c]{ll}%
\multicolumn{2}{c}{$\Omega_{j\delta}=\Omega_{\delta}\cap \hat{\Omega}_{j\delta}$ with the following 
$\Omega_{\delta}$ and $\hat{\Omega}_{j\delta}$} \\[8pt] \hline
Set                      & Definition: pairs $(x,z) \in \mathbb{R} \times \mathbb{C}$ s.t.  \\[8pt]  
$\Omega_{\delta}$        & $\delta \leq x-1 \leq 1/\delta$, $\left\vert z \right\vert \leq 1/\delta$,
                           and $\inf_{y \in \mathbb{R}\backslash[0,\infty)}
                            \left\vert z-y \right\vert \geq \delta$\\[6pt]
$\hat{\Omega}_{1\delta}$ & $\operatorname{Re}z \geq -2x+1$\\[6pt]
$\hat{\Omega}_{2\delta}$ & $\inf_{y \in \mathbb{R}\backslash(-\infty,1]}
\left\vert z-y \right\vert \geq \delta$ and $x$ is unconstrained.\\ \hline
\end{tabular}
\label{Table Def Omegajdelta}
\end{table}

Point-wise asymptotic approximation (\ref{ConfluentApprox}) was established in
%Passemier et al (2014)
\cite{pamcch14} 
for $j=1,$ and in 
%Paris (2013a,b) 
\cite{par13a, par13b}
for $j=2.$ However,
those papers do not study the uniformity of the approximation error, which is
important for our analysis. 
Lemma \ref{1F1approximation} is proved at
length in SM \ref{Approximations to 1F1 and 2F1}. It is fair to say that 
the corresponding derivations constitute the technically most challenging
part of our analysis. This further highlights the technical 
difficulties that occur when going from SMD, PCA, and SigD cases to REG$_{0}$, REG, 
and CCA.  

Using Lemmas \ref{uniform1} and \ref{1F1approximation}, and Stirling's
approximation%
\begin{equation}
C_{m}=\frac{\sqrt{\pi p\left(  1-c_{1}\right)  }}{r}\exp\left\{  m\left(
\kappa-1\right)  \ln\left(  \kappa-1\right)  -m\kappa
\ln\kappa\right\}  \left(  1+o(1)\right)  \label{StirlingCm}%
\end{equation}
we set the components of the \textquotedblleft Laplace form\textquotedblright%
\ (\ref{LaplaceForm}) of $\left.  _{\mathsf{p}}F_{\mathsf{q}}\right.  $ 
for the $\mathsf{q}=1$ cases as follows%
\begin{equation}
f_{\mathrm{h}}(z)=
\begin{cases}
\frac{1-c_{1}}{c_{1}}\varphi_{0}\left(  t_{0}\right)  & \text{REG}_{0}\\
\frac{1-c_{1}}{c_{1}}\left(  \varphi_{j}\left(  t_{j}\right)  +\kappa
\ln\kappa-\left(  \kappa-1\right)  \ln\left(  \kappa-1\right)
\right)  \text{ } & \text{REG, CCA}%
\end{cases}
\label{f2}%
\end{equation}
% \begin{equation}
% f_{\mathrm{h}}(z)=\left\{
% \begin{array}
% [c]{ll}%
% \frac{1-c_{1}}{c_{1}}\varphi_{0}\left(  t_{0}\right)  & \text{for REG}_{0}\\
% \frac{1-c_{1}}{c_{1}}\left(  \varphi_{j}\left(  t_{j}\right)  +\kappa
% \ln\kappa-\left(  \kappa-1\right)  \ln\left(  \kappa-1\right)
% \right)  \text{ } & \text{for REG and CCA}%
% \end{array}
% \right.  \label{f2}%
% \end{equation}
and%
\begin{equation}
g_{\mathrm{h}}(z)=
\begin{cases}
\left(  1+4\eta_{0}\right)  ^{-1/4}\left(  1+o\left(  1\right)  \right)  &
\text{REG}_{0}\\
\sqrt{c_{1}/r^{2}}e^{-\mathrm{i}\omega_{j}/2}\left\vert \varphi_{j}%
^{\prime\prime}\left(  t_{j}\right)  \right\vert ^{-1/2}\psi_{j}\left(
t_{j}\right)  \left(  1+o\left(  1\right)  \right)  & \text{REG, CCA}%
\end{cases}
\label{g2}%
\end{equation}
% \begin{equation}
% g_{\mathrm{h}}(z)=\left\{
% \begin{array}
% [c]{ll}%
% \left(  1+4\eta_{0}\right)  ^{-1/4}\left(  1+o\left(  1\right)  \right)  &
% \text{for REG}_{0}\\
% \sqrt{c_{1}/r^{2}}e^{-\mathrm{i}\omega_{j}/2}\left\vert \varphi_{j}%
% ^{\prime\prime}\left(  t_{j}\right)  \right\vert ^{-1/2}\psi_{j}\left(
% t_{j}\right)  \left(  1+o\left(  1\right)  \right)  & \text{for REG and CCA}%
% \end{array}
% \right.  . \label{g2}%
% \end{equation}
To express $t_{j}$ and $\eta_{j}$ in terms of $z,$ one should use (\ref{tj})
and (\ref{etas}). We do not need to know how exactly the  $o\left(
  1\right)  $ in 
(\ref{g2}) depend on $z$. For our purposes, the knowledge of the fact that
$o\left(  1\right)  $ are analytic functions of $\eta_{j}$ that converge to
zero uniformly with respect to $\left(  \kappa,\eta_{j}\right)  \in
\Omega_{j\delta}$ is sufficient. The analyticity of $o(1)$ follows from the
analyticity of the functions on the left hand sides, and of the factors of
$1+o(1)$ on the right hand sides of the equations (\ref{BesselApprox}) and
(\ref{ConfluentApprox}).

% Using the definitions of $\varphi_{j}$ and $t_{j},$ it is straightforward to
% verify that $f_{\mathrm{h}}^{(\rm{REG})}(z)$ and $f_{\mathrm{h}}^{(\rm{CCA})}(z)$ converge to
% $f_{\mathrm{h}}^{(\rm{REG_{0}})}(z)$ as $c_{2}\rightarrow0$. Since, as has been shown
% above, a similar convergence holds for $f_{\mathrm{c}}$ and $f_{\mathrm{e}},$ we have
% $f^{(\rm{REG})}(z),f^{(\rm{CCA})}(z)\rightarrow f^{(\rm{REG_{0}})}(z)$ as $c_{2}\rightarrow
% 0$. Elementary calculations that use equations (\ref{f1Case}), (\ref{phi0}),
% (\ref{f2}) together with the explicit forms of $f_{\mathrm{c}}^{(\rm{REG_{0}})}$ and
% $f_{\mathrm{c}}^{(\rm{SMD})}$ given in Table \ref{Table 2a} show that $f^{(\rm{REG_{0}}%
% )}(z)\rightarrow f^{(\rm{SMD})}(z)$ as $c_{1}\rightarrow0$ after transformations
% $\theta\mapsto\sqrt{c_{1}}\theta$ and $z\mapsto\sqrt{c_{1}}z+1$.

\smallskip
\textit{Confluences of functions $f$.} \ 
As $c_2 \to 0$ with $c_1$ held fixed, we have
\begin{equation}
  \label{eq:c2to0}
  \begin{split}
  f^{\mathrm{SigD}}(z) & \to f^{\mathrm{PCA}}(z), \\
  f^{\mathrm{REG}}(z), f^{\mathrm{CCA}}(z) & \to f^{\mathrm{REG}_0}(z). 
  \end{split}
\end{equation}
Also, as $c_1 \to 0$,
\begin{equation}
  \label{eq:c1to0}
  f^{\mathrm{PCA}}(z), f^{\mathrm{REG}_0}(z)  \to f^{\mathrm{SMD}}(z),
\end{equation}
after making the substitutions
$\theta \to \sqrt c_1 \theta$ and $z \to \sqrt c_1 z+1$ on the left
hand side. Some details appear in SM \ref{sec:proof-confluences}.
% \texttt{[currently Appendix E.]}

\subsection{Saddlepoints and Contours of steep descent}

We shall now show how to deform contours $\mathcal{K}$ in (\ref{Laplace form})
into the contours of steep descent. First, we find saddle points of functions
$f(z)$ for each of the six cases. Note that 
\begin{equation*}
%  - f_{\mathrm{e}}'(z) 
 - \mathrm{d}f_{\mathrm{e}}(z)/\mathrm{d}z
     = \int(\lambda-z)^{-1}\mathrm{d}F_{\mathbf{c}}(\lambda)
     = m_{\mathbf{c}}\left(  z\right),
\end{equation*}
% the negative of the derivative of $f_{\mathrm{e}}(z)$
% equals $\int(\lambda-z)^{-1}\mathrm{d}F_{\mathbf{c}}(\lambda)$, which is
the Stieltjes transform 
%$m_{\mathbf{c}}\left(  z\right)$ 
of $F_{\mathbf{c}}$.
Although the Stieltjes transform is formally defined on $\mathbb{C}^{+}$, the
definition remains valid on the part of the real line outside the support
$\left[  b_{-},b_{+}\right]  $ of $F_{\mathbf{c}}$. Since we assume that
%$\gamma_{1}\leq1,$ 
$p \leq n_{1}$,
$F_{\mathbf{c}}$ does not have any non-trivial mass at $0$.
%for sufficiently large $\mathbf{n}$ and~$p$.

To find saddle points $z_{0}$ of $f(z)$ we therefore solve the equation%
\begin{equation}
m_{\mathbf{c}}\left(  z\right)  =\mathrm{d}f_{\mathrm{h}}(z)/\mathrm{d}z.
\label{extremum points}%
\end{equation}
A proof of the following lemma appears in SM \ref{saddle points z0}. 

\begin{lemma}
\label{critical}
The saddle points $z_{0}(\theta,\mathbf{c})$ of $f(z)$ satisfy
\begin{equation}
z_{0}(\theta,\mathbf{c})=\left\{
\begin{array}
[c]{ll}%
\theta+1/\theta & \text{for SMD}\\
\left(  1+\theta\right)  \left(  \theta+c_{1}\right)  /\theta & \text{for PCA
and REG}_{0}\\
\left(  1+\theta\right)  \left(  \theta+c_{1}\right)  /\left[  \theta l\left(
\theta\right)  \right]  & \text{for SigD, REG, and CCA.}%
\end{array}
\right.  \label{zcase}%
\end{equation}
For $\theta\in\left(  0,\bar{\theta
}_{\mathbf{c}}\right)  $, 
%and sufficiently large $\mathbf{n}$ and $p$ as $\mathbf{n}%
%,p\rightarrow_{\boldsymbol{\gamma}}\infty,$
$z_{0}>b_{+}$, where $\bar{\theta}_{\mathbf{c}}$ is the threshold corresponding to $F_{\mathbf{c}}$,
which is an analogue of the threshold $\bar{\theta}_{\boldsymbol{\gamma}} \equiv \bar{\theta}$ corresponding 
to $F_{\boldsymbol{\gamma}}$ given in Table \ref{Table 3}.
% and $b_{+}$ is the upper boundary of support of $F_{\mathbf{c}%
%}$. 
\end{lemma}

As $c_{2}\rightarrow0$ while $c_{1}$ stays constant, the value of $z_{0}$ for
SigD, REG, and CCA converges to that for PCA and REG$_{0}.$ The latter value
in its turn converges to the value of $z_{0}$ for SMD when $c_{1}\rightarrow
0$, after the transformations $\theta\mapsto\sqrt{c_{1}}\theta$ and
$z_{0}\mapsto\sqrt{c_{1}}z_{0}+1$. Precisely, solving equation%
\[
\sqrt{c_{1}}z_{0}+1=\left(  1+\sqrt{c_{1}}\theta\right)  \left(  \sqrt{c_{1}%
}\theta+c_{1}\right)  /\left(  \sqrt{c_{1}}\theta\right)
\]
for $z_{0}$ and taking limit as $c_{1}\rightarrow0$ yields $z_{0}%
=\theta+1/\theta$.

\begin{remark}
For all the six cases that we study, $f(z_{0})$ equals zero. 
SM \ref{sec: verification of remark 5} has a verification of this important fact.
\label{fz00}
\end{remark}

\begin{remark}
As $\mathbf{n},p \rightarrow_{\boldsymbol{\gamma}} \infty$, 
$z_{0}(\theta,\mathbf{c}) \rightarrow z_{0}(\theta,\boldsymbol{\gamma})>\beta_{+}$, where
the latter inequality holds for any $\theta \in \left( 0, \bar{\theta} \right)$.
Since $\lambda_{1} \overset{a.s.}{\rightarrow} \beta_{+}$ the inequality
$z_{0}(\theta,\mathbf{c})>\lambda_{1}$ must hold with probability approaching one as
$\mathbf{n},p \rightarrow_{\boldsymbol{\gamma}} \infty$. 
\label{z0>lambda1}
\end{remark}

\begin{figure}[ptb]
\centering
\includegraphics[height=2.82in,width=3.242in]{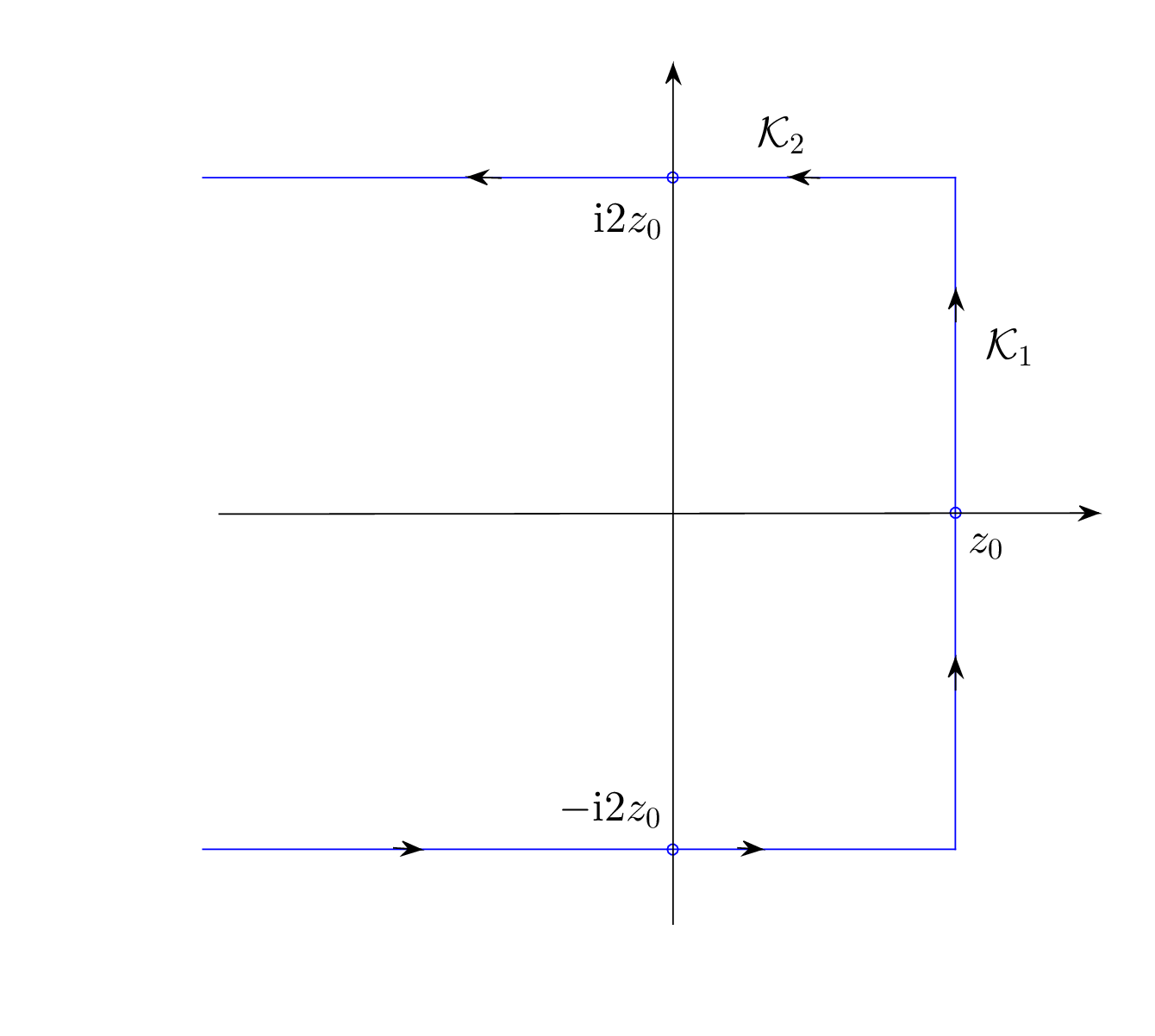}
\caption{Deformed contour $\mathcal{K}$ for SMD, PCA, and SigD.}
\label{Case 0 contour}
\end{figure}

For the rest of the paper, assume that $\theta\in\left(  0,\bar{\theta
}\right)  $. We deform contour $\mathcal{K}$ in (\ref{Laplace form}) so that
it passes through the saddle point $z_{0}$ as follows. Let $\mathcal{K}%
=\mathcal{K}_{+}\cup\mathcal{K}_{-},$ where $\mathcal{K}_{-}$ is the complex
conjugate of $\mathcal{K}_{+}$ and $\mathcal{K}_{+}=\mathcal{K}_{1}%
\cup\mathcal{K}_{2}.$ For SMD, PCA, and SigD, let%
\begin{align}
\mathcal{K}_{1}  &  =\left\{  z_{0}+\mathrm{i}t:0\leq t\leq2z_{0}\right\}
\text{ and}\label{contours}\\
\mathcal{K}_{2}  &  =\left\{  x+\mathrm{i}2z_{0}:-\infty<x\leq z_{0}\right\}
. \label{contours 1}%
\end{align}
The deformed contour is shown on Figure \ref{Case 0 contour}.

Note that the singularities of the integrand in (\ref{Laplace form}) are
situated at $z=\lambda_{j}$ (plus an additional singularity at $z=c_{1}%
(1+\theta)/\left(  \theta c_{2}\right)  <z_{0}$ for SigD). 
Since 
%$\lambda
%_{1}\overset{a.s.}{\rightarrow}\beta_{+}$ and $z_{0}>b_{+},$ inequality
$z_{0}>\lambda_{1}$ holds with probability approaching one as
$\mathbf{n},p\rightarrow_{\boldsymbol{\gamma}}\infty$, Cauchy's
theorem ensures that the deformation of the contour does not change the value of $L\left(
\theta;\Lambda\right)  $ with probability approaching one as $\mathbf{n}%
,p\rightarrow_{\boldsymbol{\gamma}}\infty.$

Strictly speaking, the deformation of the contour is not continuous because
$\mathcal{K}_{+}$ does not approach\ $\mathcal{K}_{-}$ at $-\infty$. In
particular, in contrast to the original contour, the deformed one is not
\textquotedblleft closed\textquotedblright\ at $-\infty$. Nevertheless, such
an \textquotedblleft opening up\textquotedblright\ at $-\infty$ does not lead
to the change of the value of the integral because the integrand converges
fast to zero in absolute value as $\operatorname{Re}z\rightarrow-\infty$.

\begin{remark}  \label{rem:neglig}
In the event of asymptotically negligible probability that the deformed
contour $\mathcal{K}$ does not encircle all $\lambda_{j},$ we not only lose
the equality (\ref{Laplace form}) but also face the difficulty that function
$g(z)$ ceases to be well defined as the definition of
$g_{\mathrm{e}}(z)$ contains a 
logarithm of a non-positive number. To eliminate any ambiguity, if such an
event holds we shall redefine $g_{\mathrm{e}}(z)$ as unity.
\end{remark}

\begin{figure}[ptb]
\centering
\includegraphics[height=2.98in,width=3.238in]{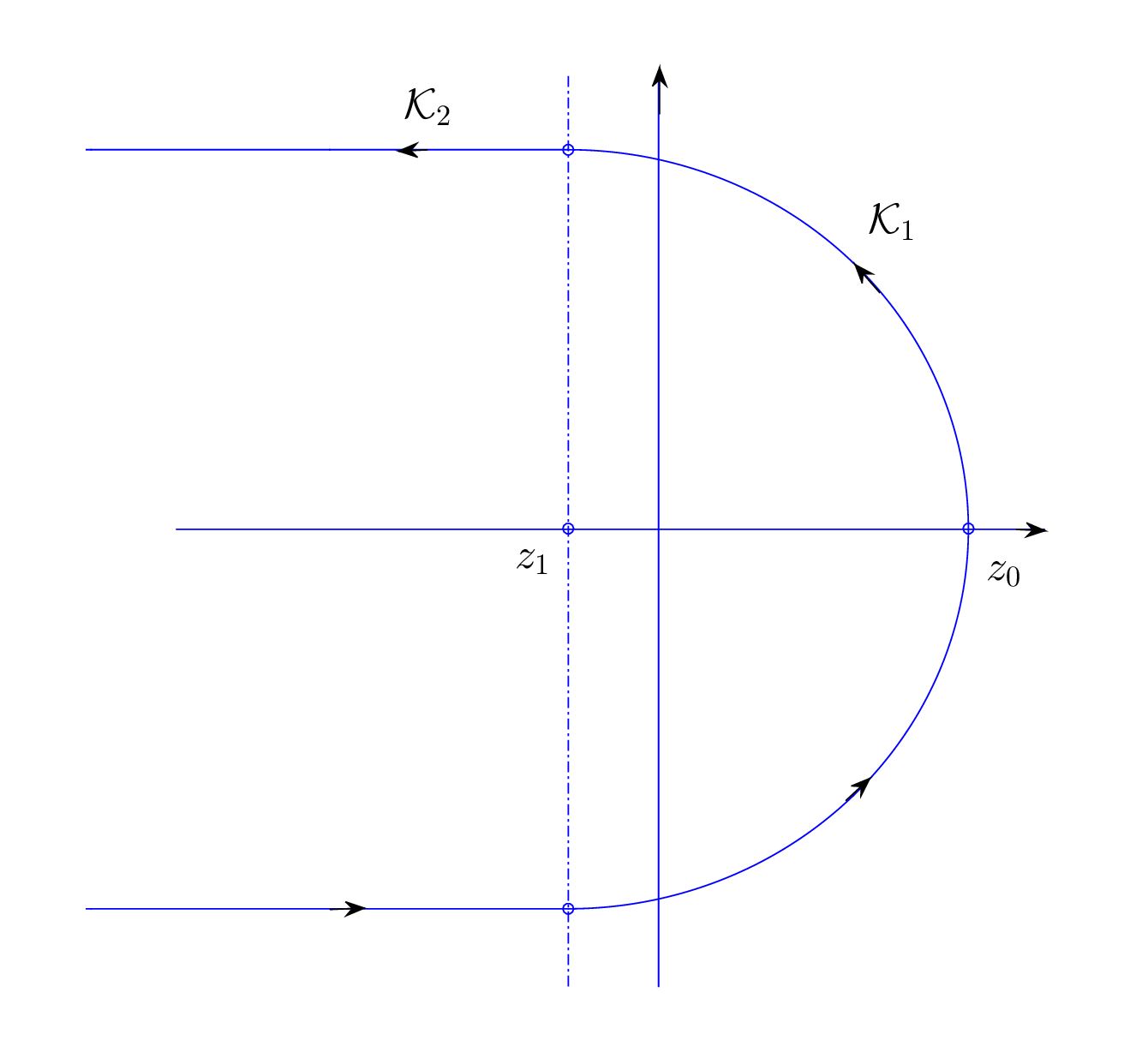}
\caption{Deformed contour $\mathcal{K}$ for REG$_{0}$ and CCA.}
\label{contourgeneral}
\end{figure}

For REG$_{0}$ and CCA, let%
\[
z_{1}=\left\{
\begin{array}
[c]{ll}%
-\left(  1-c_{1}\right)  ^{2}/\left[  4\theta\right]  & \text{for REG}_{0}\\
-c_{1}\left(  1-c_{1}\right)  ^{2}l\left(  \theta\right)  /\left[  4\theta
r^{2}\right]  & \text{for CCA}%
\end{array}
\right.  ,
\]
and let%
\begin{align*}
\mathcal{K}_{1}  &  =\left\{  z_{1}+\left\vert z_{0}-z_{1}\right\vert
\exp\left\{  \mathrm{i}\gamma\right\}  :\gamma\in\left[  0,\pi/2\right]
\right\}  \text{ and}\\
\mathcal{K}_{2}  &  =\left\{  z_{1}-x+\left\vert z_{0}-z_{1}\right\vert
\exp\left\{  \mathrm{i}\pi/2\right\}  :x\geq0\right\}  .
\end{align*}
The corresponding contour $\mathcal{K}$ is shown on Figure
\ref{contourgeneral}. Similarly to the SMD, PCA and SigD cases, the
deformation of the contour in (\ref{Laplace form}) to $\mathcal{K}$ does not
change the value of $L\left(  \theta;\Lambda\right)  $ with probability
approaching one as $\mathbf{n},p\rightarrow_{\boldsymbol{\gamma}}\infty.$%

For REG, deformed contour $\mathcal{K}$ in $z$-plane is simpler to describe as
an image of a contour $\mathcal{C}$ in $\tau$-plane, where $\tau=\eta_{1}%
t_{1}$ with%
\begin{equation}
\eta_{1}=z\theta c_{2}/\left[  c_{1}\left(  1-c_{1}\right)  \right]
\label{eta1}%
\end{equation}
and $t_{1}$ as defined in (\ref{tj}). Let $\mathcal{C}=\mathcal{C}_{+}%
\cup\mathcal{C}_{-},$ where $\mathcal{C}_{-}$ is the complex conjugate of
$\mathcal{C}_{+}$ and $\mathcal{C}_{+}=\mathcal{C}_{1}\cup\mathcal{C}_{2},$
and let%
\begin{align*}
\mathcal{C}_{1}  &  =\left\{  -\kappa+\left\vert \tau_{0}+\kappa
\right\vert \exp\left\{  \mathrm{i}\gamma\right\}  :\gamma\in\left[
0,\pi/2\right]  \right\}  \text{ and}\\
\mathcal{C}_{2}  &  =\left\{  -\kappa-x+\left\vert \tau_{0}%
+\kappa\right\vert \exp\left\{  \mathrm{i}\pi/2\right\}  :x\geq0\right\}
,
\end{align*}
where $\tau_{0}=\left(  \theta+c_{1}\right)  /\left(  1-c_{1}\right)  .$

Using (\ref{eta1}) and the identity
\begin{equation}
\eta_{1}=\tau\left(  \tau+1\right)  /(\tau+\kappa),
\label{identity for eta1}%
\end{equation}
we obtain%
\begin{equation}
z=\frac{c_{1}\left(  1-c_{1}\right)  }{\theta c_{2}}\frac{\tau\left(
\tau+1\right)  }{\tau+\kappa}. \label{transformation}%
\end{equation}
We define the deformed contour $\mathcal{K}$ in $z$-plane as the image of
$\mathcal{C}$ under the transformation $\tau\rightarrow z$ given by
(\ref{transformation}). The parts $\mathcal{K}_{+},\mathcal{K}_{-}%
,\mathcal{K}_{1}$ and $\mathcal{K}_{2}$ of $\mathcal{K}$ are defined as the
images of the corresponding parts of $\mathcal{C}$. Note that $\tau_{0}$ is
transformed to $z_{0}$ so that $\mathcal{K}$ passes through the saddle point
$z_{0}.$

The next lemma, proven in SM \ref{sec: contours fo steep descent}, shows that $\mathcal{K}_{1}$
are contours of steep descent of $-\operatorname{Re}f\left(  z\right)  $ for
all the six cases, SMD, PCA, SigD, REG$_{0}$, REG, and CCA.

\begin{lemma}
\label{ContourProp1F1} For any of the six cases that we study, as $z$ moves
along the corresponding $\mathcal{K}_{1}$ away from $z_{0}$,
$-\operatorname{Re}f\left(  z\right)  $ is strictly decreasing.
\end{lemma}

\section{Laplace approximation}

The goal of this section is to derive Laplace approximations to the
integral \eqref{integral} 
% \[
% L\left(  \theta;\Lambda\right)  =\sqrt{\pi p}\frac{1}{2\pi\mathrm{i}}%
% \int_{\mathcal{K}}\exp\left\{  -(p/2)f(z)\right\}  g(z)\mathrm{d}z
% \]
for the six cases that we study. First, consider a general integral%
\[
I_{p,\omega}=\int_{\mathcal{K}_{p,\omega}}e^{-p\phi_{p,\omega}(z)}%
\chi_{p,\omega}(z)\mathrm{d}z,
\]
where $p$ is large, $\omega \in \Omega \subset \mathbb{R}^k$ is a
$k$-dimensional parameter, and $\mathcal{K}_{p,\omega}$ is
a path in $\mathbb{C}$ that starts at
$a_{p,\omega}$ and ends at $b_{p,\omega}$.
We allow $\chi_{p,\omega}(z)$ to be a random element of the normed space of
continuous functions on $\mathcal{K}_{p,\omega}$ with the supremum norm.
Assume that there is a domain $T_{p,\omega} \supset
\mathcal{K}_{p,\omega}$ on which
for sufficiently large $p$,
$\phi_{p,\omega}(z)$ and  $\chi_{p,\omega}(z)$ are single-valued
holomorphic functions of $z$,
%in a domain that contains $\mathcal{K}_{p,\omega}$, 
in the case of $\chi_{p,\omega}$
with probability increasing to $1$. 

% where $p\rightarrow\infty,$ $\omega$ is a $k$-dimensional parameter that
% belongs to a subset $\Omega$ of $\mathbb{R}^{k},$ $\mathcal{K}_{p,\omega}$ is
% a path in $\mathbb{C}$ that starts at $a_{p,\omega}$ and ends at $b_{p,\omega
% }$, and for sufficiently large $p,$ $\phi_{p,\omega}(z)$ is a single-valued
% holomorphic function of $z$ in a domain $T_{p,\omega}$ that contains
% $\mathcal{K}_{p,\omega}$.

% Furthermore, we suppose that for any $\delta>0$, there exists $\bar{p}$ such
% that for any $p>\bar{p},$ $\chi_{p,\omega}(z)$ is a single-valued holomorphic
% function of $z$ in the domain $T_{p,\omega}$ with probability larger than
% $1-\delta$. 

We describe an extension of the Laplace approximation detailed by Olver
\cite[p. 127]{olve97} 
%  is a fairly straightforward extension of Theorem 7.1 of
% Olver (1997), p. 127 
to a situation in which functions $\phi,$ $\chi$ and
contour $\mathcal{K}$ depend on $p$ and $\omega$ 
and in addition $\chi$ is random.
In Olver's original theorem, 
%in different notation, 
both functions and contour are fixed. 
In what follows, however, we omit subscripts $p$ and $\omega$ from 
$\phi_{p,\omega},$ $\chi_{p,\omega},$ $\mathcal{K}_{p,\omega},$
etc. to lighten notation.

Suppose that $\phi^{\prime} (z)  =0$ at $z_{0}$ which is an
interior point of $\mathcal{K},$ and suppose that $\operatorname{Re}\phi(z)$
is strictly increasing as $z$ moves away from $z_{0}$ along the path. In other
words, the path $\mathcal{K}$ is a contour of steep descent\ of
$-\operatorname{Re}\phi(z)$. Denote a closed segment of $\mathcal{K}$
contained between $z_{1}$ and $z_{2}$ as $\left[  z_{1},z_{2}\right]
_{\mathcal{K}}$. Similarly denote the segments that exclude one or both
endpoints as $\left[  z_{1},z_{2}\right)  _{\mathcal{K}},$ $\left(
z_{1},z_{2}\right]  _{\mathcal{K}},$ and $\left(  z_{1},z_{2}\right)
_{\mathcal{K}}$. Let $\beta$ be the limiting value of $\arg\left(
z-z_{0}\right)  $ on the principal branch as $z\rightarrow z_{0}$ along
$\left(  z_{0},b\right)  _{\mathcal{K}}$. Finally, let $\phi_{s}$ and
$\chi_{s}$ be the coefficients in the power series
representations%
\begin{equation}
\phi\left(  z\right)  =\sum_{s=0}^{\infty}\phi_{s}\left(  z-z_{0}\right)
^{s}, \qquad \chi(z)=\sum_{s=0}^{\infty}\chi_{s}\left(  z-z_{0}\right)  ^{s}.
\label{f and g series}%
\end{equation}

We assume that there exist positive constants $C_{1},...,C_{4}$ that do not
depend on $p$ or $\omega,$ such that for all $\omega\in\Omega$, for
sufficiently large $p:$

\begin{enumerate}
\item[A0] The length of the path $\mathcal{K}$ is bounded, uniformly over
$\omega\in\Omega$ and all sufficiently large $p.$ Furthermore,%
\[
\sup_{z\in\left(  z_{0},b\right)  _{\mathcal{K}}}\left\vert z-z_{0}\right\vert
>C_{1},\text{ and }\sup_{z\in\left(  a,z_{0}\right)  _{\mathcal{K}}}\left\vert
z-z_{0}\right\vert >C_{1}%
\]

\item[A1] Functions $\phi\left(  z\right)  $ and $\chi(z)$ are holomorphic in
the ball $\left\vert z-z_{0}\right\vert \leq C_{1}$

\item[A2] The coefficient $\phi_{2}$ satisfies $C_{2}\leq\left\vert \phi
_{2}\right\vert \leq C_{3}$

\item[A3] The third derivative of $\phi\left(  z\right)  $ satisfies
inequality
\[
\sup_{\left\vert z-z_{0}\right\vert \leq C_{1}}\left\vert \mathrm{d}^{3}%
\phi\left(  z\right)  /\mathrm{d}z^{3}\right\vert \leq C_{4}%
\]

\item[A4] For any positive $\varepsilon<C_{1},$ which does not depend on $p$
and $\omega$, and for all $z_{1}\in\mathcal{K}$ such that $\left\vert
z_{1}-z_{0}\right\vert =\varepsilon,$ there exist positive constants
$C_{5},C_{6}$, such that%
\[
\operatorname{Re}\left(  \phi\left(  z_{1}\right)  -\phi_{0}\right)
>C_{5}\text{ and }\left\vert \operatorname{Im}\left(  \phi\left(
z_{1}\right)  -\phi_{0}\right)  \right\vert <C_{6}%
\]

\item[A5] For a subset $\Theta$ of $\mathbb{C}$ that consists of all points
whose Euclidean distance from $\mathcal{K}$ is no larger than $C_{1},$%
\[
\sup_{z\in\Theta}\left\vert \chi(z)\right\vert =O_{\mathrm{P}}(1)
\]
as $p\rightarrow\infty,$ where $O_{\mathrm{P}}(1)$ is uniform in $\omega
\in\Omega$.
\end{enumerate}

Assumptions A0--A5 ensure that Olver's proof of the Laplace
approximation theorem (Theorem 7.1 on p. 127 of Olver (1997)) can be
extended to cases where functions $\phi (z)$ and $\chi (z)$, as well as the
contour $\mathcal{K}$, depend on $p$ and $\omega $. Note that in Olver's
notations, $\phi (z),$ $\chi (z),$ and $p$ are, respectively $p(t),$ $q(t),$
and $z$. 

The first part of A0, which requires the boundedness of $\left\vert \mathcal{%
K}\right\vert$, taken together with A5 and the assumption that $\mathcal{K}$
is a contour of steep descent guarantee the absolute convergence of the
integral $\int_{\mathcal{K}}e^{-p\left( \phi (z)-\phi _{0}\right) }\chi (z)%
\mathrm{d}z,$ in probability. The second part of A0 ensures that as $%
p\rightarrow \infty$, $\mathcal{K}$ does not collapse to a point. 

Assumption A1 excludes situations where $z_{0}$ approaches singular points
of $\phi (z)$ or $\chi (z)$ as $p\rightarrow \infty$.
Assumption A2 guarantees that the second derivative of $\phi (z)$ at $z_{0}$
does not degenerate to $0$ or infinity as $p\rightarrow \infty$.
Assumption A3 implies that $\left\vert \phi (z)-\phi (z_{0})\right\vert $
can be bounded from below by a fixed quadratic function of $z$ in a vicinity
of $z_{0}$ as $p\rightarrow \infty $. This ensures a regular behavior of
function $\left( \phi (z)-\phi (z_{0})\right) ^{1/2}$.
Assumption A4 implies that $\left\vert \arg \left( \phi (z)-\phi
(z_{0})\right) \right\vert <\pi /2$ is some neighborhood of $z_{0}$ as $%
p\rightarrow \infty $. We need this condition to be able to use an
asymptotic expansion of an incomplete Gamma function in our proofs 
(Section \ref{sec: extension of Olver} of SM).
Assumption A5 ensures that $\left\vert \chi \left( z\right) \right\vert $
remains bounded in probability as $p\rightarrow \infty $.

% The following lemma, proved in SM 5.1,
%  is a fairly straightforward extension of Theorem 7.1 of
% Olver (1997), p. 127 to the situation where functions $\phi(z),$ $\chi(z)$ and
% the contour $\mathcal{K}$ depend on $p$ and $\omega.$ In Olver's original
% theorem, which uses different notation, both functions and contour are
% fixed. 
% %A proof of the extension can be found in the SM.

\begin{lemma}
Under assumptions A0-A5, for any positive integer $k$, as
$p\rightarrow\infty,$ we have 
\[
I_{p,\omega}=2e^{-p\phi_{0}}\left[  \sum_{s=0}^{k-1}\Gamma\left(
s+\frac{1}{2}\right)  \frac{a_{2s}}{p^{s+1/2}}+\frac{O_{\mathrm{P}}\left(
1\right)  }{p^{k+1/2}}\right]  ,
\]
where $O_{\mathrm{P}}\left(  1\right)  $ is uniform in
$\omega\in\Omega$ and 
the coefficients $a_{2s}$ can be
expressed through $\phi_s$ and $\chi_{s}$ defined above.
In particular we have $a_{0}=\chi_{0}/[2\phi_{2}^{1/2}],$ where
$\phi_{2}^{1/2}=\exp\left\{  \left(  \log\left\vert \phi_{2}\right\vert
+ \mathrm{i} \arg\phi_{2}\right)  /2\right\}  $ with the branch of
$\arg\phi_{2}$ chosen 
so that $\left\vert \arg\phi_{2}+2\beta\right\vert \leq\pi/2.$
\label{Olver}
\end{lemma}

Lemma \ref{Olver} is proved in SM \ref{sec: extension of Olver}. 
We use it to obtain the Laplace approximation to
\begin{equation}
L_{1}\left(  \theta;\Lambda\right)  =\sqrt{\pi p}\frac{1}{2\pi\mathrm{i}}%
\int_{\mathcal{K}_{1}\cup\mathcal{\bar{K}}_{1}}e^{-(p/2)f(z)}g(z)\mathrm{d}z.
\label{I1case}%
\end{equation}
Then we show that $L_{1}\left(  \theta;\Lambda\right)  $ asymptotically
dominates the \textquotedblleft residual\textquotedblright\ $L\left(
\theta;\Lambda\right)  -L_{1}\left(  \theta;\Lambda\right)  $. For this
analysis, it is important to know the values of $f(z_{0})$ and $\mathrm{d}%
^{2}f(z_{0})/\mathrm{d}z^{2}$. As was mentioned
in Remark \ref{fz00}, $f(z_{0})=0$ for
all the six cases that we study. The values of $\mathrm{d}^{2}f(z_{0}%
)/\mathrm{d}z^{2}$ are derived in SM \ref{sec: evaluation of the second derivative}. 
All of them are negative. The explicit form of $D_{2}\equiv
\theta^{2}\left(  -\mathrm{d}^{2}f(z_{0})/\mathrm{d}z^{2}\right)  ^{-1},$
which is somewhat shorter than that for $\mathrm{d}^{2}f(z_{0})/\mathrm{d}%
z^{2}$, is reported in Table \ref{Table 5}. We formulate the main result of
this section in the following theorem, proven in SM \ref{sec: proof of theorem 9}.%

\begin{table}[tbp] 
\centering
\caption{The values of $D_2\equiv\theta^{2}(-\mathrm{d}^{2}f(z_0)/\mathrm{d}z^{2})^{-1}$ for the different cases.}
\vspace{.1in} 
\begin{tabular}
[c]{llll}%
Case & Value of $D_{2}$ & Case & Value of $D_{2}$\\[12pt]
SMD & $1-\theta^{2}$ & REG$_{0}$ & $c_{1}\left(  1+c_{1}%
+2\theta\right)  \left(  c_{1}-\theta^{2}\right)$\\[12pt]
PCA & $c_{1}\left(  c_{1}-\theta^{2}\right)  \left(  1+\theta\right)
^{2}$ & REG & $c_{1}h\left(  c_{1}+\theta+\left(  1+\theta\right)
l\right)  /l^{4}$\\[12pt]
SigD & $r^{2}h\left(  1+\theta\right)  ^{2}/l^{4}$ & CCA & $c_{1}^{2}h\left(
2\left(  c_{1}+\theta\right)  +l\left(  1-c_{1}\right)  \right)  /\left(
l^{3}\left(  c_{1}+c_{2}\right)  \right)  $\\[12pt]
 &$l \equiv l(\theta)=1+(1+\theta)c_{2}/c_{1}$ & & $h \equiv h(\theta)=c_1+c_2(1+\theta)^2-\theta^2$\\
\end{tabular}
\label{Table 5}
\end{table}

\begin{theorem}
Suppose that the null hypothesis holds, that is, $\theta_{0}=0$.
Let $\bar{\theta}$ be the threshold corresponding to $F_{\boldsymbol{\gamma}}$ as
given in Table \ref{Table 3}, and let $\varepsilon$ be an arbitrarily small
fixed positive number. Then for \textit{any }$\theta\in\left(  0,\bar{\theta
}-\varepsilon\right]  ,$\textit{\ as }$\mathbf{n},p\rightarrow_{\boldsymbol{\gamma
}}\infty$,\textit{\ we have}%
\begin{equation}
L\left(  \theta;\Lambda\right)  =\frac{g(z_{0})}{\sqrt{-\mathrm{d}^{2}%
f(z_{0})/\mathrm{d}z^{2}}}+O_{\mathrm{P}}\left(  p^{-1}\right)  ,
\label{Icase equality}%
\end{equation}
\textit{where }$O_{\mathrm{P}}\left(  p^{-1}\right)  $\textit{\ is uniform in
}$\theta\in\left(  0,\bar{\theta}-\varepsilon\right]  $ and the principal
branch of the square root is taken.
\label{Icase}
\end{theorem}

\section{Asymptotics of LR}

Combining the results of Theorem \ref{Icase} with the definitions of $g(z)$
and the values of $-\mathrm{d}^{2}f(z_{0})/\mathrm{d}z^{2}$, given in Table
\ref{Table 5}, it is straightforward to establish the following
theorem,
details in SM \ref{sec: derivations for theorem 10}. Let
%(see SM for a detailed derivation). Let
\[
\Delta_{p}(\theta)=p\int\ln\left(  z_{0}(\theta)-\lambda\right)  \mathrm{d}\left(
\hat{F}\left(  \lambda\right)  -F_{\mathbf{c}}\left(  \lambda\right)  \right)
.
\]
In accordance with Remark \ref{rem:neglig}, we define $\Delta_{p}(\theta)$ as
zero in the event of asymptotically negligible probability that $z_{0}%
\leq\lambda_{1}$.

\begin{theorem}
\label{LR asymptotics}Suppose that the null hypothesis holds, that is
$\theta_{0}=0$. Let $\bar{\theta}$ be the threshold corresponding to
$F_{\boldsymbol{\gamma}}$ as given in Table \ref{Table 3}, and let $\varepsilon$
be an arbitrarily small fixed positive number. Then for \textit{any }%
$\theta\in\left(  0,\bar{\theta}-\varepsilon\right]  ,$\textit{\ as
}$\mathbf{n},p\rightarrow_{\boldsymbol{\gamma}}\infty$,\textit{\ we have}%
\[
L\left(  \theta;\Lambda\right)  =\exp\left\{  -\frac{1}{2}\Delta_{p}%
(\theta)+\frac{1}{2}\ln\left(  1-\left[  \delta_{p}\left(  \theta\right)
\right]  ^{2}\right)  \right\}  \left(  1+o_{\mathrm{P}}(1)\right)  ,
\]
where
\[
\delta_{p}\left(  \theta\right)  =\left\{
\begin{array}
[c]{ll}%
\theta & \text{for SMD}\\
\theta/\sqrt{c_{1}} & \text{for PCA and REG}_{0}\\
\theta r/\left(  c_{1}l\left(  \theta\right)  \right)  & \text{for SigD, REG,
and CCA,}%
\end{array}
\right.
\]
$r^2 = c_1 + c_2 - c_1 c_2$
and $o_{\mathrm{P}}(1)$ is \textit{uniform in }$\theta\in\left(  0,\bar
{\theta}-\varepsilon\right]  $.
\end{theorem}

Statistic $\Delta_{p}(\theta)$ is a linear spectral statistic. As follows from
the CLT for such statistics derived by 
%Bai and Yao (2005), Bai and Silverstein (2004), and Zheng (2012)
\cite{bayo05}, \cite{basi04}, and \cite{zheng12}
for the Semi-circle, Marchenko-Pastur, and Wachter limiting
distributions $F_{\mathbf{c}}\,,$ respectively, statistic $\Delta_{p}(\theta)$
weakly converges to a Gaussian process indexed by $\theta\in\left(
0,\bar{\theta}-\varepsilon\right]  .$ The explicit form of the mean and the
covariance structure can be obtained from the general formulae for the
asymptotic mean and covariance of linear spectral statistics given in 
%Theorem 1.1 of Bai and Yao (2005) 
\cite[Th. 1.1]{bayo05}
for SMD, in 
%Theorem 1.1 of Bai and Silverstein (2004) 
\cite[Th. 1.1]{basi04}
for PCA and REG$_{0}$, and in 
%Theorem 4.1 and Example 4.1 of Zheng (2012) 
\cite[Th. 4.1 and Exmpl. 4.1]{zheng12}
for the remaining cases. 
%For PCA, the corresponding calculations have
%been done in 
%%Onatski et al (2013).
%\cite{omh13}. 
%The convergence of $\Delta_{p}(\theta)$
%and Theorem \ref{LR asymptotics} imply the following theorem, proven
%in SM \ref{sec: proof of theorem 11}.
% Its proof is given in the SM.
SM \ref{sec: proof of theorem 11} provides details on the use of 
\cite{bayo05, basi04, zheng12} to establish convergence of
$\Delta_{p}(\theta)$, and the use of Theorem \ref{LR asymptotics} to
obtain the following theorem.

\begin{theorem}
Suppose that the null hypothesis holds, that is
$\theta_{0}=0$. Let $\bar{\theta}$ be the threshold corresponding to
$F_{\boldsymbol{\gamma}}$ as given in Table \ref{Table 3}, and let $\varepsilon$
be an arbitrarily small fixed positive number. Further, let $C\left[
0,\bar{\theta}-\varepsilon\right]  $ be the space of continuous functions on
$\left[  0,\bar{\theta}-\varepsilon\right]  $ equipped with the supremum norm.
Then $\ln L\left(  \theta;\Lambda\right)  $ viewed as random elements of
$C\left[  0,\bar{\theta}-\varepsilon\right]  $ converge weakly to
$\mathcal{L}\left(  \theta\right)  $ with Gaussian finite dimensional
distributions such that%

\[
\mathbb{E}\mathcal{L}\left(  \theta\right)  =\tfrac{1}{4}\ln ( 1-
\delta^2(\theta))
% \mathbb{E}\mathcal{L}\left(  \theta\right)  =\frac{1}{4}\ln\left(  1-\left[
% \delta\left(  \theta\right)  \right]  ^{2}\right)
\]
and%
\[
\mathbb{C}ov\left(  \mathcal{L}\left(  \theta_{1}\right)  ,\mathcal{L}\left(
\theta_{2}\right)  \right)  =-\tfrac{1}{2}\ln\left(  1-\delta\left(  \theta
_{1}\right)  \delta\left(  \theta_{2}\right)  \right)
\]
with%
\[
\delta\left(  \theta\right)  =\left\{
\begin{array}
[c]{ll}%
\theta & \text{for SMD}\\
\theta/\sqrt{\gamma_{1}} & \text{for PCA and REG}_{0}\\
\theta\rho/\left(  \gamma_{1}+\gamma_{2}+\theta\gamma_{2}\right)  & \text{for
SigD, REG, and CCA}%
\end{array}
\right.  .
\]
Here $\rho,\gamma_{1},\gamma_{2}$ are the limits of $r,c_{1},c_{2}$ as
$\mathbf{n},p\rightarrow_{\boldsymbol{\gamma}}\infty$.
\label{AsymptoticNormality}
\end{theorem}

Let $\left\{  \mathbb{P}_{p,\theta}\right\}  $ and $\left\{  \mathbb{P}%
_{p,0}\right\}  $ be the sequences of measures corresponding to the joint
distributions of $\lambda_{1},...,\lambda_{p}$ when $\theta_{0}=\theta$ and
when $\theta_{0}=0$ respectively. Then Theorem \ref{AsymptoticNormality}
implies, via Le Cam's first lemma, the mutual contiguity of $\left\{
\mathbb{P}_{p,\theta}\right\}  $ and $\left\{  \mathbb{P}_{p,0}\right\}  $ as
$\mathbf{n},p\rightarrow_{\boldsymbol{\gamma}}\infty$,
for each $\theta < \bar \theta$.
This reveals the
statistical meaning of the phase transition thresholds as the upper boundaries
of the contiguity regions for spiked models.

The precise form of the autocovariance of $\mathcal{L}\left(  \theta\right)  $
shows that,\footnote{
%Fyodorov et al (2013) have 
\cite{fykhsi13} has
an interesting discussion of
ubiquity of random processes with logarithmic covariance structure in physics
and engineering applications. 
%In their paper, 
In that paper,
such processes appear as
limiting objects related to the behavior of the characteristic polynomials of
large matrices from Gaussian Unitary Ensemble.} although the experiment of
observing $\lambda_{1},...,\lambda_{p}$ is asymptotically normal, it does not
converge to a Gaussian shift experiment. In particular, the optimality results
available for Gaussian shifts cannot be used in our framework. To analyze
asymptotic risks of various statistical problems related to the experiment of
observing $\lambda_{1},...,\lambda_{p},$ one should directly use Theorem
\ref{AsymptoticNormality}.

Here we use it to derive the
asymptotic power envelopes for tests of the null hypothesis $\theta_{0}=0$
against the point alternative $\theta_{0}>0.$ 
%[Such a power envelope has been derived 
%by Onatski et al (2013) for the case of PCA.] 
By the Neyman-Pearson lemma, the most powerful test 
%of the null against a point alternative $\theta_{0}=\theta$
would reject the null when $\ln L\left(  \theta;\Lambda\right)  $ is above a
critical value. By Theorem \ref{AsymptoticNormality} and Le Cam's third lemma
%(see van der Vaart (1998), chapter 6),%
(see \cite[Ch. 6]{vdva98}),
\[
\ln L\left(  \theta;\Lambda\right)  \overset{d}{\rightarrow}N\left(  
\pm \tfrac{1}{4} \ln(  1-  \delta^2 (\theta) ),
- \tfrac{1}{2} \ln(  1-  \delta^2 (\theta) ) \right)
% ,-\tfrac{1}{2}\ln\left(  1-\left[  \delta\left(  \theta\right)  \right]
% ^{2}\right)  \right)
\]
with the plus sign holding under the null, and
the minus under the alternative.
% \[
% \ln L\left(  \theta;\Lambda\right)  \overset{d}{\rightarrow}N\left(  -\frac
% {1}{4}\ln\left(  1-\left[  \delta\left(  \theta\right)  \right]  ^{2}\right)
% ,-\frac{1}{2}\ln\left(  1-\left[  \delta\left(  \theta\right)  \right]
% ^{2}\right)  \right)
% \]
% under the alternative. 
This implies the following theorem.

\begin{theorem}
\label{th:PE}
Let $\bar{\theta}$ be the threshold corresponding to $F_{\boldsymbol{\gamma}}$ as
given in Table \ref{Table 3}. For any $\theta\in\left[  0,\bar{\theta}\right)
,$ the value of the asymptotic power envelope for the tests of the null
$\theta_{0}=0$ against the alternative $\theta_{0}>0$ which are based on
$\lambda_{1},...,\lambda_{p}$ and have asymptotic size $\alpha$ is given by%
\[
PE\left(  \theta\right)  =1-\Phi\left[  \Phi^{-1}\left(  1-\alpha\right)
- \sigma(\theta) \right], \qquad 
\sigma(\theta) = 
\sqrt{-\tfrac{1}{2}\ln(  1-  \delta^2(\theta) )  }  .
\]
% \[
% PE\left(  \theta\right)  =1-\Phi\left[  \Phi^{-1}\left(  1-\alpha\right)
% -\sqrt{-\tfrac{1}{2}\ln(  1-  \delta^2(\theta) )  }\right]  .
% \]
Here $\Phi$ denotes the standard normal cumulative distribution function. For
$\theta\geq\bar{\theta}$ the value of the asymptotic power envelope equals one.
\end{theorem}

\begin{figure}[h]
\centering
\includegraphics[height=3.686in,width=3.58in]{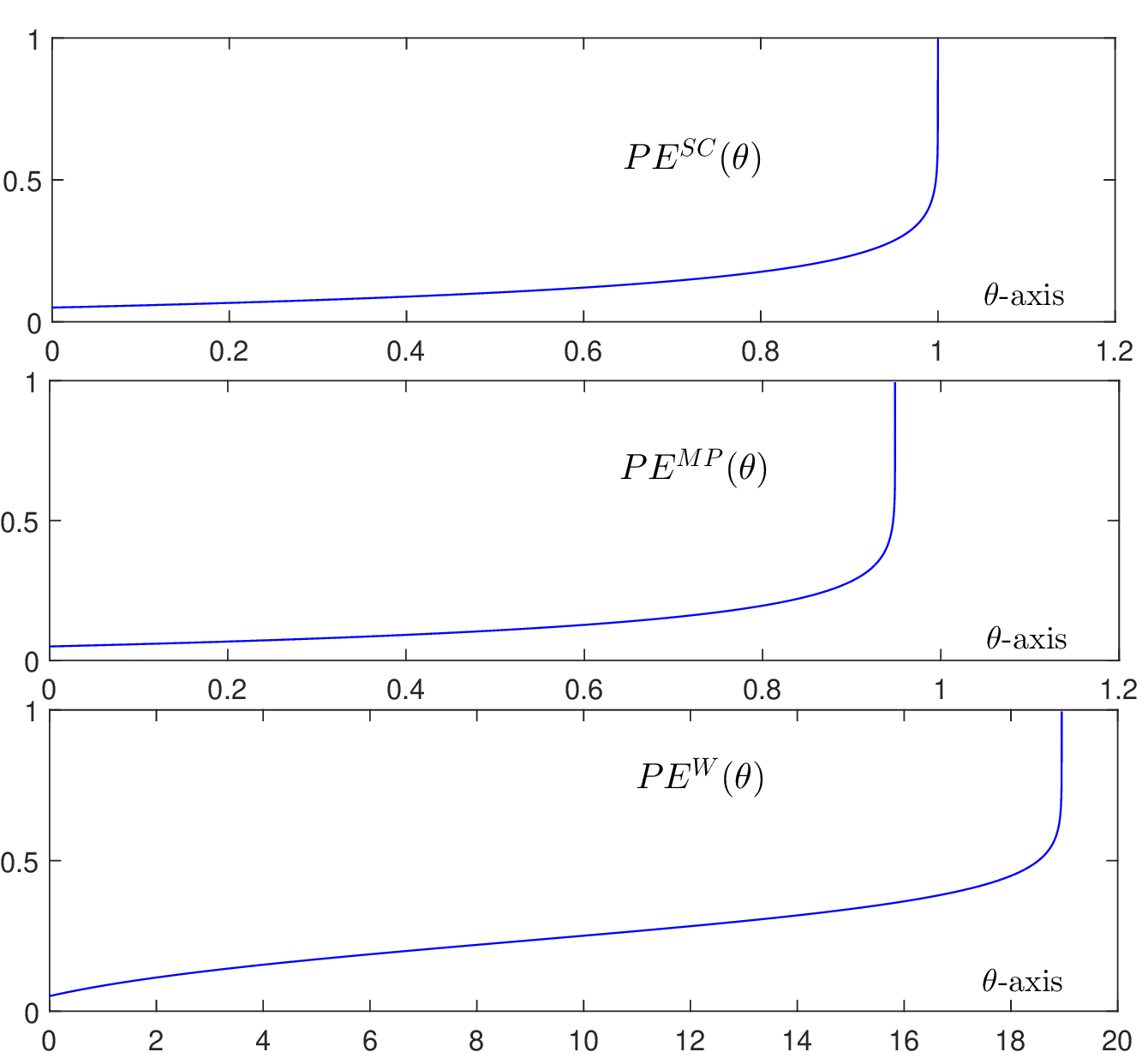}
\caption{The asymptotic power envelopes $PE^{\rm{SC}}(\theta),$ $PE^{\rm{MP}}(\theta,\gamma_{1}),$
and $PE^{\rm{W}}(\theta,\boldsymbol{\gamma})$ for $\alpha=0.05,$ $\gamma_{1}=\gamma_{2}=0.9.$}%
\label{envelopes}
\end{figure}

The envelopes differ only for cases with different
limiting spectral distributions:  Semi-circle,  Marchenko-Pastur, and
Wachter, denoted
$PE^{\rm SC}(\theta)$, $PE^{\rm MP}(\theta,\gamma_1)$ and
$PE^{\rm W}(\theta,\boldsymbol{\gamma})$ respectively.
% Therefore, we can denote $PE\left(  \theta\right)  $
% as $PE^{\rm{SC}}\left(  \theta\right)  $ for SMD, as $PE^{\rm{MP}}\left(  \theta\right)
% $ for PCA and REG$_{0}$, and as $PE^{\rm{W}}\left(  \theta\right)  $ for the
% remaining cases. 
Figure \ref{envelopes} shows the graphs of the envelopes for
$\alpha=0.05$ and $\gamma_{1}=\gamma_{2}=0.9.$ Such large values of
$\gamma_{1}$ and $\gamma_{2}$ correspond to situations where the
dimensionality $p$ is not very different from the 
degrees of freedom 
% \textquotedblleft sample sizes\textquotedblright\ 
$n_{1}$ and $n_{2}.$ 
% Of course, the values of
% $\gamma_{1}$ and $\gamma_{2}$ are irrelevant for $PE^{\rm{SC}}\left(
% \theta\right)  $, and the value of $\gamma_{2}$ is irrelevant for
% $PE^{\rm{MP}}\left(  \theta\right)  .$%

Envelope $PE^{\rm{MP}}\left(  \theta, \gamma_1
\right)  $ can 
be obtained from $PE^{\rm{W}}\left(  \theta, \boldsymbol{\gamma} \right)  $ by
sending $\gamma_{2}$ to 
zero. Further, $PE^{\rm{SC}}\left(  \theta\right)  $ can be obtained from
$PE^{\rm{MP}}\left(  \theta, \gamma_1 \right)  $ by transformation
$\theta\longmapsto 
\sqrt{\gamma_{1}}\theta.$ Further, note the difference in the horizontal scale
of the bottom panel of Figure \ref{envelopes} relative to the two other
panels. For $\gamma_{1}=\gamma_{2}=0.9$ the phase transition threshold
corresponding to the Wachter distribution is relatively large. It equals $\left(
\gamma_{2}+\rho\right)  /(1-\gamma_{2})\approx18.9$. Moreover, the value of
$PE^{\rm{W}}\left(  \theta\right)  $ becomes substantially larger than the nominal
size $\alpha=0.05$ for $\theta$ that are situated far below this threshold.
This suggests that the information in all the eigenvalues $\lambda
_{1},...,\lambda_{p}$ might be effectively used to detect spikes that are
small relative to the phase transition threshold in two sample problems. We
leave a confirmation or rejection of this speculation for future research.

\section{Concluding remarks} \label{sec-Concluding remarks}

Note that Theorem \ref{AsymptoticNormality} establishes the weak convergence 
of the log likelihood
ratio viewed as a random element of the space of continuous functions. This is
much stronger than simply the convergence of the finite dimensional
distributions of the log likelihood process. In particular, the theorem can be
used to find the asymptotic distribution of the supremum of the likelihood
ratio, and thus, to find the asymptotic critical values of the likelihood
ratio test. It also can be used to construct asymptotic confidence 
intervals for a sub-critical spike as well as to describe the asymptotic properties of 
its maximum likelihood estimator. We do not pursue this line of research here, but
provide a general outline.

Consider the log likelihood ratio $\ln L\left( \theta ;\Lambda
\right) -\ln L\left( \theta _{0};\Lambda \right) .$ According to Theorem 12,
this ratio converges to $X\left( \theta \right) \equiv \mathcal{L}\left(
\theta \right) -\mathcal{L}\left( \theta _{0}\right) .$ By Le Cam's third
lemma, under the null hypothesis that the true value of the spike equals $%
\theta _{0},$ $X\left( \theta \right) $ is a Gaussian process with mean%
\[
\mathbb{E}X\left( \theta \right) =\frac{1}{4}\ln \frac{\left( 1-\delta
^{2}\left( \theta \right) \right) \left( 1-\delta ^{2}\left( \theta
_{0}\right) \right) }{\left( 1-\delta \left( \theta \right) \delta \left(
\theta _{0}\right) \right) ^{2}}
\]%
and covariance function%
\[
\mathbb{C}ov\left( X\left( \theta _{1}\right) ,X\left( \theta _{2}\right)
\right) =-\frac{1}{2}\ln \frac{\left( 1-\delta \left( \theta _{1}\right)
\delta \left( \theta _{2}\right) \right) \left( 1-\delta \left( \theta
_{0}\right) \delta \left( \theta _{0}\right) \right) }{\left( 1-\delta
\left( \theta _{1}\right) \delta \left( \theta _{0}\right) \right) \left(
1-\delta \left( \theta _{2}\right) \delta \left( \theta _{0}\right) \right) }%
.
\]

An approximation to the distribution of the supremum of such a process over $%
\theta \in \left[ 0,\bar{\theta}-\varepsilon \right] $ can be obtained via
simulation. Alternatively, it might be expressed analytically in the form of
converging Rice series (see e.g. \cite{aw02}). Quantiles of the distribution 
can be used as asymptotic critical values for the likelihood ratio
test of the hypothesis $\theta=\theta_{0}$. Inverting the test, we obtain 
asymptotic confidence intervals for the true value of a sub-critical spike.

The maximum likelihood estimator for the spike, $\hat{\theta}_{ML}$, 
equals the $\arg \max $ of $\ln L\left( \theta ;\Lambda \right)-
\ln L\left( \theta_{0} ;\Lambda \right)$ over $%
\theta \in \left[ 0,\bar{\theta}-\varepsilon \right] $. By Lemma 2.6 of \cite{kp90},
the limiting process $X\left( \theta \right) $
achieves maximum at a unique point with probability one. Therefore by the argmax 
continuous mapping theorem, $\hat{\theta}_{ML}$ converges in distribution to 
the $\arg \max $ of 
$X\left( \theta \right) $. The distribution of such an $\arg \max $ can be
approximated using simulations.

Unfortunately, the quality of the estimator $\hat{\theta}%
_{ML}$ cannot be \textquotedblleft good\textquotedblright. For PCA, we 
were able to prove that no estimator of $\theta $ has
root mean squared error better than the order of magnitude of the
sub-critical parameter $\theta$. This result will appear in another work.

Our asymptotic discussion of James' framework can likely be extended
to a fixed number of sub-critical spikes. Such an extension would 
require developing Laplace approximations to multiple contour 
integrals, and uniform approximations to hypergeometric functions of two 
matrix arguments in terms of elementary functions.
Alternatively, one may employ large deviation analysis of spherical integrals
as in \cite{omh14}, which covers the PCA case. As this paper is already
long, the extension will appear separately.

Addressing the case of slowly increasing number of spikes 
may require new techniques, perhaps, similar to those developed in \cite{d17}. 
In such a case, relatively little is known even 
about the phase transition phenomenon. For sample covariance matrices, Theorem 1.1
of \cite{basi98} can be used to show that the phase transition still happens
at the usual threshold $\bar{\theta}=\sqrt{\gamma_{1}}$.
However, it is not clear whether the experiments of observing sample
covariance eigenvalues corresponding to the null case and an alternative
with a growing number of sub-critical eigenvalues remain mutually contiguous.

Note that, intuitively, the asymptotic power of the likelihood ratio test of 
the null hypothesis of no spikes against the alternative of one spike should not
decrease if the rank-one assumption on the alternative is wrong and there are 
additional spikes. In
SM \ref{power-lr-test},
we confirm this intuition 
for SMD and PCA cases. A confirmation or refutation of the intuition for the
other James' cases requires further analysis and is left for future research. 

% This paper derives the asymptotics of the likelihood ratio processes
% corresponding to the null hypothesis of no spikes and the alternative of a
% single spike in various high-dimensional multivariate models. We cover all the
% five classes of multivariate statistical problems identified by James (1964).
% In addition, we consider a symmetric matrix denoising problem that does not
% fit in James' classification. We find that, as the dimensionality and the
% number of observations go to infinity proportionally, the log likelihood
% processes converge to Gaussian limits as long as the value of the spike
% parameter is below corresponding phase transition thresholds. We derive
% explicit formulae for the autocovariance and the mean of the limiting
% processes and use them to obtain asymptotic power envelopes for tests for the
% presence of a spike.

In this paper, we make the assumption that $n_{2}\geq p$ to ensure the
invertibility of matrix $E$ in (\ref{basicJO}) with probability one. However, we
also make the assumption $n_{1}\geq p,$
mostly to simplify our exposition.
It can probably be lifted without a
substantial reformulation of the problem.
%We make the latter assumption mostly to simplify our exposition. 
Indeed,
for SMD the assumption is irrelevant. For PCA the case
$p>n_{1}$ was explicitly covered in \cite{omh13}.
%. The PCA
%results are obtained in 
%%Onatski et al (2013) 
%\cite{omh13}
%without using this assumption.
For REG$_{0}$ the assumption can be relaxed using the symmetry of the problem.
Specifically, the canonical REG$_{0}$ problem tests restriction $M=0$ in the
model $Y=M+\varepsilon$, where each matrix is
$n_{1}\times p$ and $% 
\varepsilon $ has i.i.d. standard normal components. Clearly, interchanging
roles of $n_{1}$ and $p$ yields essentially the same problem. 

For CCA, the
sample canonical correlations are not well defined for $p>n_{1}.$
%so we are not interested in such a case.
% This leaves us with SigD and REG cases, which
% we mark as more difficult from the point of view of relaxing $n_{1}\geq p$
% assumption.
For SigD, our derivations (not reported here) show that the equivalent of
(\ref{LRgeneral}) for $p>n_{1}$ involves the hypergeometric function
$_{2}F_{1}.$ Therefore, SigD with $p>n_{1}$ represents the fifth, rather than
the second, group of multivariate statistical problems according to James'
(1964) classification. For REG, an equivalent of
(\ref{LRgeneral}) for $p>n_{1}$ can be obtained using 
%equation (74) of James (1964). 
\cite[eq. (74)]{jame64}.
However, further analysis of SigD and REG in the situation where $p>n_{1}$
needs more
substantial changes to our derivations. We leave such an analysis for future research.  

%Our asymptotic discussion of James' framework can likely be extended
%to a fixed number of subcritical spikes. Indeed \cite{omh14} discussed the
%PCA case, albeit via a different technique. As this paper is already
%long, the extension will appear separately.

Finally, many existing results in the random matrix literature do not require that the
data are Gaussian. This suggests that some results about tests for the
presence of the spikes in the data may remain valid without the
Gaussian assumption.
We hope that the results of this paper
might provide 
a benchmark for such future studies.

\addtocontents{toc}{\protect\setcounter{tocdepth}{3}}
\newpage
\setcounter{section}{0}
\begin{frontmatter}
\title{Supplementary Material for \textquotedblleft Testing in
high-dimensional spiked models.\textquotedblright}
\runtitle{Supplementary material}
\author{\fnms{Iain M.} \snm{Johnstone}\ead[label=e1]{imj@stanford.edu}}
\address{\printead{e1}}
\and
\author{\fnms{Alexei} \snm{Onatski}\ead[label=e2]{ao319@cam.ac.uk}}
\address{\printead{e2}}
%\affiliation{???}

\runauthor{I. M. Johnstone and A. Onatski}

\begin{abstract}

This note contains supplementary material for Johnstone and Onatski (2018)
(JO in what follows). It is lined up with sections in the main text to make 
it relatively easy to see how and where the proof details fit in.

\end{abstract}

\end{frontmatter}
\tableofcontents
%\section{Introduction}

% AOS,AOAS: If there are supplements please fill:
%\begin{supplement}[id=suppA]
%  \sname{Supplement A}
%  \stitle{Title}
%  \slink[doi]{10.1214/00-AOASXXXXSUPP}
%  \sdatatype{.pdf}" 
%  \sdescription{Some text}
%\end{supplement}

%\setcounter{equation}{1}

\section{Introduction}
%\addtocontents{toc}{\bigskip text\par}
There is no supplementary material for the Introduction section of JO.

\section{Links to classical statistical problems}

\subsection{Sufficiency and invariance considerations} \label{Sufficiency and invariance considerations}

In this subsection, we clarify which sufficiency and invariance arguments
lead us to consider tests based on the solutions of 
\begin{equation}
\det \left( H-\lambda E\right) =0  \label{basic}\tag{SM1}
\end{equation}%
and%
\begin{equation}
\det \left( H-\lambda \left( E+\frac{n_{1}}{n_{2}}H\right) \right) =0
\label{basic beta form}\tag{SM2}
\end{equation}%
for SMD, PCA, SigD, RED$_{0}$, REG, and CCA problems. Most of this discussion
is standard and can be found, for example, in Muirhead
(1982). \vspace{3 mm}

\textbf{\small{SMD:}}
Consider the group of transformations 
\begin{equation}
G=\left\{ U:U\in \mathcal{O(}p\mathcal{)}\right\} ,  \label{group}\tag{SM3}
\end{equation}
where $\mathcal{O(}p\mathcal{)}$ is the group of $p\times p$ orthogonal
matrices, acting on the space of $p\times p$ symmetric matrices $X=\theta\mathbb{%
\psi } \mathbb{\psi }^{\prime }+GOE/\sqrt{p}$ by%
\begin{equation*}
U\circ X=UXU^{\prime }.
\end{equation*}%
The corresponding induced group of transformations on the parameter space of
points $\left( \mathbb{\psi },\theta \right) $ is given by%
\begin{equation*}
U\circ \left( \mathbb{\psi },\theta \right) =\left( U\mathbb{\psi },\theta
\right) .
\end{equation*}%
A maximal invariant in the parameter space is $\theta ,$ whereas that in the
sample space is given by the ordered eigenvalues $\lambda _{1}\geq ...\geq
\lambda _{p}$ of $X.$ Since neither the null nor the alternative hypothesis,%
\begin{equation}
H_{0}:\theta _{0}=0\text{ and }H_{1}:\theta _{0}>0,  \label{Basic hypotheses}\tag{SM4}
\end{equation}%
is affected by the transformations, it is natural to base the test on the
maximal invariant in the sample space. \vspace{3 mm}

\textbf{\small{PCA:}}
In this case, the data are given by $X\sim N\left( 0,\Omega \otimes
I_{n_{1}}\right) ,$ where $\Omega =\Sigma +\theta \mathbb{\psi } \mathbb{\psi 
}^{\prime },$ where $\Sigma $ is a known positive definite symmetric matrix
and $\left\Vert \Sigma ^{-1/2}\mathbb{\psi }\right\Vert =1$. Without loss of
generality, we can set $\Sigma =I_{p}$. A sufficient statistic is $%
H=XX^{\prime }/n_{1}.$ Consider the group of transformations (\ref{group})
that acts on the sample space of the sufficient statistic by%
\begin{equation*}
U\circ H=UHU^{\prime },
\end{equation*}%
and on the parameter space by%
\begin{equation*}
U\circ \left( \mathbb{\psi },\theta \right) =\left( U\mathbb{\psi },\theta
\right) .
\end{equation*}%
The maximal invariant in the parameter space is $\theta ,$ and we base the
test of (\ref{Basic hypotheses}) on a maximal invariant in the sample space
of the sufficient statistic, which is given by the ordered eigenvalues of
the sample covariance matrix $H.$ \vspace{3 mm}

\textbf{\small{SigD:}}
The data are given by independent matrices 
\begin{equation*}
X\sim N\left( 0,\Omega \otimes I_{n_{1}}\right) \text{ and }Y\sim N\left(
0,\Sigma \otimes I_{n_{2}}\right) ,
\end{equation*}
where $\Omega =\Sigma +\theta \mathbb{\psi } \mathbb{\psi }^{\prime },$ $%
\Sigma $ is an unknown positive definite symmetric matrix, and $\left\Vert
\Sigma ^{-1/2}\mathbb{\psi }\right\Vert =1$. A sufficient statistic consists
of the sample covariance matrices $H=XX^{\prime }/n_{1}$ and $E=YY^{\prime
}/n_{2}.$ Let $\mathcal{GL}\left( p\right) $ be the group of non-singular $%
p\times p$ matrices. Consider the group of transformations $G=\left\{ B:B\in 
\mathcal{GL}\left( p\right) \right\} $ that acts on the space of points $%
\left( H,E\right) \in \mathcal{S}_{p}\times \mathcal{S}_{p}$, where $%
\mathcal{S}_{p}$ is the space of positive definite symmetric $p\times p$
matrices, by%
\begin{equation*}
B\circ \left( H,E\right) =\left( BHB^{\prime },BEB^{\prime }\right)
\end{equation*}%
and on the parameter space by%
\begin{equation*}
B\circ \left( \Sigma ,\mathbb{\psi },\theta \right) =\left( B\Sigma
B^{\prime },B\mathbb{\psi },\theta \right) .
\end{equation*}%
Note that we restrict the sample space to $\mathcal{S}_{p}\times \mathcal{S}%
_{p},$ that is we exclude from consideration zero-probability event where
the matrix $HE$ is singular. The maximal invariant in the parameter space is 
$\theta $ and we base the test of (%
\ref{Basic hypotheses}) on a maximal invariant in the sample space of the
sufficient statistic, which is given by the ordered solutions to (\ref{basic}%
) or to (\ref{basic beta form}) (see Theorem 8.2.2 of Muirhead (1982)). 
The links between SigD and PCA become
particularly clear when we work with the solutions to (\ref{basic beta form}%
). 

Note that we can assume that $\Sigma=I_{p}$ wlog. 
It is because $\lambda_{1} \geq ... \geq \lambda_{p}$ that solve equation
(\ref{basic beta form}) are invariant with respect to the transformation
\begin{equation*}
(H,E) \mapsto (\Sigma^{-1/2}H\Sigma^{-1/2},\Sigma^{-1/2}E\Sigma^{-1/2}).
\end{equation*}
In particular, the joint distribution of $\lambda_{1} \geq ... \geq \lambda_{p}$
under the null hypothesis $H_{0}:\Omega=\Sigma$ is the same as in the case where
$\Omega=\Sigma=I_{p}$. Similarly, the joint distribution of
$\lambda_{1} \geq ... \geq \lambda_{p}$ under the alternative 
$H_{1}:\Omega=\Sigma+\theta\psi\psi^{\prime}$ with 
$\left\Vert {\Sigma^{-1/2}\psi}\right\Vert=1 $ is the same as in the situation where
$\Omega=I_{p}+\theta\psi\psi^{\prime}$ with $\left\Vert {\psi}\right\Vert=1$
and $\Sigma=I_{p}$.
\vspace{3 mm}

\textbf{\small{REG$_{0}$:}}
Consider linear regression $Y=X\beta +\varepsilon $,
where $Y$ is $T\times p$, $X$ is $T\times q$, $\beta $ is $q\times p$, and $%
\varepsilon $ has i.i.d. $N(0,\Sigma )$ rows. 
For REG$_{0}$,
$\Sigma$ is a know symmetric positive definite matrix, which can
be set to $I_{p}$ wlog. We would like to
test a general linear hypothesis $C\beta =0$, where $C$ is a known $%
n_{1}\times q$ matrix of rank $n_{1}$. 

As explained in Muirhead (1982, pp
433-434), the problem can be cast in the canonical form, where the matrix of
transformed response variables is split into three parts: an $n_{1}\times p$
matrix $Y_{1}$, a $\left(  q-n_{1}\right)  \times p$ matrix $Y_{2}$,
and an $n_{2} \times p$ matrix $Y_{3}$ with $n_{2}=T-q$.
Under the null hypothesis, $M \equiv \mathbb{E}Y_{1}=0,$ whereas
under the alternative,
\begin{equation}
M=\sqrt{n_{1}\theta}\varphi\psi^{\prime},
\label{REG0mean1}\tag{SM5}%
\end{equation}
where $\theta>0,$ $\| \Sigma^{-1/2}\psi\| =1,$ and
$\left\Vert \varphi\right\Vert =1$.  Matrices $Y_{2}$ and $Y_{3}$
have, respectively, unrestricted and zero means
under both the null and the alternative.

More explicitly, if $X = QR$ is the thin $QR$
decomposition of $X$, so that $Q \in \mathcal{O}(T)$ and
$R$ is $q \times q$ (upper triangular and) invertible, and
with a corresponding decomposition 
$C R^{-1} = \check{R} \check{Q}$,
then Muirhead shows that $M = \check{R}^{-1}C\beta$.
% In terms of the original regression model, matrix $M$ 
% can be expressed as the product of an invertible
% matrix, which depends only on $C$ and $X$, and matrix $C\beta$.
In particular, $M=0$ if and only if $C\beta =0$.
Alternative (\ref{REG0mean1}) corresponds to a rank-one
alternative $C\beta =\sqrt{n_{1}\theta }\tilde{\varphi}\psi ^{\prime }$ in
the original model,
where vector $\tilde{\varphi} = \check{R} \varphi$.
% is obtained from vector $%
% \varphi $ via a linear transformation that depends on matrices $C$ and $X$.

%The canonical form of the problem corresponds to the situation where the
%data are represented by independent matrices $Y_{1},$ $Y_{2},$ and $Y_{3}.$ We have 
%\begin{equation*}
%Y_{1}\sim N\left( \sqrt{n_{1}\theta }\varphi \psi ^{\prime
%},I_{n_{1}}\otimes \Sigma \right) ,
%\end{equation*}
%where $\theta >0,$ $\left\Vert \Sigma ^{-1/2}\psi \right\Vert =1,$ and $%
%\left\Vert \varphi \right\Vert =1;$ $Y_{2}\sim N\left(
%M,I_{T-n_{1}-n_{2}}\otimes \Sigma \right) ,$ where $M$ is an unrestricted $%
%\left( T-n_{1}-n_{2}\right) \times p$ matrix; and $Y_{3}\sim N\left(
%0,I_{n_{2}}\otimes \Sigma \right) .$ 
A sufficient statistic for $\sqrt{n_{1}\theta }\varphi \psi
^{\prime }$ is $Y_{1}$. Consider a group of transformations%
\begin{equation*}
G=\left\{ \left( U,V\right) :U\in \mathcal{O}\left( p\right) ,V\in \mathcal{O%
}\left( n_{1}\right) \right\}
\end{equation*}%
that acts on the points $Y_{1}$ of the sample space $\mathbb{R}%
^{n_{1}\times p}$ by%
\begin{equation*}
\left( U,V\right) \circ Y_{1}=VY_{1}U^{\prime }
\end{equation*}%
and on the parameter space by%
\begin{equation*}
\left( U,V\right) \circ \left( \varphi ,\psi ,\theta \right) =\left(
V\varphi ,U\psi ,\theta \right) .
\end{equation*}%
A maximal invariant in the parameter space is $\theta ,$ whereas the maximal
invariant statistic consists of the ordered eigenvalues of $H=Y_{1}
Y_{1}^{\prime }/n_{1}.$ \vspace{3 mm}

\textbf{\small{REG:}}
The difference between the cases REG and REG$_{0}$ is that in REG $\Sigma $
is assumed to be an unknown matrix from $\mathcal{S}_{p}$. The sufficient
statistic now is $\left( Y_{1},Y_{2},E\right) ,$ where $%
E=Y_{3}^{\prime }Y_{3}/n_{2}$. Consider a group of
transformations%
\begin{equation*}
G=\left\{ \left( B,V,A\right) :B\in \mathcal{GL}\left( p\right) ,V\in 
\mathcal{O}\left( n_{1}\right) ,A\in \mathbb{R}^{\left( T-n_{1}-n_{2}\right)
\times p}\right\}
\end{equation*}%
that acts on the points $\left( Y_{1},Y_{2},E\right) $ of
the sample space $\mathbb{R}^{n_{1}\times p}\times \mathbb{R}%
^{T-n_{1}-n_{2}\times p}\times \mathcal{S}_{p}$ by%
\begin{equation*}
\left( B,V,A\right) \circ \left( Y_{1},Y_{2},E\right)
=\left( VY_{1}B^{\prime },Y_{2}B^{\prime }+A,BEB^{\prime
}\right)
\end{equation*}%
and on the parameter space by 
\begin{equation*}
\left( B,V,A\right) \circ \left( \varphi ,\psi ,\theta ,M,\Sigma \right)
=\left( V\varphi ,B\psi ,\theta ,MB^{\prime }+A,B\Sigma B^{\prime }\right) .
\end{equation*}
A maximal invariant in the parameter space is $\theta ,$ whereas
the maximal invariant in the sample space consists of the ordered roots of
equation (\ref{basic beta form}), where $H=Y_{1}{\prime }Y_{1}/n_{1}$ 
and $E=Y_{3}^{\prime }Y_{3}/n_{2}$ (see Theorem 10.2.1
on page 437 of Muirhead (1982)). \vspace{3 mm}

\textbf{\small{CCA:}}
For this case, the sufficient statistic is $S=\left( 
\begin{array}{cc}
S_{xx} & S_{xy} \\ 
S_{yx} & S_{yy}%
\end{array}%
\right) $. Consider a group of transformations%
\begin{equation*}
G=\left\{ B:B=diag\left\{ B_{1},B_{2}\right\} ,B_{1}\in \mathcal{GL}\left(
p\right) ,\text{ }B_{2}\in \mathcal{GL}\left( n_{1}\right) \right\}
\end{equation*}%
acting on the sample space, restricted so that $S_{xx}$ and $S_{yy}$ are
invertible, by 
\begin{equation*}
B\circ S=BSB^{\prime }.
\end{equation*}
On the parameter space, the group acts by%
\begin{equation*}
B\circ \left( \Sigma _{xx},\Sigma _{yy},\psi ,\varphi ,\theta \right)
=\left( B_{1}\Sigma _{xx}B_{1}^{\prime },B_{2}\Sigma _{yy}B_{2}^{\prime
},B_{1}\psi ,B_{2}\varphi ,\theta \right) .
\end{equation*}%
As follows from Muirhead's (1982) Theorem 11.2.2, a maximal invariant in the
parameter space is $\theta $ and that in the sample space consists of the
solutions to (\ref{basic}) with 
\begin{equation*}
H=S_{xy}S_{yy}^{-1}S_{yx}\text{ and }E=S_{xx}.
\end{equation*}

\subsection{Sequential asymptotic links between the cases} \label{Sequential asymptotic links between the cases}

\textbf{\small{PCA$\rightarrow $SMD:}}
Recall that the relevant data for PCA case are represented by the solutions
to equation (\ref{basic}) with $E=\Sigma $ and $n_{1}H\sim W_{p}\left(
n_{1},\Omega \right) $. Let 
\begin{equation}
\Omega =\Sigma +\sqrt{p/n_{1}}\theta \psi  \psi ^{\prime }
\label{Omega modified}\tag{SM6}
\end{equation}%
with $\left\Vert \Sigma ^{-1/2}\psi \right\Vert =1.$ That is, let the value
of the spike in the original version of PCA be scaled by $\sqrt{p/n_{1}}$.
Equation (\ref{Omega modified}) implies that%
\begin{equation*}
\Sigma ^{-1}=\Omega ^{-1}+\frac{\sqrt{p/n_{1}}\theta }{1+\sqrt{p/n_{1}}%
\theta }\Sigma ^{-1}\psi \psi ^{\prime }\Sigma ^{-1},
\end{equation*}%
and therefore, equation (\ref{basic}) is equivalent to%
\begin{equation*}
\det \left( \left( \Omega ^{-1}+\frac{\sqrt{p/n_{1}}\theta }{1+\sqrt{p/n_{1}}%
\theta }\Sigma ^{-1}\psi \psi ^{\prime }\Sigma ^{-1}\right) H-\lambda
I_{p}\right) =0,
\end{equation*}%
which, in its turn, is equivalent to%
\begin{equation}
\det \left( \Omega ^{-1/2}H\Omega ^{-1/2}+\sqrt{p/n_{1}}\theta \eta \eta
^{\prime }\Omega ^{-1/2}H\Omega ^{-1/2}-\lambda I_{p}\right) =0,
\label{modified equation 1}\tag{SM7}
\end{equation}%
where 
\begin{equation*}
\eta =\Omega ^{1/2}\Sigma ^{-1}\psi /\left( 1+\sqrt{p/n_{1}}\theta \right)
^{1/2}
\end{equation*}%
is such that $\left\Vert \eta \right\Vert =1.$ The latter equality follows
from the fact $\psi ^{\prime }\Sigma ^{-1}\Omega \Sigma ^{-1}\psi =1+\sqrt{%
p/n_{1}}\theta ,$ which is a consequence of (\ref{Omega modified}) and of
the normalization $\left\Vert \Sigma ^{-1/2}\psi \right\Vert =1$.

Now assume that $n_{1}$ diverges to infinity while $p$ is held constant.
Then, by a CLT%
\begin{equation}
\Omega ^{-1/2}H\Omega ^{-1/2}=I_{p}+Z/\sqrt{n_{1}}+o_{\mathrm{P}}\left(
n_{1}^{-1/2}\right) ,  \label{CLT Case1}\tag{SM8}
\end{equation}%
where $Z$ belongs to GOE. Multiplying (\ref{modified equation 1}) by $(n_{1}/p)^{p/2}$ 
and using (\ref{CLT Case1}), we see that, as $n_{1}\rightarrow
\infty ,$ equation (\ref{modified equation 1}) degenerates to%
\begin{equation*}
\det \left( Z/\sqrt{p}+\theta\eta  \eta ^{\prime }-\mu I_{p}\right) =0\text{
with }\mu =\sqrt{n_{1}/p}\left( \lambda -1\right) .
\end{equation*}%
Hence, PCA degenerates to SMD. \vspace{3 mm}

\textbf{\small{SigD$\rightarrow $PCA:}}
As shown in JO, SigD degenerates to PCA as $%
n_{2}\rightarrow \infty $ while $n_{1}$ and $p$ are held constant.
Therefore, SigD can be linked to SMD via PCA. \vspace{3 mm}

\textbf{\small{REG$_{0}\rightarrow $SMD:}}
Consider REG$_{0}$ with%
\[
\mathbb{E}Y_{1}^{}=\sqrt{\left(  p/n_{1}\right)  ^{1/2}n_{1}\theta}%
\varphi\psi^{\prime},
\]
so that the original value of the spike $\theta$ (see equation (JO\ref{REG0mean}%
)) is scaled by $\left(  p/n_{1}\right)  ^{1/2}$. Suppose now that $n_{1}$
diverges to infinity while $p$ is held constant. Then, by a CLT,%
\begin{equation}
\Sigma^{-1/2}H\Sigma^{-1/2}-I_{p}=Z/\sqrt{n_{1}}+\sqrt{p/n_{1}}\theta\eta
\eta^{\prime}+o_{\mathrm{P}}\left(  n_{1}^{-1/2}\right)  ,
\label{equality link}\tag{SM9}%
\end{equation}
where $Z$ belongs to GOE and $\eta=\Sigma^{-1/2}\psi$. On the other hand,
equation (\ref{basic}) is equivalent to%
\begin{equation}
\det\left(  \Sigma^{-1/2}H\Sigma^{-1/2}-\lambda I_{p}\right)  =0.
\label{modified equation 1REG0}\tag{SM10}%
\end{equation}
Multiplying it by $(n_{1}/p)^{p/2}$ and using (\ref{equality link}), we see
that equation (\ref{modified equation 1REG0}) degenerates to%
\[
\det\left(  Z/\sqrt{p}+\theta\eta\eta^{\prime}-\mu I_{p}\right)  =0\text{ with
}\mu=\sqrt{n_{1}/p}\left(  \lambda-1\right)  .
\]
Hence, REG$_{0}$ degenerates to SMD. \vspace{3 mm}

\textbf{\small{REG$\rightarrow $REG$_{0}$:}}
The REG case degenerates to REG$_{0}$ as $n_{2}\rightarrow \infty $ while $%
n_{1}$ and $p$ are held constant. Therefore, REG can be linked to SMD via REG%
$_{0}$. \vspace{3 mm}

\textbf{\small{CCA$\rightarrow $REG$_{0}$:}}
Recall that the CCA case is based on the solutions to equation (\ref{basic})
with%
\begin{equation*}
H=S_{xy}S_{yy}^{-1}S_{yx}\text{ and }E=S_{xx},
\end{equation*}%
where $S_{xx}$ and $S_{yy}$ are sample covariance matrices corresponding to
i.i.d. $N\left( 0,\Sigma _{xx}\right) $ sample $x_{t}\in \mathbb{R}^{p},$ $%
t=1,...,n_{1}+n_{2},$ and i.i.d. $N\left( 0,\Sigma _{yy}\right) $ sample $%
y_{t}\in \mathbb{R}^{n_{1}},$ $t=1,...,n_{1}+n_{2},$ respectively. Matrices $%
S_{xy}$ and $S_{yx}$ are the corresponding sample cross-covariance matrices.
Since the transformations $x_{t}\mapsto \Sigma _{xx}^{-1/2}x_{t}$ and $%
y_{t}\mapsto \Sigma _{yy}^{-1/2}y_{t}$ do not affect the roots of (\ref%
{basic}), we shall assume without loss of generality that $\Sigma
_{xx}=I_{p} $ and $\Sigma _{yy}=I_{n_{1}}.$ Recall that, by assumption,%
\begin{equation*}
\Sigma _{xy}=\sqrt{\frac{n_{1}\theta }{n_{1}\theta +n_{1}+n_{2}}}\psi
\varphi ^{\prime }.
\end{equation*}

Suppose that $n_{2}$ diverges to infinity while $n_{1}$ and $p$ are held
constant. Then, by a CLT,

\begin{equation*}
S_{xx}=I_{p}+o_{\mathrm{P}}\left( 1\right) ,\text{ }S_{yy}=I_{n_{1}}+o_{%
\mathrm{P}}\left( 1\right) ,
\end{equation*}%
whereas%
\begin{equation*}
S_{xy}=\Sigma _{xy}+Z_{xy}/\sqrt{n_{1}+n_{2}}+o_{\mathrm{P}}\left( \left(
n_{1}+n_{2}\right) ^{-1/2}\right) ,
\end{equation*}%
where $Z_{xy}$ is a $p\times n_{1}$ matrix with i.i.d. $N(0,1)$ entries.
Therefore, equation (\ref{basic}) degenerates to%
\begin{equation}
\det \left( \frac{1}{n_{1}}\left( \tilde{\Sigma}_{xy}+Z_{xy}\right) \left( 
\tilde{\Sigma}_{xy}+Z_{xy}\right) ^{\prime }-\nu I_{p}\right) =0
\label{modified equation 2}\tag{SM11}
\end{equation}%
with%
\begin{equation*}
\tilde{\Sigma}_{xy}=\sqrt{n_{1}\theta }\psi \varphi ^{\prime }
\end{equation*}%
and%
\begin{equation*}
\nu =\left( 1+n_{2}/n_{1}\right) \lambda .
\end{equation*}%
Hence, CCA degenerates to REG$_{0}$. It can further be linked to SMD via REG$%
_{0}$.

\section{The likelihood ratios}

\subsection{SMD entry of Table JO\ref{Table 2}} \label{SMD entry}

The explicit expression for $L^{(\rm{SMD})}\left(  \theta;\Lambda\right)  $ given
in Table JO\ref{Table 2} follows from the following lemma.

\begin{lemma}
\label{Lemma 1}For SMD case, the joint density of the diagonal elements of
$\Lambda$ evaluated at the diagonal elements of $x=\operatorname*{diag}%
\left\{  x_{1},...,x_{p}\right\}  $ with $x_{1}\geq...\geq x_{p}$ equals
\begin{equation}
c_{p}\left(  x\right)  \exp\left\{  -p\theta^{2}/4\right\}  \left.  _{0}%
F_{0}\left(  \Psi,x\right)  \right.  , \label{density0}\tag{SM12}%
\end{equation}
where $c_{p}\left(  x\right)  $ is a quantity that depends on $p$ and $x,$ but
not on $\theta$, and $\Psi=\mathrm{diag}\left\{  \theta p/2,0,...,0\right\}
$. The density under the null hypothesis is obtained from the above expression
by setting $\theta=0$.
\end{lemma}

\textbf{\small{Proof:}} The proof is based on the \textquotedblleft symmetrization
trick\textquotedblright\ used by James (1955) to derive the density of
non-central Wishart distribution. Let $Y=U^{\prime}XU,$ where $U$ is a random
matrix from $\mathcal{O(}p\mathcal{)}$ and $X=Z/\sqrt{p}+\eta\theta
\eta^{\prime}$ with $Z$ from GOE, $\theta\geq0,$ and $\left\Vert
\eta\right\Vert =1.$ Note that the eigenvalues of $X$ and $Y$ are the same.
The joint density of the functionally independent elements of $Y$ evaluated at
$y$ is
\[
\left(  2\pi/p\right)  ^{-p(p+1)/4}2^{-p/2}\int_{\mathcal{O(}p\mathcal{)}%
}\operatorname*{etr}\left\{  -\frac{p}{4}\left(  uyu^{\prime}-\eta\theta
\eta^{\prime}\right)  ^{2}\right\}  \mathcal{(}\mathrm{d}u),
\]
where $\operatorname*{etr}\{\cdot\}$ denotes the exponential trace function, and
$\mathcal{(}\mathrm{d}u)$ is the normalized uniform measure over
$\mathcal{O(}p\mathcal{)}$. Taking the square under $\operatorname*{etr}$ and
factorizing, we obtain an equivalent expression%
\[
\left(  2\pi/p\right)  ^{-p(p+1)/4}2^{-p/2}\exp\left\{  -\frac{p}{4}\theta
^{2}\right\}  \operatorname*{etr}\left\{  -\frac{p}{4}y^{2}\right\}
\int_{\mathcal{O(}p\mathcal{)}}\operatorname*{etr}\left\{  \frac{p\theta}%
{2}uyu^{\prime}\eta\eta^{\prime}\right\}  \mathcal{(}\mathrm{d}u).
\]
Now change the variables from $y$ to $\left(  H,x\right)  ,$ where
$y=HxH^{\prime}$ is the spectral decomposition of $y$. Using the strategy of 
the proof of Muirhead's (1982) Theorem 3.2.17, integrate $H$ out
to obtain (\ref{density0}) with
\[
c_{p}\left(  x\right)  =\frac{p^{p(p+1)/4}\pi^{p\left(  p-1\right)  /4}%
}{2^{p\left(  p-1\right)  /4+p}\Gamma_{p}\left(  p/2\right)  }%
\operatorname*{etr}\left(  -\frac{p}{4}x^{2}\right)
%TCIMACRO{\dprod _{i<j}^{p}}%
%BeginExpansion
{\displaystyle\prod_{i<j}^{p}}
%EndExpansion
\left(  x_{i}-x_{j}\right)  ,
\]
where $\Gamma_{p}\left(  p/2\right)  $ is the multivariate Gamma function. $\square$

\subsection{Identification of the parameters of expression (JO\ref{LRgeneral})} \label{Identification of the parameters}

For the reader's convenience, we provide some extra detail on the
identification of the parameters of expression (JO\ref{LRgeneral}) for the likelihood
ratio $L(\theta ;\Lambda )$ summarized in Table JO\ref{Table 2}. To have a
self-contained source for derivations, we refer below to Muirhead (1982),
henceforth [M], in addition to James (1964), [J] below. \vspace{3 mm}

\textbf{\small{Some Notational conventions.}}
$|A|=\det (A)$, and $c_{pn}$ for a constant depending only on $p,n$. The
hypergeometric function 
\begin{equation*}
_{\mathsf{p}}F_{\mathsf{q}}(a,b;A,B)=\int_{\mathcal{O}(p)}\,_{\mathsf{p}}F_{\mathsf{q}%
}(a,b;AHBH^{\prime })(\mathrm{d}H),
\end{equation*}%
We sometimes drop explicit mention of the parameter vectors $a,b$, and write 
$_{\mathsf{p}}F_{\mathsf{q}}[A;B]$. In particular, we have 
\begin{equation}
%TCIMACRO{%
%\TeXButton{pFq}{_{\mathsf{p}}F_{\mathsf{q}}[A;B]=\,_{\mathsf{p}}F_{\mathsf{q}}[B;A]}}%
%BeginExpansion
_{\mathsf{p}}F_{\mathsf{q}}[A;B]=\,_{\mathsf{p}}F_{\mathsf{q}}[B;A]%
%EndExpansion
\qquad \text{and}\qquad _{\mathsf{p}}F_{\mathsf{q}}[cA;B]=\,_{\mathsf{p}}F_{%
\mathsf{q}}[A;cB],  \label{eq:rules}\tag{SM13}
\end{equation}%
and $%
%TCIMACRO{%
%\TeXButton{pFq}{$_{\mathsf{p}}F_{\mathsf{q}}[A;0]=\,_{\mathsf{p}}F_{\mathsf{q}}[0]=1$}}%
%BeginExpansion
$_{\mathsf{p}}F_{\mathsf{q}}[A;0]=\,_{\mathsf{p}}F_{\mathsf{q}}[0]=1$%
%EndExpansion
$
For $_{0}F_{0}(A)=\etr(A)$ we also have 
\begin{equation}
_{0}F_{0}(A,I+C)=\etr(A)\,_{0}F_{0}(A,C).  \label{eq:0F0}\tag{SM14}
\end{equation}%
To indicate the extension to rank $r$ perturbations, we write $\psi =[\psi
_{1}\cdots \psi _{r}]$ for a $p\times r$ matrix with $\psi ^{\prime }\Sigma^{-1}
\psi =I_{r}$, $\theta =\diag(\theta _{1},\ldots ,\theta _{r})$, $%
1+\theta $ for $I_{r}+\theta $ and $\sqrt{\theta }$ for $\diag(%
\sqrt{\theta }_{1},\ldots ,\sqrt{\theta }_{r})$. \vspace{3 mm}

\textbf{\small{PCA. \ [J, eq. (58)], [M, Th. 9.4.1].}}
We assume a $p\times n_{1}$ matrix $X\sim N(0,\Omega \otimes I_{n_{1}})$
with $\Omega =\Sigma +\psi \theta \psi ^{\prime }$ for $\Sigma >0$ and $\psi
^{\prime }\Sigma ^{-1}\psi =I_{r}$. Without loss of generality we can set $%
\Sigma =I_{p}$. The matrix $n_{1}H=XX^{\prime }$ has eigenvalues $\Lambda =%
\diag(\lambda _{i})$. Using the dictionary

\begin{center}
\begin{tabular}[h]{cccccc}
M: & $S $ & $m $ & $n $ & $\Sigma $ & $L $ \\ \hline
JO: & $H $ & $p $ & $n_1 $ & $\Omega $ & $\Lambda $%
\end{tabular},
\end{center}
\lbrack M, Th. 9.4.1] gives the joint density of $\Lambda $ as 
\begin{equation*}
p(\Lambda |\Omega )=c_{pn_{1}}|\Omega |^{-n_{1}/2}|\Lambda
|^{(n_{1}-p-1)/2}v(\Lambda )\,_{0}F_{0}(-\tfrac{1}{2}n_{1}\Lambda ,\Omega
^{-1}),
\end{equation*}
where
\begin{equation*}
v(\Lambda)=\prod_{i<j}^{p}(\lambda_{i}-\lambda_{j}).
\end{equation*}

Since $_{0}F_{0}(A,I)=\etr(A)$, the likelihood ratio 
\begin{equation*}
L(\theta ;\Lambda )=\frac{p(\Lambda |\Omega )}{p(\Lambda |I)}=|\Omega
|^{-n_{1}/2}\etr(\tfrac{1}{2}n_{1}\Lambda )\,_{0}F_{0}(-\tfrac{1}{2}%
n_{1}\Lambda ,\Omega ^{-1}).
\end{equation*}

We have $|\Omega |=|I+\theta |$, and $\Omega ^{-1}=I-\psi \theta (1+\theta
)^{-1}\psi ^{\prime }$, and referring to \eqref{eq:0F0}, we obtain 
\begin{equation*}
_{0}F_{0}(-\tfrac{1}{2}n_{1}\Lambda ,\Omega ^{-1})=\etr(-\tfrac{1}{2}%
n_{1}\Lambda )\,_{0}F_{0}(\tfrac{1}{2}n_{1}\Lambda ,\psi \theta (1+\theta
)^{-1}\psi ^{\prime }),
\end{equation*}
and arrive at 
\begin{equation*}
L(\theta ;\Lambda )=\left\vert 1+\theta \right\vert ^{-n_{1}/2}\,_{0}F_{0}(\tfrac{1}{2}n_{1}\psi
\theta (1+\theta )^{-1}\psi ^{\prime },\Lambda ).
\end{equation*} \vspace{3 mm}

\textbf{\small{SigD. \ [J, eq. (65), citing Constantine (unpublished)],
[M, Th. 8.2.8].}}
Now assume independent matrices 
\begin{equation*}
X\sim N(0,\Omega \otimes I_{n_{1}})\qquad \text{and}\qquad Y\sim N(0,\Sigma
\otimes I_{n_{2}}),
\end{equation*}%
with dimensions $p\times n_{1}$ and $p\times n_{2}$, and $\Omega =\Sigma
+\psi \theta \psi ^{\prime }$ for $\Sigma >0$ unknown and $\psi ^{\prime
}\Sigma ^{-1}\psi =I_{r}$. Without loss of generality (wlog) we can again set $%
\Sigma =I_{p}$. The sample covariance matrices are given by%
\begin{equation*}
H=XX^{\prime }/n_{1}\qquad \text{and}\qquad E=YY^{\prime }/n_{2}.
\end{equation*}%
Using now the dictionary

\begin{center}
\begin{tabular}[h]{ccccccccc}
M: & $A_1 $ & $A_2$ & $m $ & $n_1 $ & $n_2 $ & $\Sigma_1 $ & $\Sigma_2 $ & $%
\Delta $ \\ \hline
JO: & $n_1H $ & $n_2E $ & $p $ & $n_1 $ & $n_2 $ & $\Omega $ & $I_{p} $ & $%
\Omega $%
\end{tabular},
\end{center}
[M, Th. 8.2.8] gives the joint density of the eigenvalues $F=\diag%
(f_{1},\ldots ,f_{p})$ of 
\begin{equation}
\det (n_{1}H-f_{i}n_{2}E)=0  \label{eq:F-eq}\tag{SM15}
\end{equation}%
as 
\begin{equation*}
p(F|\Omega )=c_{pn_{1}n_{2}}|\Omega
|^{-n_{1}/2}|F|^{(n_{1}-p-1)/2}v(F)\,_{1}F_{0}(\tfrac{1}{2}n;-\Omega ^{-1},F),
\end{equation*}%
where $n=n_{1}+n_{2}$. It is helpful to transform the hypergeometric
function using [M, Lemma 8.2.10], due to Khatri (1967), which says here that 
\begin{equation*}
_{1}F_{0}[-\Omega ^{-1},F]=|I+F|^{-n/2}\,_{1}F_{0}[I-\Omega
^{-1},F(I+F)^{-1}]
\end{equation*}%
Note that, as for PCA, $I-\Omega ^{-1}=\psi \theta (1+\theta )^{-1}\psi
^{\prime }$. The (generalized) eigenvalues $\Lambda =\text{diag}(\lambda
_{1},\ldots ,\lambda _{p})$ of (JO\ref{basic beta form}) are seen to be related to those of (\ref%
{eq:F-eq}) via the transformation $\Lambda =(n_{2}/n_{1})F(I+F)^{-1}$. In
forming the likelihood ratio, terms not depending on $\theta $ cancel,
including the Jacobian of this transformation. Hence we arrive at 
\begin{equation*}
L(\theta ;\Lambda )=\frac{p(\Lambda |\Omega )}{p(\Lambda |I)}=\frac{%
p(F|\Omega )}{p(F|I)}=|1+\theta |^{-n_{1}/2}\,_{1}F_{0}[(n_{1}/n_{2})\psi
\theta (1+\theta )^{-1}\psi ^{\prime },\Lambda ].
\end{equation*} \vspace{3 mm}

\textbf{\small{REG$_{0}$. \ [J, eq. (68)], [M, Exer. 10.9].}}
After reduction to canonical form, we assume that we observe an $n_1 \times
p $ matrix $Y_{1} \sim N(M, I_{n_1} \otimes \Sigma).$ The unnormalized sample
covariance matrix $n_1 H = Y_{1}^{\prime }Y_{1}$ has a non-central Wishart
distribution with non-centrality matrix $\Omega = \Sigma^{-1} M^{\prime }M$.
Without loss of generality we can set $\Sigma = I_p$. Using the dictionary

\begin{center}
\begin{tabular}[h]{cccc}
M: & $A $ & $m $ & $n $ \\ \hline
JO: & $n_1 H $ & $p $ & $n_1 $%
\end{tabular},
\end{center}
[M, Exer. 10.9] gives the joint density of the eigenvalues $W=\diag(w_{i})$ of $%
n_{1}H$ as 
\begin{equation*}
p(W|\Omega )=c_{pn_{1}}\etr(-\tfrac{1}{2}\Omega )\etr(-%
\tfrac{1}{2}W)|W|^{(n_{1}-p-1)/2}v(W)\,_{0}F_{1}(\tfrac{1}{2}n_{1};\tfrac{1}{4}%
\Omega ,W).
\end{equation*}

The low rank assumption (JO\ref{REG0mean}) posits $M=\sqrt{n_{1}}\varphi \sqrt{\theta }\psi
^{\prime }$ with $\varphi ^{\prime }\varphi =\psi ^{\prime }\Sigma ^{-1}\psi
=I_{r} $, so that with $\Sigma =I_{p}$, we have $\Omega =\Omega _{\theta
}=n_{1}\psi \theta \psi ^{\prime }$. Note that $\mathbb{E}H=I+\psi \theta
\psi ^{\prime } $, which explains the normalization chosen for $M$.

The eigenvalues $\Lambda $ of $H$ are related to the eigenvalues $W$ of $%
n_{1}H$ by $\Lambda =W/n_{1}$ and so 
\begin{align*}
L(\theta ;\Lambda )=\frac{p(\Lambda |\Omega _{\theta })}{p(\Lambda |\Omega
_{0})}=\frac{p(W|\Omega _{\theta })}{p(W|\Omega _{0})}& =\etr\{-%
\tfrac{1}{2}n_{1}\theta \}\ _{0}F_{1}[\tfrac{1}{4}n_{1}\psi \theta \psi
^{\prime },n_{1}\Lambda ] \\
& =\etr\{-\tfrac{1}{2}n_{1}\theta \}\ _{0}F_{1}[\tfrac{1}{4}%
n_{1}^{2}\psi \theta \psi ^{\prime },\Lambda ],
\end{align*}%
where we used \eqref{eq:rules}. \vspace{3 mm}

\textbf{\small{REG. \ [J, eq. (73), citing Constantine (1963)], [M, Th. 10.4.2].}}
We are in the situation of REG$_{0}$, but with $\Sigma $ unknown and
estimated by an independent Wishart matrix $n_{2}E\sim W_{p}(n_{2},\Sigma )$%
. [M, Th. 10.4.2] gives the joint density of the eigenvalues $F$ of equation %
\eqref{eq:F-eq}. Using the dictionary

\begin{center}
\begin{tabular}[h]{cccccccc}
M: & $A $ & $B$ & $m $ & $r $ & $n-p $ & $\Sigma $ & $\Omega$ \\ \hline
JO: & $n_1H $ & $n_2E $ & $p $ & $n_1 $ & $n_2 $ & $I_{p} $ & $\Omega $%
\end{tabular},
\end{center}
this may be written as 
\begin{equation*}
p(F|\Omega )=c_{pn_{1}n_{2}}\etr(-\tfrac{1}{2}\Omega
)w(F)\,_{1}F_{1}(\tfrac{1}{2}n,\tfrac{1}{2}n_{1};\tfrac{1}{2}\Omega
,F(I+F)^{-1}).
\end{equation*}%
where $w(F)=|F|^{(n_{1}-p-1)/2}|I+F|^{-(n_{1}+n_{2})/2}v(F)$ does not depend
on $\theta $.

As for SigD, we make the transformation $\Lambda =(n_{2}/n_{1})F(I+F)^{-1}$
to the generalized eigenvalues of (JO\ref{basic beta form}). So, as in previous cases, 
\begin{align*}
L(\theta ;\Lambda )=\frac{p(\Lambda |\Omega _{\theta })}{p(\Lambda |\Omega
_{0})}=\frac{p(F|\Omega _{\theta })}{p(F|\Omega _{0})}& =\etr\{-%
\tfrac{1}{2}n_{1}\theta \}\ _{1}F_{1}[\tfrac{1}{2}n_{1}\psi \theta \psi
^{\prime },(n_{1}/n_{2})\Lambda ] \\
& =\etr\{-\tfrac{1}{2}n_{1}\theta \}\ _{1}F_{1}[\tfrac{1}{2}%
(n_{1}^{2}/n_{2})\psi \theta \psi ^{\prime },\Lambda ].
\end{align*} \vspace{3 mm}

\textbf{\small{CCA. \ [J, eq. (76), citing Constantine (1963)], [M, Th. 11.3.2].}}
We recall some of the steps from [M, Th. 11.2.6], borrowing some text from
Johnstone and Nadler (2015). The canonical correlation problem is invariant
under change of basis for each of the two sets of variables, e.g. [M, Th.
11.2.2]. We may therefore assume that the matrix $\Sigma $ takes the
canonical form 
\begin{equation*}
\Sigma =%
\begin{pmatrix}
I_{p} & \tilde{P} \\ 
\tilde{P}^{\prime} & I_{n_{1}}%
\end{pmatrix}%
,\quad \tilde{P}=[P\ 0],\quad P=\text{diag}(\rho _{1},\ldots ,\rho
_{r},0,\ldots ,0)
\end{equation*}%
where $\tilde{P}$ is $p\times n_{1}$ and the matrix $P$ is of size $p\times
p $ with r non-zero population canonical correlations $\rho_{1},...\rho_{r}$.
Furthermore, in this new basis, we decompose the sample covariance matrix as
follows, 
\begin{equation}
nS=%
\begin{pmatrix}
X^{\prime }X & X^{\prime }Y \\ 
Y^{\prime }X & Y^{\prime }Y%
\end{pmatrix}
\label{eq:DSDT1}\tag{SM16}
\end{equation}%
where the columns of the $n\times p$ matrix $X$ contain the first $p$
variables of the $n\equiv n_{1}+n_{2}$ samples, now assumed to have mean $0$, represented in
the transformed basis. Similarly, the columns of $n\times n_{1}$ matrix $Y$
contain the remaining $n_{1}$ variables. For future use, we note that the
matrix $Y^{\prime }Y\sim W_{n_{1}}(n,I_{n_{1}})$.

The squared canonical correlations $\{r_{i}^{2}\}$ are the eigenvalues of $%
S_{xx}^{-1}S_{xy}S_{yy}^{-1}S_{yx}$. According to [M, Th. 11.3.2], the joint
density of $R^{2}=\diag(r_{1}^{2},\ldots ,r_{p}^{2})$ is given by 
\begin{equation*}
p(R^{2}|P^{2})=c_{pn_{1}n_{2}}|I_{p}-P^{2}|^{n/2}w(R^{2})\,_{2}F_{1}(\tfrac{1}{2}%
n,\tfrac{1}{2}n;\tfrac{1}{2}n_{1};P^{2},R^{2}),
\end{equation*}%
where $w(R^{2})=|R^{2}|^{(n_{1}-p-1)/2}|I_{p}-R^{2}|^{(n_{2}-p-1)/2}v(R^{2})$
does not depend on $P^{2}$. Below, we abbreviate the hypergeometric function
as $_{2}F_{1}[P^{2},R^{2}]$ since the parameters $(\tfrac{1}{2}n,\tfrac{1}{2}n;%
\tfrac{1}{2}n_{1})$ don't change.

If we set $P_{Y}=Y(Y^{\prime }Y)^{-1}Y^{\prime }$ the canonical correlations 
$r_{i}^{2}$ can be rewritten as the roots of $\det (r^{2}X^{\prime
}X-X^{\prime }P_{Y}X)=0.$ Now set $n_{1}H=X^{\prime }P_{Y}X$ and $%
n_{2}E=X^{\prime }(I-P_{Y})X$: the previous equation becomes 
\begin{equation}
\det (n_{1}H-r^{2}(n_{1}H+n_{2}E))=0.  \label{eq:canc}\tag{SM17}
\end{equation}%
We now recall a standard partitioned Wishart argument. Conditional on $Y$,
matrix $X$ is Gaussian with independent rows, and mean and covariance
matrices 
\begin{align*}
M(Y)& =Y\Sigma _{yy}^{-1}\Sigma _{yx}=Y\tilde{P}^{\prime } \\
\Sigma _{xx\cdot y}& =\Sigma _{xx}-\Sigma _{xy}\Sigma _{yy}^{-1}\Sigma
_{yx}=I-P^{2}.
%:=\Phi .
\end{align*}%
Conditional on $Y$, and using Cochran's theorem, the matrices 
\begin{equation*}
n_{1}H\sim W_{p}(n,\Sigma _{xx\cdot y},\Phi (Y))\quad \text{and}\quad
n_{2}E\sim W_{p}(n_{2},\Sigma _{xx\cdot y})
\end{equation*}%
are independent, where the noncentrality matrix 
\begin{equation*}
\Phi (Y)=\Sigma _{xx\cdot y}^{-1}M(Y)^{\prime }M(Y).
%=\Phi ^{-1}\tilde{P}%
%X^{\prime }X\tilde{P}^{\prime }.
\end{equation*}%

The generalized eigenvalues $\lambda _{i}$ of (JO\ref{basic beta form}) are related to the
canonical correlations $r_{i}^{2}$, the generalized eigenvalues of %
\eqref{eq:canc}, by $\lambda _{i}=(n_{2}/n_{1})r_{i}^{2}$. Thus we obtain
the interpretation of the roots of (JO\ref{basic beta form}) in terms of a pair of matrices $%
n_{1}H$ and $n_{2}E$ which are conditionally independently Wishart given
(part of the data) $Y$. Further, as for the previous cases, we can write the
likelihood ratio as 
\begin{align*}
\frac{p(\Lambda |P^{2})}{p(\Lambda |0)}& =\frac{p(R^{2}|P^{2})}{p(R^{2}|0)}%
=|I_{p}-P^{2}|^{n/2}\,_{2}F_{1}[P^{2},R^{2}] \\
& =|I_{p}-P^{2}|^{n/2}\,_{2}F_{1}[(n_{1}/n_{2})P^{2},\Lambda ].
\end{align*}%
Now in our rank $r$ setting, $P^{2}=\sum_{1}^{r}\rho
_{i}^{2}e_{i}e_{i}^{\prime }$ with $\rho _{i}^{2}=n_{1}\theta
_{i}/(n_{1}\theta _{i}+n_{1}+n_{2})$. From the previous display we obtain,
after setting $\psi =[I_{r}\ 0_{r\times (p-r)}]^{\prime }$, 
\begin{equation*}
L(\theta ,\Lambda )=\frac{p(\Lambda |P^{2})}{p(\Lambda |0)}=|I_{r}+n_{1}\theta
/n|^{-n/2}\,_{2}F_{1}[ n_{1}^{2}\psi \theta(n_{2}^{2}I_{r}+n_{2}n_{1}(I_{r}+\theta
))^{-1}\psi ^{\prime },\Lambda ].
\end{equation*}

\section{Contour integral representation}

\subsection{Derivations for Table JO\ref{Table 2a}} \label{Derivations for Table 4}

In this subsection, we obtain decomposition (JO\ref{decomposition 0})%
\begin{equation*}
\mathcal{A}\equiv \frac{\Gamma \left( s+1\right) \alpha \left( \theta
\right) q_{s}}{\sqrt{\pi p}\Psi _{11}^{s}}=\exp \left\{ -\left( p/2\right)
f_{\mathrm{c}}\right\} g_{\mathrm{c}},
\end{equation*}%
where $s=p/2-1,$ and $g_{\mathrm{c}}$ and $f_{\mathrm{c}}$ remain bounded as 
$\mathbf{n},p\rightarrow _{\boldsymbol{\gamma }}\infty ,$ for SMD, PCA, SigD, RED%
$_{0}$, REG, and CCA. The values of $g_{\mathrm{c}}$ and $f_{\mathrm{c}}$
for the different cases are given in Table JO\ref{Table 2a}. \vspace{3 mm}

\textbf{\small{Structure of the prefactor $\mathcal{A}$.}}
Let us rewrite 
\begin{equation}
\mathcal{A}=\alpha (\theta )\frac{\Gamma (s+1)}{(p/2)^{s}\sqrt{\pi p}}\left[ 
\frac{p/2}{\Psi _{11}}\right] ^{s}q_{s}  \label{eq:1}\tag{SM18}
\end{equation}%
as a product of terms $\mathcal{A}_{k}=g_{k}e^{-(p/2)f_{k}}(1+o(1))$ where $%
f_{k},g_{k}$ depend only on $(c_{1},c_{2},\theta ),$ $\mathsf{p}$, and $%
\mathsf{q}$, see \eqref{eq:6} below. The idea is to show the dependence on $%
\mathsf{p}$, $\mathsf{q}$. Referring to Table JO\ref{Table 2}, we have $a_{j}=n/2$ and $%
b_{j}=n_{1}/2$ whenever they are present, and so 
\begin{equation}
q_{s}=\left[ \frac{\Gamma (n/2-s)}{\Gamma (n/2)}\right] ^{\mathsf{p}}\left[ 
\frac{\Gamma (n_{1}/2)}{\Gamma (n_{1}/2-s)}\right] ^{\mathsf{q}}.
\label{eq:2}\tag{SM19}
\end{equation}%
Table JO\ref{Table 2} also shows that 
\begin{equation}
\alpha (\theta )=A(\theta )^{-n_{1}/2},\qquad \qquad \frac{p/2}{\Psi _{11}}=%
\frac{p}{n_{1}}\frac{1}{B(\theta )}\frac{(n_{2}/2)^{\mathsf{p}}}{(n_{1}/2)^{%
\mathsf{q}}}  \label{eq:3}\tag{SM20}
\end{equation}%
where $A(\theta )$ and $B(\theta )$ depend on the particular case in James'
classification. This dependence is shown in Table \ref{tab:AB} below.

\setcounter{table}{0}
\setcounter{figure}{0}

\begin{table}[h]
\caption{Terms $A(\protect\theta), B(\protect\theta)$ in the prefactor ${}_ 
\mathsf{p}\mathcal{A}_\mathsf{q}$, formula \eqref{eq:6}.}
\label{tab:AB}\centering
\vspace{.1in} 
\begin{tabular}[h]{p{.5in}ccc}
Case & $_{\mathsf{p}}F_{\mathsf{q}}$ & $A(\theta)$ & $B(\theta)$ \\ \hline
&  &  &  \\ 
SMD$^*$ & $_0F_0$ & $e^{\theta^2/2}$ & $\theta$ \\[12pt] 
PCA & $_0F_0$ & $1 + \theta$ & $\theta/(1+\theta)$ \\ 
SigD & $_1F_0$ & $1 + \theta$ & $\theta/(1+\theta)$ \\[12pt] 
REG$_0$ & $_0F_1$ & $e^\theta$ & $\theta$ \\ 
REG & $_1F_1$ & $e^\theta$ & $\theta$ \\[12pt] 
CCA & $_2F_1$ & $(1+n_1 \theta/n)^{n/n_1}$ & $\theta/l(\theta)$ \\[20pt] 
\multicolumn{4}{l}{$(^*)$ replace $n_1$ by $p$, \quad $l(\theta) = 1 + \dfrac{%
c_2}{c_1}(1+\theta)$}%
\end{tabular}%
%\tag{SM1}
\end{table}

Combine like terms in \eqref{eq:2} and \eqref{eq:3} to get 
\begin{equation}
\frac{(n_{2}/2)^{\mathsf{p}s}}{(n_{1}/2)^{\mathsf{q}s}}q_{s}=\left( \frac{%
n_{2}}{n}\right) ^{\mathsf{p}s}\frac{\Theta ^{\mathsf{p}}(n/2,p/2)}{\Theta ^{%
\mathsf{q}}(n_{1}/2,p/2)},  \label{eq:4}\tag{SM21}
\end{equation}%
where we define 
\begin{equation}
\Theta (N,M)=\frac{N^{M-1}\Gamma (N-M+1)}{\Gamma (N)}\sim e^{M}\left( 1-%
\frac{M}{N}\right) ^{N-M+1/2}.  \label{eq:5}\tag{SM22}
\end{equation}%
[This is verified at the end of this section.] Finally, define 
\begin{equation}
E(M)=\frac{\Gamma (M)}{M^{M-1}\sqrt{2\pi M}}\sim e^{-M}.  \label{eq:7}\tag{SM23}
\end{equation}

Combining \eqref{eq:1}, \eqref{eq:3}--\eqref{eq:5}, we obtain the desired
form 
\begin{equation}
\mathcal{A}={}_{\mathsf{p}}\mathcal{A}_{\mathsf{q}}=E(p/2)A(\theta
)^{-n_{1}/2}B(\theta )^{-s}\left( \frac{p}{n_{1}}\right) ^{s}\left( \frac{%
n_{2}}{n}\right) ^{\mathsf{p}s}\frac{\Theta ^{\mathsf{p}}(n/2,p/2)}{\Theta ^{%
\mathsf{q}}(n_{1}/2,p/2)}.  \label{eq:6}\tag{SM24}
\end{equation}
Each factor in this product is easily cast in the form $g_k e^{-(p/2)
f_k}(1+o(1))$, with the resulting values of $f_k$ and $g_k$ shown in Table %
\ref{tab:fg}. When needed, we factorize $g_k = \check{g}_k \tilde g_k$ to
show the leading term $\check{g}_k$ and the term $\tilde g_k = 1 + o(1)$,
with the specific dependence of the $o(1)$ term (which comes from the error
bound in Stirling's formula) shown in the final column of Table \ref{tab:fg}.

\begin{table}[h]
\caption{Form of each term in \eqref{eq:6}, when expressed as $%
g_{k}e^{-(p/2)f_{k}}$, with $g_{k}=\check{g}_{k}\tilde{g}_{k}$. Here $%
\protect\vartheta _{m}$ denotes a term that is $O(m^{-1})$.}
\label{tab:fg}\centering
\vspace{.1in} 
\begin{tabular}[h]{c|ccc}
& $f_k$ & $\check{g}_k$ & $\tilde{g}_k$ \\ \hline
&  &  &  \\ 
$E(p/2)$ & $1$ & $1$ & $1 + \vartheta_p$ \\[6pt] 
$A(\theta)^{-n_1/2}$ & $c_1^{-1} \log A(\theta) $ & $1$ & $1$ \\[6pt] 
$B(\theta)^{-s}$ & $\log B(\theta)$ & $B(\theta)$ & $1$ \\[6pt] 
$(p/n_1)^s$ & $-\log c_1$ & $1/c_1$ & $1$ \\[6pt] 
$\left(\dfrac{n_2}{n} \right)^{s}$ & $\log \Big(1+\dfrac{c_2}{c_1} \Big) $ & 
$1+\dfrac{c_2}{c_1} $ & $1$ \\[10pt] 
$\Theta (n/2,p/2) $ & $-1 - \dfrac{r^2}{c_1 c_2} \log \dfrac{r^2}{c_1 + c_2}$
& $\dfrac{r}{(c_1+c_2)^{1/2}}$ & $1 + \vartheta_n + \vartheta_{n-p}$ \\%
[10pt] 
$\Theta^{-1} (n_1/2,p/2)$ \quad & $\ \ 1 + \dfrac{1 - c_1}{c_1} \log(1- c_1) 
$ \quad & \quad $(1 - c_1)^{-1/2}$ \quad & \quad $1 + \vartheta_{n_1} +
\vartheta_{n_1-p}$%
\end{tabular}%
\end{table}

\begin{table}[h]
\caption{Table JO\ref{Table 2a}. Values of $f_{c}$ and $\check g_{c} = g_{c}/(1+o(1))$
for the different cases. The terms $o(1)$ do not depend on $\protect\theta$
and converge to zero as $\mathbf{n},p\rightarrow _{\mathbf{\protect\gamma }%
}\infty $. In the table, $l(\protect\theta)=1+(1+\protect\theta)c_2/c_1$ and 
$r^{2}=c_{1}+c_{2}-c_{1}c_{2}$.}
\label{tab:JO3}\centering
\vspace{.1in} 
\begin{tabular}[h]{lll}
Case & $f_c$ & $\check g_c = g_c/(1+o(1))$ \\ 
&  &  \\ 
SMD & $1 + \theta^2/2 + \log \theta$ & $\theta$ \\[8pt] 
PCA & $1 + \dfrac{1-c_1}{c_1} \log(1+\theta) + \log \dfrac{\theta}{c_1}$ & $%
\theta (1+\theta)^{-1} c_1^{-1}$ \\[8pt] 
SigD & $f_{\mathrm{c}}^{\mathrm{PCA}} + f_{10}$ & $\check g_c^{\mathrm{PCA}}
\check g_{10}$ \\[8pt] 
REG$_0$ & $1 + \dfrac{\theta + c_1}{c_1} + \log \dfrac{\theta}{c_1} + \dfrac{%
1-c_1}{c_1} \log(1-c_1)$ & $\theta c_1^{-1} (1-c_1)^{-1/2}$ \\[5pt] 
REG & $f_{\mathrm{c}}^{\mathrm{REG}_0} + f_{10}$ & $\check g_c^{\mathrm{REG}%
_0} \check g_{10}$ \\[8pt] 
CCA & $f_{\mathrm{c}}^{\mathrm{REG}} + f_{21}$ & $\check g_c^{\mathrm{REG}}
\check g_{10}/l(\theta)$ \\[12pt] 
& $f_{10} = -1 - \dfrac{r^2}{c_1 c_2} \log \dfrac{r^2}{c_1 + c_2} + \log \dfrac{%
c_1 + c_2}{c_1}$ \qquad & $\check g_{10} = c_1^{-1} r (c_1 + c_2)^{1/2}$ \\%
[5pt] 
& $f_{21} = -1 - \dfrac{\theta}{c_1} - \dfrac{r^2}{c_1 c_2} \log \dfrac{r^2}{%
c_1 l(\theta)}$ &  \\[5pt] 
&  & 
\end{tabular}%
\end{table}
\vspace{3 mm}

%\begin{table}[h]
%\caption{Table JO\ref{Table 2a}.}
%\label{tab:JO3}\centering
%\vspace{.1in}  
%\begin{tabular}[h]{lll}
%Case & $f_c$ & $\check g_c = g_c/(1+o(1))$ \\ \hline
%&  &  \\ 
%SMD & $1 + \theta^2/2 + \log \theta$ & $\theta$ \\[12pt] 
%PCA & $1 + \frac{1-c_1}{c_1} \log(1+\theta) + \log \frac{\theta}{c_1}$ & $%
%\theta (1+\theta)^{-1} c_1^{-1}$ \\[5pt] 
%SigD & $f_{\mathrm{c}}^{\mathrm{PCA}} + f_{10}$ & $\check g_c^{\mathrm{PCA}}
%\check g_{10}$ \\[12pt] 
%REG$_0$ & $1 + \frac{\theta + c_1}{c_1} + \log \frac{\theta}{c_1} + \frac{%
%1-c_1}{c_1} \log(1-c_1)$ & $\theta c_1^{-1} (1-c_1)^{-1/2}$ \\[5pt] 
%REG & $f_{\mathrm{c}}^{\mathrm{REG}_0} + f_{10}$ & $\check g_c^{\mathrm{REG}%
%_0} \check g_{10}$ \\[12pt] 
%CCA & $f_{\mathrm{c}}^{\mathrm{REG}} + f_{21}$ & $\check g_c^{\mathrm{REG}}
%\check g_{10}/l(\theta)$ \\[5pt] \hline
%&  &  \\ 
%& $f_{10} = -1 - \frac{r^2}{c_1 c_2} \log \frac{r^2}{c_1 + c_2} + \log \frac{%
%c_1 + c_2}{c_1}$ \qquad & $\check g_{10} = c_1^{-1} r (c_1 + c_2)^{1/2}$ \\%
%[5pt] 
%& $f_{21} = -1 - \frac{\theta}{c_1} - \frac{r^2}{c_1 c_2} \log \frac{r^2}{%
%c_1 l(\theta)}$ &  \\[5pt] \hline
%\end{tabular}
%\end{table}

\textbf{\small{Verification of Table JO\ref{Table 2a}.}}
We write $g_{\mathrm{c}}=\check{g}_{\mathrm{c}}(1+o(1))$ and $g_{10}=\check{g%
}_{10}(1+o(1))$ when we seek to be explicit about the leading term. The $o(1)
$ term differs from row to row, but depends only on $p,n_{1},n_{2}$ (and not 
$\theta $). The explicit dependence can be constructed from the rows of
Table \ref{tab:fg}, from which it is seen in fact always to be $O(m^{-1})$,
where $m=\min (p,n_{1}-p)$.

The lines for SMD, PCA and REG$_{0}$ in Table JO\ref{Table 2a} -- reproduced here as
Table \ref{tab:JO3} below -- are immediately verified from Table \ref{tab:fg}%
. Next, we consider ratios in which the $\mathsf{p}$ index decreases by one
from numerator to denominator. We then have from \eqref{eq:6} 
\begin{equation*}
\frac{\mathcal{A}^{\mathrm{SigD}}}{\mathcal{A}^{\mathrm{PCA}}}=\frac{{}_{1}%
\mathcal{A}_{0}}{{}_{0}\mathcal{A}_{0}}=\frac{\mathcal{A}^{\mathrm{REG}}}{%
\mathcal{A}^{\mathrm{REG}_{0}}}=\frac{{}_{1}\mathcal{A}_{1}}{{}_{0}\mathcal{A%
}_{1}}=\left( \frac{n_{2}}{n}\right) ^{s}\Theta
(n/2,p/2)=g_{10}e^{-(p/2)f_{10}}.
\end{equation*}%
%
%
%
%
% after noting that the numerator and denominator of each ratio have the
% same $\bq$ and differ by one in $\mathsf{p}$.
Referring to Table \ref{tab:fg}, we recover the terms $f_{10}$ and $g_{10}$
and hence the lines for SigD and REG in Table JO\ref{Table 2a}. For future reference, it
is useful to decompose 
\begin{align}
f_{10}& =k_{1}+k_{0},  \notag \\
k_{1}& =-1-\frac{r^{2}}{c_{1}c_{2}}\log r^{2},\qquad k_{0}=\frac{r^{2}}{%
c_{1}c_{2}}\log (c_{1}+c_{2})+\log \frac{c_{1}+c_{2}}{c_{1}}
\label{eq:k01def}\tag{SM25}
\end{align}%
%
%
%
%
% \begin{equation*}
%   f_{10} = k_1 + k_0, \qquad
%   k_1 = -1 - \frac{r^2}{c_1 c_2} \log r^2, \qquad 
%   k_0 = \log \frac{c_1 + c_2}{c_1} + \frac{r^2}{c_1 c_2} \log (c_1 + c_2).
% \end{equation*}

Using \eqref{eq:6} we have, in an obvious notation, 
\begin{align*}
\frac{\mathcal{A}^{\mathrm{CCA}}}{\mathcal{A}^{\mathrm{REG}}} = \frac{{}_2%
\mathcal{A}_1}{{}_1\mathcal{A}_1} & = \left( \frac{A^C}{A^R}
\right)^{-n_1/2} \left( \frac{B^C}{B^R} \right)^{-s} \cdot \left( \frac{n_2}{%
n} \right)^s \Theta(n/2,p/2) = \mathcal{R} \cdot \frac{\mathcal{A}^{\mathrm{%
REG}}}{\mathcal{A}^{\mathrm{REG}_0}}, \\
& = g_{21} e^{-(p/2) f_{21}}.
\end{align*}
and referring to Table \ref{tab:fg}, $\mathcal{R} = l^{-1}(\theta) \exp \{
-(p/2) f_{20} \}$ , where 
\begin{align}
f_{20} = \frac{n_1}{p} \log \frac{A^C}{A^R} + \log \frac{B^C}{B^R} & = \frac{%
c_1 + c_2}{c_1 c_2} \log \frac{c_1 l(\theta)}{c_1 + c_2} - \frac{\theta}{c_1}
- \log l(\theta)  \label{eq:f20def}\tag{SM26} \\
& = k_2 - k_0,  \notag
\end{align}
and 
\begin{equation}  \label{eq:k2def}\tag{SM27}
k_2 = - \frac{\theta}{c_1} + \frac{r^2}{c_1 c_2} \log c_1 l(\theta).
\end{equation}
This establishes the CCA line of Table JO\ref{Table 2a} after we note that 
\begin{equation}  \label{eq:f21}\tag{SM28}
f_{21} = f_{20} + f_{10} = k_2 - k_0 + k_1 + k_0 = k_2 + k_1.
\end{equation}
\vspace{2 mm}
%\newpage

\textbf{\small{Verification of \eqref{eq:5}.}}
 Use Stirling's formula \eqref{eq:7}
twice: 
\begin{align*}
\Gamma(N) & \sim \sqrt{2 \pi N} \ N^{N-1} e^{-N}, \qquad \text{and } \\
\Gamma(N-M+1) & = (N-M) \Gamma(N-M) \\
& \sim \sqrt{2 \pi (N-M)} (N-M)^{N-M} e^{-N+M}
\end{align*}
to arrive at 
\begin{equation*}
\frac{N^{-1} \Gamma(N-M+1)}{\Gamma(N)} \sim \left(\frac{N-M}{N}\right)^{1/2} 
\frac{(N-M)^{N-M}}{N^N} e^M
\end{equation*}
and hence formula \eqref{eq:5}.

\subsection{Proof of Lemma JO\ref{uniform1} (approximation to $_{0}F_{1}$)} \label{Approximation to 0F1}

By equation 9.6.47 in Abramowitz and Stegun (1964), we have%
\begin{equation}
F_{0}=\Gamma \left( m+1\right) \left( m^{2}\eta _{0}\right)
^{-m/2}I_{m}\left( 2m\eta _{0}^{1/2}\right) ,  \label{hyperBessel}\tag{SM29}
\end{equation}%
where the principal branches of the fractional powers are taken, and $%
I_{m}\left( \cdot \right) $ is the modified Bessel function of the first
kind. Using equation 9.7.7 in Abramowitz and Stegun (1964), we obtain%
\begin{equation}
I_{m}\left( 2m\eta _{0}^{1/2}\right) =\frac{\eta _{0}^{m/2}}{\left( 1+4\eta
_{0}\right) ^{1/4}\sqrt{2\pi m}}e^{m\left( 2t_{0}-1-\ln t_{0}\right) }\left(
1+o(1)\right) ,  \label{uniformAsy}\tag{SM30}
\end{equation}%
where $o(1)\rightarrow 0$ as $m\rightarrow \infty $ uniformly with respect
to $\eta _{0}\in \Omega _{0\delta }$ for any $\delta >0.$ Using (\ref%
{uniformAsy}) in (\ref{hyperBessel}), and invoking Stirling's approximation%
\begin{equation*}
\Gamma \left( m+1\right) =m^{m}e^{-m}\sqrt{2\pi m}\left( 1+o(1)\right) ,
\end{equation*}%
we obtain%
\begin{equation*}
F_{0}=\left( 1+4\eta _{0}\right) ^{-1/4}e^{-m\left( -2t_{0}+2+\ln
t_{0}\right) }\left( 1+o(1)\right) .
\end{equation*}%
Since $1-t_{0}=-\eta _{0}/t_{0},$ we obtain $-2t_{0}+2+\ln t_{0}=\varphi
_{0}\left( t_{0}\right) $ and thus,%
\begin{equation*}
F_{0}=\left( 1+4\eta _{0}\right) ^{-1/4}e^{-m\varphi _{0}\left( t_{0}\right)
}\left( 1+o(1)\right) .
\end{equation*}

\subsection{Proof of Lemma JO\ref{1F1approximation} (approximation to $_{1}F_{1}$, $_{2}F_{1}$)} \label{Approximations to 1F1 and 2F1}

First, let us change variable of integration in%
\begin{equation*}
F_{j}=\frac{C_{m}}{2\pi \mathrm{i}}\int_{0}^{(1+)}\exp \left\{ -m\varphi
_{j}\left( t\right) \right\} \psi _{j}\left( t\right) \mathrm{d}t
\end{equation*}%
from $t$ to $\tau =t\eta _{j}.$ We obtain%
\begin{equation}
F_{j}=\frac{C_{m}\eta _{j}^{-m}}{2\pi \mathrm{i}}\int_{0}^{(\eta _{j}+)}\exp
\left\{ -m\phi _{j}\left( \tau \right) \right\} \chi _{j}\left( \tau \right) 
\mathrm{d}\tau ,  \label{changed variables}\tag{SM31}
\end{equation}%
where%
\begin{equation}
\phi _{j}(\tau )=\left\{ 
\begin{array}{ll}
-\tau -\kappa \ln \tau +\left( \kappa -1\right) \ln \left( \tau -\eta
_{j}\right) & \text{for }j=1 \\ 
-\kappa \ln \left( \tau /\left( 1-\tau \right) \right) +\left( \kappa
-1\right) \ln \left( \tau -\eta _{j}\right) & \text{for }j=2%
\end{array}%
\right.  \label{faij}\tag{SM32}
\end{equation}%
and%
\begin{equation*}
\chi _{j}\left( \tau \right) =\left\{ 
\begin{array}{ll}
\left( \tau -\eta _{j}\right) ^{-1} & \text{for }j=1 \\ 
\left( \tau -\eta _{j}\right) ^{-1}\left( 1-\tau \right) ^{-1} & \text{for }%
j=2%
\end{array}%
\right. .
\end{equation*}%
Note that, for $j=2$, the contour in (\ref{changed variables}) does not
encircle $1$.

To obtain point-wise asymptotic approximation to (\ref{changed variables}),
the method of the steepest descent (ascent) is very convenient. However,
establishing the uniformity of the approximation requires the knowledge of
details of the structure of the steepest descent paths. For example, this
knowledge becomes essential when some of the steepest descent paths contain
two saddle points. Unfortunately, for our problem, the steepest descent
paths are relatively complicated. Therefore, we will consider very simple
paths that are steep (but not the steepest) in a neighborhood of a saddle
point. This strategy allows us to rigorously establish the uniformity for
relatively large sets of parameters $\kappa $ and $\eta _{j}$. A downside of
this approach is that we need to explicitly characterize the behavior of $%
\phi _{j}\left( \tau \right) $ on the simple paths, which requires some
relatively lengthy but elementary calculus.

We shall prove Lemma JO\ref{1F1approximation} separately for $j=1$ (REG) and for $j=2$ (CCA).
Therefore, we shall omit subscripts $j$ from the notation below.
\vspace{3 mm}

\textit{Proof of Lemma JO\ref{1F1approximation} for REG.}

\textbf{\small{Saddle points, REG.}}
The saddle points satisfy%
\begin{equation*}
\frac{\mathrm{d}}{\mathrm{d}\tau }\phi \left( \tau \right) =-1-\frac{\kappa 
}{\tau }+\frac{\kappa -1}{\tau -\eta }=-\frac{\tau ^{2}+(1-\eta )\tau
-\kappa \eta }{\tau \left( \tau -\eta \right) }=0.
\end{equation*}%
There are two solutions to this equation%
\begin{equation}
\tau _{\pm }=\tfrac{1}{2}\left\{ \eta +2\kappa -1\pm \sqrt{\left( \eta
+2\kappa -1\right) ^{2}-4\kappa \left( \kappa -1\right) }\right\} -\kappa ,
\label{tauplusminus}\tag{SM33}
\end{equation}%
where we choose the principal branch of the square root (cut along $(-\infty,0]$) when $\newre\eta
\geq -2\kappa +1,$ and the other branch when $\newre\eta
<-2\kappa +1$. The following lemma collects facts about the behavior of $%
\tau _{+}$ for various $\left( \kappa ,\eta \right) $. Suppose that $\kappa
>1$ (which is certainly true if $0<p<\min \left\{ n_{1},n_{2}\right\} $).
Let $\beta =\arg \eta .$ Here and in what follows the principal branch of $%
\arg $ (cut along $(-\infty,0]$) is considered, unless stated otherwise.

\begin{lemma}
\label{ts1facts}(i) If $\newim\eta >0$, then $0<\arg \left( \tau
_{+}-\eta \right) <\beta ;$ if $\newim\eta <0,$ then $\beta <\arg \left(
\tau _{+}-\eta \right) <0$. For real $\eta >0,$ $\tau _{+}$ is real and $%
\tau _{+}>\eta .\medskip $\newline
(ii) For $\eta \in \mathbb{C}\backslash \left( -\infty ,0\right] ,$ function 
$\newre\phi \left( \tau \right) $ is strictly increasing as $\tau $ moves
away from $\tau _{+}$ (in any direction) along the circle with center $\eta $
and radius $\left\vert \tau _{+}-\eta \right\vert $ until it reaches a point 
$B$ on the circle. If $\newim\eta >0$, then $-\pi \leq \arg \left( B-\eta
\right) \leq \beta -\pi $. If $\newim\eta <0,$ then $\pi +\beta \leq \arg
\left( B-\eta \right) \leq \pi $. If $\eta >0,$ then $B=2\eta -\tau _{+}$.
\end{lemma}

\textbf{\small{Proof:} } \textit{(i)} For $\newim\eta >0$ and the branch of the square root
chosen as described above, we have%
\begin{equation*}
\newim\sqrt{\left( \eta +2\kappa -1\right) ^{2}-4\kappa (\kappa -1)}>%
\newim\left( \eta +2\kappa -1\right) =\newim\eta .
\end{equation*}%
Since%
\begin{equation*}
2\newim\left( \tau _{+}-\eta \right) =-\newim\eta +\newim\sqrt{%
\left( \eta +2\kappa -1\right) ^{2}-4\kappa (\kappa -1)},
\end{equation*}%
we have $\newim\left( \tau _{+}-\eta \right) >0.$ Therefore,%
\begin{equation}
\text{if }\newim\eta >0,\text{ then }0<\arg \left( \tau _{+}-\eta \right)
<\pi .  \label{i1}\tag{SM34}
\end{equation}%
Similarly, we can show that if $\newim\eta <0,$ then $-\pi <\arg \left(
\tau _{+}-\eta \right) <0.$

Now let $\rho =\left\vert \tau _{+}-\eta \right\vert $. Then, for $\tau
=\eta +\rho e^{\mathrm{i}x}$ we have%
\begin{eqnarray*}
\newre\phi \left( \tau \right) &=&\left( \kappa -1\right) \ln \rho -\newre\eta -\rho \cos x \\
&&-\frac{\kappa }{2}\ln \left( \rho ^{2}+\left\vert \eta \right\vert
^{2}+2\rho \left\vert \eta \right\vert \cos \left( x-\beta \right) \right) ,
\end{eqnarray*}%
and therefore%
\begin{equation}
\frac{\mathrm{d}}{\mathrm{d}x}\newre\phi \left( \tau \right) =\rho
\left\{ \sin x+k_{\eta }\left( x\right) \sin \left( x-\beta \right) \right\}
,  \label{first derivative on circle}\tag{SM35}
\end{equation}%
where%
\begin{equation}
k_{\eta }\left( x\right) =\frac{\kappa \left\vert \eta \right\vert }{\rho
^{2}+\left\vert \eta \right\vert ^{2}+2\rho \left\vert \eta \right\vert \cos
\left( x-\beta \right) }>0,  \label{ketax}\tag{SM36}
\end{equation}%
unless $\cos \left( x-\beta \right) =-1$ and $\rho =\left\vert \eta
\right\vert $, in which case $\tau =\eta +\rho e^{\mathrm{i}x}=0$ and $\frac{%
\mathrm{d}}{\mathrm{d}x}\newre\phi \left( \tau \right) \rightarrow
-\infty $ as $x\downarrow \beta -\pi $ and $\frac{\mathrm{d}}{\mathrm{d}x}%
\newre\phi \left( \tau \right) \rightarrow +\infty $ as $x\uparrow \beta
+\pi $.

For $\newim\eta >0$, (\ref{first derivative on circle}) implies that 
\begin{align}
\frac{\mathrm{d}}{\mathrm{d}x}\newre\phi \left( \tau \right) &>0\text{
for }x\in \left[ \beta ,\pi \right] \text{, and}  \label{first inequality}\tag{SM37} \\
\frac{\mathrm{d}}{\mathrm{d}x}\newre\phi \left( \tau \right) &<0\text{
for }x\in \left[ \beta -\pi ,0\right] \text{.}  \label{second inequality}\tag{SM38}
\end{align}%
But, since $\tau _{+}$ is a saddle point of $\phi \left( \tau \right) $, 
\begin{equation*}
\frac{\mathrm{d}}{\mathrm{d}x}\newre\phi \left( \tau \right) =0\text{ for 
}x=\arg \left( \tau _{+}-\eta \right) .
\end{equation*}%
Therefore, inequalities (\ref{i1}) and (\ref{first inequality}) guarantee
that $\arg \left( \tau _{+}-\eta \right) \in \left( 0,\beta \right) .$
Similarly, we can show that $\newim\eta <0$ implies that $\arg \left(
\tau _{+}-\eta \right) \in \left( \beta ,0\right) $. The part of \textit{(i)} that
deals with real $\eta >0$ holds by inspection.

\textit{(ii)} Consider the case $\newim\eta >0$. Let us show that there are no
zeros of $\frac{\mathrm{d}}{\mathrm{d}x}\newre\phi \left( \tau \right) $
on $\left( 0,\beta \right) $ other than $\arg \left( \tau _{+}-\eta \right) $%
. First, suppose that $k_{\eta }\left( \beta /2\right) <1$, where $k_{\eta
}\left( \cdot \right) $ is as defined in (\ref{ketax}). Then, since $k_{\eta
}\left( x\right) $ is a decreasing function of $x\in \left( 0,\beta \right)
, $ equation (\ref{first derivative on circle}) implies that all zeros of $%
\frac{\mathrm{d}}{\mathrm{d}x}\newre\phi \left( \tau \right) $ on $x\in
\left( 0,\beta \right) $ must belong to $\left( 0,\beta /2\right) .$ Furthermore,
at any zero $x$ of $\frac{\mathrm{d}}{\mathrm{d}x}\newre\phi \left( \tau
\right) $ we must have $k_{\eta }\left( x\right) <1$. 

Indeed, let $x=\beta/2+y$. Then, 
\begin{equation*}
\sin x+k_{\eta}(x)\sin(x-\beta)=\sin(\beta/2+y)-k_{\eta}(x)\sin(\beta/2-y)
\end{equation*}
On the other hand,
\begin{equation*}
\sin(\beta/2+y)-\sin(\beta/2-y)=2\sin y \cos(\beta/2),
\end{equation*}
which is positive for $0<y<\beta/2$ and negative for $0>y>-\beta/2$. Therefore,
$\sin x+k_{\eta}(x)\sin(x-\beta)$ (and $\frac{\mathrm{d}}{\mathrm{d}x}\newre\phi \left( \tau \right)$) 
is positive for $x\in \left( \beta/2,\beta \right) $, and
it may equal zero for some $x\in \left( 0, \beta/2 \right) $ only if
$k_{\eta }\left( x\right) <1$.

If there are more than
one zero for $x\in \left( 0, \beta/2 \right) $, then by the mean-value theorem there must 
exist $x_{1}\in \left( 0,\beta /2\right) $ such that,
at $x=x_{1}$, $\frac{\mathrm{d}^{2}}{\mathrm{d}x^{2}}\newre\phi \left(
\tau \right) \leq 0$ and $\frac{\mathrm{d}}{\mathrm{d}x}\newre\phi \left(
\tau \right) \geq 0$. The latter inequality and the fact that
$\frac{\mathrm{d}}{\mathrm{d}x}\newre\phi \left(
\tau \right)<0$ at $x=0$
implies that some zeros of $%
\frac{\mathrm{d}}{\mathrm{d}x}\newre\phi \left( \tau \right) $ must be
less than or equal to $x_{1},$ and hence, $k_{\eta }\left( x_{1}\right) <1.$

To summarize, if there are more than one zero of $\frac{\mathrm{d}}{\mathrm{d%
}x}\newre\phi \left( \tau \right) $ on $\left( 0,\beta /2\right) ,$ we
must have 
\begin{equation}
\rho \left\{ \cos x_{1}+k_{\eta }\left( x_{1}\right) \cos \left( x_{1}-\beta
\right) +k_{\eta }^{\prime }\left( x_{1}\right) \sin \left( x_{1}-\beta
\right) \right\} \leq 0  \label{derivative at x1}\tag{SM39}
\end{equation}%
for some $x_{1}\in \left( 0,\beta /2\right) $ with $k_{\eta }\left(
x_{1}\right) <1$. Since 
\begin{equation*}
k_{\eta }^{\prime }\left( x_{1}\right) \sin \left( x_{1}-\beta \right) >0 \text{and}%
\text{ } (1-k_{\eta}(x_{1}))\cos x_{1}>0,
\end{equation*}%
we must have%
\begin{equation*}
\cos x_{1}+\cos \left( x_{1}-\beta \right) <0.
\end{equation*}%
Therefore,%
\begin{equation}
2\cos \left( x_{1}-\beta /2\right) \cos \left( \beta /2\right) <0,
\label{impossible inequality}\tag{SM40}
\end{equation}%
which is impossible for $x_{1}\in \left( 0,\beta \right) $ and $0<\beta <\pi 
$.

Now, suppose that $k_{\eta }\left( \beta /2\right) >1$. Then, all zeros of $%
\frac{\mathrm{d}}{\mathrm{d}x}\newre\phi \left( \tau \right) $ on $x\in
\left( 0,\beta \right) $ belong to $\left( \beta /2,\beta \right) $. If
there are more than one zero, there must exist $x_{1}\in \left( \beta
/2,\beta \right) $ such that $\frac{\mathrm{d}^{2}}{\mathrm{d}x^{2}}\newre%
\phi \left( \tau \right) \leq 0$ at $x=x_{1}$ with $k_{\eta }\left(
x_{1}\right) >1.$ That is, (\ref{derivative at x1}) holds. Since 
\begin{equation*}
k_{\eta }^{\prime }\left( x_{1}\right) \sin \left( x_{1}-\beta \right) >0,%
\text{ }\cos \left( x_{1}-\beta \right) >0,\text{ and }k_{\eta }\left(
x_{1}\right) >1,
\end{equation*}%
we still must have (\ref{impossible inequality}), which is impossible.

Finally, if $k_{\eta }\left( \beta /2\right) =1$ then, since $k_{\eta
}\left( x\right) $ is decreasing, (\ref{first derivative on circle}) implies
that there is only one zero of $\frac{\mathrm{d}}{\mathrm{d}x}\newre\phi
\left( \tau \right) $ on $x\in \left( 0,\beta \right) ,$ which is $x=\beta
/2 $. To summarize, we have shown that 
\begin{equation*}
x_{+}\equiv \arg \left( \tau _{+}-\eta \right)
\end{equation*}%
is the only zero of $\frac{\mathrm{d}}{\mathrm{d}x}\newre\phi \left( \tau
\right) $ on $\left( 0,\beta \right) $. Similar arguments show that there
exists only one zero, say $x_{-}$, of $\frac{\mathrm{d}}{\mathrm{d}x}\newre%
\phi \left( \tau \right) $ on $\left( -\pi ,\beta -\pi \right) $. (If $%
\left\vert \eta \right\vert =\rho $ so that $\newre\phi \left( \tau
\right) $ is singular at $x=\beta -\pi $ with $\frac{\mathrm{d}}{\mathrm{d}x}%
\newre\phi \left( \tau \right) \rightarrow -\infty $ as $x\downarrow
\beta -\pi $ and $\frac{\mathrm{d}}{\mathrm{d}x}\newre\phi \left( \tau
\right) \rightarrow +\infty $ as $x\uparrow \beta +\pi ,$ we formally define $%
\frac{\mathrm{d}}{\mathrm{d}x}\newre\phi \left( \tau \right) $ at $%
x=\beta -\pi $ as zero). 

We will set 
\begin{equation*}
B=\eta +\rho \exp \left\{ \mathrm{i}x_{-}\right\} .
\end{equation*}%
The uniqueness of the zeros of $\frac{\mathrm{d}}{\mathrm{d}x}\newre\phi
\left( \tau \right) $ on $\left( 0,\beta \right) $ and on $\left( -\pi
,\beta -\pi \right) ,$ and inequalities (\ref{first inequality},\ref{second
inequality}) imply \textit{(ii)} for the situation where $\newim\eta >0.$ The
analysis for the cases with $\newim\eta <0$ is similar to the above, and
we omit it.

It remains to note that for real $\eta $ such that $\eta >0$, we have%
\begin{equation*}
\frac{\mathrm{d}}{\mathrm{d}x}\newre\phi \left( \tau \right) =\rho  \left\{ 1+k_{\eta }\left( x\right) \right\} \sin x
\end{equation*}%
which implies the validity of \textit{(ii)} for $\eta >0$. $\square $ \vspace{3 mm}

\textbf{\small{Contours of steep descent, REG.}}
We shall choose the contour of integration in (\ref{changed variables}) so
that it passes through $\tau _{+}$, and $\newre\phi \left( \tau \right) $
increases as $\tau $ moves away from $\tau _{+}$ along the contour, at least
in a neighborhood of $\tau _{+}$. Such contours are called
contours of steep descent\ (of $-\newre\phi \left( \tau \right) $). The
contour consists of a circle with center $\eta $ and radius $\rho
=\left\vert \tau _{+}-\eta \right\vert $ (which, in what follows, we refer
to as the circle) and two overlapping straight segments of opposite
orientations.

\begin{figure}[h]
\centering
\includegraphics[height=4.5247in,width=5.1603in]{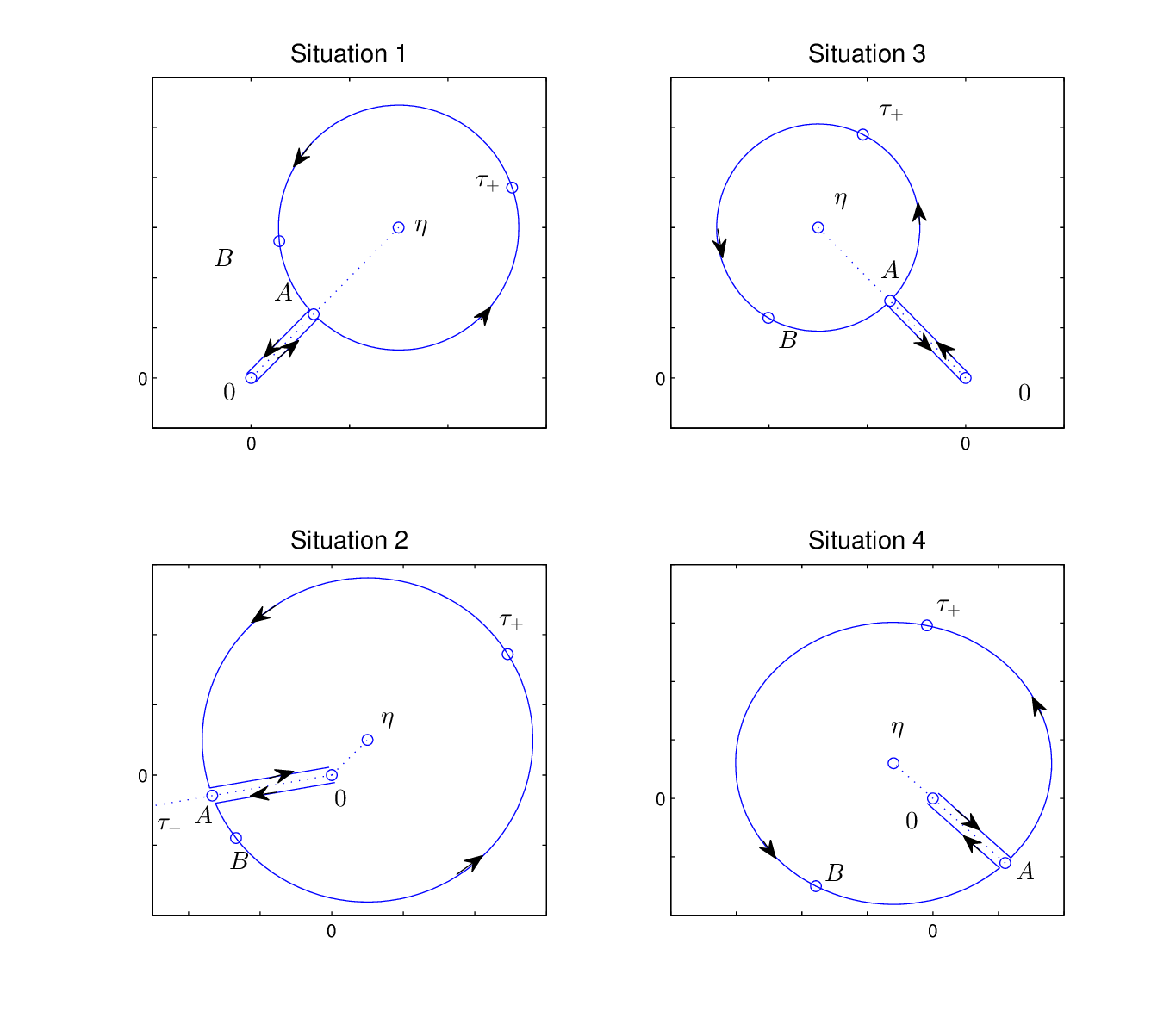}
\caption{Contours of steep descent, $j=1$.}
\label{contourKbeta}
\end{figure}

We consider four situations. The first and the second ones correspond to $%
\newre\eta >0$ and to $\rho <\left\vert \eta \right\vert $ and $\rho \geq
\left\vert \eta \right\vert $, respectively. The third and the fourth ones
correspond to $\newre\eta \leq 0$ and to $\rho <\left\vert \eta
\right\vert $ and $\rho \geq \left\vert \eta \right\vert $, respectively. In
situations 1, 3, and 4, the two straight segments of opposite orientation
connect zero and the point $A$ where the circle is intersected by a
half-line that starts at $\eta $ and passes through zero. In situation 2,
the point $A$ is the intersection point of the circle and a half-line that
starts at $\tau _{-}$ and passes through zero. Figure \ref{contourKbeta}
illustrates the choice of the contour. The points $B$ on the circles are as
defined in Lemma \ref{ts1facts}.

%\FRAME{ftbpFU}{5.1603in}{4.5247in}{0pt}{\Qcb{Contours of steep descent, $j=1$%
%.}}{\Qlb{contourKbeta}}{contourf11all1.eps}{\special{language "Scientific
%Word";type "GRAPHIC";maintain-aspect-ratio TRUE;display "USEDEF";valid_file
%"F";width 5.1603in;height 4.5247in;depth 0pt;original-width
%9.0952in;original-height 7.9667in;cropleft "0";croptop "1";cropright
%"1";cropbottom "0";filename 'contourF11ALL1.eps';file-properties "XNPEU";}}

Let us show that in situation 2, that is when $\newre\eta >0$ and $\rho
>\left\vert \eta \right\vert ,$ the circle intersects the straight
segment $\left[ \tau _{-},0\right) $, as shown in Figure \ref{contourKbeta}. 
Indeed, by definition of $\tau_{\pm}$ we have
\begin{equation}
-\left( \tau _{-}-\eta \right) =\left( \tau _{+}-\eta \right) +\eta +1.
\label{minus plus}\tag{SM41}
\end{equation}%
Since, by Lemma \ref{ts1facts} (i), $\newim\left( \tau _{+}-\eta \right) $
has the same sign as $\newim\eta $ and $\newre\left( \tau _{+}-\eta
\right) \geq 0,$ and since $\newre\left( \eta +1\right) >0$ and $\newim%
\left( \eta +1\right) =\newim\eta ,$ we have%
\begin{equation*}
\left\vert \newre\left\{ -\left( \tau _{-}-\eta \right) \right\}
\right\vert >\left\vert \newre\left( \tau _{+}-\eta \right) \right\vert 
\text{ and }\left\vert \newim\left\{ -\left( \tau _{-}-\eta \right)
\right\} \right\vert \geq \left\vert \newim\left( \tau _{+}-\eta \right)
\right\vert ,
\end{equation*}%
which implies that the circle must intersect the straight segment $\left[
\tau _{-},0\right) $.

We shall split the contour, which we shall call $\mathcal{K}$, in three
parts. In situations 1, 3, and 4, the splitting is%
\begin{equation}
\mathcal{K}=\mathcal{K}_{\left[ 0,A\right] }+\mathcal{K}_{\left[ A,\tau
_{+},B\right] }+\mathcal{K}_{\left[ B,A,0\right] }.  \label{splitting 1}\tag{SM42}
\end{equation}%
This decomposition assumes that $\newim\eta \geq 0$. If the sign of $%
\newim\eta $ changes to the negative, so that $\eta \longmapsto \bar{\eta}
$, then $\mathcal{K}$ is transformed to its complex conjugate, and the
orientation of such a complex conjugate must be changed to the opposite one.
The decomposition then becomes%
\begin{equation}
\mathcal{K}=\mathcal{K}_{\left[ 0,A,B\right] }+\mathcal{K}_{\left[ B,\tau
_{+},A\right] }+\mathcal{K}_{\left[ A,0\right] }.  \label{splitting 2}\tag{SM43}
\end{equation}%
In situation 2, when $\newim \eta \geq 0,$ the splitting is (\ref%
{splitting 2}) because $\arg \left( B-\eta \right) \geq \arg \left( A-\eta
\right) $. (We will verify the latter inequality shortly.) In what follows, we consider only
the case $\newim \eta \geq 0$. The complex conjugate case is analyzed
similarly, and we omit details of such an analysis.

As follows from the proof of Lemma \ref{ts1facts}, $\newre\phi \left(
\tau \right) $ is strictly increasing as $\tau $ is going along $\mathcal{K}%
_{\left[ A,\tau _{+},B\right] }$ away from $\tau _{+}$. In other words, $%
\mathcal{K}_{\left[ A,\tau _{+},B\right] }$ is a contour of steep descent.
Below, we shall use Lemma JO\ref{Olver} to analyze 
\begin{equation*}
\mathcal{I}_{_{\left[ A,\tau _{+},B\right] }}=\int_{\mathcal{K}_{\left[
A,\tau _{+},B\right] }}e^{-m\phi \left( \tau \right) }\chi \left( \tau
\right) \mathrm{d}\tau .
\end{equation*}%
We shall then show that $\mathcal{I}_{\left[ 0,A\right] }$ and $\mathcal{I}_{%
\left[ B,A,0\right] }$, which are defined similarly to $\mathcal{I}_{\left[
A,\tau _{+},B\right] },$ are asymptotically dominated by $\mathcal{I}_{\left[
A,\tau _{+},B\right] }$. However, before we embark on this agenda, let us 
show that $\arg \left( B-\eta \right) \geq \arg \left( A-\eta
\right) $ as was claimed above.

\begin{figure}[h]
\centering
\includegraphics[height=3.5in,width=4.6in]{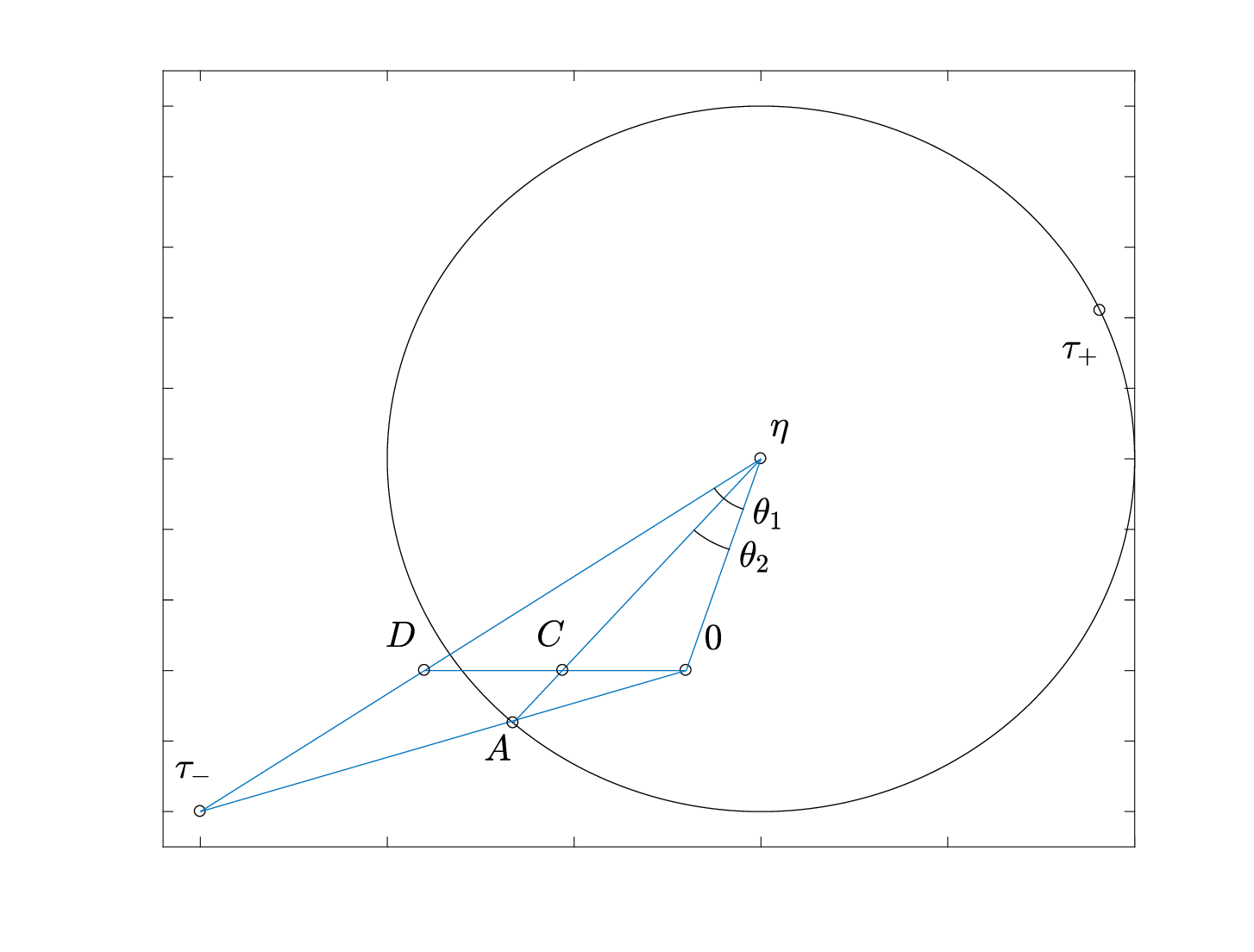}
\caption{An illustration to the argument that $\frac{\mathrm{d}}{\mathrm{d}x}\newre\phi \left( \tau \right)$
is positive at $\tau=A$.}
\label{Sit2Expl}
\end{figure}

As follows from Lemma \ref{ts1facts}, to see that the latter inequality holds, 
it is sufficient to verify that
$\frac{\mathrm{d}}{\mathrm{d}x}\newre\phi \left( \tau \right)$
is positive at $\tau=A$. For such a verification, we will refer to Figure \ref{Sit2Expl}.

First, note that $\tau _{+}\tau_{-}=-\kappa \eta $ and $(\tau _{+}-\eta)(\tau
_{-}-\eta)=(1-\kappa) \eta $, where by assumption, $\kappa>1$. The first of these equalities
implies that $\arg \tau_{-}=-\pi+\arg \eta-\arg \tau_{+}>-\pi,$ so that point
$C$ on Figure \ref{Sit2Expl} rightly belongs to $[A,\eta]$ (the line passing through $D,C,0$
 is a horizontal line). The second of the equalities
implies that the angle $\angle D \eta 0 \equiv \theta_{1}$ equals $\arg (\tau_{+}-\eta)$. Furthermore, we have
\begin{equation*}
\frac{\mathrm{d}}{\mathrm{d}x}\newre\phi \left( \tau \right)=\rho\{\sin x+\frac{|\tau_{+}||\tau_{-}|}{|\tau|^{2}} \sin (x-\beta)\}. 
\end{equation*}

For $\tau=A$, we have $x=\arg (A-\eta)$ and $\beta-x-\pi$ equals $\angle C \eta 0 \equiv \theta_{2} $. This implies
\begin{equation}
\label{theta2eq}\tag{SM44}
\frac{\mathrm{d}}{\mathrm{d}x}\newre\phi \left( \tau \right)=\rho\{-\sin (\beta-\theta_{2})+\frac{|\tau_{+}||\tau_{-}|}{|A|^{2}} \sin \theta_{2}\}.
\end{equation}
For $\tau=\tau_{+}$, the derivative is zero, and hence
\begin{equation}
\label{theta1eq}\tag{SM45}
0=\sin \theta_{1}+\frac{|\tau_{-}|}{|\tau_{+}|} \sin (\theta_{1}-\beta).
\end{equation}

Now, by the law of sines applied to the triangle $\eta C 0$, we have
\begin{equation}
\label{sintheta2}\tag{SM46}
\frac{\sin \theta_{2}}{|C|}=\frac{\sin (\beta-\theta_{2})}{|\eta|}.
\end{equation}
Similarly, for the triangle $\eta D 0$, we have
\begin{equation}
\label{sintheta1}\tag{SM47}
\frac{\sin \theta_{1}}{|D|}=\frac{\sin (\beta-\theta_{1})}{|\eta|}.
\end{equation}
Combining (\ref{sintheta2}) with (\ref{theta2eq}), we obtain
\begin{equation}
\label{deriv}\tag{SM48}
\frac{\mathrm{d}}{\mathrm{d}x}\newre\phi \left( \tau \right)=\rho \sin \theta_{2} \left \{ \frac{|\tau_{+}||\tau_{-}|}{|A|^{2}} -\frac{|\eta|}{|C|} \right \}.
\end{equation}
Combining (\ref{sintheta1}) with (\ref{theta1eq}), we obtain
\begin{equation}
\label{etaeq}\tag{SM49}
|\eta|=\frac{|\tau_{+}||D|}{|\tau_{-}|}
\end{equation}
Using (\ref{etaeq}) in (\ref{deriv}), we get
\begin{equation*}
\frac{\mathrm{d}}{\mathrm{d}x}\newre\phi \left( \tau \right)=\rho \sin \theta_{2} \frac{|\tau_{+}|}{|\tau_{-}|} 
\left \{ \frac{|\tau_{-}|^{2}}{|A|^{2}} -\frac{|D|}{|C|} \right \}>0. \square
\end{equation*}
\vspace{3 mm}

\textbf{\small{Saddle point approximation for $\mathcal{I}_{\left[ A,\protect%
\tau _{+},B\right] },$ REG.}}
We shall now derive a saddle point approximation to the integral $\mathcal{I}%
_{\left[ A,\tau _{+},B\right] }$ which is uniform with respect $\left(
\kappa ,\eta \right) \in \Omega _{\delta },$ where 
\begin{equation}
\Omega _{\delta }=\left\{ \left( \kappa ,\eta \right) \in \mathbb{R} \times \mathbb{C} :\delta \leq \kappa
-1\leq \delta ^{-1},\text{ }\mathrm{dist}\left( \eta ,\mathbb{R}^{-}\right)
\geq \delta ,\text{ and }\left\vert \eta \right\vert \leq \delta
^{-1}\right\} ,  \label{Sigmadelta}\tag{SM50}
\end{equation}%
$\delta $ is an arbitrary fixed number that satisfies inequalities $0<\delta
<1$, $\mathbb{R}^{-}=\left( -\infty ,0\right) $ and, for any $A\subseteq 
\mathbb{C}$ and $B\subseteq \mathbb{C}$, 
\begin{equation*}
\mathrm{dist}\left( A,B\right) =\inf_{a\in A,b\in B}\left\vert
a-b\right\vert .
\end{equation*}%

Let us show that assumptions A0-A5 of Lemma JO\ref{Olver} hold. For this, 
we shall need the following lemma.

\begin{lemma}
\label{Corollary} The quantities $\left\vert \tau _{+}-\eta \right\vert $
and $\left\vert \tau _{+}\right\vert $ are bounded away from zero and
infinity, uniformly with respect to $\left( \kappa ,\eta \right) \in \Omega
_{\delta }.$
\end{lemma}

\textbf{\small{Proof:}} Note that $\left\vert \tau _{+}-\eta \right\vert $ and $%
\left\vert \tau _{+}\right\vert $ are continuous functions of $\left( \kappa
,\eta \right) \in \Omega _{\delta }.$ On the other hand, the definitions (%
\ref{tauplusminus}, \ref{Sigmadelta}) of $\tau _{+}$ and $\Omega _{\delta }$
together with Lemma \ref{ts1facts} imply that $\left\vert \tau _{+}-\eta
\right\vert \neq 0$ and $\left\vert \tau _{+}\right\vert \neq 0$ for any $%
\left( \kappa ,\eta \right) \in \Omega _{\delta }.$ The lemma follows from
these observations and the compactness of $\Omega _{\delta }.\square $

Lemma \ref{Corollary} implies that the length of $\mathcal{K}_{\left[ A,\tau
_{+},B\right] }$ is bounded uniformly with respect to $\left( \kappa ,\eta
\right) \in \Omega _{\delta }$. Further,%
\begin{equation*}
\sup_{\tau \in \mathcal{K}_{_{\left[ A,\tau _{+}\right] }}}\left\vert \tau
-\tau _{+}\right\vert \geq \left\vert A-\tau _{+}\right\vert \text{ and }%
\sup_{\tau \in \mathcal{K}_{_{\left[ \tau _{+},B\right] }}}\left\vert \tau
-\tau _{+}\right\vert \geq \left\vert B-\tau _{+}\right\vert 
\end{equation*}%
with $\left\vert A-\tau _{+}\right\vert $ and $\left\vert
B-\tau _{+}\right\vert $ being continuous functions of $\left( \kappa ,\eta
\right) \in \Omega _{\delta },$ which are not equal to zero for any $\left(
\kappa ,\eta \right) \in \Omega _{\delta }.$ Therefore, $\left\vert A-\tau
_{+}\right\vert $ and $\left\vert B-\tau _{+}\right\vert $ are bounded away
from zero, uniformly with respect to $\left( \kappa ,\eta \right) \in \Omega
_{\delta }$ and assumption A0 holds.

Assumptions A1, A2, A3 and A5 follow from Lemma \ref{Corollary}. Finally,
let $\tau _{1}$ and $\tau _{2}$ be the points of intersection of $\mathcal{K}
$ with a circle with center at $\tau _{+}$ and a sufficiently small fixed
radius $\varepsilon _{1}$. The validity of Assumption A4 follows from the
fact that $\newre\left( \phi \left( \tau _{s}\right) -\phi \left( \tau
_{+}\right) \right) ,$ $s=1,2,$ are positive continuous functions of $\left(
\kappa ,\eta \right) \in \Omega _{\delta }$ (the positivity being a
consequence of Lemma \ref{ts1facts} (ii)) and $\newim\left( \phi \left(
\tau _{s}\right) -\phi \left( \tau _{+}\right) \right) ,$ $s=1,2,$ are
continuous functions of $\left( \kappa ,\eta \right) \in \Omega _{\delta }$.

Since assumptions A0-A5 hold, we have by Lemma JO\ref{Olver}
\begin{equation*}
\mathcal{I}_{\left[ A,\tau _{+},B\right] }=2e^{-m\phi _{0}}\left[ \sqrt{\pi }%
\frac{a_{0}}{m^{1/2}}+\frac{O\left( 1\right) }{m^{3/2}}\right] ,
\end{equation*}%
where $O\left( 1\right) $ is uniform with respect to $\left( \kappa ,\eta
\right) \in \Omega _{\delta }$, 
\begin{equation}
\phi _{0}=-\tau _{+}-\kappa \ln \tau _{+}+\left( \kappa -1\right) \ln \left(
\tau _{+}-\eta \right)  \label{phi0}\tag{SM51}
\end{equation}%
and%
\begin{equation}
a_{0}=\frac{\left( \tau _{+}-\eta \right) ^{-1}}{\sqrt{2\kappa /\tau
_{+}^{2}-2\left( \kappa -1\right) /\left( \tau _{+}-\eta \right) ^{2}}}
\label{a0}\tag{SM52}
\end{equation}%
with the branch of the square root chosen as described in Lemma JO\ref{Olver}.

Precisely, let%
\begin{equation*}
\alpha =\pi /2+\arg \left( \tau _{+}-\eta \right) ,
\end{equation*}%
where the principal branch of $\arg \left( \cdot \right) $ is taken, and let 
\begin{equation}
w=\arg \left( 2\kappa /\tau _{+}^{2}-2\left( \kappa -1\right) /\left( \tau
_{+}-\eta \right) ^{2}\right) ,
\label{omega SM}\tag{SM53}
\end{equation}%
where the branch of $\arg \left( \cdot \right) $ is chosen so that
$\left\vert w+2\alpha \right\vert \leq \pi /2.$
Then%
\begin{equation}
a_{0}=\frac{e^{-\mathrm{i}w/2}\left( \tau _{+}-\eta \right) ^{-1}}{\sqrt{%
\left\vert 2\kappa /\tau _{+}^{2}-2\left( \kappa -1\right) /\left( \tau
_{+}-\eta \right) ^{2}\right\vert }}.  \label{a0adjusted}\tag{SM54}
\end{equation}
\vspace{3 mm}

\textbf{\small{Analysis of $\mathcal{I}_{\left[ 0,A\right] }$ and $\mathcal{I}_{%
\left[ B,A,0\right] },$ REG.}}
Let us show that $\mathcal{I}_{\left[ A,\tau _{+},B\right] }$ asymptotically
dominates $\mathcal{I}_{\left[ 0,A\right] }$ and $\mathcal{I}_{\left[ B,A,0%
\right] }$ uniformly with respect to $\left( \kappa ,\eta \right) \in \Omega
_{1\delta }$, where%
\begin{equation*}
\Omega _{1\delta }=\Omega _{\delta }\cap \left\{ \left( \kappa ,\eta \right) \in \mathbb{R} \times \mathbb{C}
:\newre\eta \geq -2\kappa +1\right\} .
\end{equation*}%
It is sufficient to prove that there exists a positive constant $S,$ such
that, for $\tau $ on $\mathcal{K}_{\left[ 0,A\right] }$ and $\mathcal{K}_{%
\left[ B,A,0\right] },$ we have $\newre\phi \left( \tau \right) \geq 
\newre\phi _{0}+S,$ uniformly with respect to $\left( \kappa ,\eta
\right) \in \Omega _{1\delta }.$ For concreteness, we again assume that $%
\newim\eta \geq 0$. The complex conjugate case is very similar, and we
omit its analysis.

Note that, by Lemma \ref{ts1facts} (ii), for any $\tau \in \mathcal{K}_{%
\left[ B,A\right] },$ $\newre\phi \left( \tau \right) \geq \newre\phi
\left( A\right) .$ Hence, it is sufficient to prove that $\newre\phi
\left( \tau \right) \geq \newre\phi _{0}+S$ for $\tau $ from $\mathcal{K}%
_{\left[ 0,A\right] }$. Moreover, for situations 1, 2 and 4, shown on Figure %
\ref{contourKbeta}, it is sufficient to establish the fact that $\newre%
\phi \left( A\right) \geq \newre\phi _{0}+S.$ Indeed, let $\tau \in 
\mathcal{K}_{\left[ 0,A\right] },$ and let $x=\left\vert \tau \right\vert .$
For situations 1 and 4, using the definition of $\phi(\tau)$ we have, respectively,%
\begin{equation*}
\frac{\mathrm{d}}{\mathrm{d}x}\newre\phi \left( \tau \right) =-\cos \beta
-\kappa /x-\left( \kappa -1\right) /\left( \left\vert \eta \right\vert
-x\right) <0,
\end{equation*}%
and%
\begin{equation*}
\frac{\mathrm{d}}{\mathrm{d}x}\newre\phi \left( \tau \right) =-\cos
\left( \pi -\beta \right) -\kappa /x+\left( \kappa -1\right) /\left(
\left\vert \eta \right\vert +x\right) <0,
\end{equation*}%
where $\beta=\arg \eta$.
Therefore, 
\begin{equation}
\newre\phi \left( A\right) =\inf_{\tau \in \mathcal{K}_{\left[ 0,A\right]
}}\newre\phi \left( \tau \right) .  \label{infOA}\tag{SM55}
\end{equation}%
For situation 2, we have%
\begin{equation*}
\newre\phi \left( \tau \right) =-x\cos \left( \arg \tau _{-}\right)
-\kappa \ln x+\frac{\kappa -1}{2}\ln \left( x^{2}+\left\vert \eta
\right\vert ^{2}-2x\left\vert \eta \right\vert \cos \gamma \right) ,
\end{equation*}%
where $\gamma =2\pi +\arg \tau _{-}-\beta ,$ and thus,%
\begin{eqnarray*}
\frac{\mathrm{d}^{2}}{\mathrm{d}x^{2}}\newre\phi \left( \tau \right)
&=&\kappa /x^{2}+\left( \kappa -1\right) \frac{-x^{2}-\left\vert \eta
\right\vert ^{2}\cos \left( 2\gamma \right) +2x\left\vert \eta \right\vert
\cos \gamma }{\left( x^{2}+\left\vert \eta \right\vert ^{2}-2x\left\vert
\eta \right\vert \cos \gamma \right) ^{2}} \\
&\geq &\kappa /\left\vert \tau \right\vert ^{2}-\left( \kappa -1\right)
/\left\vert \tau -\eta \right\vert ^{2}.
\end{eqnarray*}%
On the other hand, using the fact that $\tau _{+}\tau _{-}=-\kappa \eta $
and Lemma \ref{ts1facts} (i), it is straightforward to verify that $\gamma
>\pi /2$ and therefore, $\left\vert \tau \right\vert ^{2}<\left\vert \tau
-\eta \right\vert ^{2}$ for any $\tau \in \mathcal{K}_{\left[ 0,A\right] }.$
Hence, $\frac{\mathrm{d}^{2}}{\mathrm{d}x^{2}}\newre\phi \left( \tau
\right) >0$. But the first derivative of $\newre\phi \left( \tau \right) $
with respect to $x$ must become positive for $x\rightarrow \infty $,
negative for $x\rightarrow 0$, and zero for $x=\left\vert \tau
_{-}\right\vert ,$ where $x$ is any point on the ray connecting $0$ with $%
\tau _{-}$. Hence, $\frac{\mathrm{d}}{\mathrm{d}x}\newre\phi \left( \tau
\right) $ must be negative for $\tau \in \mathcal{K}_{\left[ 0,A\right] },$
and (\ref{infOA}) again holds.

For situation 3, let $\tau \in \mathcal{K}_{\left[ 0,A\right] }.$ There are
two possibilities. First, there exists $\tau _{1}$ on the circle, such that $%
\newre\tau =\newre\tau _{1}$ and $\left\vert \newim\tau
_{1}\right\vert \leq \left\vert \newim\eta \right\vert $. In such a case, 
$\newre\phi \left( \tau \right) \geq \newre\phi \left( \tau
_{1}\right) $. Furthermore, by Lemma \ref{ts1facts} (ii), $\newre\phi
\left( \tau _{1}\right) >\newre\phi \left( C\right) ,$ where
$C=\eta +\left\vert \tau _{+}-\eta \right\vert $.
Hence,%
\begin{equation}
\newre\phi \left( \tau \right) >\newre\phi \left( C\right) .
\label{more than right}\tag{SM56}
\end{equation}%
Second, $\newre\tau >\newre\eta +\left\vert \tau _{+}-\eta \right\vert 
$. Assuming that $\left( \kappa ,\eta \right) \in \Omega _{1\delta },$ the
latter inequality implies that $\newre\tau \geq -\kappa $. Indeed, for $%
\left( \kappa ,\eta \right) \in \Omega _{1\delta }$, the definition (\ref%
{tauplusminus}) of $\tau _{+}$ implies that $\newre\tau _{+}\geq -\kappa
. $ Therefore, 
\begin{equation*}
\newre\tau >\newre\eta +\left\vert \tau _{+}-\eta \right\vert \geq 
\newre\tau _{+}\geq -\kappa .
\end{equation*}%
Let $x=\left\vert \tau \right\vert ,$ then%
\begin{equation*}
\frac{\mathrm{d}}{\mathrm{d}x}\newre\phi \left( \tau \right) =-\cos \beta
-\kappa /x-\left( \kappa -1\right) /\left( \left\vert \eta \right\vert
-x\right) ,
\end{equation*}%
and%
\begin{equation*}
\frac{1}{\cos \beta }\frac{\mathrm{d}}{\mathrm{d}x}\newre\phi \left( \tau
\right) =-1-\frac{\kappa }{\newre\tau }-\frac{\kappa -1}{\newre\eta -%
\newre\tau }.
\end{equation*}%
But 
\begin{equation*}
-1-\frac{\kappa }{\newre\tau }\geq 0\text{ and }-\frac{\kappa -1}{\newre%
\eta -\newre\tau }>0.
\end{equation*}%
Therefore $\frac{\mathrm{d}}{\mathrm{d}x}\newre\phi \left( \tau \right)
<0,$ and (\ref{more than right}) holds. Note that in the analysis of
situation 3 we used the assumption $\left( \kappa ,\eta \right) \in \Omega
_{1\delta },$ and in particular that $\newre\eta \geq -2\kappa +1$. If
the latter inequality is not satisfied, the minimum of $\newre\phi \left(
\tau \right) $ on $\mathcal{K}$ can be achieved at some point on $\mathcal{K}%
_{[0,A]}$. This fact will be used later, in our proof of Theorem JO\ref{Icase}.

It remains to show that, for some positive $S,$%
\begin{equation}
\newre\phi \left( A\right) \geq \newre\phi _{0}+S  \label{ReA}\tag{SM57}
\end{equation}%
uniformly with respect to $\left( \kappa ,\eta \right) \in \Omega _{1\delta
},$ and 
\begin{equation}
\newre\phi \left( C\right) \geq \newre\phi _{0}+S  \label{ReC}\tag{SM58}
\end{equation}%
uniformly with respect to $\left( \kappa ,\eta \right) \in \tilde{\Omega}%
_{1\delta },$ where%
\begin{equation*}
\tilde{\Omega}_{1\delta }=\Omega _{1\delta }\cap \left\{ \kappa \geq 1,\newre%
\eta \leq 0\right\} .
\end{equation*}%
Inequality (\ref{ReC}) follows from the fact that function $\newre\phi
\left( C\right) -\newre\phi _{0}$ is continuous and positive for $\left(
\kappa ,\eta \right) \in \tilde{\Omega}_{1\delta }$ and from the compactness
of $\tilde{\Omega}_{1\delta }.$

We cannot use a similar argument to establish inequality (\ref{ReA}) because
$\newre\phi \left( A\right) -\newre\phi _{0}$ is not a
continuous function of $\left( \kappa ,\eta \right) \in \Omega _{1\delta }$,
as we may have $A=B=0$ and $\newre\phi \left( A\right) =+\infty $
for some $\left( \kappa ,\eta \right) \in \Omega _{1\delta }$. However, we
can bound $\newre\phi \left( A\right) -\newre\phi _{0}$ from below by
the minimum of two positive continuous functions $\newre\phi \left( \tau
_{1}\right) -\newre\phi _{0}$ and $\newre\phi \left( \tau _{2}\right) -%
\newre\phi _{0},$ where $\tau _{1}$ and $\tau _{2}$ are the points of the
intersection of the circle with center $\eta $ and radius $\left\vert \tau
_{+}-\eta \right\vert $ and a circle with center $\tau _{+}$ and a fixed
radius, which is smaller than $\left\vert A-\tau _{+}\right\vert $,
uniformly with respect to $\left( \kappa ,\eta \right) \in \Omega _{1\delta
} $. Therefore, there exists $S>0$ such that (\ref{ReA}) holds uniformly
with respect to $\left( \kappa ,\eta \right) \in \Omega _{1\delta }.$
\vspace{3 mm}

\textbf{\small{Asymptotics in terms of $\protect\varphi $ and $\protect\psi $,
REG.}}
The above analysis implies the following asymptotic representation%
\begin{equation*}
F_{1}=\frac{C_{m}\eta ^{-m}}{2\pi \mathrm{i}}2e^{-m\phi _{0}}\left[ \sqrt{%
\pi }\frac{e^{-\mathrm{i}w/2}\left( \tau _{+}-\eta \right) ^{-1}}{m^{1/2}%
\sqrt{\left\vert 2\kappa /\tau _{+}^{2}-2\left( \kappa -1\right) /\left(
\tau _{+}-\eta \right) ^{2}\right\vert }}+\frac{O\left( 1\right) }{m^{3/2}}%
\right] ,
\end{equation*}%
where $O\left( 1\right) $ is uniform with respect to $\left( \kappa ,\eta
\right) \in \Omega _{1\delta }.$ We would like to express this formula in
terms of $t_{1},\varphi _{1}\left( \cdot \right) ,$ and $\psi _{1}\left(
\cdot \right)$. 
As follows from the definition of $\varphi_{1}$ (see equation (JO\ref{phij})) and 
the fact that $\tau_{+}=t_{1}\eta$, 
\begin{equation*}
\left\vert 2\kappa /\tau _{+}^{2}-2\left( \kappa -1\right) /\left( \tau
_{+}-\eta \right) ^{2}\right\vert =\left\vert 2\varphi _{1}^{\prime \prime
}\left( t_{1}\right) /\eta ^{2}\right\vert ,
\end{equation*}%
Furthermore, 
\begin{equation*}
\varphi _{1}\left( t_{1}\right) =\phi _{0}-\ln \eta ,
\end{equation*}%
and by (JO\ref{gammaj})%
\begin{equation*}
\left( \tau _{+}-\eta \right) ^{-1}=\psi _{1}\left( t_{1}\right) \eta ^{-1}.
\end{equation*}%
Therefore, we have%
\begin{equation*}
F_{1}=C_{m}e^{-m\varphi _{1}\left( t_{1}\right) }\left[ \frac{e^{-\mathrm{i}%
w/2}e^{-\mathrm{i\arg }\eta }}{\mathrm{i}}\frac{\psi _{1}\left( t_{1}\right) 
}{\sqrt{\left\vert 2\pi m\varphi _{1}^{\prime \prime }\left( t_{1}\right)
\right\vert }}+\frac{O\left( 1\right) }{m^{3/2}}\right] .
\end{equation*}%
On the other hand, by definition (\ref{omega SM}),%
\begin{equation*}
w=\arg \left( \varphi _{1}^{\prime \prime }\left( t_{1}\right) \right)
-2\arg \eta =\omega _{1}-\pi -2\arg \eta ,
\end{equation*}%
where $\omega _{1}$ is as defined in equation (JO\ref{omegaj}). Hence,%
\begin{eqnarray*}
F_{1} &=&C_{m}e^{-m\varphi _{1}\left( t_{1}\right) }e^{-\mathrm{i}\omega
_{1}/2}\left[ \frac{\psi _{1}\left( t_{1}\right) }{\sqrt{\left\vert 2\pi
m\varphi _{1}^{\prime \prime }\left( t_{1}\right) \right\vert }}+\frac{%
O\left( 1\right) }{m^{3/2}}\right] \\
&=&C_{m}\psi _{1}\left( t_{1}\right) e^{-\mathrm{i}\omega _{1}/2}\left\vert
2\pi m\varphi _{1}^{\prime \prime }\left( t_{1}\right) \right\vert
^{-1/2}\exp \left\{ -m\varphi _{1}\left( t_{1}\right) \right\} \left(
1+o(1)\right) .
\end{eqnarray*}
\vspace{3 mm}

\textit{Proof of Lemma JO\ref{1F1approximation} for CCA.}

\textbf{\small{Saddle points, CCA.}}
From equation (\ref{faij}) with $j=2$, we see that the saddle points satisfy%
\begin{equation*}
\frac{\mathrm{d}}{\mathrm{d}\tau }\phi \left( \tau \right) =-\frac{\kappa }{%
\tau }-\frac{\kappa }{1-\tau }+\frac{\kappa -1}{\tau -\eta }=\frac{\tau
^{2}\left( \kappa -1\right) +\tau -\kappa \eta }{\tau \left( \tau -1\right)
\left( \tau -\eta \right) }=0.
\end{equation*}%
There are two solutions to this equation%
\begin{equation}
\tau _{\pm }=\frac{-1\pm \sqrt{1+4\kappa \left( \kappa -1\right) \eta }}{%
2\left( \kappa -1\right) },  \label{tauplusminusj2}\tag{SM59}
\end{equation}%
where we choose the principal branch of the square root cut along $(-\infty,0]$. 

The following lemma collects facts about the behavior of $\tau _{+}$ for
various $\left( \kappa ,\eta \right) $. As usual, we assume that $\kappa
>1.$ In addition, we assume that $\eta \notin \left( -\infty ,\eta _{\ast
}\right) \cup \left[ 1,\infty \right)$, where 
\begin{equation*}
\eta _{\ast }=-\frac{1}{4\kappa \left( \kappa -1\right) }.
\end{equation*}
Note that set $\left( -\infty ,\eta _{\ast
}\right) \cup \left[ 1,\infty \right)$ does not intersect with $\Omega_{2\delta}$
for any $\delta>0$.

\begin{lemma}
\label{ts1facts2F1}(i) $\left\vert \tau _{+}-\eta \right\vert <\left\vert
1-\eta \right\vert ,$ and $\tau _{+}=0$ if and only if $\eta =0.\medskip $%
\newline
(ii) If $\newim\eta >0$ and $\newre\tau _{+}>1/2,$ then $0<\newim%
\tau _{+}<\newim\eta .$ If $\newim\eta >0$ and $\newre\tau
_{+}<1/2, $ then $\newim\tau _{+}>\newim\eta .$ Similarly, if $\newim%
\eta <0$ and $\newre\tau _{+}>1/2,$ then $0>\newim\tau _{+}>\newim%
\eta .$ If $\newim\eta <0$ and $\newre\tau _{+}<1/2,$ then $\newim%
\tau _{+}<\newim\eta .\medskip $\newline
(iii) For $\eta \notin \left( -\infty ,\eta _{\ast }\right) \cup \left[
1,\infty \right) ,$ function $\newre\phi \left( \tau \right) $ is
strictly increasing as $\tau $ moves away from $\tau _{+}$ (in any
direction) along the circle with center $\eta $ and radius $\left\vert \tau
_{+}-\eta \right\vert $ until it reaches a point $B$ on the circle.
\end{lemma}

\textbf{\small{Proof:} } \textit{(i)} Let 
\begin{equation}
-\eta _{\ast }^{-1}\left( \eta -\eta _{\ast }\right) =\rho ^{2}\exp \left\{ 
\mathrm{i}2\theta \right\}  \label{polar}\tag{SM60}
\end{equation}%
with $\theta \in \left( -\pi /2,\pi /2\right) $. Then 
\begin{equation}
\tau _{+}=\frac{-1+\rho \exp \left\{ \mathrm{i}\theta \right\} }{2\left(
\kappa -1\right) }  \label{tauplusj2}\tag{SM61}
\end{equation}%
and a direct calculation (we perform it using Maple's symbolic algebra
software) shows that%
\begin{equation*}
\eta _{\ast }^{-2}\left( \left\vert \tau _{+}-\eta \right\vert
^{2}-\left\vert 1-\eta \right\vert ^{2}\right) =-4\kappa \left( \kappa +\rho
\cos \theta -1\right) \left( \left( 2\kappa -1\right) ^{2}+\rho ^{2}-2\rho
\left( 2\kappa -1\right) \cos \theta \right) .
\end{equation*}%
Since $\theta \in \left( -\pi /2,\pi /2\right) $ and $\kappa>1$, the latter expression is
less than zero. Further, equation (\ref{tauplusj2}) implies
that $\tau _{+}=0$ if and only if $\theta =0$ and $\rho
=1. $ The latter two equalities are equivalent to $\eta =0$.

\textit{(ii)} From (\ref{tauplusj2}), we see that $\newre\tau _{+}>1/2$ if and
only if%
\begin{equation}
\rho \cos \theta >\kappa .  \label{ReMore}\tag{SM62}
\end{equation}%
On the other hand,%
\begin{equation}
-\eta _{\ast }^{-1}\newim\left( \tau _{+}-\eta \right) =2\rho \sin \theta
\left( \kappa -\rho \cos \theta \right) .  \label{ImMore}\tag{SM63}
\end{equation}%
Combining (\ref{ReMore}) and (\ref{ImMore}), we obtain \textit{(ii)}.

\begin{figure}[h]
\centering
\includegraphics[height=2.8184in,width=5.77in]{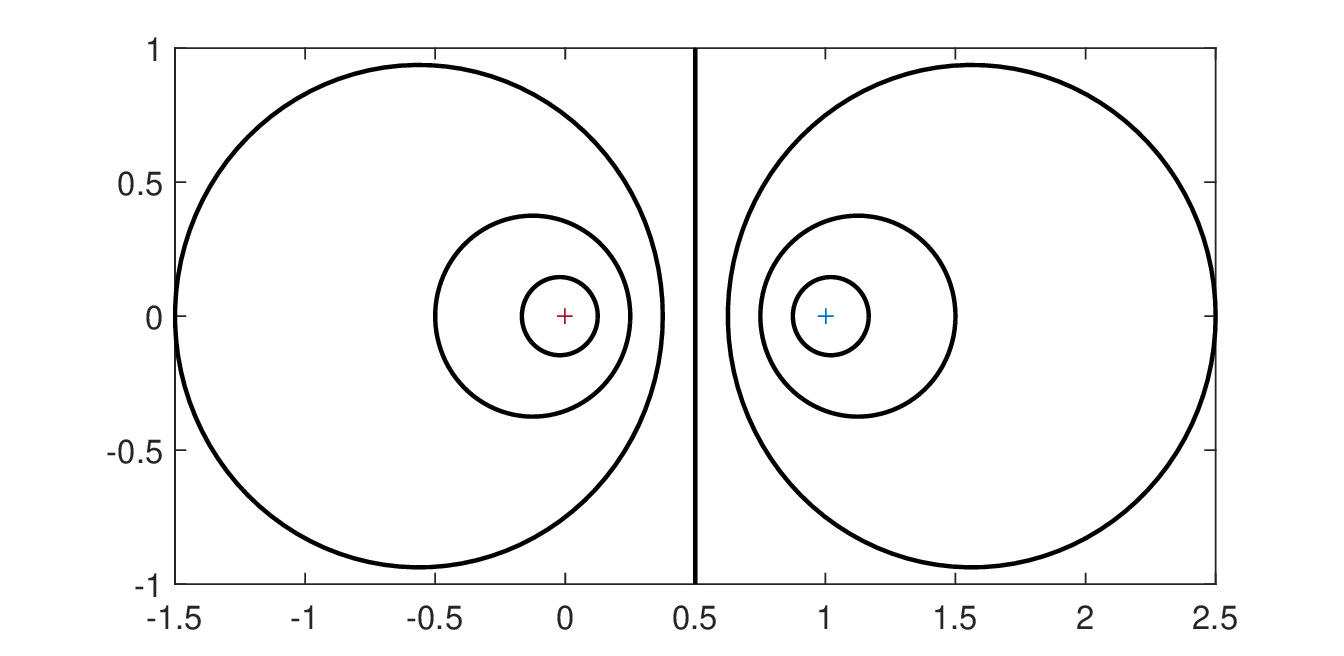}
\caption{Isolines of the function $\left\vert \protect\tau /\left( 1-\protect\tau \right) \right\vert .$}
\label{isolines}
\end{figure}

\textit{(iii)} Recall that%
\begin{equation*}
\phi \left( \tau \right) =-\kappa \ln \frac{\tau }{1-\tau }+\left( \kappa
-1\right) \ln \left( \tau -\eta \right) .
\end{equation*}%
Therefore, on the circle with center $\eta $ and radius $\left\vert \tau
_{+}-\eta \right\vert$, $\newre\phi \left( \tau \right) $ equals 
$-\kappa \ln \left\vert \tau /\left( 1-\tau \right)
\right\vert $ plus a constant. Further, for $c>0$ such that $c\neq 1,$ the set of $\tau $
that satisfy equality $\left\vert \tau /\left( 1-\tau \right) \right\vert =c$
is a circle with center $c^{2}/\left( c^{2}-1\right) $ and radius $%
c/\left\vert c^{2}-1\right\vert .$ For $c=1,$ $\left\vert \tau /\left(
1-\tau \right) \right\vert =c$ along the line $\newre\tau =1/2.$ Figure 
\ref{isolines} shows the iso-lines of $\left\vert \tau /\left( 1-\tau
\right) \right\vert $. For $c<1,$ the isolines are encircling $0,$ for $c>1,$
they are encircling $1$.

%\FRAME{ftbpFU}{5.77in}{2.8184in}{0pt}{\Qcb{Isolines of the function $%
%\left\vert \protect\tau /\left( 1-\protect\tau \right) \right\vert .$}}{\Qlb{%
%isolines}}{isolines.eps}{\special{language "Scientific Word";type
%"GRAPHIC";maintain-aspect-ratio TRUE;display "USEDEF";valid_file "F";width
%5.77in;height 2.8184in;depth 0pt;original-width 8.9577in;original-height
%4.3613in;cropleft "0";croptop "1";cropright "1";cropbottom "0";filename
%'isolines.eps';file-properties "XNPEU";}}

Since $\tau_{+}$ is a critical point of $\newre\phi \left( \tau \right)$,
the circle with center $\eta $ and radius $%
\left\vert \tau _{+}-\eta \right\vert $ must have a common tangent with one
of the isolines at $\tau =\tau _{+}.$ Therefore, $\newre\phi \left( \tau
\right) $ must be strictly monotone as $\tau $ moves away from $\tau _{+}$
along the circle with center $\eta $ and radius $\left\vert \tau _{+}-\eta
\right\vert $ until it reaches a point $B$ on the circle. Part \textit{(ii)} of the
lemma implies that $\newre\phi \left( \tau \right) $ is strictly
increasing. $\square $
\vspace{3 mm}

\textbf{\small{Contours of steep descent, CCA.}}
We shall choose the contour of integration in (\ref{changed variables}),
which we shall call $\mathcal{K}$, so that it passes through $\tau _{+}$,
and $\newre\phi \left( \tau \right) $ increases as $\tau $ moves away
from $\tau _{+}$ along the contour, at least in a neighborhood of $\tau
_{+}$. The contour consists of a circle with center $\eta $ and radius $%
r=\left\vert \tau _{+}-\eta \right\vert $, which, in what follows, we refer
to as $C_{1}$, and two overlapping circular segments of opposite
orientations, which we will refer to as $C_{2}$.

We consider four situations. The first and the second ones correspond to $%
r<\left\vert \eta \right\vert $ and to $\newre\eta <0$ and $\newre\eta
\geq 0$, respectively. The third and the fourth ones correspond to $r\geq
\left\vert \eta \right\vert $ and to $\newre\eta <0$ and $\newre\eta
\geq 0$, respectively. Using (\ref{tauplusj2}), we obtain%
\begin{equation*}
\eta _{\ast }^{-2}\left\vert \tau _{+}-\eta \right\vert ^{2}-\eta _{\ast
}^{-2}\left\vert \eta \right\vert ^{2}=4\kappa \left( \rho ^{2}-2\rho \cos
\theta +1\right) \left( \kappa -\rho \cos \theta -1\right) .
\end{equation*}%
Therefore, situations 3 or 4 are realized whenever%
\begin{equation}
\rho \cos \theta \leq \kappa -1.  \label{34}\tag{SM64}
\end{equation}%
In particular, the corresponding $\tau _{+}$ must be such that $\newre%
\tau _{+}<1/2$ (compare to (\ref{ReMore})).

\begin{figure}[h]
\centering
\includegraphics[height=2.3981in,width=4.6in]{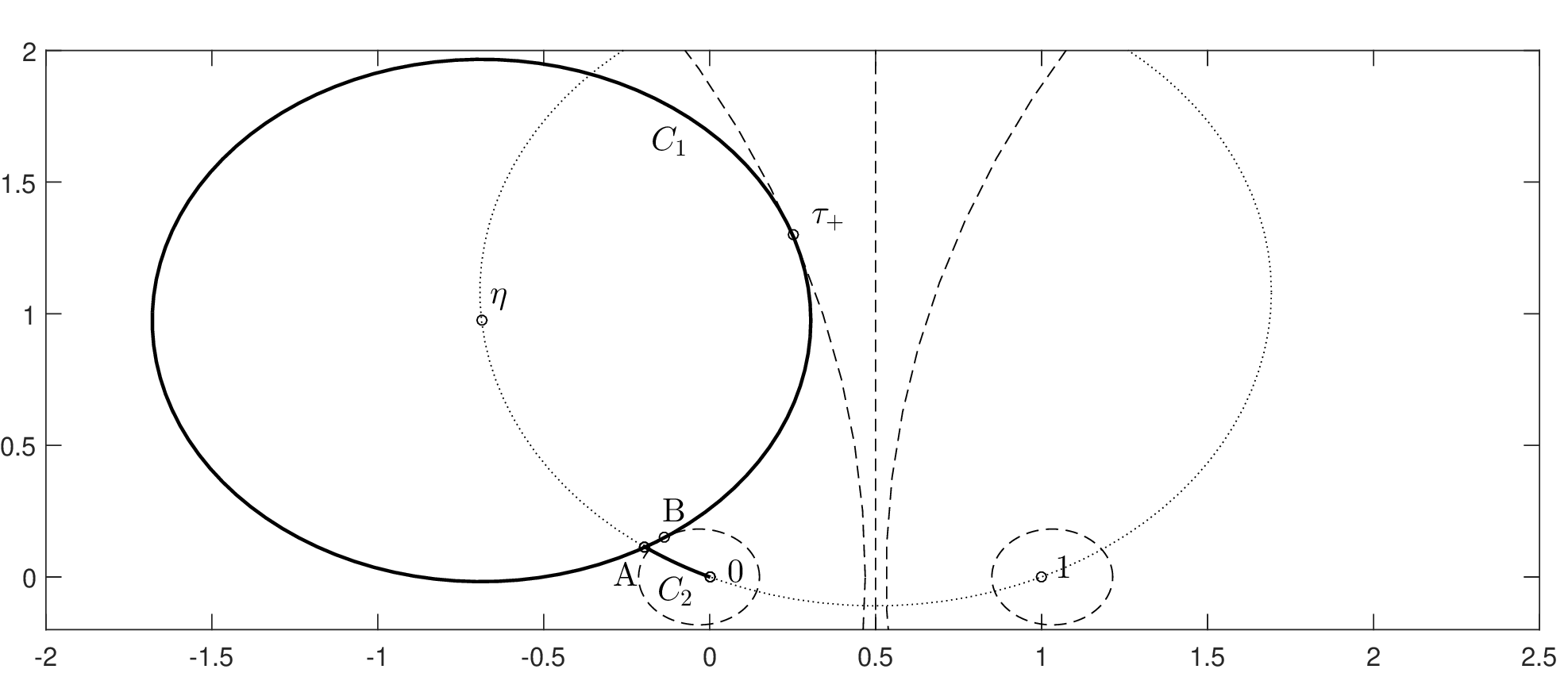}
\caption{Choice of contour $\mathcal{K}$
in situation 1. The contour is represented by the dark black circle and the
circle segment ending at $0$. The dashed lines are iso-lines of function $%
\left\vert \protect\tau /\left( \protect\tau -1\right) \right\vert .$}
\label{situation0}
\end{figure}

\begin{figure}[h]
\centering
\includegraphics[height=2.4111in,width=4.6in]{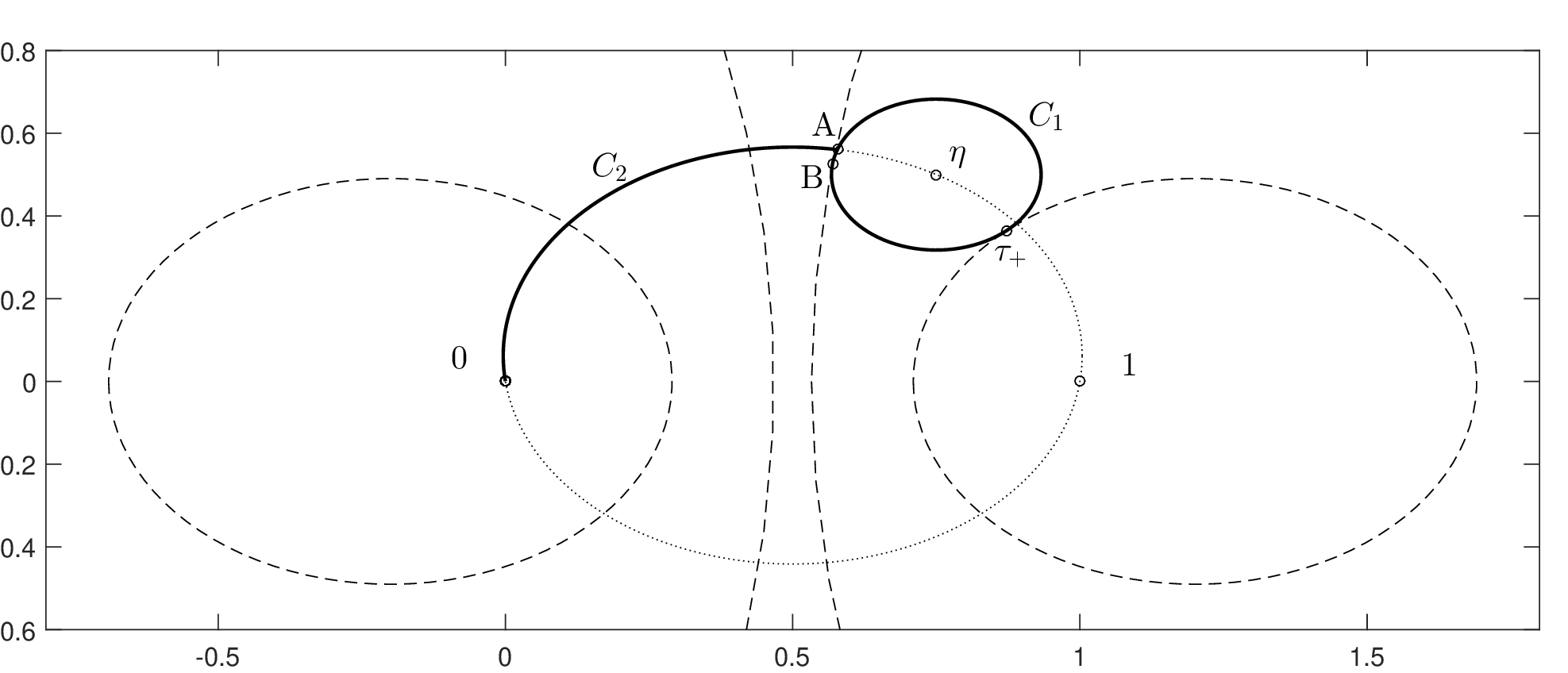}
\caption{Choice of contour $\mathcal{K}$
in situation 2. The contour is represented by the dark black circle and the
circle segment ending at $0$. The dashed lines are iso-lines of function $%
\left\vert \protect\tau /\left( \protect\tau -1\right) \right\vert .$}
\label{situationnew}
\end{figure}

\begin{figure}[h]
\centering
\includegraphics[height=2.2563in,width=4.8in]{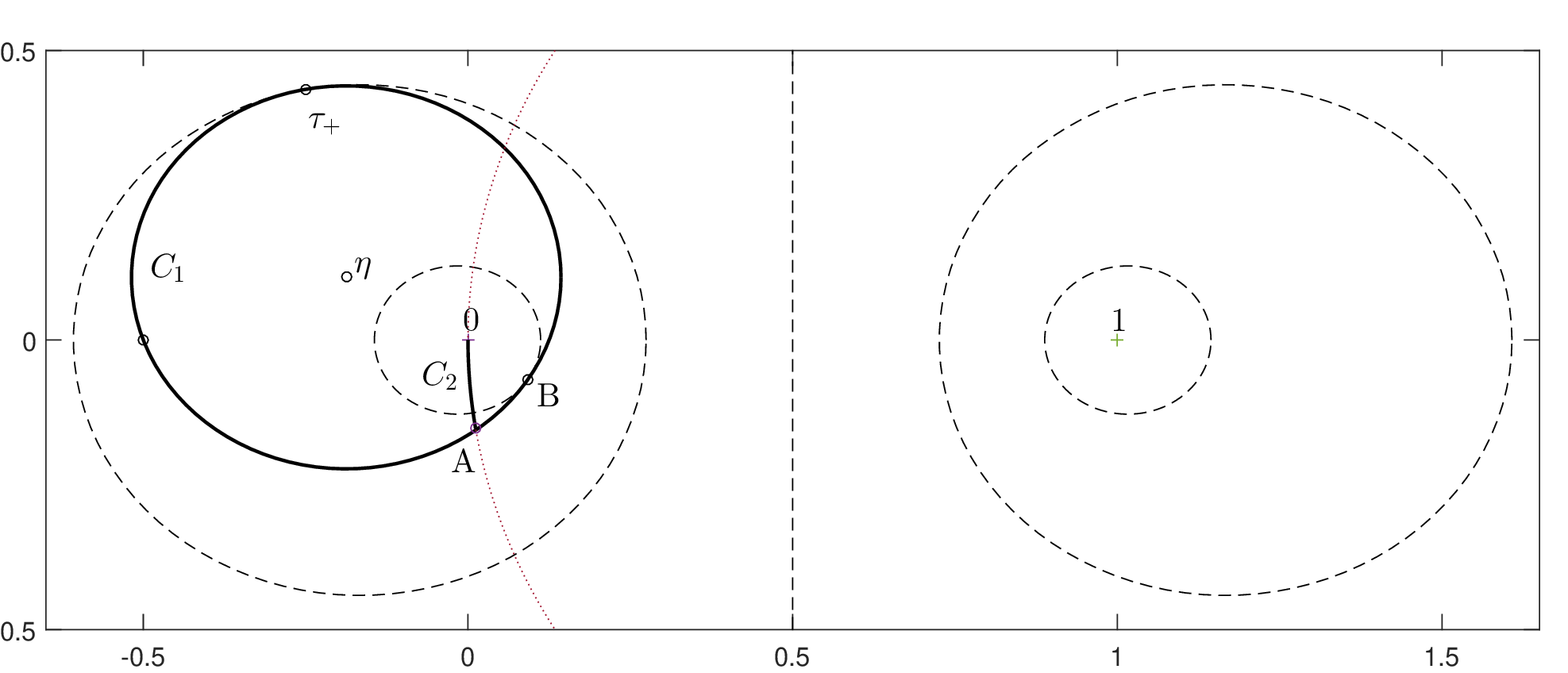}
\caption{Choice of contour $\mathcal{K}$
in situation 3. The contour is represented by the dark black circle and the
circle segment ending at $0$. The dashed lines are iso-lines of function $%
\left\vert \protect\tau /\left( \protect\tau -1\right) \right\vert .$}
\label{situation2}
\end{figure}

\begin{figure}[h]
\centering
\includegraphics[height=2.4189in,width=5.0in]{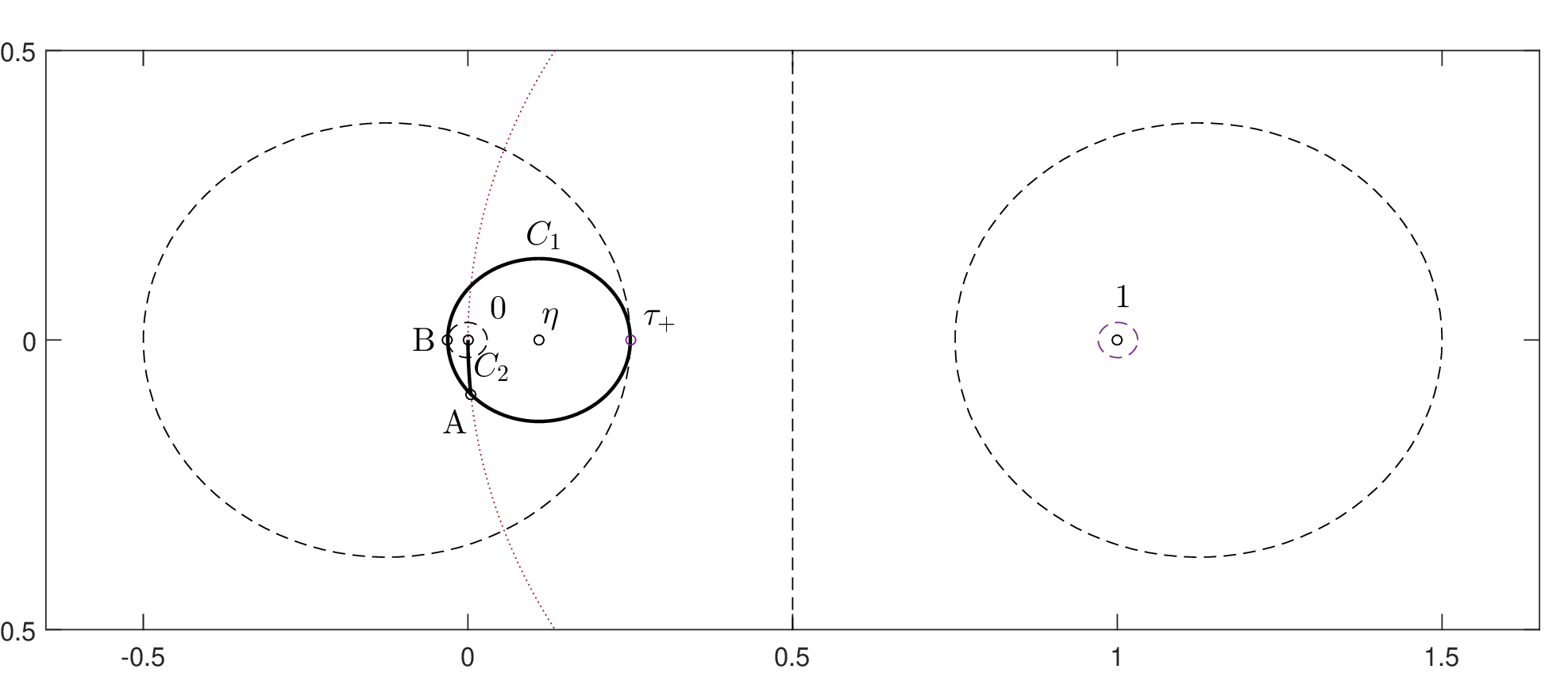}
\caption{Choice of contour $\mathcal{K}$
in situation 4. The contour is represented by the dark black circle and the
circle segment ending at $0$. The dashed lines are iso-lines of function $%
\left\vert \protect\tau /\left( \protect\tau -1\right) \right\vert .$}
\label{situation4}
\end{figure}

For situation 1 and 2, $C_{2}$ consists of a segment of the circle that
passes through $0,$ $1,$ and $\eta .$ The segment starts at the closest to $%
0 $ intersection of the latter circle with $C_{1}$ and ends at $0$. It does
not pass through $1$ or $\eta $. For situation 3 and 4, $C_{2}$ consists of
the segment of the circle with center at $1$ and radius $1$ that connects $0$
with the point $A$ of the intersection of this circle with $C_{1},$ and lies
inside $C_{1}.$ Out of the two intersection points we choose the one with
the imaginary part of the opposite sign to that of $\newim\eta $. Figures %
\ref{situation0}, \ref{situationnew}, \ref{situation2}, and \ref{situation4}
illustrate the choice of $\mathcal{K}$ for situations 1, 2, 3, and 4,
respectively.

%\FRAME{ftbpFU}{5.5858in}{2.3981in}{0pt}{\Qcb{Choice of contour $\mathcal{K}$
%in Situation 1. The contour is represented by the dark black circle and the
%circle segment ending at $0$. The dashed lines are iso-lines of function $%
%\left\vert \protect\tau /\left( \protect\tau -1\right) \right\vert .$}}{\Qlb{%
%situation0}}{situation1.eps}{\special{language "Scientific Word";type
%"GRAPHIC";maintain-aspect-ratio TRUE;display "USEDEF";valid_file "F";width
%5.5858in;height 2.3981in;depth 0pt;original-width 13.1702in;original-height
%5.636in;cropleft "0";croptop "1";cropright "1";cropbottom "0";filename
%'situation1.eps';file-properties "XNPEU";}}

%\FRAME{ftbpFU}{5.6135in}{2.4111in}{0pt}{\Qcb{Choice of contour $\mathcal{K}$
%in Situation 2. The contour is represented by the dark black circle and the
%circle segment ending at $0$. The dashed lines are iso-lines of function $%
%\left\vert \protect\tau /\left( \protect\tau -1\right) \right\vert .$}}{\Qlb{%
%situationnew}}{situationnew.eps}{\special{language "Scientific Word";type
%"GRAPHIC";maintain-aspect-ratio TRUE;display "USEDEF";valid_file "F";width
%5.6135in;height 2.4111in;depth 0pt;original-width 13.1702in;original-height
%5.636in;cropleft "0";croptop "1";cropright "1";cropbottom "0";filename
%'situationnew.eps';file-properties "XNPEU";}}

%\FRAME{ftbpFU}{5.604in}{2.2563in}{0pt}{\Qcb{Choice of contour $\mathcal{K}$
%in Situation 3. The contour is represented by the dark black circle and the
%circle segment ending at $0$. The dashed lines are iso-lines of function $%
%\left\vert \protect\tau /\left( \protect\tau -1\right) \right\vert .$}}{\Qlb{%
%situation2}}{situation2.eps}{\special{language "Scientific Word";type
%"GRAPHIC";maintain-aspect-ratio TRUE;display "USEDEF";valid_file "F";width
%5.604in;height 2.2563in;depth 0pt;original-width 14.0428in;original-height
%5.6317in;cropleft "0";croptop "1";cropright "1";cropbottom "0";filename
%'situation2.eps';file-properties "XNPEU";}}

%\FRAME{ftbpFU}{5.6342in}{2.4189in}{0pt}{\Qcb{Choice of contour $\mathcal{K}$
%in Situation 4. The contour is represented by the dark black circle and the
%circle segment ending at $0$. The dashed lines are iso-lines of function $%
%\left\vert \protect\tau /\left( \protect\tau -1\right) \right\vert .$}}{\Qlb{%
%situation4}}{situation4.eps}{\special{language "Scientific Word";type
%"GRAPHIC";maintain-aspect-ratio TRUE;display "USEDEF";valid_file "F";width
%5.6342in;height 2.4189in;depth 0pt;original-width 13.1702in;original-height
%5.636in;cropleft "0";croptop "1";cropright "1";cropbottom "0";filename
%'situation4.eps';file-properties "XNPEU";}}

We split the contour in three parts%
\begin{equation}
\mathcal{K}=\mathcal{K}_{\left[ 0,A\right] }+\mathcal{K}_{\left[ A,\tau
_{+},B\right] }+\mathcal{K}_{\left[ B,A,0\right] },  \label{splitting2j2}\tag{SM65}
\end{equation}%
or%
\begin{equation}
\mathcal{K}=\mathcal{K}_{\left[ 0,A,B\right] }+\mathcal{K}_{\left[ B,\tau
_{+},A\right] }+\mathcal{K}_{\left[ A,0\right] }  \label{splitting1j2}\tag{SM66}
\end{equation}%
depending on whether moving counter-clockwise along $C_{1}$ from $A$ to $B$
reaches $\tau _{+}$ or not. In the rest of this note, we shall refer to (\ref%
{splitting2j2}) for concreteness. Our arguments do not depend on the
specific form of the splitting.

As follows from Lemma \ref{ts1facts2F1}, $\newre\phi \left( \tau \right) $
is strictly increasing as $\tau $ is going along $\mathcal{K}_{\left[ A,\tau
_{+},B\right] }$ away from $\tau _{+}$. In other words, $\mathcal{K}_{\left[
A,\tau _{+},B\right] }$ is a contour of steep descent. Below, we shall use
Lemma JO\ref{Olver} to analyze 
\begin{equation*}
\mathcal{I}_{_{\left[ A,\tau _{+},B\right] }}=\int_{\mathcal{K}_{\left[
A,\tau _{+},B\right] }}e^{-m\phi \left( \tau \right) }\chi \left( \tau
\right) \mathrm{d}\tau .
\end{equation*}%
We shall then show that $\mathcal{I}_{\left[ 0,A\right] }$ and $\mathcal{I}_{%
\left[ B,A,0\right] }$, which are defined similarly to $\mathcal{I}_{\left[
A,\tau _{+},B\right] },$ are asymptotically dominated by $\mathcal{I}_{\left[
A,\tau _{+},B\right] }$.
\vspace{3 mm}

\textbf{\small{Saddle point approximation for $\mathcal{I}_{\left[ A,\protect%
\tau _{+},B\right] },$ CCA.}}
We now derive a saddle point approximation to the integral $\mathcal{I}%
_{\left[ A,\tau _{+},B\right] }$ which is uniform with respect $\left(
\kappa ,\eta \right) \in \Omega _{2\delta },$ where 
\begin{equation}
\Omega _{2\delta }=\left\{ \left( \kappa ,\eta \right) :\delta \leq \kappa
-1\leq \delta ^{-1},\text{ }\mathrm{dist}\left( \eta ,\mathbb{R}\backslash %
\left[ 0,1\right] \right) \geq \delta ,\text{ and }\left\vert \eta
\right\vert \leq \delta ^{-1}\right\} ,  \label{Omega2delta}\tag{SM67}
\end{equation}%
and $\delta $ is an arbitrary fixed number that satisfies inequalities $%
0<\delta <1$. Let us verify assumptions A0-A5 of Lemma JO\ref{Olver}.
For this verification, we need the following lemma.

\begin{lemma}
\label{Corollaryj2} The quantities $\left\vert \tau _{+}-\eta \right\vert $
and $\left\vert \tau _{+}\right\vert $ are bounded away from zero and
infinity, uniformly with respect to $\left( \kappa ,\eta \right) \in \Omega
_{2\delta }.$
\end{lemma}

\textbf{\small{Proof:}} The lemma follows from Lemma \ref{ts1facts2F1} (i,ii), the
fact that $\tau _{+}\neq \eta $ for $\left( \kappa ,\eta \right) \in \Omega
_{2\delta }$, and the compactness of $\Omega _{2\delta }$. $\square $

Lemma \ref{Corollaryj2} implies that the length of $\mathcal{I}_{\left[
A,\tau _{+},B\right] }$ is bounded uniformly with respect to $\left( \kappa
,\eta \right) \in \Omega _{2\delta }$. Further,%
\begin{equation*}
\sup_{\tau \in \mathcal{K}_{_{\left[ A,\tau _{+}\right] }}}\left\vert \tau
-\tau _{+}\right\vert \geq \left\vert A-\tau _{+}\right\vert \text{ and }%
\sup_{\tau \in \mathcal{K}_{_{\left[ \tau _{+},B\right] }}}\left\vert \tau
-\tau _{+}\right\vert \geq \left\vert B-\tau _{+}\right\vert ,
\end{equation*}%
where $\left\vert A-\tau _{+}\right\vert $ and $\left\vert
B-\tau _{+}\right\vert $ are continuous functions of $\left( \kappa ,\eta
\right) \in \Omega _{2\delta },$ which are not equal to zero for any $\left(
\kappa ,\eta \right) \in \Omega _{2\delta }.$ Therefore, $\left\vert A-\tau
_{+}\right\vert $ and $\left\vert B-\tau _{+}\right\vert $ are bounded away
from zero, uniformly with respect to $\left( \kappa ,\eta \right) \in \Omega
_{2\delta }$ and assumption A0 holds.

Assumptions A1, A2, A3 and A5 follow from Lemma \ref{Corollaryj2}. Finally,
let $\tau _{1}$ and $\tau _{2}$ be the points of intersection of $\mathcal{K}
$ with a circle with center at $\tau _{+}$ and a sufficiently small fixed
radius $\varepsilon _{1}$. The validity of Assumption A4 follows from the
fact that $\newre\left( \phi \left( \tau _{s}\right) -\phi \left( \tau
_{+}\right) \right) ,$ $s=1,2,$ are positive continuous functions of $\left(
\kappa ,\eta \right) \in \Omega _{2\delta }$ (the positivity being a
consequence of Lemma \ref{ts1facts2F1} (iii)) and $\newim\left( \phi
\left( \tau _{s}\right) -\phi \left( \tau _{+}\right) \right) ,$ $s=1,2,$
are continuous functions of $\left( \kappa ,\eta \right) \in \Omega
_{2\delta }$.

Since assumptions A0-A5 hold, by Lemma JO\ref{Olver}, we have%
\begin{equation*}
\mathcal{I}_{\left[ A,\tau _{+},B\right] }=2e^{-m\phi _{0}}\left[ \sqrt{\pi }%
\frac{a_{0}}{m^{1/2}}+\frac{O\left( 1\right) }{m^{3/2}}\right] ,
\end{equation*}%
where $O\left( 1\right) $ is uniform with respect to $\left( \kappa ,\eta
\right) \in \Omega _{2\delta }$, 
\begin{equation}
\phi _{0}=-\kappa \ln \frac{\tau _{+}}{1-\tau _{+}}+\left( \kappa -1\right)
\ln \left( \tau _{+}-\eta \right)  \label{phi0j2}\tag{SM68}
\end{equation}%
and%
\begin{equation}
a_{0}=\frac{\left( \tau _{+}-\eta \right) ^{-1}\left( 1-\tau _{+}\right)
^{-1}}{\sqrt{2\kappa \left( 1-2\tau _{+}\right) /\left( \left( 1-\tau
_{+}\right) ^{2}\tau _{+}^{2}\right) -2\left( \kappa -1\right) /\left( \tau
_{+}-\eta \right) ^{2}}}  \label{aoj2}\tag{SM69}
\end{equation}%
with the branch of the square root chosen as described in Lemma JO\ref{Olver}.

Precisely, let%
\begin{equation*}
\alpha =\pi /2+\arg \left( \tau _{+}-\eta \right) ,
\end{equation*}%
where the principal branch of $\arg \left( \cdot \right) $ is taken, and let 
\begin{equation*}
w=\arg \left( 2\kappa \left( 1-2\tau _{+}\right) /\left( \left( 1-\tau
_{+}\right) ^{2}\tau _{+}^{2}\right) -2\left( \kappa -1\right) /\left( \tau
_{+}-\eta \right) ^{2}\right) ,
\end{equation*}%
where the branch of $\arg \left( \cdot \right) $ is chosen so that%
\begin{equation*}
\left\vert w+2\alpha \right\vert \leq \pi /2.
\end{equation*}%
Then%
\begin{equation}
a_{0}=\frac{e^{-\mathrm{i}w/2}\left( \tau _{+}-\eta \right) ^{-1}\left(
1-\tau _{+}\right) ^{-1}}{\sqrt{\left\vert 2\kappa \left( 1-2\tau
_{+}\right) /\left( \left( 1-\tau _{+}\right) ^{2}\tau _{+}^{2}\right)
-2\left( \kappa -1\right) /\left( \tau _{+}-\eta \right) ^{2}\right\vert }}.
\label{a0adjj2}\tag{SM70}
\end{equation}
\vspace{3 mm}

\textbf{\small{Analysis of $\mathcal{I}_{\left[ 0,A\right] }$ and $\mathcal{I}_{%
\left[ B,A,0\right] },$ CCA.}}
Let us show that $\mathcal{I}_{\left[ A,\tau _{+},B\right] }$ asymptotically
dominates $\mathcal{I}_{\left[ 0,A\right] }$ and $\mathcal{I}_{\left[ B,A,0%
\right] }$ uniformly with respect to $\left( \kappa ,\eta \right) \in \Omega
_{2\delta }$. It is sufficient to prove that there exists a positive
constant $S,$ such that, for $\tau $ on $\mathcal{I}_{\left[ 0,A\right] }$
or on $\mathcal{I}_{\left[ B,A,0\right] },$ we have $\newre\phi \left(
\tau \right) \geq \newre\phi _{0}+S,$ uniformly with respect to $\left(
\kappa ,\eta \right) \in \Omega _{2\delta }.$

Note that, by Lemma \ref{ts1facts2F1} (iii), for any $\tau \in \mathcal{K}_{%
\left[ A,B\right] },$ $\newre\phi \left( \tau \right) \geq \newre\phi
\left( A\right) .$ Hence, it is sufficient to prove that $\newre\phi
\left( \tau \right) \geq \newre\phi _{0}+S$ for $\tau $ from $\mathcal{K}%
_{\left[ 0,A\right] }$. Moreover, it is sufficient to establish the fact that%
\begin{equation}
\newre\phi \left( A\right) \geq \newre\phi _{0}+S  \label{ReAj2}\tag{SM71}
\end{equation}%
uniformly with respect to $\left( \kappa ,\eta \right) \in \Omega _{2\delta
}.$ It is because for any $\tau \in \mathcal{K}_{\left[ 0,A\right] },$ $%
\newre\phi \left( \tau \right) \geq \newre\phi \left( A\right) .$

Indeed, for situations 1 and 2 this property of $\newre\phi \left( \tau
\right) $ follows from the fact that $\left\vert \tau /\left( 1-\tau \right)
\right\vert $ is strictly decreasing and $\left\vert \tau -\eta \right\vert $
is strictly increasing as $\tau $ moves along $\mathcal{K}_{\left[ 0,A\right]
}$ away from $A.$ For situation 3, we have%
\begin{eqnarray*}
\newre\phi \left( \tau \right) -\newre\phi \left( A\right) &=&-\kappa
\log \frac{\left\vert \tau \right\vert }{\left\vert A\right\vert }+\left(
\kappa -1\right) \log \frac{\left\vert \tau -\eta \right\vert }{\left\vert
A-\eta \right\vert } \\
&>&-\kappa \log \frac{\left\vert \tau \right\vert }{\left\vert \tau -\eta
\right\vert }+\kappa \log \frac{\left\vert A\right\vert }{\left\vert A-\eta
\right\vert },
\end{eqnarray*}%
where the latter inequality holds because $\left\vert \tau -\eta \right\vert
<\left\vert A-\eta \right\vert .$ The iso-lines of function $\left\vert \tau
\right\vert /\left\vert \tau -\eta \right\vert $ are similar to those shown
on Figure \ref{isolines} with the concentration points $0$ and $\eta $
instead of $0$ and $1$. As $\tau $ moves along $\mathcal{K}_{\left[ 0,A%
\right] }$ away from $A,$ the isolines are crossed so that $\left\vert \tau
\right\vert /\left\vert \tau -\eta \right\vert $ is decreasing. Therefore,%
\begin{equation}
\newre\phi \left( \tau \right) -\newre\phi \left( A\right) \geq 0.
\label{AlargerTau}\tag{SM72}
\end{equation}%
For situation 4, the analysis is more involved. We have the
following lemma.

\begin{lemma}
\label{sit2}Inequality (\ref{AlargerTau}) holds for situation 4.
\end{lemma}

\textbf{\small{Proof:}} The analysis is similar to that of situation 3. However, in
contrast to situation 3, we cannot immediately claim that as $\tau $ moves
along $\mathcal{K}_{\left[ 0,A\right] }$ away from $A,$ the isolines of the
function $\left\vert \tau \right\vert /\left\vert \tau -\eta \right\vert $
are crossed so that the function is decreasing. For this claim to be valid,
we must verify that 
\begin{equation}
\left\vert A\right\vert /\left\vert A-\eta \right\vert <1,  \label{need}\tag{SM73}
\end{equation}
so that $A$ and $0$ lie on the same side of the iso-line $\left\vert \tau \right\vert /\left\vert
\tau-\eta \right\vert =1$. 

Let $x$ be the point on $C_{1}$ where $\left\vert x\right\vert /\left\vert
x-\eta \right\vert =1,$ such that $\newim\left( \eta -x\right) $ has the
same sign as $\newim\eta $. To establish (\ref{need}), it is sufficient
to show that $x$ lies inside the circle with center at $1$ and radius $1$
(circumference of which contains $\mathcal{K}_{\left[ 0,A\right] }$). That is,
it is sufficient to show that $\left\vert 1-x\right\vert ^{2}<1.$

For concreteness, let us focus on the case $\newim\eta >0.$ Then, we have%
\begin{equation*}
x=\eta /2-\mathrm{i}\eta a,\text{ where }a=\sqrt{\left( r^{2}-\left\vert
\eta /2\right\vert ^{2}\right) /\left\vert \eta \right\vert ^{2}}
\end{equation*}%
and $r$ is the radius of $C_{1}.$ A straightforward algebra shows that%
\begin{equation*}
\left\vert 1-x\right\vert ^{2}=r^{2}+1-\newre\eta -2a\newim\eta .
\end{equation*}%
Furthermore, since $r^{2}\geq \left\vert \eta \right\vert ^{2},$ we have $a>%
\sqrt{3}/2>1.$ Therefore, the inequality $\left\vert 1-x\right\vert ^{2}<1$ 
would follow from the inequality 
$r^{2}\leq \left\vert \eta \right\vert <\newre\eta +\newim\eta$ .
Let us now show that in situation 4, $r^{2}\leq \left\vert \eta \right\vert$.

Let $z=\rho \exp \left\{ \mathrm{i}\theta \right\} ,$ where $\rho $ and $%
\theta $ are as in (\ref{polar}). Situation 4 imposes the following
constraints on $z$: 1) $\newre z\geq 0,$ 2) $\newre z\leq \kappa -1,$ 3) 
$({\newre z})^{2}-({\newim z})^{2}\geq 1$. The first one is equivalent to $%
\theta \in \left[ -\pi /2,\pi /2\right] $, which must be true by definition of $\theta$.
The second one is equivalent to (\ref{34}), and the last one ensures that $%
\newre\eta \geq 0$. We have%
\begin{equation*}
\tau _{+}=\frac{-1+z}{2\left( \kappa -1\right) }\text{ and }\eta =\frac{%
-1+z^{2}}{4\kappa \left( \kappa -1\right) }.
\end{equation*}%
Therefore,%
\begin{equation*}
r^{2}\equiv \left\vert \tau _{+}-\eta \right\vert ^{2}=\frac{\left\vert
z-1\right\vert ^{2}\left\vert z-\left( 2\kappa -1\right) \right\vert ^{2}}{%
16\kappa ^{2}\left( \kappa -1\right) ^{2}}
\end{equation*}%
and%
\begin{equation*}
\left\vert \eta \right\vert =\frac{\left\vert z-1\right\vert \left\vert
z+1\right\vert }{4\kappa \left( \kappa -1\right) }.
\end{equation*}%
For $z$ that satisfies the above three constraints, we must have $\left\vert
z+1\right\vert >\left\vert z-1\right\vert $. Therefore, to establish
inequality $r^{2}\leq \left\vert \eta \right\vert ,$ it is sufficient to
show that%
\begin{equation*}
\frac{\left\vert z-\left( 2\kappa -1\right) \right\vert ^{2}}{4\kappa \left(
\kappa -1\right) }\leq 1.
\end{equation*}%
The latter inequality is equivalent to%
\begin{equation*}
({\newim z})^{2}+1\leq 2\newre z\left( 2\kappa -1\right) -({\newre z})^{2}.
\end{equation*}%
In view of the third constraint, it is sufficient to show that%
\begin{equation*}
2\newre z\left( 2\kappa -1\right) -({\newre z})^{2}\geq ({\newre z})^{2}.
\end{equation*}%
But the second constraint implies this inequality. The situation where $%
\newim\eta \leq 0$ is analyzed similarly. $\square $

To summarize, in all the four situations we only need to show that (\ref{ReAj2})
holds. Note that $\newre\phi \left( A\right) -\newre\phi _{0}$ is not
a continuous function of $\left( \kappa ,\eta \right) \in \Omega _{2\delta }$
because we may have $A=B=0$ and $\newre\phi \left( A\right) =+\infty $
for some $\left( \kappa ,\eta \right) \in \Omega _{2\delta }$. However, we
can bound $\newre\phi \left( A\right) -\newre\phi _{0}$ from below by
the minimum of two positive continuous functions $\newre\phi \left( \tau
_{1}\right) -\newre\phi _{0}$ and $\newre\phi \left( \tau _{2}\right) -%
\newre\phi _{0},$ where $\tau _{1}$ and $\tau _{2}$ are points of the
intersection of the circle with center $\eta $ and radius $\left\vert \tau
_{+}-\eta \right\vert $ and a circle with center $\tau _{+}$ and a fixed
radius, which is smaller than $\left\vert A-\tau _{+}\right\vert $,
uniformly with respect to $\left( \kappa ,\eta \right) \in \Omega _{2\delta
} $. Therefore, there exists $S>0$ such that (\ref{ReAj2}) holds uniformly
with respect to $\left( \kappa ,\eta \right) \in \Omega _{2\delta }.$
\vspace{3 mm}

\textbf{\small{Asymptotics in terms of $\protect\varphi $ and $\protect\psi $,
CCA.}}

The above analysis implies the following asymptotic representation%
\begin{equation*}
F_{2} =\frac{C_{m}\eta ^{-m}}{2\pi \mathrm{i}}2e^{-m\phi _{0}} 
\left[ \frac{\sqrt{\pi } e^{-\mathrm{i}w/2}\left( \tau _{+}-\eta
\right) ^{-1}\left( 1-\tau _{+}\right) ^{-1}}{m^{1/2}\sqrt{\left\vert
2\kappa \left( 1-2\tau _{+}\right) /\left( \left( 1-\tau _{+}\right)
^{2}\tau _{+}^{2}\right) -2\left( \kappa -1\right) /\left( \tau _{+}-\eta
\right) ^{2}\right\vert }}+\frac{O\left( 1\right) }{m^{3/2}}\right] ,
\end{equation*}%
where $O\left( 1\right) $ is uniform with respect to $\left( \kappa ,\eta
\right) \in \Omega _{2\delta }.$ We would like to express this formula in
terms of $t_{2},\varphi _{2}\left( \cdot \right) ,$ and $\psi _{2}\left(
\cdot \right) .$ Since%
\begin{equation*}
\left\vert 2\kappa \left( 1-2\tau _{+}\right) /\left( \left( 1-\tau
_{+}\right) ^{2}\tau _{+}^{2}\right) -2\left( \kappa -1\right) /\left( \tau
_{+}-\eta \right) ^{2}\right\vert =\left\vert 2\varphi _{2}^{\prime \prime
}\left( t_{2}\right) /\eta ^{2}\right\vert ,
\end{equation*}%
\begin{equation*}
\varphi _{2}\left( t_{2}\right) =\phi _{0}-\ln \eta ,
\end{equation*}%
and%
\begin{equation*}
\left( \tau _{+}-\eta \right) ^{-1}\left( 1-\tau _{+}\right) ^{-1}=\psi
_{2}\left( t_{2}\right) \eta ^{-1},
\end{equation*}%
we have%
\begin{equation*}
F_{2}=C_{m}e^{-m\varphi _{2}\left( t_{2}\right) }\left[ \frac{e^{-\mathrm{i}%
w/2}e^{-\mathrm{i\arg }\eta }}{\mathrm{i}}\frac{\psi _{2}\left( t_{2}\right) 
}{\sqrt{\left\vert 2\pi m\varphi _{2}^{\prime \prime }\left( t_{2}\right)
\right\vert }}+\frac{O\left( 1\right) }{m^{3/2}}\right] .
\end{equation*}%
On the other hand, by definition,%
\begin{equation*}
w=\arg \left( \varphi _{2}^{\prime \prime }\left( t_{2}\right) \right)
-2\arg \eta =\omega _{2}-\pi -2\arg \eta ,
\end{equation*}%
where $\omega _{2}$ is as defined in equation (JO\ref{omegaj}). Therefore,%
\begin{eqnarray*}
F_{2} &=&C_{m}e^{-m\varphi _{2}\left( t_{2}\right) }e^{-\mathrm{i}\omega
_{2}/2}\left[ \frac{\psi _{2}\left( t_{2}\right) }{\sqrt{\left\vert 2\pi
m\varphi _{2}^{\prime \prime }\left( t_{2}\right) \right\vert }}+\frac{%
O\left( 1\right) }{m^{3/2}}\right] \\
&=&C_{m}\psi _{2}\left( t_{2}\right) e^{-\mathrm{i}\omega _{2}/2}\left\vert
2\pi m\varphi _{2}^{\prime \prime }\left( t_{2}\right) \right\vert
^{-1/2}\exp \left\{ -m\varphi _{2}\left( t_{2}\right) \right\} \left(
1+o(1)\right) .
\end{eqnarray*}

\subsection{Proof of Confluences} \label{sec:proof-confluences}

% \section{Proof of confluences}
% \label{sec:proof-confluences}

The confluences (JO\ref{eq:c2to0}) are established by showing convergence of each of
the components in (JO\ref{eq:fdecomp}).
For the $f_\mathrm{c}$ and $f_\mathrm{e}$ components, this follows from inspection of
Tables JO\ref{Table 2a} and JO\ref{Table 3} respectively, while for 
$f_\mathrm{h}^{\mathrm{SigD}}(z)$, this follows from
  (JO\ref{elementary}).
For $f_\mathrm{h}^{\mathrm{REG}}(z)$ and
$f_\mathrm{h}^{\mathrm{CCA}}(z)$, one uses the definitions of
$\varphi_{j}$ and $t_{j}$ and calculation, though one can also appeal to the
confluences
\begin{align*}
%  _{1}F_{1} \, (\frac{a}{\epsilon};b;\epsilon x) 
  _{1}F_{1} \, (a/\epsilon;b;\epsilon x) 
     & \to \ _{0}F_{1} \, (b; ax) \\
  _{2}F_{1} \, (a_1/\epsilon, a_2/\epsilon ;b;\epsilon^2 x) 
     & \to \ _{0}F_{1} \, (b; a_1 a_2 x)
\end{align*}
as $\epsilon \rightarrow 0$,
and observe in (JO\ref{pFq}) that with $\kappa \sim c_1/((1-c_1) c_2)$ and 
as $c_2 \to 0$, we have
\begin{align*}
  (m \kappa + 1) m \eta_1 & \to m^2 \eta_0  \\
   (m \kappa + 1)^2 \eta_2 & \to m^2 \eta_0.
\end{align*}

For the confluences (JO\ref{eq:c1to0}), there is some crosstalk
between the components. 
We write $f_\mathrm{c}[\theta]$ to show the dependence on $\theta$ explicitly.
Writing
$\theta = \sqrt c_1 \xi$, it is
direct to verify that 
\begin{align*}
  f_\mathrm{c}^{\mathrm{PCA}}[ \sqrt c_1 \xi]
      & = f_\mathrm{c}^{\mathrm{SMD}}[\xi] + \xi/\sqrt c_1 
           - \xi^2 - \log \sqrt c_1 + O(\sqrt c_1) \\
  f_\mathrm{c}^{\mathrm{REG}_0}[ \sqrt c_1 \xi]
      & = f_\mathrm{c}^{\mathrm{SMD}}[\xi] + \xi/\sqrt c_1 
           - \xi^2/2 - \log \sqrt c_1 + O(\sqrt c_1).
\end{align*}
From (JO\ref{f1Case}) and the MP entry in Table JO\ref{Table 3}, and
writing $z = 1 + \sqrt c_1 w$, we have
\begin{equation*}
  f_\mathrm{e}^{\mathrm{PCA}}(1+\sqrt c_1 w)
 =f_\mathrm{e}^{\mathrm{REG}_0}(1+\sqrt c_1 w)
 =f_\mathrm{e}^{\mathrm{SMD}}(w) + \log \sqrt c_1 + o(1).
\end{equation*}

For the $\mathrm{h}$ term, we write $f_\mathrm{h}(z; \theta)$ to show
the dependence on $\theta$ explicitly.
From (JO\ref{elementary}), one quickly has
\begin{equation*}
  f_\mathrm{h}^{\mathrm{PCA}}(1+\sqrt c_1 w; \sqrt c_1 \xi)
   = f_\mathrm{h}^{\mathrm{SMD}}(w; \xi) - \xi/\sqrt c_1 
           + \xi^2 + O(\sqrt c_1).
\end{equation*}
For $f_\mathrm{h}^{\mathrm{REG}_0}$, we have
$\varphi_0(t_0) = \log t_0 - 2(t_0 -1)$ 
and that $t_0 = 1 + \eta_0 - \eta_0^2 + O(\eta_0^3)$ for small
$\eta_0$. 
This leads to $f_{\mathrm{h}}(z;\theta) = (1-c_1)c_1^{-1}[-\eta_0 + \tfrac{1}{2}
\eta_0^2 + O(\eta_0^3)]$ and thence by elementary evaluation to 
\begin{equation*}
  f_\mathrm{h}^{\mathrm{REG}_0}(1+\sqrt c_1 w; \sqrt c_1 \xi)
   = f_\mathrm{h}^{\mathrm{SMD}}(w; \xi) - \xi/\sqrt c_1 
           + \xi^2/2 + O(\sqrt c_1).
\end{equation*}
Combining terms from the preceding displays yields the confluences
(JO\ref{eq:c1to0}).

\subsection{Proof of Lemma JO\ref{critical} (saddle points $z_{0}$)} \label{saddle points z0} 

\textit{$\mathsf{q}=0$ cases: (SMD, PCA, SigD).}

First, note that 
\begin{equation}
f^{\prime }(z)=f_{\mathrm{e}}^{\prime }(z)+f_{\mathrm{h}}^{\prime }(z)=-m_{%
\mathbf{c}}(z)+f_{\mathrm{h}}^{\prime }(z),  \label{eq:derivatives}\tag{SM74}
\end{equation}%
where $m_{\mathbf{c}}(z)$ is the appropriate Stieltjes transform. We
proceed, then, by solving for $z$ in the equation $f_{\mathrm{h}}^{\prime
}(z)=m_{\mathbf{c}}(z)$.
\vspace{3 mm}

\textbf{\small{SMD.}}
We have $f_{\mathrm{h}}^{\prime }(z)=-\theta $, so substituting $m_{\mathbf{c%
}}(z_{0})=-\theta $ into the quadratic equation 
\begin{equation}
m^{2}+zm+1=0  \label{eq:mq-SMD}\tag{SM75}
\end{equation}%
satisfied by $m=m_{\mathbf{c}}^{\mathrm{SC}}(z)$, we get 
\begin{equation*}
z_{0}(\theta )=-\frac{m^{2}+1}{m}=\frac{\theta ^{2}+1}{\theta }=\theta
+1/\theta .
\end{equation*}%
Obviously, for any $\theta \in \left( 0,\bar{\theta}^{\mathrm{SMD}}\right)
\equiv \left( 0,1\right) ,$ $z_{0}(\theta )$ is larger than $b_{+}^{\mathrm{SMD%
}}=\beta _{+}^{\mathrm{SMD}}\equiv 2.$
\vspace{3 mm}

\textbf{\small{PCA.}}
Now $f_{\mathrm{h}}^{\prime }(z)=-\theta /[c_{1}(1+\theta )]$, so we
substitute $m_{\mathbf{c}}(z_{0})=-\theta /[c_{1}(1+\theta )]$ into the
quadratic equation 
\begin{equation}
c_{1}zm^{2}+(z+c_{1}-1)m+1=0  \label{eq:mq-PCA}\tag{SM76}
\end{equation}%
satisfied by $m=m_{\mathbf{c}}^{\mathrm{MP}}(z)$. This is a linear equation
for $z$ whose solution is 
\begin{equation*}
z_{0}(\theta )=(\theta +1)(\theta +c_{1})/\theta .
\end{equation*}%
Note that the minimum of $z_{0}(\theta )$ over $\theta >0$ equals $b_{+}^{%
\mathrm{PCA}}\equiv \left( 1+\sqrt{c_{1}}\right) ^{2}$ and is achieved at 
\begin{equation*}
\theta =\bar{\theta}_{\mathbf{c}}\equiv \sqrt{c_{1}}.
\end{equation*}%
Therefore, since $m_{\mathbf{c}}^{\mathrm{MP}}(z)$ is well defined for $%
z>b_{+}^{\mathrm{PCA}}$, 
%and since $\bar{\theta}_{p}\rightarrow \bar{\theta}$
%as $\mathbf{n},p\rightarrow _{\boldsymbol{\gamma} }\infty ,$ 
$m_{\mathbf{c}}^{\mathrm{MP}%
}(z_{0})$ must be well defined for any $\theta \in \left( 0,\bar{\theta}_{\mathbf{c}}%
\right) $.
%, for sufficiently large $\mathbf{n},p.$
\vspace{3 mm}

\textbf{\small{SigD.}}
The Stieltjes transform $m=m_{\mathbf{c}}^{\mathrm{W}}(z)$ of the Wachter
distribution, as normalized here, satisfies the quadratic equation 
\begin{equation}
c_{1}z(c_{1}-c_{2}z)m^{2}+[c_{1}(1-c_{2})z-(1-c_{1})(c_{1}-c_{2}z)]m+r^{2}=0,
\label{eq:m-W}\tag{SM77}
\end{equation}%
while 
\begin{equation}
f_{\mathrm{h}}^{\prime }(z)=\frac{ab}{1-az},\qquad a=\frac{c_{2}\theta }{%
c_{1}(1+\theta )},\qquad b=-\frac{r^{2}}{c_{1}c_{2}}.  \label{eq:fhp}\tag{SM78}
\end{equation}%
To solve $m(z)=f_{\mathrm{h}}^{\prime }(z)$, insert $m=ab/(1-az)$ into %
\eqref{eq:m-W} to obtain an apparently quadratic equation. However the
coefficient of $z^{2}$ vanishes, so that as with SMD and PCA, $z_{0}$ is the
solution of a linear equation $\beta z+\gamma =0$, where in this case 
\begin{align*}
\beta & =abc_{1}l(\theta )/(1+\theta ) \\
\gamma & =-abc_{1}(c_{1}+\theta )/\theta ,
\end{align*}%
so that 
\begin{equation}
z_{0}(\theta )=\frac{(c_{1}+\theta )(1+\theta )}{\theta l(\theta )}.
\label{z0SD}\tag{SM79}
\end{equation}%
It also follows that 
\begin{align}
az_{0}& =\frac{c_{2}(c_{1}+\theta )}{c_{1}l(\theta )},\qquad 1-az_{0}=\frac{%
r^{2}}{c_{1}l(\theta )}  \notag \\
f_{\mathrm{h,SigD}}^{\prime }(z_{0})& =m(z_{0})=-\frac{\theta l(\theta )}{%
c_{1}(1+\theta )}.  \label{eq:m-wa}\tag{SM80}
\end{align}%
Recall that $l(\theta )=1+(1+\theta )c_{2}/c_{1}.$ Therefore, (\ref{z0SD})
implies that that the minimum of $z_{0}(\theta )$ over $\theta >0$ equals 
\begin{equation*}
b_{+}^{\mathrm{SigD}}\equiv c_{1}\left( \frac{r+1}{r+c_{2}}\right) ^{2}
\end{equation*}%
and is achieved at 
\begin{equation*}
\theta =\bar{\theta}_{\mathbf{c}}\equiv \frac{c_{2}+r}{1-c_{2}}.
\end{equation*}%
Therefore, $m_{\mathbf{c}}^{\mathrm{W}}(z)$ is well defined for $z>b_{+}^{%
\mathrm{SigD}}$,
% and since $\bar{\theta}_{p}\rightarrow \bar{\theta}$
as $\mathbf{n},p\rightarrow _{\boldsymbol{\gamma} }\infty ,$ 
and $m_{\mathbf{c}}^{\mathrm{W}%
}(z_{0})$ must be well defined for any $\theta \in \left( 0,\bar{\theta}_{\mathbf{c}}%
\right) $.
%, for sufficiently large $\mathbf{n},p.$ 
\vspace{3 mm}

\textit{$\mathsf{q}=1$ cases: (REG$_{0}$, REG, CCA).}

We find the critical points $z_{0}(\theta )$ for the $\mathsf{q}=1$ cases by
showing that they are the same as for the corresponding $\mathsf{q}=0$
cases. This is cast as a verification rather than a derivation as we still
lack a good explanation for this curious fact.

We have seen, based on \eqref{eq:derivatives}, that%
\begin{equation*}
f_{\mathrm{h,PCA}}^{\prime }(z_{0})=m_{\mathbf{c}}^{\mathrm{MP}%
}(z_{0}),\qquad f_{\mathrm{h,SigD}}^{\prime }(z_{0})=m_{\mathbf{c}}^{\mathrm{%
W}}(z_{0}),
\end{equation*}%
for $z_{0}=z_{0}^{\mathrm{PCA}}$ and $z_{0}^{\mathrm{SigD}}$ respectively.
We now show that 
\begin{equation*}
f_{\mathrm{h,REG_{0}}}^{\prime }(z_{0})=f_{\mathrm{h,PCA}}^{\prime
}(z_{0}),\qquad f_{\mathrm{h,REG}}^{\prime }(z_{0})=f_{\mathrm{h,CCA}%
}^{\prime }(z_{0})=f_{\mathrm{h,SigD}}^{\prime }(z_{0})
\end{equation*}%
for $z_{0}=z_{0}^{\mathrm{PCA}}$ and $z_{0}^{\mathrm{SigD}}$ respectively.
In combination with \eqref{eq:derivatives}, this verifies that $z_{0}^{%
\mathrm{PCA}}$ and $z_{0}^{\mathrm{SigD}}$ are critical points for the $%
\mathsf{q}=1$ cases as well.

The functions defined in (JO\ref{phi0}) and (JO\ref{phij}) will sometimes be written in
the form $\varphi _{j}(t,\eta _{j})$ to show the dependence on $\eta _{j}$
explicitly. We have 
\begin{equation*}
f_{\mathrm{h}}(z)=\frac{1-c_{1}}{c_{1}}[\varphi _{j}(t_{j},\eta _{j})+\gamma
_{j}],
\end{equation*}%
where $\gamma _{j}=\gamma _{j}(\mathbf{n},p)$ and $t_{j}=t_{j}(\eta _{j})$
satisfies 
\begin{equation}
\frac{\partial }{\partial t}\varphi _{j}(t,\eta _{j})=0,
\label{eq:criticalpt}\tag{SM81}
\end{equation}%
a quadratic equation for $t_{j}$ with coefficients depending on $\eta _{j}$
and $\kappa $. We therefore have, dropping the subscript $j$ temporarily, 
\begin{equation}
f_{\mathrm{h}}^{\prime }(z)=\frac{1-c_{1}}{c_{1}}\frac{\mathrm{d}}{\mathrm{d}%
\eta }\varphi (t(\eta ),\eta )\frac{\mathrm{d}\eta }{\mathrm{d}z}=\frac{%
1-c_{1}}{c_{1}}\frac{\partial }{\partial \eta }\varphi (t(\eta ),\eta )\frac{%
\mathrm{d}\eta }{\mathrm{d}z}  \label{eq:partials}\tag{SM82}
\end{equation}%
From definitions (JO\ref{phi0}) and (JO\ref{phij}), again with $t_{j}=t_{j}(\eta _{j})$, and 
$\kappa =r^{2}/[c_{2}(1-c_{1})]$, 
\begin{equation}
\frac{\partial }{\partial \eta _{j}}\varphi (t_{j},\eta _{j})=%
\begin{cases}
-1/t_{0} \\ 
-t_{1} \\ 
-\kappa t_{2}/(1-\eta _{2}t_{2}).%
\end{cases}
\label{eq:t_j}\tag{SM83}
\end{equation}%
We now turn to the specifics of the three cases.
\vspace{3 mm}

\textbf{\small{REG$_{0}$.}}
We show that $z=z_{0}^{\mathrm{PCA}}=(\theta +1)(\theta +c_{1})/\theta $
solves 
\begin{equation*}
f_{\mathrm{h,REG_{0}}}^{\prime }(z)=m(z_{0})=-\frac{\theta }{c_{1}(\theta +1)%
}.
\end{equation*}%
From \eqref{eq:partials} and \eqref{eq:t_j}, 
\begin{equation*}
f_{\mathrm{h,REG_{0}}}^{\prime }(z)=-\frac{\theta }{c_{1}(1-c_{1})}\frac{1}{%
t_{0}(\eta )},
\end{equation*}%
so that we should solve $t_{0}(\eta _{0}(z))=(\theta +1)/(1-c_{1})$ for $z$.
Since $t_{0}$ satisfies a quadratic equation, the equation for $z$ becomes 
\begin{equation*}
\eta _{0}(z)=t_{0}^{2}-t_{0}=\frac{(\theta +1)(\theta +c_{1})}{(1-c_{1})^{2}}%
,
\end{equation*}%
which implies that $z_{0}^{\mathrm{REG_{0}}}=\theta ^{-1}(1-c_{1})^{2}\eta
_{0}=(\theta +1)(\theta +c_{1})/\theta =z_{0}^{\mathrm{PCA}}(\theta )$.
\vspace{3 mm}

\textbf{\small{REG.}}
This time we solve for $z$ in 
\begin{equation*}
f_{\mathrm{h,REG}}^{\prime }(z)=f_{\mathrm{h,SigD}}^{\prime }(z_{0})=-\frac{%
\theta l(\theta )}{c_{1}(1+\theta )},
\end{equation*}%
where the second equality uses \eqref{eq:m-wa}. From \eqref{eq:partials} and %
\eqref{eq:t_j} we have $f_{\mathrm{h,REG}}^{\prime }(z)=-c_{2}\theta
t_{1}(\eta )/c_{1}^{2}$, and so 
\begin{equation}
t_{1}(\eta )=\frac{c_{1}l(\theta )}{c_{2}(1+\theta )},\qquad t_{1}(\eta )-1=%
\frac{c_{1}}{c_{2}(1+\theta )}.  \label{eq:t1}\tag{SM84}
\end{equation}%
The quadratic equation for $t_{1}$ is $\eta _{1}t_{1}^{2}+(1-\eta
_{1})t_{1}-\kappa =0$, so that 
\begin{equation}
\eta _{1}=\frac{\kappa -t_{1}}{t_{1}(t_{1}-1)}=\frac{c_{2}(\theta +1)(\theta
+c_{1})}{(1-c_{1})c_{1}l(\theta )}  \label{eq:eta1}\tag{SM85}
\end{equation}%
which implies that $z_{0}^{\mathrm{REG}}=c_{1}(1-c_{1})\eta _{1}/(\theta
c_{2})=z_{0}^{\mathrm{SigD}}(\theta )$.
\vspace{3 mm}

\textbf{\small{CCA.}}
Treat this as a modification of REG. Thus 
\begin{align}
\varphi _{2}(t)& =\kappa \log (1-\eta _{2}t)+\eta _{1}t+\varphi
_{1}(t),\qquad \text{and}  \notag \\
\varphi _{2}^{\prime }(t)& =-\frac{\kappa \eta _{2}}{1-\eta _{2}t}+\eta
_{1}+\varphi _{1}^{\prime }(t).  \label{eq:phi2p}\tag{SM86}
\end{align}%
We verify that at $z=z_{0}^{\mathrm{SigD}}$, 
\begin{equation}
t_{2}=\frac{c_{1}l(\theta )}{c_{2}(1+\theta )}=t_{1}(\eta _{1}(z_{0}))
\label{eq:t2form}\tag{SM87}
\end{equation}%
satisfies $\varphi _{2}^{\prime }(t_{2})=0$ for $\eta _{2}=\eta _{2}(z_{0})$%
. Indeed, writing $L(\theta )=c_{1}l(\theta )$, we have 
\begin{equation}
\eta _{2}=\frac{c_{2}^{2}(c_{1}+\theta )(1+\theta )}{L^{2}(\theta )},\qquad
t_{2}\eta _{2}=\frac{c_{2}(c_{1}+\theta )}{L(\theta )},\qquad 1-t_{2}\eta
_{2}=\frac{r^{2}}{L(\theta )},  \label{eq:t2eta2}\tag{SM88}
\end{equation}%
and 
\begin{equation*}
\frac{\kappa }{1-t_{2}\eta _{2}}=\frac{L(\theta )}{(1-c_{1})c_{2}}=\frac{%
\eta _{1}}{\eta _{2}},
\end{equation*}%
so that from \eqref{eq:phi2p}, $\varphi _{2}^{\prime }(t_{2})=\varphi
_{1}^{\prime }(t_{1})=0$. But now we can see that, at $z=z_{0}$, 
\begin{align*}
f_{\mathrm{h,CCA}}^{\prime }(z_{0})& =-\frac{1-c_{1}}{c_{1}}\frac{\kappa }{%
1-\eta _{2}t_{2}}\cdot t_{2}\cdot \frac{\mathrm{d}\eta _{2}}{\mathrm{d}z} \\
& =-\frac{1-c_{1}}{c_{1}}\frac{\eta _{1}}{\eta _{2}}\cdot t_{1}\cdot \frac{%
\mathrm{d}\eta _{2}}{\mathrm{d}z}=-\frac{1-c_{1}}{c_{1}}\cdot t_{1}\cdot 
\frac{\mathrm{d}\eta _{1}}{\mathrm{d}z}=f_{\mathrm{h,REG}}^{\prime }(z_{0}),
\end{align*}%
so that $z_{0}$ also satisfies $f_{\mathrm{CCA}}^{\prime }(z_{0})=0$.

\subsection{Verification of Remark JO\ref{fz00}: that $f(z_{0})=0$} \label{sec: verification of remark 5}

Recall that $f(z_{0})=f_{\mathrm{c}}+f_{\mathrm{e}}(z_{0})+f_{\mathrm{h}}(z_{0})$. The term 
$f_{\mathrm{c}}$ is given in Table \ref{tab:JO3}. The next term,
\[
f_{\mathrm{e}}\left(  z_{0}\right)  =\int_{b_{-}}^{b_{+}}\ln\left(  z_{0}%
-\lambda\right)  \mathrm{d}F_{\mathbf{c}}\left(  \lambda\right),
\]
takes on three different values: one for SMD, another for PCA and REG$_{0},$
and the third one for SigD, REG, and CCA.

\begin{lemma}
\label{Lemmaf0}For SigD, REG, and CCA, for any $\theta\in\left(  0,\bar
{\theta}\right)  $ and for sufficiently large $\mathbf{n},p,$ we have%
\begin{equation}
f_{\mathrm{e}}\left(  z_{0}\right)  =2\ln c_{1}-\ln\theta-\frac{1-c_{1}}{c_{1}}%
\ln\left(  1+\theta\right)  -\frac{c_{1}+c_{2}}{c_{1}c_{2}}\ln\left(
c_{1}+c_{2}\right)  +\frac{r^{2}}{c_{1}c_{2}}\ln\left[  c_{1}l\left(
\theta\right)  \right]  . \label{formula1}\tag{SM89}%
\end{equation}

\end{lemma}

\textbf{\small{Proof:}} We follow the usual strategy of reduction to a contour
integral. First make the change of variables $\lambda=\alpha-\beta\cos
\varphi.$ In order to arrange that $\lambda=b_{-}$ and $b_{+}$ at $\varphi=0$
and $\pi$ respectively, we set%
\begin{equation}
\alpha=\frac{b_{+}+b_{-}}{2}=\frac{c_{1}\left(  r^{2}+c_{1}^{2}\right)
}{\left(  c_{1}+c_{2}\right)  ^{2}},\qquad\beta=\frac{b_{+}-b_{-}}{2}%
=\frac{2rc_{1}^{2}}{\left(  c_{1}+c_{2}\right)  ^{2}}. \label{alphabeta}\tag{SM90}%
\end{equation}
We obtain%
\[
f_{\mathrm{e}}\left(  z_{0}\right)  =\frac{c_{1}+c_{2}}{4\pi c_{1}}\int_{0}^{2\pi
}\frac{\beta^{2}\sin^{2}\varphi\ln\left(  z_{0}-\alpha+\beta\cos
\varphi\right)  }{\left(  \alpha-\beta\cos\varphi\right)  \left(  c_{1}%
-c_{2}\alpha+c_{2}\beta\cos\varphi\right)  }\mathrm{d}\varphi
\]
after extending the integral from $\left[  0,\pi\right]  $ to $\left[
0,2\pi\right]  $ using the symmetry of the integrand about $\varphi=\pi$. Now
introduce $z=e^{\mathrm{i}\varphi}.$ Since $\cos\varphi=\left(  z+z^{-1}%
\right)  /2,$ we have from (\ref{alphabeta}) the factorizations%
\begin{align*}
c_{1}\left(  \alpha-\beta\cos\varphi\right)   &  =\frac{\beta}{2r}\left(
r-c_{1}z\right)  \left(  r-c_{1}z^{-1}\right)  ,\\
c_{1}-c_{2}\alpha+c_{2}\beta\cos\varphi &  =\frac{\beta}{2r}\left(
r+c_{2}z\right)  \left(  r+c_{2}z^{-1}\right)  ,\\
z_{0}-\alpha+\beta\cos\varphi &  =q(z)q\left(  z^{-1}\right)  \text{ with}\\
q\left(  z\right)   &  =\frac{c_{1}}{c_{1}+c_{2}}\left(  \sqrt{c_{1}l\left(
\theta\right)  /\theta}+rz\sqrt{\theta/\left[  c_{1}l\left(  \theta\right)
\right]  }\right)  .
\end{align*}
Our integral becomes
\[
f_{\mathrm{e}}\left(  z_{0}\right)  =\frac{-\left(  c_{1}+c_{2}\right)  r^{2}}%
{4\pi\mathrm{i}}\int_{\mathcal{C}}\frac{\left(  z-z^{-1}\right)  ^{2}%
\ln\left(  q(z)q\left(  z^{-1}\right)  \right)  }{\left(  r-c_{1}z\right)
\left(  r-c_{1}z^{-1}\right)  \left(  r+c_{2}z\right)  \left(  r+c_{2}%
z^{-1}\right)  }\frac{\mathrm{d}z}{z}.
\]
The integral has form $I=
{\displaystyle\oint}
\ln\left(  q(z)q\left(  z^{-1}\right)  \right)  H\left(  z\right)
z^{-1}\mathrm{d}z$ with $H(z)=H\left(  z^{-1}\right)  $. Hence, expanding the
logarithm yields two identical terms, so that%
\[
f_{\mathrm{e}}\left(  z_{0}\right)  =\frac{-\left(  c_{1}+c_{2}\right)  }%
{2\pi\mathrm{i}}\int_{\mathcal{C}}\frac{\left(  z^{2}-1\right)  ^{2}\ln
q(z)}{\left(  r-c_{1}z\right)  \left(  z-c_{1}/r\right)  \left(
r+c_{2}z\right)  \left(  z+c_{2}/r\right)  }\frac{\mathrm{d}z}{z}.
\]
For $\theta\in\left(  0,\bar{\theta}\right)  $ and sufficiently large
$\mathbf{n},p,$ we have $\theta\in\left(  0,\bar{\theta}_{p}\right)  $ with
$\bar{\theta}_{p}=\left(  c_{2}+r\right)  /\left(  1-c_{2}\right)  .$ On the
other hand, for $\theta\in\left(  0,\bar{\theta}_{p}\right)  ,$ the function
$\ln q\left(  z\right)  $ is analytic inside the circle $\left\vert
z\right\vert =1,$ and so the whole integrand is analytic inside the circle
except for simple poles at $z=0,c_{1}/r$ and $-c_{2}/r.$ The residues at these
poles are respectively%
\[
\frac{c_{1}+c_{2}}{c_{1}c_{2}}\ln\frac{c_{1}\sqrt{c_{1}l/\theta}}{c_{1}+c_{2}%
},-\frac{1-c_{1}}{c_{1}}\ln\frac{c_{1}\left(  1+\theta\right)  }{\sqrt{\theta
c_{1}l}},\text{ and }-\frac{1-c_{2}}{c_{2}}\ln\frac{c_{1}}{\sqrt{\theta
c_{1}l}}%
\]
and their sum, after collecting terms, yields formula (\ref{formula1}%
). $\square$

\begin{corollary}
For PCA and REG$_{0},$ for any $\theta\in\left(
0,\bar{\theta}\right)  $ and for sufficiently large $\mathbf{n},p,$ we have%
\begin{equation}
f_{\mathrm{e}}\left(  z_{0}\right)  =\ln c_{1}-\ln\theta-\frac{1-c_{1}}{c_{1}}%
\ln\left(  1+\theta\right)  +\theta/c_{1}. \label{formula2}\tag{SM91}%
\end{equation}
\label{Corollaryf0}
\end{corollary}

\textbf{\small{Proof:}} The corollary is obtained from Lemma \ref{Lemmaf0} by taking
the limit as $c_{2}\rightarrow0$. $\square$

\begin{corollary}
For SMD, for any $\theta\in\left(  0,\bar{\theta}\right)  $ and for
sufficiently large $\mathbf{n},p,$ we have%
\begin{equation}
f_{\mathrm{e}}\left(  z_{0}\right)  =-\ln\theta+\theta^{2}/2. \label{formula3}\tag{SM92}%
\end{equation}
\label{Corollaryf0SMD}
\end{corollary}

\textbf{\small{Proof:}} We remarked earlier that SMD is a limit of PCA and REG$_{0}$
as $c_{1}\rightarrow0$ after the transformations $\theta\mapsto\sqrt{c_{1}%
}\theta$ and $z\mapsto\sqrt{c_{1}}z+1.$ In particular,%
\[
z_{0}^{\rm{SMD}}=\lim_{c_{1}\rightarrow0,\theta\mapsto\sqrt{c_{1}}\theta}%
(z_{0}^{\rm{PCA}}-1)/\sqrt{c_{1}}\text{ and }F^{\rm{SC}}\left(  \lambda\right)
=\lim_{c_{1}\rightarrow0}F_{\mathbf{c}}^{\rm{MP}}\left(  \sqrt{c_{1}}%
\lambda+1\right)  .
\]
These equations imply that%
\[
f_{\mathrm{e}}^{\rm{SMD}}\left(  z_{0}^{\rm{SMD}}\right)  =\lim_{c_{1}\rightarrow
0,\theta\mapsto\sqrt{c_{1}}\theta}\left[  f_{\mathrm{e}}^{\rm{PCA}}\left(  z_{0}%
^{\rm{PCA}}\right)  -\ln\sqrt{c_{1}}\right]  .
\]
Using this relationship together with Corollary \ref{Corollaryf0} yields
$f_{\mathrm{e}}\left(  z_{0}\right)  =-\ln\theta+\theta^{2}/2$ for SMD. $\square$

Observe from Lemma \ref{Lemmaf0} and Corollary \ref{Corollaryf0} that 
\begin{align*}
f_{\mathrm{e}}^{\mathrm{SigD}}& =f_{\mathrm{e}}^{\mathrm{PCA}}+\log c_{1}-\frac{c_{1}+c_{2}}{%
c_{1}c_{2}}\log (c_{1}+c_{2})+\frac{r^{2}}{c_{1}c_{2}}\log [c_{1}l(\theta
)]-\theta /c_{1} \\
& =f_{\mathrm{e}}^{\mathrm{PCA}}+f_{20},\qquad \quad \text{and } \\
f_{\mathrm{e}}^{\mathrm{REG}}& =f_{\mathrm{e}}^{\mathrm{REG_{0}}}+f_{20},
\end{align*}%
where $f_{20}$ is defined at \eqref{eq:f20def}. Combining Table \ref{tab:JO3}
%\ref{tab:AB} 
for $f_{\mathrm{c}}$ with this display and Corollaries \ref{Corollaryf0}
and \ref{Corollaryf0SMD} for $f_{%
\mathrm{e}}(z_{0})$, we can summarize the results for $f_{\mathrm{c}}+f_{\mathrm{e}}(z_{0})
$ by case in Table \ref{tab:fc+fe} below. For the SigD and REG lines we use %
\eqref{eq:f21}, namely $f_{21}=f_{10}+f_{20}$, while for CCA we recall that $%
f_{\mathrm{e}}^{\mathrm{REG}}=f_{\mathrm{e}}^{\mathrm{CCA}}$.

\begin{table}[h]
\centering
\begin{tabular}[h]{ll}
Case \qquad \quad & $F = f_{\mathrm{c}} + f_{\mathrm{e}}(z_{0})$ \\ \hline
&  \\ 
SMD & $1 + \theta^2$ \\[5pt] 
PCA & $1 + \theta/c_1$ \\[5pt] 
SigD & $F^{\mathrm{PCA}} + f_{21}$ \\[5pt] 
REG$_0$ & $2(1+\theta/c_1) + \dfrac{1-c_1}{c_1} \log \dfrac{1-c_1}{1+\theta}$
\\[5pt] 
REG & $F^{\mathrm{REG_0}} + f_{21}$ \\[5pt] 
CCA & $F^{\mathrm{REG}} + f_{21}$%
\end{tabular}%
\caption{Explicit form of $f_{\mathrm{c}}+f_{\mathrm{e}}(z_{0})$ for the different
cases.}
\label{tab:fc+fe}
\end{table}

We turn to the evaluation of $f_{\mathrm{h}}(z_{0})$: in each case it will
turn out to equal $-F=-f_{\mathrm{c}}-f_{\mathrm{e}}(z_{0})$ as shown in
Table \ref{tab:fc+fe}. Again we start with the $\mathsf{q}=0$ cases, in
which $f_{\mathrm{h}}(z)$ is an elementary function.
\vspace{3 mm}

\textbf{\small{SMD.}}
We immediately have $f_{\mathrm{h}}(z_{0})=-z_{0}\theta =-\theta ^{2}-1$. 
% To summarize, 
% we have
% \begin{equation*}
%   2f_{j,{\rm SMD}}(z_0) =
%   \begin{cases}
%     1 + \theta^2/2 + \log \theta  \qquad & j = \text{c} \\
%     -\log \theta + \theta^2/2  &    j = \text{e} \\
%     -\theta^2 - 1 &                 j = \text{h}
%   \end{cases}
% \end{equation*}
% and 
% $f_{\rm SMD}(z_0) = 0$ is immediate.
\vspace{3 mm}

\textbf{\small{PCA.}}
Now $f_{\mathrm{h}}(z_{0})=-z_{0}\theta /[c_{1}(1+\theta )]=-1-\theta /c_{1}$%
. % and again $f_{\rm PCA}(z_0) = 0$. 
\vspace{3 mm}

\textbf{\small{SigD.}}
This time, referring to the definition of $f_{21}$ in Table \ref{tab:JO3}, 
% We organize this by reference to PCA. Thus from Table
% \ref{tab:AB} and some algebra with Lemma 13,
% \begin{align*}
%   2f_{\rm c,SigD} \ \  \quad & = 2f_{\rm c,PCA} \quad \ \  + k_1 + k_0, \\
%   2f_{\rm e,SigD}(z_0) & = 2f_{\rm e,PCA}(z_0) + k_2 - k_0,
% \end{align*}
% where $k_0, k_1$ and $k_2$ are defined at \eqref{eq:k01def} and
% \eqref{eq:k2def} respectively. 
% Finally
\begin{align*}
f_{\mathrm{h}}^{\mathrm{SigD}}(z_{0})& =\frac{r^{2}}{c_{1}c_{2}}\log \left[ 1-\frac{%
c_{2}}{c_{1}}\frac{c_{1}+\theta }{l(\theta )}\right] =\frac{r^{2}}{c_{1}c_{2}%
}\log \left[ \frac{r^{2}}{c_{1}l(\theta )}\right]  \\
& =-f_{21}-F^{\mathrm{PCA}}.
\end{align*}%
%Combining the last two displays, we get $f_{\rm SigD} = 0$.
\vspace{3 mm} 

\textbf{\small{REG$_{0}$.}}
Since $t_{0}=\tfrac{1}{2}(1+\sqrt{1+4\eta _{0}})$ satisfies $%
t_{0}^{2}-t_{0}-\eta _{0}=0$, we have 
\begin{equation*}
\varphi _{0}(t_{0})=\log t_{0}-t_{0}-\eta _{0}/t_{0}+1=\log t_{0}-2\eta
_{0}/t_{0}.
\end{equation*}%
Since $\eta _{0}(z_{0})=(1+\theta )(c_{1}+\theta )/(1-c_{1})^{2}$, we find
after algebra that 
\begin{equation}
\sqrt{1+4\eta _{0}}=\frac{1+c_{1}+2\theta }{1-c_{1}},  \label{eq:eta0-D0}\tag{SM93}
\end{equation}%
so that 
\begin{equation*}
t_{0}=\frac{1+\theta }{1-c_{1}},\qquad \frac{\eta _{0}}{t_{0}}=\frac{%
c_{1}+\theta }{1-c_{1}},
\end{equation*}%
and 
\begin{equation*}
f_{\mathrm{h}}(z_{0})=\frac{1-c_{1}}{c_{1}}\varphi _{0}(t_{0})=\frac{1-c_{1}%
}{c_{1}}\log \frac{1+\theta }{1-c_{1}}-\frac{2(c_{1}+\theta )}{c_{1}}=-F^{%
\mathrm{REG_{0}}}.
\end{equation*}
\vspace{3 mm}

\textbf{\small{REG.}}
Combining the definitions of $f_{\mathrm{h}}$ and $\varphi _{1}(t_{1})$ we
have 
\begin{equation*}
f_{\mathrm{h}}^{\mathrm{REG}}(z)=\frac{1-c_{1}}{c_{1}}\left\{ -\eta _{1}t_{1}+\kappa
\log \frac{\kappa }{t_{1}}+(\kappa -1)\log \left( \frac{t_{1}-1}{\kappa -1}%
\right) \right\} .
\end{equation*}%
Combining \eqref{eq:t1} and \eqref{eq:eta1} gives 
\begin{equation}
\eta _{1}t_{1}=\frac{c_{1}+\theta }{1-c_{1}},\qquad \frac{\kappa }{t_{1}}=%
\frac{r^{2}}{1-c_{1}}\frac{1+\theta }{L(\theta )},\qquad \frac{t_{1}-1}{%
\kappa -1}=\frac{1-c_{1}}{1+\theta },  \label{eq:t1eta1}\tag{SM94}
\end{equation}%
so that 
\begin{equation*}
f_{\mathrm{h}}^{\mathrm{REG}}(z_{0})=-\frac{c_{1}+\theta }{c_{1}}+\frac{r^{2}}{%
c_{1}c_{2}}\log \frac{r^{2}}{L(\theta )}+\frac{1-c_{1}}{c_{1}}\log \frac{%
1+\theta }{1-c_{1}}
\end{equation*}%
We can now compare REG with REG$_{0}$ just as SigD was compared with PCA:
thus 
\begin{equation*}
f_{\mathrm{h}}^{\mathrm{REG}}(z_{0})-f_{\mathrm{h}}^{\mathrm{REG_{0}}}(z_{0})=\frac{r^{2}}{%
c_{1}c_{2}}\log \frac{r^{2}}{c_{1}l(\theta )}+\frac{c_{1}+\theta }{c_{1}}%
=-f_{21},
\end{equation*}%
and so 
\begin{equation}
f_{\mathrm{h}}^{\mathrm{REG}}(z_{0})=-F^{\mathrm{REG_{0}}}-f_{21}=-F^{\mathrm{REG}}.
\label{eq:fh-REG}\tag{SM95}
\end{equation}
\vspace{3 mm}

\textbf{\small{CCA.}}
Combining the definitions of $f_{\mathrm{h}}$ and $\varphi _{2}(t_{2})$ we
have 
\begin{equation*}
f_{\mathrm{h}}^{\mathrm{CCA}}(z)=\frac{1-c_{1}}{c_{1}}\left\{ \kappa \log (1-\eta
_{2}t_{2})+\kappa \log \frac{\kappa }{t_{2}}+(\kappa -1)\log \left( \frac{%
t_{2}-1}{\kappa -1}\right) \right\} .
\end{equation*}%
In particular, recalling that $t_{2}=t_{1}$, 
\begin{align*}
f_{\mathrm{h}}^{\mathrm{CCA}}(z_{0})-f_{\mathrm{h}}^{\mathrm{REG}}(z_{0})& =\frac{1-c_{1}}{c_{1}}%
\left[ \kappa \log (1-\eta _{2}t_{2})+\eta _{1}t_{1}\right]  \\
& =\frac{r^{2}}{c_{1}c_{2}}\log \frac{r^{2}}{L(\theta )}+\frac{c_{1}+\theta 
}{c_{1}}=-f_{21},
\end{align*}%
after substitution from \eqref{eq:t2eta2} and \eqref{eq:t1eta1}. 
% \begin{equation*}
%   2f_{\rm h, CCA}(z_0) = \frac{1-c_1}{c_1} \left[\kappa \log(1-\eta_2
%     t_2) + \eta_1 t_1 \right] + 2f_{\rm h, REG}(z_0).
% \end{equation*}
% After substitution from \eqref{eq:t2eta2} and \eqref{eq:t1eta1}, we
% arrive at
% \begin{equation*}
%   2f_{\rm h, CCA}(z_0) -  2f_{\rm h, REG}(z_0)
%    = \frac{r^2}{c_1 c_2} \log \frac{r^2}{L(\theta)} +
%    \frac{c_1+\theta}{c_1} = -f_{21},
% \end{equation*}
In combination with \eqref{eq:fh-REG}, we get 
\begin{equation*}
f_{\mathrm{h}}^{\mathrm{CCA}}(z_{0})=-F^{\mathrm{REG}}-f_{21}=-F^{\mathrm{CCA}}.
\end{equation*}

\subsection{Proof of Lemma JO\ref{ContourProp1F1} (contours of steep descent)} \label{sec: contours fo steep descent}

For SMD, PCA, and SigD, $\left\vert z-\lambda\right\vert $ is obviously
strictly increasing for any $\lambda\in\mathbb{R}$ and as $z$ moves away from
$z_{0}$ along $\mathcal{K}_{1}$. Therefore,
\[
\operatorname{Re}f_{\mathrm{e}}(z)\equiv\int\ln\left\vert z-\lambda\right\vert
\mathrm{d}F_{\mathbf{c}}\left(  \lambda\right)
\]
is strictly increasing. On the other hand, the definition (JO\ref{elementary})
of $f_{\mathrm{h}}(z)$ implies that
$\operatorname{Re}f_{\mathrm{h}}\left(  z\right)  $ is non-decreasing. Hence
$\operatorname{Re}f\left(  z\right)  $ is strictly increasing.

For REG$_{0}$ and CCA, $\left\vert z-\lambda\right\vert $ is strictly
increasing for any $\lambda\geq0$ as $z$ moves away from $z_{0}$ along
$\mathcal{K}_{1}$ because the center of the circumference that includes
$\mathcal{K}_{1}$ is a negative real number. Therefore, $\operatorname{Re}%
f_{\mathrm{e}}(z)$ is strictly increasing. To show that $\operatorname{Re}%
f_{\mathrm{h}}\left(  z\right)  $ is strictly increasing too, it is sufficient to
prove that $\operatorname{Re}\varphi_{j}\left(  t_{j}\right)  $ is strictly
increasing for $j=0,2$.
\vspace{3 mm}

\textbf{\small{Proof of the monotonicity of $\operatorname{Re}\varphi_{j}\left(  t_{j}\right)  $ for $j=0,2$.}}
Let us show that $\operatorname{Re}\varphi _{j}\left( t_{j}\right) $ is
strictly increasing for $j=0,2$ as $z$ moves away from $z_{0}$ along $%
\mathcal{K}_{1}.$ Recall that for $z\in \mathcal{K}_{1},$ we have 
\begin{equation*}
z=z_{1}+\left\vert z_{0}-z_{1}\right\vert \exp \left\{ \mathrm{i}\gamma
\right\} ,\gamma \in \left[ 0,\pi /2\right] .
\end{equation*}%
Let%
\begin{equation}
R_{j}=\left\{ 
\begin{array}{ll}
\left\vert z_{0}-z_{1}\right\vert \theta /\left( 1-c_{1}\right) ^{2} & \text{%
for }j=0 \\ 
\left\vert z_{0}-z_{1}\right\vert \theta c_{2}^{2}/\left[ c_{1}^{2}l\left(
\theta \right) \right] & \text{for }j=2%
\end{array}%
\right. .  \label{Rj}\tag{SM96}
\end{equation}%
For REG$_{0}$, using 
\begin{equation}
\eta _{j}=\left\{ 
\begin{array}{ll}
z\theta /\left( 1-c_{1}\right) ^{2} & \text{for }j=0 \\ 
z\theta c_{2}^{2}/\left[ c_{1}^{2}l\left( \theta \right) \right] & \text{for 
}j=2%
\end{array}%
\right. .  \label{etas}\tag{SM97}
\end{equation}%
and the definition of $\varphi _{0}$ and $t_{0}$, we obtain%
\begin{equation*}
\operatorname{Re}\varphi _{0}\left( t_{0}\right) =\tfrac{1}{2}\ln \left(
1+4R_{0}^{1/2}\cos \left( \gamma /2\right) +4R_{0}\right) -2R_{0}^{1/2}\cos
\left( \gamma /2\right) +1-\ln 2.
\end{equation*}%
Since the derivative of the above expression with respect to $\gamma \in %
\left[ 0,\pi /2\right) $ is positive, $\operatorname{Re}\varphi _{0}\left(
t_{0}\right) $ does strictly increase as $z$ moves away from $z_{0}$ along $%
\mathcal{K}_{1}$.

For CCA, using the identity 
\begin{equation*}
1-\eta _{2}t_{2}=\frac{\kappa }{\kappa -1}\frac{t_{2}-1}{t_{2}}
\end{equation*}%
we obtain%
\begin{equation}
\operatorname{Re}\varphi _{2}\left( t_{2}\right) =-2\kappa \ln \left\vert
t_{2}\right\vert +\left( 2\kappa -1\right) \ln \left\vert t_{2}-1\right\vert
+\kappa \ln \frac{\kappa }{\kappa -1}  \label{newf2CCA}\tag{SM98}
\end{equation}%
Further, we have 
\begin{equation*}
\eta _{2}=-\frac{1}{4\kappa \left( \kappa -1\right) }+R_{2}\exp \left\{ 
\mathrm{i}\gamma \right\} 
\end{equation*}%
and 
\begin{equation}
t_{2}=\frac{2\kappa }{\left( k_{2}\exp \left\{ \mathrm{i}\gamma /2\right\}
+1\right) },  \label{tau1onC1}\tag{SM99}
\end{equation}%
where $k_{2}=\sqrt{4R_{2}\kappa \left( \kappa -1\right) }$. Taking the
derivative of $\operatorname{Re}\varphi _{2}\left( t_{2}\right) $ with respect to $%
\gamma ,$ we obtain%
\begin{equation*}
\frac{\mathrm{d}}{\mathrm{d}\gamma }\operatorname{Re}\varphi _{2}\left( t_{2}\right)
=\frac{-k_{2}\sin \gamma /2}{2\left\vert k_{2}\exp \left\{ \mathrm{i}\gamma
/2\right\} +1\right\vert ^{2}}+\frac{k_{2}\sin \gamma /2}{2\left\vert 1-%
\frac{k_{2}}{2\kappa -1}\exp \left\{ \mathrm{i}\gamma /2\right\} \right\vert
^{2}}.
\end{equation*}%
For $\gamma \in \left[ 0,\pi /2\right] ,$ the above derivative is positive if%
\begin{equation*}
\left\vert k_{2}\exp \left\{ \mathrm{i}\gamma /2\right\} +1\right\vert
>\left\vert 1-\frac{k_{2}}{2\kappa -1}\exp \left\{ \mathrm{i}\gamma
/2\right\} \right\vert .
\end{equation*}%
The latter inequality does hold because $k_{2}/\left( 2\kappa -1\right) <k_{2}.$ Hence, $%
\frac{\mathrm{d}}{\mathrm{d}\gamma }\operatorname{Re}\varphi _{2}\left( t_{2}\right)
>0$ for $\gamma \in \left[ 0,\pi /2\right] .$ $\square $
\vspace{3 mm}

It remains to prove Lemma JO\ref{ContourProp1F1} for REG. In the REG case, $z$ moves away from $z_{0}$ 
along $\mathcal{K}_{1}$ when $\tau$ moves
away from $\tau_{0}$ along $\mathcal{C}_{1}$. Using the definition of $\varphi_{j}$ (JO\ref{phij}),
%\ref{phij}
the formula (JO\ref{f2}) for $f_{\mathrm{h}}(z)$,
% (\ref{f2}),
and the expression (JO\ref{identity for eta1}) for $\eta_{1}$,
%(\ref{identity for eta1}),
 we obtain
\[
\operatorname{Re}f_{\mathrm{h}}\left(  \tau\right)  =\frac{1-c_{1}}{2c_{1}}\left(
-\operatorname{Re}\tau+\ln\left\vert \tau+1\right\vert +\kappa
\ln\left\vert \tau+\kappa\right\vert +\kappa\ln\kappa\right)
.
\]
On the other hand, $\left\vert \tau+\kappa\right\vert $ remains constant
on $\mathcal{C}_{1}$ whereas both $-\operatorname{Re}\tau$ and $\left\vert
\tau+1\right\vert $ increase as $\tau$ moves away from $\tau_{0}$ along
$\mathcal{C}_{1}.$ To see that $\left\vert \tau+1\right\vert $ indeed
increases recall that the center $-\kappa$ of the circumference that
represents $\mathcal{C}_{1}$ is to the left of the point $-1$. Hence,
$\operatorname{Re}f_{\mathrm{h}}\left(  \tau\right)  $ is strictly increasing.

To show that $\operatorname{Re}f_{\mathrm{e}}\left(  \tau\right)  $ is strictly
increasing too it is sufficient to verify that%
\[
\left\vert z-\lambda\right\vert \equiv\left\vert \frac{c_{1}\left(
1-c_{1}\right)  }{\theta c_{2}}\frac{\tau\left(  \tau+1\right)  }%
{\tau+\kappa}-\lambda\right\vert
\]
is strictly increasing for any $\lambda$ from the support of $F_{\mathbf{c}}.$
Since $\left\vert \tau+\kappa\right\vert $ remains constant, it is
sufficient to show that%
\[
\gamma\left(  \tau,x\right)  \equiv\left\vert \tau\left(  \tau+1\right)
-x\left(  \tau+\kappa\right)  \right\vert ^{2}%
\]
increases as $\tau$ moves away from $\tau_{0}$ along $\mathcal{C}_{1}$ for any
$x=\lambda\theta c_{2}/\left[  c_{1}\left(  1-c_{1}\right)  \right]  .$

Parameterize $\tau\in\mathcal{C}_{1}$ as $-\kappa+\rho e^{\mathrm{i}%
\alpha},\alpha\in\left[  0,\pi/2\right]  .$ Then elementary calculations
yield
\begin{align*}
\gamma\left(  \tau,x\right)   &  =\rho^{4}+\left(  2\kappa-1+x\right)
^{2}\rho^{2}-2\rho^{3}\left(  2\kappa-1+x\right)  \cos\alpha\\
&  +\kappa^{2}\left(  \kappa-1\right)  ^{2}+2\left(  \rho^{2}%
\cos2\alpha-\left(  2\kappa-1+x\right)  \rho\cos\alpha\right)
\kappa\left(  \kappa-1\right)
\end{align*}
so that%
\begin{equation}
\frac{\mathrm{d}\gamma\left(  \tau,x\right)  }{\mathrm{d}\cos\alpha}%
=2\rho\left\{  -\left(  2\kappa-1+x\right)  \left[  \rho^{2}%
+\kappa\left(  \kappa-1\right)  \right]  +4\rho\kappa\left(
\kappa-1\right)  \cos\alpha\right\}  . \label{derivative cosinus1}\tag{SM100}%
\end{equation}
We would like to prove that the derivative $\mathrm{d}\gamma\left(
\tau,x\right)  /\mathrm{d}\cos\alpha$ is negative. As is seen from
(\ref{derivative cosinus1}), the derivative is decreasing in $x$ and
increasing in $\cos\alpha$. Since $x\geq0$ and $\cos\alpha\leq1$, it is
sufficient to show that $\mathrm{d}\gamma\left(  \tau,0\right)  /\mathrm{d}%
\cos\alpha$ is negative for $\cos\alpha=1.$ We have%
\begin{equation*}
\left. \frac{\mathrm{d}\gamma\left(  \tau,0\right)  }{\mathrm{d}\cos\alpha
}\right\vert _{\cos\alpha=1}  =-2\rho\left(  2\kappa-1\right)
\left\{  \left(  \rho-\frac{2\kappa\left(  \kappa-1\right)
}{2\kappa-1}\right)  ^{2}+
\kappa\left(  \kappa-1\right)  -\left(  \frac
{2\kappa\left(  \kappa-1\right)  }{2\kappa-1}\right)
^{2}\right\}  .
\end{equation*}
This is negative because the expression in the figure brackets is positive. 
The positivity follows from the observation that%
\[
\kappa\left(  \kappa-1\right)  \left(  2\kappa-1\right)
^{2}-4\kappa^{2}\left(  \kappa-1\right)  ^{2}=\kappa\left(
\kappa-1\right)  >0.
\]
To summarize, both $\operatorname{Re}f_{\mathrm{e}}\left(  \tau\right)  $ and
$\operatorname{Re}f_{\mathrm{h}}\left(  \tau\right)  $ are strictly increasing as
$\tau$ moves away from $\tau_{0}$ along $\mathcal{C}_{1}.$ Hence, the image of
$\mathcal{C}_{1},$ $\mathcal{K}_{1},$ is a contour of steep descent of
$-\operatorname{Re}f(z)$ in $z$-plane. $\square$

\section{Laplace approximation}

\subsection{Proof of Lemma JO\ref{Olver} (extends Olver's asy. expansion)} \label{sec: extension of Olver}

We closely follow Olver's (1997, pp. 121-125) derivation of an approximation
to a similar integral, augmenting Olver's proof by explicit uniform bounds
on the approximation errors. First, focus on the integral 
\begin{equation*}
I^{+}=\int_{\left[ z_{0},b\right] _{\mathcal{K}}}e^{-p\phi (z)}\chi (z)%
\mathrm{d}z.
\end{equation*}%
Let us introduce new variables $v$ and $w$ by the equations%
\begin{equation}
w^{2}=v=\phi \left( z\right) -\phi _{0},  \label{w and v}\tag{SM101}
\end{equation}%
where the branch of $w$ is determined by $\lim \left\{ \arg v\right\} =\arg
\phi _{2}+2\beta $ as $z\rightarrow z_{0}$ along $\left( z_{0},b\right) _{%
\mathcal{K}}$, and by continuity elsewhere. Here $\beta=\lim \arg (z-z_{0})$
as $z\rightarrow z_{0}$ along $\left( z_{0},b\right) _{%
\mathcal{K}}$. 
%Since, by assumption, $\newre \phi(z)$ is strictly increasing
%as $z$ moves away from $z_{0}$ along $\left( z_{0},b\right] _{\mathcal{K}}$, 
%$w$ cannot equal zero on 
%$\left( z_{0},b\right] _{\mathcal{K}}$, and thus, the branch of $w$ is uniquely specified there.

Consider $w$ as a function of $z.$ A proof of the following auxiliary lemma
is given in the next subsection of this note.

\begin{lemma}
\label{Auxiliary A}Let $B\left( \alpha ,R\right) $\ and $\overline{B}\left(
\alpha ,R\right) $ denote, respectively, the open and closed balls in $%
\mathbb{C}$ with center at $\alpha $\ and radius $R.$ Suppose that
assumptions A0-A4 hold. Then, there exist $\rho _{1},\rho _{2}>0$ with $\rho
_{2}<\rho _{1}$, which do not depend on $p$ and $\omega $, such that, for
sufficiently large $p,$ \vspace{2 mm} \newline
(i) $w(z)$\textit{\ is holomorphic in }$\overline{B}\left( z_{0},\rho
_{1}\right) .$\textit{\ Furthermore, for any }$\zeta _{1},\zeta _{2}$\textit{%
\ from }$\overline{B}\left( z_{0},\rho _{1}\right) $,\textit{\ we have } \\ $%
\left\vert w\left( \zeta _{2}\right) -w\left( \zeta _{1}\right) \right\vert
\geq \tfrac{1}{2}\left\vert \phi _{2}^{1/2}\right\vert \left\vert \zeta
_{2}-\zeta _{1}\right\vert $. \vspace{2 mm} \newline
(ii) $w(z)$\textit{\ maps }$B\left( z_{0},\rho _{1}\right) $\textit{\ on
an open set }$W$ that contains $0.$\textit{\ The inverse function }$z(w)$%
\textit{\ is holomorphic in }$W$. \vspace{2 mm} \newline
(iii) For any $z_{1}\in \left[ z_{0},b\right] _{\mathcal{K}}$ such that $%
\left\vert z_{1}-z_{0}\right\vert =\rho _{2}$, $\overline{B}\left(
0,2\left\vert w\left( z_{1}\right) \right\vert \right) $\textit{\ is
contained in }$W.$
\end{lemma}

Let $z_{1}$ be a point of $\left[ z_{0},b\right] _{\mathcal{K}}$ satisfying
Lemma \ref{Auxiliary A} (iii). Then the portion $\left[ z_{0},z_{1}\right] _{%
\mathcal{K}}$ of $\mathcal{K}$ can be deformed, without changing the value
of the integral 
\begin{equation*}
\overline{I^{+}}=\int_{\left[ z_{0},z_{1}\right] _{\mathcal{K}}}e^{-p\phi
(z)}\chi (z)\mathrm{d}z,
\end{equation*}%
to make its $w(z)$ map a straight line. Since $\chi (z)$ may be random, the
latter statement is only true under qualification: \textquotedblleft with
probability arbitrarily close to one (w.p.a.c.1) for sufficiently large $p.$%
\textquotedblright\ Transformation to the variable $v$ gives%
\begin{equation}
\overline{I^{+}}=e^{-p\phi _{0}}\int_{_{[0,\tau ]}}e^{-pv}\varphi (v)\mathrm{%
d}v,  \label{change of variables}\tag{SM102}
\end{equation}%
where 
\begin{equation}
\tau =\phi \left( z_{1}\right) -\phi _{0},\text{ }\varphi (v)=\chi (z)/\phi ^{\prime }(z),
\label{tau and fi}\tag{SM103}
\end{equation}%
and the path for the integral on the right-hand side of (\ref{change of
variables}) is also a straight line.

For $\left\vert v\right\vert \leq \tau $ with $\left\vert v\right\vert \neq
0 $, $\varphi (v)$ has a convergent expansion of the form%
\begin{equation}
\varphi (v)=\sum_{s=0}^{\infty }a_{s}v^{\left( s-1\right) /2},
\label{fi_v_w as series}\tag{SM104}
\end{equation}%
w.p.a.c.1 for sufficiently large $p.$ Indeed,  
it is sufficient to show that expansion%
\begin{equation}
w\varphi(v) \equiv w \chi (z)/\phi'
(z)=\sum_{s=0}^{\infty }a_{s}w^{s}  \label{fi_v_w as series detailed}\tag{SM105}
\end{equation}%
converges for $w\in W$, w.p.a.c.1 for sufficiently large $p$. 
But by Lemma \ref{Auxiliary A}, $w \chi (z)$ and $\phi' (z)$ viewed as functions of $w$,
are holomorphic
in $W$, w.p.a.c.1 for sufficiently large $p$. Furthermore, since 
\begin{equation*}
\phi' (z)\frac{\mathrm{d}}{\mathrm{d}w}%
z(w)=2w,
\end{equation*}%
$\phi' (z)$ is not equal to zero 
for $w \in \overline{B}\left( 0,2\tau ^{1/2}\right) \backslash \left\{ 0\right\} $,
and, since $\phi _{2}\neq 0$, $\phi' (z)$
has a simple zero at $w=0$. Therefore, the desired convergence holds,
w.p.a.c.1 for sufficiently large $p.$ 

The coefficients $a_{s}$ in (\ref%
{fi_v_w as series}) can be computed from the coefficients $\phi _{j}$ and $%
\chi _{j}$ defined by equation (JO\ref{f and g series}). The formulae for $%
a_{0},a_{1},$ and $a_{2}$ are given, for example, on p. 86 of Olver (1997).
We use the formula for $a_{0}$ in the statement of Lemma JO\ref{Olver}.

Define $\varphi _{k}(v)$, $k=0,1,2,...$ by the relations $\varphi
_{k}(0)=a_{k}$ and%
\begin{equation}
\varphi (v)=\sum_{s=0}^{k-1}a_{s}v^{\left( s-1\right) /2}+v^{\left(
k-1\right) /2}\varphi _{k}\left( v\right) \text{ for }v\neq 0.
\label{remainders}\tag{SM106}
\end{equation}%
Then the integral on the right-hand side of (\ref{change of variables}) can
be rearranged in the form%
\begin{equation}
\int_{_{\lbrack 0,\tau ]}}e^{-pv}\varphi (v)\mathrm{d}v=\sum_{s=0}^{k-1}%
\Gamma \left( \frac{s+1}{2}\right) \frac{a_{s}}{p^{\left( s+1\right) /2}}%
-\varepsilon _{k,1}\left( p,\omega \right) +\varepsilon _{k,2}\left(
p,\omega \right) ,  \label{epsilons introduced}\tag{SM107}
\end{equation}%
where%
\begin{align}
\varepsilon _{k,1}\left( p,\omega \right) & =\sum_{s=0}^{k-1}\Gamma \left( 
\frac{s+1}{2},\tau p\right) \frac{a_{s}}{p^{\left( s+1\right) /2}},
\label{epsilon1}\tag{SM108} \\
\varepsilon _{k,2}\left( p,\omega \right) & =\int_{[0,\tau
]}e^{-pv}v^{\left( k-1\right) /2}\varphi _{k}\left( v\right) \mathrm{d}v,
\label{epsilon2}\tag{SM109}
\end{align}%
and 
$\Gamma \left( \alpha ,x\right) = \int_{x}^{\infty
}e^{-t} t^{\alpha -1}\mathrm{d}t$
is the incomplete Gamma function. Keep in mind that $\tau
,$ $a_{s},$ and $\varphi _{k}$ depend on $p$ and $\omega .$ 

Note that $\arg
v$ is a continuous function of $z,$ and as mentioned above, $\lim  \left\vert
\arg v \right\vert =\left\vert \arg \phi _{2}+2\beta
\right\vert $ as $z\rightarrow z_{0}$ along $\left( z_{0},b\right) _{%
\mathcal{K}}.$ On the other hand, Lemma JO\ref{Olver} requires that $%
\left\vert \arg \phi _{2}+2\beta \right\vert \leq \pi /2.$ Therefore, $%
\lim \left\vert  \arg v \right\vert \leq \pi /2$ as $%
z\rightarrow z_{0}$ along $\left( z_{0},b\right) _{\mathcal{K}}.$ But since 
$\mathcal{K}$ is a path of steep descent (of $-\newre \phi(z)$), 
$\operatorname{Re}\left( v\right) $ must be
positive for $z\in \left( z_{0},b\right] _{\mathcal{K}}$. Hence, by
continuity, $\left\vert \arg v\right\vert <\pi /2$ for $z\in \left( z_{0},b%
\right] _{\mathcal{K}}$. In particular, $\left\vert \arg \tau \right\vert
=\left\vert \arg \left( \phi \left( z_{1}\right) -\phi _{0}\right)
\right\vert <\pi /2$. Therefore, each incomplete Gamma function in (\ref%
{epsilon1}) takes its principal value.

Consider $\varphi (v)w$ as a function of $w.$ 
%According to (\ref{fi_v_w as
%series detailed}), $\varphi (v)w$ has a convergent series representation%
%\begin{equation}
%\varphi (v)w=\sum_{s=0}^{\infty }a_{s}w^{s}  \label{back to old}
%\end{equation}%
%for $w\in W.$ 
Since $\phi' \left( z\right)
=2w(z)w'(z),$ we have 
\begin{equation}
\varphi (v)w=\chi \left( z \right) /(2w'(z)) .  \label{fi w through z}\tag{SM110}
\end{equation}%
By Lemma \ref{Auxiliary A} (i),%
\begin{equation}
\left\vert w'(z)\right\vert >\tfrac{1}{2}%
\left\vert \phi _{2}^{1/2}\right\vert  \label{bound on w prime}\tag{SM111}
\end{equation}%
for $z\in \overline{B}\left( z_{0},\rho _{1}\right) .$ Equation (\ref{fi w
through z}), inequality (\ref{bound on w prime}), and Assumptions A2, A5
imply that 
\begin{equation}
\sup_{w\in W}\left\vert \varphi (v)w\right\vert =\sup_{z\in \overline{B}%
\left( z_{0},\rho _{1}\right) }\left\vert \chi (z)/(2w'(z)) 
\right\vert =O_{\mathrm{P}}(1)
\label{fiw is Op1}\tag{SM112}
\end{equation}
as $p \rightarrow \infty$,
where $O_{\mathrm{P}}(1)$ is uniform in $\omega \in \Omega .$

Further, by Assumption A4, there exist positive constants $\tau _{1}$ and $\tau _{2}$
(that may depend on $\rho_{2} \equiv \left\vert z_{1}-z_{0} \right\vert$) such that 
for all $\omega \in \Omega $ and sufficiently large $p,$ $\operatorname{Re}%
\tau >\tau _{1}$ and $\left\vert \operatorname{Im}\tau \right\vert <\tau _{2}$.
Since $\left\vert \tau \right\vert \geq \left\vert \operatorname{Re}\tau \right\vert
>\tau _{1},$ $B\left( 0,\left\vert \tau _{1}\right\vert ^{1/2}\right) $ is
contained in $W$, where $\varphi (v)w$ is analytic. Using Cauchy's estimates
for the derivatives of an analytic function (see Theorem 10.26 in Rudin
(1987)), (\ref{fi_v_w as series detailed}) and (\ref{fiw is Op1}), we get%
\begin{equation}
\left\vert a_{s}\right\vert \leq \left\vert \tau _{1}\right\vert
^{-s/2}\sup_{w\in B\left( 0,\left\vert \tau _{1}\right\vert ^{1/2}\right)
}\left\vert \varphi w\right\vert =O_{\mathrm{P}}(1).  \label{as are Op1}\tag{SM113}
\end{equation}

Next, Olver (1997, ch. 4, pp.109-110) shows that $\Gamma \left( \alpha
,\zeta \right) =O\left( e^{-\zeta }\zeta ^{\alpha -1}\right) $ as $%
\left\vert \zeta \right\vert \rightarrow \infty ,$ uniformly in the sector $%
\left\vert \arg \left( \zeta \right) \right\vert \leq \pi/2-\delta $
for an arbitrary positive $\delta .$ Let us take $\alpha =(s+1)/2$ and 
$\zeta =\tau p.$ Since $\operatorname{Re}\tau >\tau _{1}$ and $\left\vert \operatorname{Im}%
\tau \right\vert <\tau _{2}$, we have 
\begin{equation*}
\left\vert \tau p\right\vert >\tau _{1}p\rightarrow \infty
\end{equation*}%
and 
\begin{equation*}
\left\vert \arg \left( \tau p\right) \right\vert =\left\vert \arctan (%
\operatorname{Im}\tau / \operatorname{Re}\tau )\right\vert <\arctan (\tau _{2}/\tau
_{1})<\pi/2,
\end{equation*}%
uniformly in $\omega \in \Omega $ for sufficiently large $p.$ Therefore, 
\begin{equation}
\Gamma \left( \frac{s+1}{2},\tau p\right) =O\left( e^{-\tau p}\left( \tau
p\right) ^{\frac{s-1}{2}}\right) =O\left( e^{-\tfrac{1}{2}\tau
_{1}p}\right)  \label{incomplete gamma 1}\tag{SM114}
\end{equation}%
for any integer $s$, uniformly in $\omega \in \Omega .$ Equality (\ref%
{incomplete gamma 1}), the definition (\ref{epsilon1}) of $\varepsilon
_{k,1}\left( p,\omega \right) $, and inequality (\ref{as are Op1}) imply
that 
\begin{equation}
\varepsilon _{k,1}\left( p,\omega \right) =O_{\mathrm{P}}\left( e^{-\frac{1}{%
2}\tau _{1}p}\right) ,  \label{epsilon 1 bound}\tag{SM115}
\end{equation}%
where $O_{\mathrm{P}}$ is uniform in $\omega \in \Omega .$

Now consider $w^{k}\varphi _{k}\left( v\right) $ as a function of $w.$
Since, by definition, 
\begin{equation*}
w^{k}\varphi _{k}\left( v\right) =\varphi \left( v\right)
w-\sum_{s=0}^{k-1}a_{s}w^{s},
\end{equation*}%
it can be interpreted as a remainder in the Taylor expansion of $\varphi
\left( v\right) w.$ As explained above, such an expansion is valid in $W$,
which includes the ball $B\left( 0,2\left\vert \tau \right\vert
^{1/2}\right) $ by Lemma \ref{Auxiliary A} (iii). By a general formula for
remainders in Taylor expansions, for any $w\in B\left( 0,\left\vert \tau
\right\vert ^{1/2}\right) $,%
\begin{equation}
\left\vert w^{k}\varphi _{k}\left( v\right) \right\vert \leq \frac{%
\left\vert w\right\vert ^{k}}{k!}\max_{w\in B\left( 0,\left\vert \tau
\right\vert ^{1/2}\right) }\left\vert \frac{\mathrm{d}^{k}}{\mathrm{d}w^{k}}%
\left( w\varphi \left( v\right) \right) \right\vert .
\label{remainder general}\tag{SM116}
\end{equation}

Further, for any $w\in B\left( 0,\left\vert \tau \right\vert ^{1/2}\right) $%
, a ball with radius $\left\vert \tau _{1}\right\vert ^{1/2}$ centered in $w$
is contained in the ball $B\left( 0,2\left\vert \tau \right\vert
^{1/2}\right) \subset W$. Therefore, using (\ref{fiw is Op1}) and Cauchy's
estimates for the derivatives of an analytic function (see Theorem 10.26 in
Rudin~(1987)), we get%
\begin{equation}
\max_{w\in B\left( 0,\left\vert \tau \right\vert ^{1/2}\right) }\left\vert 
\frac{\mathrm{d}^{k}}{\mathrm{d}w^{k}}\left( w\varphi \left( v\right)
\right) \right\vert \leq k!\left\vert \tau _{1}\right\vert ^{-k/2}\sup_{w\in
W}\left\vert w\varphi \left( v\right) \right\vert =O_{\mathrm{P}}(1).
\label{Cauchy again}\tag{SM117}
\end{equation}

Combining (\ref{remainder general}) and (\ref{Cauchy again}), we have 
\begin{equation*}
\sup_{v\in \left( 0,\tau \right] }\left\vert \varphi _{k}\left( v\right)
\right\vert =O_{\mathrm{P}}(1).
\end{equation*}%
This equality together with (\ref{as are Op1}) and the fact that, by
definition, $\varphi _{k}\left( 0\right) =a_{k}$ imply that%
\begin{equation}
\max_{v\in \left[ 0,\tau \right] }\left\vert \varphi _{k}\left( v\right)
\right\vert =O_{\mathrm{P}}(1),  \label{Op1 everywhere}\tag{SM118}
\end{equation}%
where $O_{\mathrm{P}}(1)$ is uniform in $\omega \in \Omega .$

For $\varepsilon _{k,2}\left( p,\omega \right) $, the substitution of
variable $v=\tau x/p$ in the integral (\ref{epsilon2}) yields%
\begin{equation*}
\varepsilon _{k,2}\left( p,\omega \right) =p^{-\left( k+1\right)
/2}\int_{0}^{p}e^{-\tau x}x^{\frac{k-1}{2}}\tau ^{\frac{k+1}{2}}\varphi
_{k}\left( v\right) \mathrm{d}x.
\end{equation*}%
Therefore,%
\begin{align}
\left\vert \varepsilon _{k,2}\left( p,\omega \right) p^{\left( k+1\right)
/2}\right\vert & <\max_{v\in \left[ 0,\tau \right] }\left\vert \varphi
_{k}\left( v\right) \right\vert \int_{0}^{p}e^{-\operatorname{Re}\tau x}x^{\frac{k-1%
}{2}}\left\vert \tau \right\vert ^{\frac{k+1}{2}}\mathrm{d}x
\label{eps2 as an integral}\tag{SM119} \\
& <\max_{v\in \left[ 0,\tau \right] }\left\vert \varphi _{k}\left( v\right)
\right\vert \int_{0}^{\infty }e^{-\frac{\operatorname{Re}\tau }{\left\vert \tau
\right\vert }y}y^{\frac{k-1}{2}}\mathrm{d}y.  \notag
\end{align}%
Since $\operatorname{Re}\tau >\tau _{1}$ and $\left\vert \operatorname{Im}\tau \right\vert
<\tau _{2}$, we have 
\begin{equation*}
\frac{\operatorname{Re}\tau }{\left\vert \tau \right\vert }\geq \frac{\operatorname{Re}\tau 
}{\left\vert \operatorname{Re}\tau \right\vert +\left\vert \operatorname{Im}\tau \right\vert 
}>\frac{\tau _{1}}{\tau _{1}+\tau _{2}}
\end{equation*}%
for all $\omega \in \Omega $ and sufficiently large $p.$ Therefore, the
integral in (\ref{eps2 as an integral}) is bounded uniformly in $\omega \in
\Omega .$ Using (\ref{Op1 everywhere}), we conclude that%
\begin{equation}
\varepsilon _{k,2}\left( p,\omega \right) =O_{\mathrm{P}}\left( p^{-\left(
k+1\right) /2}\right) .  \label{epsilon 2 final}\tag{SM120}
\end{equation}

Combining (\ref{change of variables}), (\ref{epsilons introduced}), (\ref%
{epsilon 1 bound}), and (\ref{epsilon 2 final}), we obtain%
\begin{equation}
\overline{I^{+}}=e^{-p\phi _{0}}\left( \sum_{s=0}^{k-1}\Gamma \left( \frac{%
s+1}{2}\right) \frac{a_{s}}{p^{\left( s+1\right) /2}}+\frac{O_{\mathrm{P}%
}\left( 1\right) }{p^{\left( k+1\right) /2}}\right) ,  \label{main part}\tag{SM121}
\end{equation}%
where $O_{\mathrm{P}}\left( 1\right) $ is uniform in $\omega \in \Omega .$

Let us now consider the contribution of $\left[ z_{1},b\right] _{\mathcal{K}%
} $ to the contour integral%
\begin{equation*}
I^{+}=\int_{\left[ z_{0},b\right] _{\mathcal{K}}}e^{-p\phi (z)}\chi (z)%
\mathrm{d}z.
\end{equation*}%
Since $\mathcal{K}$ is a contour of steep descent,%
\begin{equation*}
\inf_{z\in \left[ z_{1},b\right] _{\mathcal{K}}}\operatorname{Re}\left( \phi \left(
z\right) -\phi _{0}\right) \geq \operatorname{Re}\tau >\tau _{1}.
\end{equation*}%
Therefore, by assumptions A5 and A0, we have 
\begin{align}
\left\vert I^{+}-\overline{I^{+}}\right\vert & \leq e^{-p\phi _{0}}e^{-p\tau
_{1}}\int_{\left[ z_{1},b\right] _{\mathcal{K}}}\left\vert \chi (z)\mathrm{d}%
z\right\vert  \label{vertical remainder}\tag{SM122} \\
& \leq e^{-p\phi _{0}}e^{-p\tau _{1}}\left\vert \mathcal{K}\right\vert O_{%
\mathrm{P}}\left( 1\right) =e^{-p\phi _{0}}e^{-p\tau _{1}}O_{\mathrm{P}%
}\left( 1\right) ,  \notag
\end{align}%
where $O_{\mathrm{P}}\left( 1\right) $ is uniform in $\omega \in \Omega $.

Combining (\ref{main part}) and (\ref{vertical remainder}), we obtain%
\begin{equation}
I^{+}=e^{-p\phi _{0}}\left( \sum_{s=0}^{k-1}\Gamma \left( \frac{s+1}{2}%
\right) \frac{a_{s}}{p^{\left( s+1\right) /2}}+\frac{O_{\mathrm{P}}\left(
1\right) }{p^{\left( k+1\right) /2}}\right) .  \label{Kplus integral}\tag{SM123}
\end{equation}

Finally, note that 
\begin{equation*}
I_{p,\omega }=I^{+}-I^{-},
\end{equation*}%
where%
\begin{equation*}
I^{-}=\int_{\left[ z_{0},a\right] _{\mathcal{K}}}e^{-p\phi (z)}\chi (z)%
\mathrm{d}z
\end{equation*}%
where $\left[ z_{0},a\right] _{\mathcal{K}}$ is a contour that coincides
with $\left[ a,z_{0}\right] _{\mathcal{K}}$ but has the opposite
orientation. The integral $I^{-}$ can be analyzed similarly to $I^{+}$. As
explained in Olver (1997, pp.121--122), $a_{s}$ with odd $s$ in the
asymptotic expansion for $I^{-}$ coincides with the corresponding $a_{s}$ in
the asymptotic expansion for $I^{-}.$ However, $a_{s}$ with even $s$ in the
two expansions differ by the sign. Therefore, coefficients $a_{s}$ with odd $%
s$ cancel out, but those with even $s$ double in the difference of the two
expansions. Setting $k=2m$, we have%
\begin{equation*}
I_{p,\omega }=2e^{-p\phi _{0}}\left( \sum_{s=0}^{m-1}\Gamma \left( s+\frac{1%
}{2}\right) \frac{a_{2s}}{p^{s+1/2}}+\frac{O_{\mathrm{P}}\left( 1\right) }{%
p^{m+1/2}}\right) ,
\end{equation*}%
which establishes the lemma. $\square $ \vspace{3 mm}

\textbf{\small{Proof of Lemma \protect\ref{Auxiliary A}.}}

First, we show that there exists $\rho_{1}$ such that 
$w(z)$\ is holomorphic in $\overline{B}\left( z_{0},\rho
_{1}\right) $ and that $\frac{\mathrm{d}}{\mathrm{d}z}w\left( z_{0}\right)
=\phi _{2}^{1/2}$. Let $\phi ^{(j)}(z)$ denote the $j$-th order derivative
of $\phi (z).$ Consider a Taylor expansion of $\phi ^{(j)}(z)$ at $z_{0}$%
\begin{equation*}
\phi ^{(j)}(z)=\sum\limits_{s=0}^{k}\frac{1}{s!}\phi ^{\left( j+s\right)
}\left( z_{0}\right) \left( z-z_{0}\right) ^{s}+R_{j,k+1}.
\end{equation*}%
In general, for any $z\in \overline{B}\left( z_{0},x\right) $, the remainder 
$R_{j,k+1}$ satisfies 
\begin{equation}
\left\vert R_{j,k+1}\right\vert \leq \frac{\left\vert z-z_{0}\right\vert
^{k+1}}{\left( k+1\right) !}\max_{\left\vert t-z_{0}\right\vert \leq
x}\left\vert \phi ^{(j+k+1)}\left( t\right) \right\vert .
\label{remainder simple}\tag{SM124}
\end{equation}%
By assumptions A1--A3, there exist constants $C_{1},C_{2},$ and $C_{4},$
such that 
\begin{equation}
\left\vert \phi ^{(3)}\left( t\right) \right\vert \leq \frac{C_{4}}{C_{2}}%
\left\vert \phi ^{(2)}\left( z_{0}\right) \right\vert  \label{ratio 2}\tag{SM125}
\end{equation}%
for any $t\in \overline{B}\left( z_{0},C_{1}\right) .$ Let $\rho _{1}=\min
\left\{ C_{1},\frac{C_{2}}{2C_{4}}\right\} .$ Then, combining (\ref{ratio 2}%
) with (\ref{remainder simple}) and recalling that $\frac{1}{j!}\phi^{(j)}(z_{0})=\phi_{j}$,
we obtain for $z\in \overline{B}\left(
z_{0},\rho _{1}\right) $,%
\begin{equation}
\left\vert R_{0,3}\right\vert \leq \frac{\left\vert z-z_{0}\right\vert ^{2}}{%
6}\left\vert \phi _{2}\right\vert ,\text{ and }\left\vert R_{1,2}\right\vert
\leq \frac{\left\vert z-z_{0}\right\vert }{2}\left\vert \phi _{2}\right\vert
.  \label{next remainder}\tag{SM126}
\end{equation}%
Further, since 
\begin{equation*}
R_{0,2}=\phi _{2}\left( z-z_{0}\right) ^{2}+R_{0,3},
\end{equation*}%
the first of the inequalities in (\ref{next remainder}) implies that, for $z\in \overline{B}\left( z_{0},\rho
_{1}\right) $,%
\begin{equation}
\frac{5}{6}\left\vert \phi _{2}\right\vert \left\vert z-z_{0}\right\vert
^{2}\leq \left\vert R_{0,2}\right\vert \leq \frac{7}{6}\left\vert \phi
_{2}\right\vert \left\vert z-z_{0}\right\vert ^{2}.
\label{remainder inequalities}\tag{SM127}
\end{equation}

Next, since $\phi _{1}=0$, inequalities (\ref{remainder inequalities}) imply
that%
\begin{equation}
\left\vert \phi (z)-\phi _{0}\right\vert =\left\vert R_{0,2}\right\vert \geq 
\frac{5}{6}\left\vert \phi _{2}\right\vert \left\vert z-z_{0}\right\vert ^{2}
\label{no zeros of f}\tag{SM128}
\end{equation}%
for any $z\in \overline{B}\left( z_{0},\rho _{1}\right) .$ Since $\phi
_{2}\neq 0$, inequality (\ref{no zeros of f}) implies that $\phi (z)-\phi
_{0}$ does not have zeros in $\overline{B}\left( z_{0},\rho _{1}\right) $
except a zero of the second order at $z=z_{0}.$ Therefore, 
\begin{equation*}
\sqrt{\frac{\phi (z)-\phi _{0}}{\left( z-z_{0}\right) ^{2}}}=\frac{w\left(
z\right) }{z-z_{0}}
\end{equation*}%
is holomorphic inside $\overline{B}\left( z_{0},\rho _{1}\right) $, and
converges to $\phi _{2}^{1/2}$ as $z\rightarrow z_{0}.$ This implies that $%
w\left( z\right) $ is holomorphic in $\overline{B}\left( z_{0},\rho
_{1}\right) $ and $\frac{\mathrm{d}}{\mathrm{d}z}w\left( z_{0}\right) =\phi
_{2}^{1/2}$.

Now let us show that, for any $z \in \overline{B}\left( z_{0},\rho _{1}\right) $,%
\begin{equation}
\left\vert \frac{\mathrm{d}}{\mathrm{d}z}w\left( z\right) -\frac{\mathrm{d}}{%
\mathrm{d}z}w\left( z_{0}\right) \right\vert \leq \tfrac{1}{2}\left\vert 
\frac{\mathrm{d}}{\mathrm{d}z}w\left( z_{0}\right) \right\vert .
\label{continuity of derivative}\tag{SM129}
\end{equation}%
Indeed, since 
\begin{equation*}
\frac{\mathrm{d}}{\mathrm{d}z}w\left( z\right) =\frac{\phi ^{(1)}\left(
z\right) }{2w\left( z\right) }=\tfrac{1}{2}\left( \phi \left( z\right) -\phi
_{0}\right) ^{-1/2}\phi ^{(1)}\left( z\right)
\end{equation*}%
and $\frac{\mathrm{d}}{\mathrm{d}z}w\left( z_{0}\right) =\phi _{2}^{1/2}\neq
0$,%
\begin{equation}
\frac{\frac{\mathrm{d}}{\mathrm{d}z}w\left( z\right) }{\frac{\mathrm{d}}{%
\mathrm{d}z}w\left( z_{0}\right) }=\left( 1+\frac{R_{0,3}}{\phi _{2}\left(
z-z_{0}\right) ^{2}}\right) ^{-\tfrac{1}{2}}\left( 1+\frac{R_{1,2}}{2\phi
_{2}\left( z-z_{0}\right) }\right) .  \label{product form}\tag{SM130}
\end{equation}%
Note that for any $y_{1}$ and $y_{2}$ such that $\left\vert y_{2}\right\vert
<1$,%
\begin{equation}
\left\vert \frac{1+y_{1}}{\sqrt{1+y_{2}}}-1\right\vert \leq \frac{\left\vert
y_{1}\right\vert +\left\vert y_{2}\right\vert }{1-\left\vert
y_{2}\right\vert },  \label{general inequality}\tag{SM131}
\end{equation}%
where the principal branch of the square root is used. This follows from the
facts that, for $\left\vert y_{2}\right\vert <1$, $\left\vert \sqrt{1+y_{2}}%
\right\vert \geq 1-\left\vert y_{2}\right\vert $ and $\left\vert 1+y_{1}-%
\sqrt{1+y_{2}}\right\vert \leq \left\vert y_{1}\right\vert +\left\vert
y_{2}\right\vert .$ Both of these inequalities follow from 
$\left\vert 1-\sqrt{1+y_{2}} \right\vert \leq \left\vert y_{2} \right\vert$, which
can be established by denoting $\sqrt{1+y_{2}}$ as $x$ so that the inequality becomes
$\left\vert 1-x \right\vert \leq \left\vert x^{2}-1 \right\vert$ and using the fact that
$1 \leq \left\vert x+1 \right\vert$ (because $\newre x \geq 0$ 
when $\left\vert y_{2}\right\vert <1$).
Setting 
\begin{equation*}
y_{1}=\frac{R_{1,2}}{2\phi _{2}\left( z-z_{0}\right) }\text{ and }y_{2}=%
\frac{R_{0,3}}{\phi _{2}\left( z-z_{0}\right) ^{2}}
\end{equation*}%
and using (\ref{next remainder}) and (\ref{product form}), we obtain 
\begin{equation*}
\left\vert \frac{\frac{\mathrm{d}}{\mathrm{d}z}w\left( z\right) }{\frac{%
\mathrm{d}}{\mathrm{d}z}w\left( z_{0}\right) }-1\right\vert \leq \tfrac{1}{2}.
\end{equation*}%
Hence, (\ref{continuity of derivative}) holds.

Finally, let $\zeta _{1}$ and $\zeta _{2}$ be any two points in $\overline{B}%
\left( z_{0},\rho _{1}\right) $, and let $\gamma (t)=\left( 1-t\right) \zeta
_{1}+t\zeta _{2}$, where $t\in \left[ 0,1\right] .$ We have 
\begin{equation*}
\int_{0}^{1}\left( \frac{\mathrm{d}}{\mathrm{d}z}w\left( \gamma (t)\right) -%
\frac{\mathrm{d}}{\mathrm{d}z}w\left( z_{0}\right) \right) \mathrm{d}t=\frac{%
w\left( \zeta _{2}\right) -w\left( \zeta _{1}\right) }{\zeta _{2}-\zeta _{1}}%
-\frac{\mathrm{d}}{\mathrm{d}z}w\left( z_{0}\right) .
\end{equation*}%
Therefore, using (\ref{continuity of derivative}), we obtain 
\begin{equation*}
\left\vert \frac{w\left( \zeta _{2}\right) -w\left( \zeta _{1}\right) }{%
\zeta _{2}-\zeta _{1}}-\frac{\mathrm{d}}{\mathrm{d}z}w\left( z_{0}\right)
\right\vert \leq \tfrac{1}{2}\left\vert \frac{\mathrm{d}}{\mathrm{d}z}w\left(
z_{0}\right) \right\vert .
\end{equation*}%
This inequality and the fact that $\frac{\mathrm{d}}{\mathrm{d}z}w\left(
z_{0}\right) =\phi _{2}^{1/2}$ imply part (i) of the lemma.

Part (ii) of the lemma is a simple consequence of part (i). Indeed, by the
open mapping theorem, $W$ is an open set. Next, by (i), $w(z)$ is one-to-one
mapping of $B\left( z_{0},\rho _{1}\right) $ on $W$ and has a non-zero
derivative in $B\left( z_{0},\rho _{1}\right) .$ Further, let $\psi \left(
w\right) $ be defined on $W$ by $\psi \left( w\left( z\right) \right) =z.$
Fix $\tilde{w}\in W.$ Then $\psi \left( \tilde{w}\right) =\tilde{z}$ for a
unique $\tilde{z}$ in $B\left( z_{0},\rho _{1}\right) .$ If $w\in W$ and $%
\psi \left( w\right) =z$, we have 
\begin{equation*}
\frac{\psi \left( w\right) -\psi \left( \tilde{w}\right) }{w-\tilde{w}}=%
\frac{z-\tilde{z}}{w\left( z\right) -w\left( \tilde{z}\right) }.
\end{equation*}%
By (i), $w\rightarrow \tilde{w}$ as $z\rightarrow \tilde{z}$, and the latter
equality implies $\frac{\mathrm{d}}{\mathrm{d}w}\psi \left( \tilde{w}\right)
=\frac{1}{\frac{\mathrm{d}}{\mathrm{d}z}w\left( \tilde{z}\right) }.$
Therefore, $z(w)\equiv \psi \left( w\right) $ is an analytic inverse of $%
w(z) $ on~$W$.

Finally, part (iii) of the lemma can be established as follows. Note that by
part~(i),%
\begin{equation*}
\left\vert w\left( z_{0}+\rho _{1}e^{i\varphi }\right) -w\left( z_{0}\right)
\right\vert \geq \frac{\rho _{1}}{2}\left\vert \frac{\mathrm{d}}{\mathrm{d}z}%
w\left( z_{0}\right) \right\vert 
\end{equation*}%
for any $\varphi \in \left[ 0,2\pi \right] .$ Therefore, for any $w_{1}$
such that $\left\vert w_{1}-w\left( z_{0}\right) \right\vert \leq \frac{\rho
_{1}}{4}\left\vert \frac{\mathrm{d}}{\mathrm{d}z}w\left( z_{0}\right)
\right\vert $, we have 
\begin{equation*}
\min_{\varphi \in \left[ 0,2\pi \right] }\left\vert w_{1}-w\left( z_{0}+\rho
_{1}e^{i\varphi }\right) \right\vert \geq \frac{\rho _{1}}{4}\left\vert 
\frac{\mathrm{d}}{\mathrm{d}z}w\left( z_{0}\right) \right\vert .
\end{equation*}%
By a corollary to the maximum modulus theorem (see Rudin (1987),
p. 212), the latter inequality implies that the function $w\left( z\right)
-w_{1}$ has a zero in $B(z_{0},\rho _{1}).$ Thus, region $W$ includes $%
\overline{B}(0,\frac{\rho _{1}}{4}\left\vert \frac{\mathrm{d}}{\mathrm{d}z}%
w\left( z_{0}\right) \right\vert )$. On the other hand, 
\begin{equation}
\left\vert w\left( z_{1}\right) \right\vert \leq 2\rho _{2}\left\vert \frac{%
\mathrm{d}}{\mathrm{d}z}w\left( z_{0}\right) \right\vert .  \label{wlt2rho}\tag{SM132}
\end{equation}%
Indeed, consider the identity 
\begin{equation*}
w^{2}\left( z_{1}\right) =\phi _{1}\left( z_{1}-z_{0}\right) +R_{0,2}.
\end{equation*}%
Since $\phi _{1}=0$, (\ref{remainder inequalities}) imply 
\begin{equation}
\left\vert w\left( z_{1}\right) \right\vert ^{2}\leq \frac{7}{6}\left\vert
\phi _{2}\right\vert \left\vert z_{1}-z_{0}\right\vert ^{2}.  \label{absTau}\tag{SM133}
\end{equation}%
But, by definition, 
\begin{equation}
\left\vert z_{1}-z_{0}\right\vert =\rho _{2}.  \label{absZ1Z2}\tag{SM134}
\end{equation}%
Since $\frac{\mathrm{d}}{\mathrm{d}z}w\left( z_{0}\right) =\phi _{2}^{1/2},$
(\ref{absTau}) and (\ref{absZ1Z2}) imply (\ref{wlt2rho}). Setting $\rho
_{2}=\rho _{1}/16,$ we obtain that $W$ includes $\overline{B}(0,2\left\vert
w\left( z_{1}\right) \right\vert ).$

\subsection{Evaluation of \MakeLowercase{$\textrm{d}^{2}f(z_{0})/\textrm{d}z^{2}$}} \label{sec: evaluation of the second derivative}

Note that
$-$\textrm{d}$^{2}f_{\mathrm{e}}\left(  z_{0}\right)  /\mathrm{d}z^{2}=\mathrm{d}%
m_{\mathbf{c}}\left(  z_{0}\right)  /$\textrm{d}$z.$ Therefore \textrm{d}%
$^{2}f_{\mathrm{e}}\left(  z_{0}\right)  /\mathrm{d}z^{2}$ can be directly evaluated
using explicit expressions for the Stieltjes transforms of the semicircle,
Marchenko-Pastur and Wachter distributions. Further, using the definition 
of $f_{\mathrm{h}}(z),$ we directly
evaluate \textrm{d}$^{2}f_{\mathrm{h}}\left(  z_{0}\right)  /\mathrm{d}z^{2}.$
Combining the expressions for the second derivatives of $f_{\mathrm{e}}$ and
$f_{\mathrm{h}},$ we obtain values of the second derivative of $f$ reported in Table
JO\ref{Table 5}. \vspace{3 mm}
%\ref{Table 5}.
 
\textbf{\small{Evaluation of $\mathrm{d}m_{\mathbf{c}}\left(  z_{0}\right)/\mathrm{d}z$.}}
For each of the three cases, it is a little easier to evaluate 
\begin{equation}
a(\theta )=\frac{m^{\prime }(z_{0})}{m^{2}(z_{0})}=-\frac{\mathrm{d}}{%
\mathrm{d}z}\left. \left( \frac{1}{m}\right) \right\vert _{z=z_{0}}.
\label{eq:a-def}\tag{SM135}
\end{equation}%
In each case $v=-1/m$ satisfies a quadratic equation in $v=v(z)$. 
%and $z$. 
Differentiation with respect to $z$ yields an equation for $v^{\prime }$
which we write in the form 
\begin{equation}
(C+\Delta )v^{\prime }=C.  \label{eq:vp-eqn}\tag{SM136}
\end{equation}

\textbf{\small{SMD.}}
From \eqref{eq:mq-SMD}, $v=-1/m$ satisfies $1-zv+v^{2}=0$, 
% \begin{equation*}
%   1 - zv + v^2 = 0,
% \end{equation*}
and so, differentiating w.r.t. $z$, 
\begin{equation*}
(2v-z)v^{\prime }=v.
\end{equation*}%
At $z=z_{0}=\theta +1/\theta $, with $m(z_{0})=-\theta $, we get $%
C=v_{0}=1/\theta $ and $\Delta =v_{0}-z_{0}=-\theta $, and 
% \begin{equation*}
%   C = v_0 = 1/\theta, \qquad \quad
%   \Delta = v_0 - z_0 = -\theta,
% \end{equation*}
%and
\begin{equation*}
a(\theta )=v^{\prime }(z_{0})=\frac{C}{C+\Delta }=\frac{1}{1-\theta ^{2}}.
\end{equation*}

\textbf{\small{PCA.}}
From \eqref{eq:mq-PCA}, $v=-1/m$ satisfies $c_{1}z-(z+c_{1}-1)v+v^{2}=0$,
and so, differentiating, 
\begin{equation*}
(2v-z-c_{1}+1)v^{\prime }=v-c_{1}.
\end{equation*}%
At $z_{0}=(1+\theta )(c_{1}+\theta )/\theta $ and $v_{0}=c_{1}(1+\theta
)/\theta $, we have $C=v_{0}-c_{1}=c_{1}/\theta $ and $\Delta
=v_{0}-z_{0}+1=-\theta $, so that 
\begin{equation*}
a(\theta )=v^{\prime }(z_{0})=\frac{C}{C+\Delta }=\frac{c_{1}}{c_{1}-\theta
^{2}}.
\end{equation*}

\textbf{\small{SigD.}}
From \eqref{eq:m-W}, $v=-1/m$ satisfies 
\begin{equation*}
c_{1}z(c_{1}-c_{2}z)-[c_{1}^{2}-c_{1}+(c_{1}+c_{2}-2c_{1}c_{2})z]v+r^{2}v^{2}=0,
\end{equation*}%
and so $v^{\prime }$ satisfies \eqref{eq:vp-eqn} with 
\begin{equation*}
C=c_{1}^{2}-2c_{1}c_{2}(z-v)-(c_{1}+c_{2})v,\qquad \quad \Delta
=-c_{1}+(c_{1}+c_{2})(z-v).
\end{equation*}%
At $z_{0}=(1+\theta )(c_{1}+\theta )/\theta l(\theta )$ and $%
v_{0}=c_{1}(1+\theta )/\theta l(\theta )$, we find $z_{0}-v_{0}=(1+\theta
)/l(\theta )$, and eventually 
\begin{equation*}
C=-\frac{c_{1}}{\theta l(\theta )}[h(\theta )+\theta ^{2}],\qquad \quad
\Delta =\frac{c_{1}}{\theta l(\theta )}\theta ^{2},
\end{equation*}%
with $h(\theta )=c_{1}+c_{2}(1+\theta )^{2}-\theta ^{2}$, and hence 
\begin{equation*}
a(\theta )=v^{\prime }(z_{0})=\frac{h(\theta )+\theta ^{2}}{h(\theta )}.
\end{equation*}%
The results are summarized for later reference in Table \ref{tab:stieltjes}.

\begin{table}[h]
\centering
%  \begin{tabular}[h]{llll}
\begin{tabular}[h]{lccc}
& $m(z_0)$ & $a(\theta)$ & $m^{\prime }(z_0)$ \\ \hline
&  &  &  \\[2pt] 
SMD & $-\theta$ & $\dfrac{1}{1-\theta^2}$ & $\dfrac{\theta^2}{1-\theta^2}$ \\%
[14pt] 
PCA, REG$_0$ & $-\dfrac{\theta}{c_1(1+\theta)}$ & $\dfrac{c_1}{c_1-\theta^2}$
& $\dfrac{\theta^2}{c_1 (1+\theta)^2 (c_1-\theta^2)}$ \\[14pt] 
SigD, REG, CCA & $-\dfrac{\theta l(\theta)}{c_1(1+\theta)}$ & $\dfrac{%
h(\theta)+\theta^2}{h(\theta)}$ & $\dfrac{\theta^2 l^2(\theta)}{c_1^2
(1+\theta)^2} \dfrac{h(\theta)+\theta^2}{h(\theta)}$%
\end{tabular}%
\caption{Summary of Stieltjes transform quantities. $a(\protect\theta)$ is
defined at \eqref{eq:a-def}, $h(\protect\theta) = c_1 + c_2(1+\protect\theta%
)^2 - \protect\theta^2$.}
\label{tab:stieltjes}
\end{table}
\vspace{3 mm}

\textbf{\small{Computation of \textrm{d}$^{2}f_{\mathrm{h}}\left(  z_{0}\right)  /\mathrm{d}z^{2}$.}}
Since $f^{\prime \prime }(z)=f_{\mathrm{e}}^{\prime \prime }(z)+f_{\mathrm{h}%
}^{\prime \prime }(z)$ and $f_{\mathrm{e}}^{\prime \prime }(z)=-m^{\prime
}(z)$, we have 
\begin{equation*}
-f^{\prime \prime }(z_{0})=m^{\prime }(z_{0})-f_{\mathrm{h}}^{\prime \prime
}(z_{0}).
\end{equation*}%
We will see that in each case there is a factorization 
\begin{align*}
m^{\prime }(z_{0})& =m^{2}(z_{0})a(\theta ) \\
f_{\mathrm{h}}^{\prime \prime }(z_{0})& =m^{2}(z_{0})b(\theta ).
\end{align*}%
Note that the functions $a(\theta ),b(\theta )$ are distinct from the
constants $a,b$ in \eqref{eq:fhp}. Thus 
\begin{equation*}
-f^{\prime \prime }(z_{0})=m^{2}(z_{0})[a(\theta )-b(\theta )],
\end{equation*}%
and the entries of Table JO\ref{Table 5} are 
\begin{equation}
D_{2}=\frac{\theta ^{2}}{-f^{\prime \prime }(z_{0})}=\frac{\theta ^{2}}{%
m^{2}(z_{0})}\frac{1}{a(\theta )-b(\theta )}.  \label{eq:tabJO4}\tag{SM137}
\end{equation}

\begin{table}[h]
\centering
\begin{tabular}[h]{lccc}
& $b(\theta)$ & $\dfrac{\theta^2}{m^2(z_0)}$ & $\dfrac{1}{a(\theta)-b(\theta)%
}$ \\ \hline
&  &  &  \\ 
SMD & $0$ & $1$ & $1-\theta^2$ \\[10pt] 
PCA & $0$ & $c_1^2(1+\theta)^2$ & $\dfrac{h_0}{c_1}$ \\[10pt] 
SigD & $-\dfrac{c_1c_2}{r^2}$ & $\dfrac{c_1^2(1+\theta)^2}{l^2}$ & $\dfrac{h
r^2}{c_1^2 l^2}$ \\[10pt] 
REG$_0$ & $\dfrac{c_1}{K_0}$ & $c_1^2(1+\theta)^2$ & $\dfrac{h_0 K_0}{%
c_1(1+\theta)^2}$ \\[10pt] 
REG & $\dfrac{c_1}{K_1}$ & $\dfrac{c_1^2(1+\theta)^2}{l^2}$ & $\dfrac{hK_1}{%
c_1(1+\theta)^2 l^2}$ \\[10pt] 
CCA & $\dfrac{c_1-c_2(1+\theta)}{K_2}$ & $\dfrac{c_1^2(1+\theta)^2}{l^2}$ & $%
\dfrac{hK_2}{(c_1+c_2)(1+\theta)^2l}$%
\end{tabular}%
\caption{Remaining quantities needed for Table JO\ref{Table 5}: as shown at 
\eqref{eq:tabJO4}, the entries there are obtained by multiplying the last
two columns of this table. In the last three cases, some algebra is required
to verify that $a(\protect\theta )-b(\protect\theta )$ factorizes as shown
in the last column. Here $h_{0}=c_{1}-\protect\theta ^{2}$, $K_{0}=1+c_{1}+2%
\protect\theta $, $K_{1}=c_{1}+\protect\theta +(1+\protect\theta )l$ and $%
K_{2}=2(c_{1}+\protect\theta )+(1-c_{1})l$. As $c_{2}\rightarrow 0$, we have 
$h\rightarrow h_{0},l\rightarrow 1,r^{2}\rightarrow c_{1}$ and $%
K_{1},K_{2}\rightarrow K_{0}$.}
\end{table}
\vspace{3 mm}

\textbf{\small{Evaluation of $b(\protect\theta )$.}}
For \textbf{\small{SMD}} and \textbf{\small{PCA}}, $f_{\mathrm{h}}(z)$ is linear in $z$ so $b(\theta )=0$.

For \textbf{\small{SigD}}, from \eqref{eq:fhp} and \eqref{eq:m-wa}, %\eqref{eq:fhp}, 
we find that 
\begin{equation*}
f^{\prime \prime }_{\mathrm{h}}(z_0) = \frac{1}{b} \left( \frac{ab}{1-az_0}
\right)^2 = - \frac{c_1 c_2}{r^2} m^2(z_0).
\end{equation*}

For the $\mathsf{q}=1$ cases, we have from \eqref{eq:partials} that 
\begin{equation}
f_{\mathrm{h}}^{\prime \prime }(z)=\frac{1-c_{1}}{c_{1}}\frac{\mathrm{d}}{%
\mathrm{d}\eta }\chi (\eta )\left( \frac{\mathrm{d}\eta }{\mathrm{d}z}%
\right) ^{2}.  \label{eq:fhpp}\tag{SM138}
\end{equation}%
where $\chi (\eta )=\chi _{j}(\eta _{j})=(\partial /\partial \eta
_{j})\varphi _{j}(t_{j},\eta _{j})$ is given by \eqref{eq:t_j}. 
%\eqref{eq:criticalpt}.

\textbf{\small{REG$_{0}$.}}
Recall that $t_{0}=\tfrac{1}{2}(1+\sqrt{1+4\eta _{0}})=(1+\theta )/(1-c_{1})$%
, so that from \eqref{eq:eta0-D0} 
\begin{equation*}
\dot{t}_{0}=\frac{\mathrm{d}}{\mathrm{d}\eta }t_{0}=(1+4\eta _{0})^{-1/2}=%
\frac{1-c_{1}}{1+c_{1}+2\theta }.
\end{equation*}%
We have both 
\begin{equation*}
\frac{\mathrm{d}}{\mathrm{d}\eta }\chi _{0}(\eta )=\frac{\mathrm{d}}{\mathrm{%
d}\eta }\left( -\frac{1}{t_{0}(\eta )}\right) =\frac{\dot{t}_{0}}{t_{0}^{2}}%
,\qquad \text{and}\qquad \frac{\mathrm{d}\eta }{\mathrm{d}z}=\frac{\theta }{%
(1-c_{1})^{2}},
\end{equation*}%
and so 
\begin{equation*}
f_{\mathrm{h}}^{\prime \prime }(z_{0})=\frac{\theta ^{2}}{%
c_{1}(1+c_{1}+2\theta )(1+\theta )^{2}}=m^{2}(z_{0})\frac{c_{1}}{%
1+c_{1}+2\theta }.
\end{equation*}

\textbf{\small{REG.}}
We have $\chi _{1}(\eta _{1})=-t_{1}(\eta _{1})$ and recall that $t_{1}$
satisfies a quadratic equation $\eta _{1}t^{2}+(1-\eta _{1})t-\kappa =0$, so
that $\dot{t}_{1}=\mathrm{d}t_{1}/\mathrm{d}\eta $ 
%\frac{\diff}{\diff \eta} t_1$ 
solves 
\begin{equation*}
\lbrack 2\eta _{1}t_{1}+1-\eta _{1}]\dot{t}_{1}=t_{1}(1-t_{1}).
\end{equation*}%
Using \eqref{eq:t1}, we can evaluate 
\begin{equation*}
t_{1}(1-t_{1})=-\frac{c_{1}^{2}l(\theta )}{c_{2}^{2}(1+\theta )^{2}},
\end{equation*}%
and setting 
\begin{equation*}
K_{1}(\theta )=c_{1}+\theta +(1+\theta )l(\theta ),
\end{equation*}%
we also have from \eqref{eq:t1eta1} and \eqref{eq:eta1} 
\begin{equation*}
2\eta _{1}t_{1}+1-\eta _{1}=\frac{c_{1}K_{1}(\theta )}{(1-c_{1})L(\theta )}.
\end{equation*}%
We then find from \eqref{eq:fhpp}, the previous displays and $\mathrm{d}\eta
_{1}/\mathrm{d}z=\theta c_{2}/[c_{1}(1-c_{1})]$ that 
%$\frac{\diff \eta_1}{\diff z} = \frac{\theta c_2}{ c_1 (1-c_1)}$ that
\begin{equation*}
f_{\mathrm{h}}^{\prime \prime }(z_{0})=\frac{\theta ^{2}l^{2}(\theta )}{%
(1+\theta )^{2}}\frac{1}{c_{1}K_{1}(\theta )}=m^{2}(z_{0})\frac{c_{1}}{%
K_{1}(\theta )}.
\end{equation*}%
%
% \begin{align*}
%   2f''_{\rm h}(z_0) 
%     & = \frac{\theta^2 l^2(\theta)}{(1+\theta)^2} 
%         \frac{1}{c_1 D_1(\theta)} \\
%    & = m^2(z_0) \frac{c_1}{D_1(\theta)}.
% \end{align*}
% \begin{align*}
%   2f''_{\rm h}(z_0) 
%     & = \frac{\theta^2 L^2(\theta)}{c_1^2 (1+\theta)^2} 
%         \frac{1}{c_1(c_1+\theta) + (1+\theta)L(\theta)} \\
%    & = m^2(z_0) \frac{c_1^2}{h(\theta)+(c_1+\theta)^2}.
% \end{align*}

\textbf{\small{CCA.}}
Recall that $t_{2}(\eta _{2})$ satisfies $\eta _{2}(\kappa -1)t^{2}+t-\kappa
=0$, and hence $\dot{t}_{2}=\mathrm{d}t_{2}(\eta )/\mathrm{d}\eta $ is given
by 
\begin{equation*}
\dot{t}_{2}=\frac{-(\kappa -1)t_{2}^{2}}{1+2\eta _{2}(\kappa -1)t_{2}}.
\end{equation*}%
Since $\chi _{2}(t_{2})=-\kappa t_{2}/(1-\eta _{2}t_{2})$, we have 
\begin{equation*}
\frac{\mathrm{d}}{\mathrm{d}\eta }\chi _{2}(\eta )=\frac{-\kappa }{(1-\eta
_{2}t_{2})^{2}}(t_{2}^{2}+\dot{t}_{2}).
\end{equation*}%
We have $\kappa =r^{2}/c_{2}(1-c_{1})$ and $\kappa -1=c_{1}/c_{2}(1-c_{1})$,
and so from \eqref{eq:t2eta2}, 
\begin{equation*}
(\kappa -1)\eta _{2}t_{2}=\frac{c_{1}+\theta }{(1-c_{1})l(\theta )}
\end{equation*}%
and if we define 
\begin{equation*}
K_{2}(\theta )=(1-c_{1})l(\theta )+2(c_{1}+\theta ),
\end{equation*}%
we arrive at 
\begin{equation*}
t_{2}^{2}+\dot{t}=\frac{t_{2}^{2}}{c_{2}}\left[ c_{2}-\frac{L(\theta )}{%
K_{2}(\theta )}\right] .
\end{equation*}%
Some algebra shows that 
\begin{equation*}
c_{1}[c_{2}K_{2}(\theta )-L(\theta )]=r^{2}[c_{2}(1+\theta )-c_{1}].
\end{equation*}%
From \eqref{eq:fhpp} and the preceding displays, 
\begin{equation*}
f_{\mathrm{h}}^{\prime \prime }(z_{0})=-\frac{1-c_{1}}{\kappa c_{1}}\left[ 
\frac{\kappa t_{2}}{1-\eta _{2}t_{2}}\frac{\mathrm{d}\eta _{2}}{\mathrm{d}z}%
\right] ^{2}\frac{r^{2}}{c_{1}c_{2}}\frac{c_{2}(1+\theta )-c_{1}}{%
K_{2}(\theta )}.
\end{equation*}%
Now from \eqref{eq:t2form} and \eqref{eq:t2eta2}, 
\begin{equation*}
\frac{\kappa t_{2}}{1-\eta _{2}t_{2}}\frac{\mathrm{d}\eta }{\mathrm{d}z}=%
\frac{\theta l(\theta )}{(1+\theta )(1-c_{1})}=-\frac{c_{1}}{1-c_{1}}m(z_{0})
\end{equation*}%
and so finally 
\begin{equation*}
f_{\mathrm{h}}^{\prime \prime }(z_{0})=m^{2}(z_{0})\frac{c_{1}-c_{2}(1+%
\theta )}{K_{2}(\theta )}.
\end{equation*}

\subsection{Proof of Theorem JO\ref{Icase}} \label{sec: proof of theorem 9}
%\ref{Icase}

First, let us show that%
\begin{equation}
L_{1}\left(  \theta;\Lambda\right)  =\frac{g(z_{0})}{\sqrt{-\textrm{d}%
^{2}f(z_{0})/\textrm{d}z^{2}}}+O_{\mathrm{P}}\left(  p^{-1}\right)  ,
\label{I1case integral}\tag{SM139}%
\end{equation}
where $O_{\mathrm{P}}\left(  1\right)  $\textit{\ is uniform with respect to
}$\theta\in\left(  0,\bar{\theta}-\varepsilon\right]  .$ Changing the variable
of integration in (JO\ref{I1case})
%(\ref{I1case}) 
from $z$ to $\zeta=\theta z,$ we obtain%
\begin{equation}
L_{1}\left(  \theta;\Lambda\right)  =\sqrt{\pi p}\frac{1}{2\pi\mathrm{i}}%
\int_{\mathcal{\tilde{K}}}e^{-p\phi(\zeta)}\chi(\zeta)\mathrm{d}\zeta,
\label{olver integral}\tag{SM140}%
\end{equation}
where
\[
\phi(\zeta)=f\left(  \zeta/\theta\right)/2  \text{, }\chi(\zeta)=g(\zeta
/\theta)/\theta,
\]
and $\mathcal{\tilde{K}}$ is the image of $\mathcal{K}_{1}\cup\mathcal{\bar
{K}}_{1}$ under the transformation $z\mapsto\zeta.$ The set of possible values
of $\theta$ is $\Omega\equiv\left(  0,\bar{\theta}-\varepsilon\right]  $.

Using Table JO\ref{Table 5}
%\ref{Table 5} 
and the definitions of $\mathcal{K}_{1},$ $z_{0},$
$f(z),$ and $g(z),$ it is straightforward to verify that the assumptions A0-A4
of Lemma JO\ref{Olver}
%\ref{Olver} 
hold for the integral in (\ref{olver integral}) for all
the six cases that we consider. The validity of A5 follows from Lemma
\ref{gBound} given below and from the definitions of $g\left(  z\right)  $.
Let%
\begin{equation}
\Delta(\zeta)=p\int\ln\left(  \zeta/\theta-\lambda\right)  \mathrm{d}\left(
\hat{F}\left(  \lambda\right)  -F_{\mathbf{c}}\left(  \lambda\right)  \right)
, \label{delta definition}\tag{SM141}%
\end{equation}
so that $\Delta(\zeta)=-2\ln g_{\mathrm{e}}(\zeta/\theta).$

\begin{lemma}
\label{gBound}Suppose that the null hypothesis holds, that is $\theta_{0}%
=0$.\textit{\ Then there exists a positive constant }$C_{1},$ such that for a
subset $\Theta$ of $\mathbb{C}$ that consists of all points whose Euclidean
distance from $\mathcal{\tilde{K}}$ is no larger than $C_{1},$\textit{ we
have}%
\[
\sup_{\zeta\in\Theta}\left\vert \Delta\left(  \zeta\right)  \right\vert
=O_{\mathrm{P}}(1)
\]
\textit{\ as }$\mathbf{n},p\rightarrow_{\boldsymbol{\gamma}}\infty$\textit{, where
}$O_{\mathrm{P}}(1)$\textit{\ is uniform with respect to }$\theta\in
\Omega\equiv\left(  0,\bar{\theta}-\varepsilon\right]  .$
\end{lemma}

\textbf{\small{Proof:}} Let us rewrite (\ref{delta definition}) in the following
equivalent form%
\[
\Delta(\zeta)=p\int\ln\left(  1-\lambda\theta/\zeta\right)  \mathrm{d}\left(
\hat{F}\left(  \lambda\right)  -F_{\mathbf{c}}\left(  \lambda\right)  \right)
.
\]
Statistic $\Delta(\zeta)$ is a special form of a linear spectral statistic
\[
\Delta(\varphi)=p\int\varphi\left(  \lambda\right)  \mathrm{d}\left(  \hat
{F}(\lambda)-F_{\mathbf{c}}\left(  \lambda\right)  \right)
\]
studied by Bai and Yao (2005), Bai and Silverstein (2004), and Zheng (2012)
for the cases of the Semi-circle, Marchenko-Pastur, and Wachter limiting
distributions, respectively. These papers note that%
\[
\Delta(\varphi)=-\frac{p}{2\pi\mathrm{i}}\int_{\mathcal{P}}\varphi\left(
\xi\right)  \left(  \hat{m}\left(  \xi\right)  -m_{\mathbf{c}}\left(
\xi\right)  \right)  \mathrm{d}\xi,
\]
where
\[
\hat{m}\left(  \xi\right)  =\int\frac{1}{\lambda-\xi}\mathrm{d}\hat{F}%
(\xi),\text{ }m_{\mathbf{c}}\left(  \xi\right)  =\int\frac{1}{\lambda-\xi
}\mathrm{d}F_{\mathbf{c}}\left(  \lambda\right)
\]
are the Stieltjes transforms of $\hat{F}$ and $F_{\mathbf{c}},$ and
$\mathcal{P}$ is a positively oriented contour in an open neighborhood of the
supports of $\hat{F}$ and $F_{\mathbf{c}},$ where $\varphi\left(  \xi\right)
$ is analytic, that encloses these supports. Theorem 2.1 and equation (2.3) of
Bai and Yao (2005) for SMD case, and Lemma 1.1 of Bai and Silverstein (2004)
for the rest of the cases, imply that if the distance from $\mathcal{P}$ to
the supports of $\hat{F}$ and $F_{\mathbf{c}}$ stays away from zero with
probability approaching one as $\mathbf{n},p\rightarrow_{\boldsymbol{\gamma}%
}\infty$, then
\[
\int_{\mathcal{P}}\left\vert p\left(  \hat{m}\left(  \xi\right)
-m_{\mathbf{c}}\left(  \xi\right)  \right)  \mathrm{d}\xi\right\vert
=O_{\mathrm{P}}\left(  1\right)  .
\]
(Throughout these notes, notation
$\int_{\mathcal{P}}\left\vert f(\xi) \mathrm{d} \xi \right\vert$
should be interpreted as $\int_{\alpha}^{\beta}\left\vert f(\mathcal{P}(t))\mathcal{P}^{\prime}(t)\right\vert \mathrm{d} t$,
where $\mathcal{P}$ is parameterized as a continuously differentiable complex function on
$[\alpha,\beta] \subseteq \mathbb{R}^{1}$.
For piecewise continuously differential pathes, $[\alpha,\beta]$ should be split into 
a finite number of sub-intervals where $\mathcal{P}$ is continuously differentiable.)
Therefore, for any $\delta>0,$ there exists $B>0,$ such that%
\begin{equation}
\Pr\left(  \left\vert \Delta(\varphi)\right\vert \leq B\sup_{\xi\in
\mathcal{P}}\left\vert \varphi\left(  \xi\right)  \right\vert \right)
>1-\delta\label{boundONdelta}\tag{SM142}%
\end{equation}
for all $\mathbf{n}$ and $p,$ where constant $B$ does not depend on $\varphi$.
Now, consider a family of functions $\varphi_{\zeta,\theta}\left(  \xi\right)
$
\[
\left\{  \varphi_{\zeta,\theta}\left(  \xi\right)  =\ln\left(  1-\xi
\theta/\zeta\right)  :\zeta\in\Theta\text{ and }\theta\in\Omega\right\}  .
\]
By the definitions of $\Theta$ and $\Omega,$ there exists an open neighborhood
$\mathcal{N}$ of the supports of $\hat{F}$ and $F_{\mathbf{c}}$ and a constant
$B_{1}$, such that, with probability arbitrarily close to one, for
sufficiently large $\mathbf{n}$ and $p,$ $\varphi_{\zeta,\theta}\left(
\xi\right)  $ are analytic in $\mathcal{N}$ for all $\zeta\in\Theta$ and
$\theta\in\Omega$ and
\[
\sup_{\theta\in\Omega}\sup_{\zeta\in\Theta}\sup_{\xi\in\mathcal{N}}\left\vert
\varphi_{\zeta,\theta}\left(  \xi\right)  \right\vert \leq B_{1}.
\]
Since $\Delta(\varphi_{\zeta,\theta})=\Delta(\zeta)$, we obtain from
(\ref{boundONdelta}) that for any $\delta>0,$ there exists $B_{2}>0$ such that
for sufficiently large $\mathbf{n}$ and $p,$%
\[
\Pr\left(  \sup_{\theta\in\Omega}\sup_{\zeta\in\Theta}\left\vert \Delta
(\zeta)\right\vert \leq B_{2}\right)  >1-\delta.
\]
In other words, $\sup_{\zeta\in\Theta}\left\vert \Delta(\zeta)\right\vert
=O_{\mathrm{P}}(1)$ uniformly over $\theta\in\Omega$. $\square$

Applying Lemma JO\ref{Olver}
%\ref{Olver} 
to the integral in (\ref{olver integral}) and using
the fact that $f(z_{0})=0,$ we obtain (\ref{I1case integral}). It remains to
show that $L_{2}\left(  \theta;\Lambda\right)  $ is asymptotically dominated
by $L_{1}\left(  \theta;\Lambda\right)  ,$ where
\[
L_{2}\left(  \theta;\Lambda\right)  =L\left(  \theta;\Lambda\right)
-L_{1}\left(  \theta;\Lambda\right)  .
\]

For \textbf{\small{SMD}}, \textbf{\small{PCA}}, and \textbf{\small{SigD}} we have%
\begin{align*}
\left\vert L_{2}\left(  \theta;\Lambda\right)  \right\vert  &  =\left\vert
\frac{\sqrt{\pi p}}{2\pi\mathrm{i}}\int_{\mathcal{K}_{2}\cup\mathcal{\bar{K}%
}_{2}}e^{-(p/2)\left(  f_{\mathrm{c}}+f_{\mathrm{h}}(z)\right)  }g_{\mathrm{c}}g_{\mathrm{h}}(z)%
{\displaystyle\prod_{j=1}^{p}}
\left(  z-\lambda_{j}\right)  ^{-1/2}\mathrm{d}z\right\vert \\
&  \leq\sqrt{\frac{p}{\pi}}e^{-(p/2)f_{\mathrm{c}}}g_{\mathrm{c}}\left(  2z_{0}\right)  ^{-p/2}%
\int_{\mathcal{K}_{2}}\left\vert e^{-(p/2)f_{\mathrm{h}}(z)}g_{\mathrm{h}}(z)\mathrm{d}%
z\right\vert \\
&  \leq\sqrt{\frac{p}{\pi}}e^{-(p/2)f_{\mathrm{c}}}g_{\mathrm{c}}\left(  2z_{0}\right)  ^{-p/2}%
\int_{-\infty}^{z_{0}}e^{-(p/2)f_{\mathrm{h}}(x)}g_{\mathrm{h}}(x)\mathrm{d}x.
\end{align*}
Explicitly evaluating the latter integral and using the exact form of $g_{\mathrm{c}}$,
available from Table JO\ref{Table 2a},
%\ref{Table 2a}
we obtain%
\[
\left\vert L_{2}\left(  \theta;\Lambda\right)  \right\vert \leq\frac{2C}%
{\sqrt{\pi p}}e^{-(p/2)f_{\mathrm{c}}}\left(  2z_{0}\right)  ^{-p/2}e^{-(p/2)f_{\mathrm{h}}(z_{0}%
)}\left(  1+o(1)\right)  ,
\]
where $o(1)$ does not depend on $\theta$, $C=1$ for SMD and PCA, and
$C=\sqrt{c_{1}+c_{2}}/r$ for SigD. Therefore,%
\begin{align*}
\left\vert L_{2}\left(  \theta;\Lambda\right)  \right\vert  &  \leq\frac
{2C}{\sqrt{\pi p}}e^{-(p/2)f\left(  z_{0}\right)  }\exp\left\{  -(p/2)\left(
\ln(2z_{0})-f_{\mathrm{e}}\left(  z_{0}\right)  \right)  \right\}  \left(
1+o(1)\right) \\
&  =\frac{2C}{\sqrt{\pi p}}\exp\left\{  -\frac{p}{2}\int\ln\left(
\frac{2z_{0}}{z_{0}-\lambda}\right)  \mathrm{d}F_{\mathbf{c}}(\lambda
)\right\}  \left(  1+o(1)\right)  ,
\end{align*}
where we used the fact that $f\left(  z_{0}\right)  =0$. But $\ln\left(
2z_{0}/\left(  z_{0}-\lambda\right)  \right)  $ is positive and bounded away
from zero uniformly over $\theta\in\left(  0,\bar{\theta}-\varepsilon\right]
$ with probability arbitrarily close to one, for sufficiently large
$\mathbf{n},p$. Hence, there exists a positive constant $K$ such that%
\[
\left\vert L_{2}\left(  \theta;\Lambda\right)  \right\vert \leq\frac{2C}%
{\sqrt{\pi p}}e^{-pK}\left(  1+o(1)\right)
\]
with probability arbitrarily close to one for sufficiently large
$\mathbf{n},p$. Combining this inequality with (\ref{I1case integral}), we
establish Theorem JO\ref{Icase}
%\ref{Icase}
 for SMD, PCA, and SigD.

For \textbf{\small{REG$_{0}$}}, we shall need the following lemma.

\begin{lemma}
For sufficiently large $\mathbf{n}$ and $p$, we have%
\begin{equation}
\left\vert \left.  _{0}F_{1}\left(  b-s;\Psi_{11}z\right)  \right.
\right\vert <4\sqrt{\pi m}\left\vert \exp\left\{  -m\varphi_{0}(t_{0}%
)\right\}  \right\vert \label{Bessel inequality}\tag{SM143}%
\end{equation}
for any $z$ and any $\theta>0$.
\end{lemma}

\textbf{\small{Proof:}} We use the identity (see formula 9.6.3 in Abramowitz and
Stegun (1964))%
\[
I_{m}\left(  \zeta\right)  =e^{-m\pi\mathrm{i}/2}J_{m}\left(  \mathrm{i}%
\zeta\right)  \text{ for }-\pi<\arg\zeta\leq\pi/2,
\]
where $J_{m}\left(  \cdot\right)  $ is the Bessel function. The identity and
(JO\ref{bessel rep})
%(\ref{bessel rep}) 
imply that%
\begin{equation}
\left.  _{0}F_{1}\left(  b-s;\Psi_{11}z\right)  \right.  =\Gamma\left(
m+1\right)  \left(  m^{2}\eta_{0}\right)  ^{-m/2}e^{-m\pi\mathrm{i}/2}%
J_{m}\left(  \mathrm{i}2m\eta_{0}^{1/2}\right)  . \label{new bessel}\tag{SM144}%
\end{equation}
On the other hand, for any $\zeta$ and any positive $K$,
\begin{equation}
\left\vert J_{K}\left(  K\zeta\right)  \right\vert \leq\left\{  1+\left\vert
\frac{\sin K\pi}{K\pi}\right\vert \right\}  \left\vert \left\{  \frac
{\zeta\exp\left\{  \sqrt{1-\zeta^{2}}\right\}  }{1+\sqrt{1-\zeta^{2}}%
}\right\}  ^{K}\right\vert , \label{Watson inequality}\tag{SM145}%
\end{equation}
(see Watson (1944), p. 270). The latter inequality, equation (\ref{new bessel}%
), and the Stirling formula for $\Gamma\left(  m+1\right)  $ imply that
(\ref{Bessel inequality}) holds for sufficiently large $m$, for any $z$ and
$\theta>0$. The constant $4$ on the right hand side of
(\ref{Bessel inequality}) is not the smallest possible one, but it is
sufficient for our purposes. $\square$

Using inequality (\ref{Bessel inequality}), we obtain for REG$_{0}$%
\begin{equation}
\left\vert L_{2}\left(  \theta;\Lambda\right)  \right\vert \leq4e^{-(p/2)f_{\mathrm{c}}%
}g_{\mathrm{c}}\sqrt{pm}\int_{\mathcal{K}_{2}}\left\vert \exp\left\{  -m\varphi
_{0}(t_{0})\right\}
%TCIMACRO{\dprod _{j=1}^{p}}%
%BeginExpansion
{\displaystyle\prod_{j=1}^{p}}
%EndExpansion
\left(  z-\lambda_{j}\right)  ^{-1/2}\mathrm{d}z\right\vert .
\label{first additional}\tag{SM146}%
\end{equation}
It is straightforward to verify that $\operatorname{Re}\varphi_{0}(t_{0})$ is
strictly increasing as $z$ is moving along $\mathcal{K}_{2}$ towards $-\infty
$. Therefore, for any $z\in\mathcal{K}_{2},$
\[
\operatorname{Re}\varphi_{0}\left(  t_{0}(z)\right)  >\operatorname{Re}%
\varphi_{0}(t_{0}(\bar{z})),
\]
where $\bar{z}=z_{1}+\mathrm{i}\left(  z_{0}-z_{1}\right)  $ is the point of
$\mathcal{K}_{2}$ where $\mathcal{K}_{2}$ meets $\mathcal{K}_{1}.$ The latter
inequality together with (\ref{first additional}) yields%
\[
\left\vert L_{2}\left(  \theta;\Lambda\right)  \right\vert \leq
4e^{-(p/2)\operatorname{Re}f(\bar{z})}g_{\mathrm{c}}\left\vert g_{\mathrm{e}}\left(  \bar
{z}\right)  \right\vert \sqrt{pm}\int_{\mathcal{K}_{2}}%
%TCIMACRO{\dprod _{j=1}^{p}}%
%BeginExpansion
{\displaystyle\prod_{j=1}^{p}}
%EndExpansion
\left\vert \frac{\bar{z}-\lambda_{j}}{z-\lambda_{j}}\right\vert ^{1/2}%
\left\vert \mathrm{d}z\right\vert .
\]
Since, for some constant $\tau_{1}$, $\operatorname{Re}f\left(  \bar
{z}\right)  >f\left(  z_{0}\right)  +\tau_{1}=\tau_{1}$ and since, by Lemma
\ref{gBound}, $4g_{\mathrm{e}}\left(  \bar{z}\right)  =O_{\mathrm{P}}\left(  1\right)
$ uniformly over $\theta\in\left(  0,\bar{\theta}-\varepsilon\right]  ,$ we
obtain%
\begin{equation}
\left\vert L_{2}\left(  \theta;\Lambda\right)  \right\vert \leq e^{-(p/2)\tau_{1}%
}g_{\mathrm{c}}\sqrt{pm}\int_{\mathcal{K}_{2}}\prod_{j=1}^{p}\left\vert \frac{\bar
{z}-\lambda_{j}}{z-\lambda_{j}}\right\vert ^{1/2}\left\vert \mathrm{d}%
z\right\vert O_{\mathrm{P}}\left(  1\right)  . \label{to combine with}\tag{SM147}%
\end{equation}
Note that for any $z\in\mathcal{K}_{2}$ and any $j=1,...,p,$ $\left\vert
\left(  \bar{z}-\lambda_{j}\right)  /\left(  z-\lambda_{j}\right)  \right\vert
\leq1$ and $\left\vert z-\lambda_{j}\right\vert >\left\vert z\right\vert .$
Further, since $z_{0}<\left\vert \bar{z}\right\vert $ and with probability
arbitrary close to one, for sufficiently large $\mathbf{n}$ and $p,$
$\lambda_{1}<z_{0},$ we have $\left\vert \bar{z}-\lambda_{j}\right\vert
<\left\vert \bar{z}-z_{0}\right\vert <2\left\vert \bar{z}\right\vert .$ Thus,
for $p\geq4,$ we have%
\[
\int_{\mathcal{K}_{2}}\prod_{j=1}^{p}\left\vert \frac{\bar{z}-\lambda_{j}%
}{z-\lambda_{j}}\right\vert ^{1/2}\left\vert \mathrm{d}z\right\vert \leq
\int_{\mathcal{K}_{2}}4\left\vert z/\bar{z}\right\vert ^{-2}\left\vert
\mathrm{d}z\right\vert =\left\vert \bar{z}\right\vert O(1)
\]
Combining this with (\ref{to combine with}) and noting that $g_{\mathrm{c}}\left\vert
\bar{z}\right\vert =O\left(  1\right)  $ uniformly over $\theta\in\left(
0,\bar{\theta}-\varepsilon\right]  ,$ we obtain
\begin{equation}
\left\vert L_{2}\left(  \theta;\Lambda\right)  \right\vert \leq\sqrt
{pm}e^{-(p/2)\tau_{1}}O_{\mathrm{P}}\left(  1\right)  , \label{REG0inequality}\tag{SM148}%
\end{equation}
where $O_{\mathrm{P}}\left(  1\right)  $ is uniform with respect to $\theta
\in\left(  0,\bar{\theta}-\varepsilon\right]  $. Theorem JO\ref{Icase}
%\ref{Icase} 
for
REG$_{0}$ follows from the latter equality and (\ref{I1case integral}).

For \textbf{\small{REG}} and \textbf{\small{CCA}}, the Theorem follows from (\ref{I1case integral}) and
inequalities%
\begin{equation}
\left\vert L_{2}\left(  \theta;\Lambda\right)  \right\vert \leq pe^{-p\tau
_{2}}O_{\mathrm{P}}\left(  1\right)  , \label{restInequalities}\tag{SM149}%
\end{equation}
where $\tau_{2}$ is a positive constant. We obtain (\ref{restInequalities}) by
combining the method used to derive (\ref{REG0inequality}) with upper bounds
on $_{1}F_{1}$ and $_{2}F_{1},$ which we establish using the integral
representations (JO\ref{integralrepF11}).
%(\ref{integralrepF11})
$\square$ \vspace{3 mm}

\textbf{\small{A proof of the domination of $L_{2}\left( \theta ;\Lambda
\right) $ by $L_{1}\left( \theta ;\Lambda \right) $ (via establishing (\ref{restInequalities})).}}
By definition, we have 
\begin{equation}
L_{2}\left( \theta ;\Lambda \right) =\sqrt{\pi p}\exp \left\{ -\frac{p}{2}%
\left( f_{\mathrm{c}}+f_{\mathrm{e}}(z_{0}\right) )\right\} \frac{g_{\mathrm{%
c}}g_{\mathrm{e}}(z_{0})}{2\pi \mathrm{i}}\int_{\mathcal{K}_{2}\cup \mathcal{%
\bar{K}}_{2}}F_{j}\prod_{j=1}^{p}\left( \frac{z_{0}-\lambda _{j}}{z-\lambda
_{j}}\right) ^{1/2}\mathrm{d}z  \label{L2}\tag{SM150}
\end{equation}%
with $j=1$ for REG and $j=2$ for CCA. The idea of the proof is to use the
integral representations (\ref{changed variables}), that is%
\begin{equation*}
F_{j}=\frac{C_{m}\eta _{j}^{-m}}{2\pi \mathrm{i}}\int_{0}^{(\eta _{j}+)}\exp
\left\{ -m\phi _{j}\left( \tau \right) \right\} \chi _{j}\left( \tau \right) 
\mathrm{d}\tau ,
\end{equation*}%
to find simple upper bounds for $\left\vert F_{j}\right\vert $ corresponding
to $z\in \mathcal{K}_{2}\cup \mathcal{\bar{K}}_{2}$. Note that since $%
F_{j}\left( \bar{z}\right) =\overline{F_{j}\left( z\right) },$ it is
sufficient to establish the bounds for $z\in \mathcal{K}_{2}$. These upper
bounds will then be used to estimate the integral in (\ref{L2}) from above,
and eventually to establish the domination of $L_{2}\left( \theta ;\Lambda
\right) $ by $L_{1}\left( \theta ;\Lambda \right) $. \vspace{2 mm}

\textbf{\small{REG.}}
\begin{lemma}
\label{oldcontour}Let $\tau _{+}\in \mathcal{C}_{2}$ and $z$ be the
corresponding point of $\mathcal{K}_{2}.$ Then $\operatorname{Re}f_{\mathrm{h}%
}(z)>f_{\mathrm{h}}(z_{0})+\alpha ,$ where $\alpha >0$ does not depend on $%
\tau _{+}\in \mathcal{C}_{2}$ and does not depend on $\theta $.
\end{lemma}

\textbf{\small{Proof:}} Parametrize points $\tau _{+}\in \mathcal{C}_{2}$ as 
\begin{equation}
\tau _{+}=-\kappa -x+\left\vert \tau _{0}+\kappa \right\vert \exp \left\{ 
\mathrm{i}\pi /2\right\} ,  \label{parametrization of tau+}\tag{SM151}
\end{equation}%
$x\geq 0$. As $x$ goes from $0$ to $\infty ,$ the corresponding $z$ tracks
contour $\mathcal{K}_{2}$ from the point $\zeta ,$ where $\mathcal{K}_{2}$
and $\mathcal{K}_{1}$ meet, to $-\infty $. Recall that%
\begin{equation}
-\frac{p}{2}f_{\mathrm{h}}(z)=-m\left( \varphi _{1}(t_{1})+k\right)
=-m\left( \phi _{1}\left( \tau _{+}\right) +\ln \eta _{1}+k\right) ,
\label{connections}\tag{SM152}
\end{equation}%
where $k=\kappa \ln \kappa -\left( \kappa -1\right) \ln \left( \kappa
-1\right) .$ Using the definition (\ref{faij}) of $\phi _{1}$ and the
identity%
\begin{equation}
\eta _{1}=\frac{\tau _{+}\left( \tau _{+}+1\right) }{\tau _{+}+\kappa },
\label{identityforeta1again}\tag{SM153}
\end{equation}%
we obtain%
\begin{equation*}
\frac{p}{2m}\operatorname{Re}f_{\mathrm{h}}(z)=-\operatorname{Re}\tau _{+}+\ln \left\vert
\tau _{+}+1\right\vert -\kappa \ln \left\vert \tau _{+}+\kappa \right\vert
+\kappa \ln \kappa .
\end{equation*}%
Taking the derivative of both sides of the latter equality with respect to $%
x,$ we obtain%
\begin{equation*}
\frac{p}{2m}\frac{\mathrm{d}}{\mathrm{d}x}\operatorname{Re}f_{\mathrm{h}}(z)=1+\frac{%
x+\kappa -1}{\left\vert \tau _{+}+1\right\vert ^{2}}-\frac{\kappa x}{%
\left\vert \tau _{+}+\kappa \right\vert ^{2}}.
\end{equation*}%
For $x\geq 0,$ we have%
\begin{eqnarray*}
\left\vert \tau _{+}+\kappa \right\vert  &\equiv &\left\vert -x+\left\vert
\tau _{0}+\kappa \right\vert \exp \left\{ \mathrm{i}\pi /2\right\}
\right\vert >x\text{ and} \\
\left\vert \tau _{+}+\kappa \right\vert  &\equiv &\left\vert -x+\left\vert
\tau _{0}+\kappa \right\vert \exp \left\{ \mathrm{i}\pi /2\right\}
\right\vert >\kappa .
\end{eqnarray*}%
Therefore, $\kappa x/\left\vert \tau _{+}+\kappa \right\vert ^{2}<1$ and $%
\frac{\mathrm{d}}{\mathrm{d}x}\operatorname{Re}f_{\mathrm{h}}(z)>0$. This implies
that%
\begin{equation*}
\operatorname{Re}f_{\mathrm{h}}(z)>\operatorname{Re}f_{\mathrm{h}}(\zeta ).
\end{equation*}%
On the other hand, as shown in subsection \ref{sec: contours fo steep descent} (pp 38-40 of these notes), 
$\operatorname{Re}f_{\mathrm{h}}(z)$ strictly
increases as $z$ moves along $\mathcal{K}_{1}$ from $z_{0}$ to $\zeta $.
Hence, there exists $\alpha >0$ that does not depend on $\tau _{+}\in 
\mathcal{C}_{2}$, such that 
\begin{equation*}
\operatorname{Re}f_{\mathrm{h}}(z)>\operatorname{Re}f_{\mathrm{h}}(z_{0})+\alpha =f_{\mathrm{%
h}}(z_{0})+\alpha .
\end{equation*}%
From the definitions of $\mathcal{C}_{1}$ (the image of which under $\tau
\mapsto z$ transformation is $\mathcal{K}_{1}$) and of $f_{\mathrm{h}}(z),$
it is easy to see that $\alpha $ can be chosen so that it does not depend on 
$\theta $ as well. $\square $

\begin{lemma}
\label{BoundF1}There exist positive constants $\alpha $ and $\alpha _{1}$
that do not depend on $\theta $ such that, for any $\tau _{+}\in \mathcal{C}%
_{2}$%
\begin{equation}
\left\vert F_{1}\right\vert \leq \alpha _{1}\sqrt{p}\left\vert \eta
_{1}\right\vert \exp \left\{ -\frac{p}{2}\left( f_{\mathrm{h}}(z_{0})+\alpha
\right) \right\} .  \label{LemmasEquation}\tag{SM154}
\end{equation}
\end{lemma}

\textbf{\small{Proof:}} Let $\tau _{+}\in \mathcal{C}_{2}$ and $z$ be the
corresponding point of $\mathcal{K}_{2}.$ Choose the contour in the integral
representation (\ref{changed variables}) of $F_{1}$ as in subsection 4.3 
(that contains a proof of Lemma JO3) of this note.
We shall call such a contour $%
\mathcal{K}^{\ast }$. As explained in subsection 4.3, the minimum of $\operatorname{Re}%
\phi _{1}\left( \tau \right) $ over $\tau \in \mathcal{K}^{\ast }$ is
achieved either at $\tau _{+}$ or, in some cases corresponding to situation
3, at $\tau ^{\ast }$ that belongs to $\left[ 0,A\right] $ and is such that $%
\operatorname{Re}\tau ^{\ast }\leq -\kappa $ (see a discussion around equation (\ref%
{more than right}), which shows that points $\tau ^{\ast }\in \left[ 0,A%
\right] $ with $\operatorname{Re}\tau ^{\ast }>-\kappa $ cannot correspond to the
minimum of $\operatorname{Re}\phi _{1}\left( \tau \right) $ over $\tau \in \mathcal{K%
}^{\ast }$).

If the minimum of $\operatorname{Re}\phi _{1}\left( \tau \right) $ over $\tau \in 
\mathcal{K}^{\ast }$ is achieved at $\tau _{+}$ then using (\ref{changed
variables}), (\ref{connections}), and the Stirling's approximation%
\begin{equation}
C_{m}=\frac{\sqrt{\pi p\left( 1-c_{1}\right) }}{r}\exp \left\{ m\left(
\kappa -1\right) \ln \left( \kappa -1\right) -m\kappa \ln \kappa \right\}
\left( 1+o(1)\right) ,  \label{Stir}\tag{SM155}
\end{equation}%
we obtain, for some $\tilde{\alpha}>0$ that does not depend on $\tau _{+}$
and on $\theta ,$%
\begin{equation}
\left\vert F_{1}\right\vert \leq \tilde{\alpha}\sqrt{p}\exp \left\{ -\frac{p%
}{2}\operatorname{Re}f_{\mathrm{h}}(z)\right\} \int_{\mathcal{K}^{\ast }}\left\vert
\chi _{1}\left( \tau \right) \mathrm{d}\tau \right\vert .  \label{firstineq}\tag{SM156}
\end{equation}%
Recall that $\chi _{1}\left( \tau \right) =\left( \tau -\eta _{1}\right)
^{-1}.$ By definition of $\mathcal{K}^{\ast },$%
\begin{equation}
\sup_{\tau \in \mathcal{K}^{\ast }}\left\vert \chi _{1}\left( \tau \right)
\right\vert \leq \max \left\{ \left\vert \tau _{+}-\eta _{1}\right\vert
^{-1},\left\vert \eta _{1}\right\vert ^{-1}\right\} \text{ and }\left\vert 
\mathcal{K}^{\ast }\right\vert \leq \left\vert \eta _{1}\right\vert +2\pi
\left\vert \tau _{+}-\eta _{1}\right\vert .  \label{secondineq}\tag{SM157}
\end{equation}%
Identity (\ref{identityforeta1again}) implies that $\left\vert \tau
_{+}-\eta _{1}\right\vert =\left( \kappa -1\right) \left\vert \tau
_{+}/\left( \tau _{+}+\kappa \right) \right\vert $ is bounded away from zero
uniformly with respect to $\tau _{+}\in \mathcal{C}_{2}.$ Therefore, (\ref%
{firstineq}) and (\ref{secondineq}) imply that there exists $\alpha _{1}>0$
that does not depend on $\tau _{+}$ and on $\theta $ such that%
\begin{equation*}
\left\vert F_{1}\right\vert \leq \alpha _{1}\sqrt{p}\left\vert \eta
_{1}\right\vert \exp \left\{ -\frac{p}{2}\operatorname{Re}f_{\mathrm{h}}(z)\right\} .
\end{equation*}%
Combining this with Lemma \ref{oldcontour}, we obtain (\ref{LemmasEquation}).

If the minimum of $\operatorname{Re}\phi _{1}\left( \tau \right) $ over $\tau \in 
\mathcal{K}^{\ast }$ is achieved at $\tau ^{\ast }$ then we must be in
situation 3 so that $\left\vert \tau ^{\ast }-\eta _{1}\right\vert
>\left\vert \tau _{+}-\eta _{1}\right\vert $ and%
\begin{eqnarray*}
\operatorname{Re}\phi _{1}\left( \tau ^{\ast }\right)  &=&-\operatorname{Re}\tau ^{\ast
}-\kappa \ln \left\vert \tau ^{\ast }\right\vert +\left( \kappa -1\right)
\ln \left\vert \tau ^{\ast }-\eta _{1}\right\vert  \\
&>&-\operatorname{Re}\tau ^{\ast }-\kappa \ln \left\vert \tau ^{\ast }\right\vert
+\left( \kappa -1\right) \ln \left\vert \tau _{+}-\eta _{1}\right\vert .
\end{eqnarray*}%
Let $\tau $ be any point on the ray starting at $0$ and passing through $%
\eta _{1},$ let $\arg \eta _{1}=\beta $ (note that $\beta >\pi /2$ so that $%
\cos \beta <0$), and let $x=\left\vert \tau \right\vert .$ Then%
\begin{equation*}
-\operatorname{Re}\tau ^{\ast }-\kappa \ln \left\vert \tau ^{\ast }\right\vert \geq
-\max_{x\geq 0}\left\{ x\cos \beta +\kappa \ln x\right\} =\kappa -\kappa \ln
\left( -\kappa /\cos \beta \right) .
\end{equation*}%
Therefore,%
\begin{equation*}
\operatorname{Re}\phi _{1}\left( \tau ^{\ast }\right) >\kappa -\kappa \ln \left(
-\kappa /\cos \beta \right) +\left( \kappa -1\right) \ln \left\vert \tau
_{+}-\eta _{1}\right\vert .
\end{equation*}%
This inequality implies%
\begin{equation*}
\operatorname{Re}\phi _{1}\left( \tau ^{\ast }\right) +\ln \left\vert \eta
_{1}\right\vert >\kappa -\kappa \ln \left( -\kappa /\cos \beta \right)
+\left( \kappa -1\right) \ln \left\vert \tau _{+}-\eta _{1}\right\vert +\ln
\left\vert \eta _{1}\right\vert .
\end{equation*}%
Using (\ref{identityforeta1again}) and the fact that $\tau _{+}\in \mathcal{C%
}_{2}$, we obtain%
\begin{align}
\operatorname{Re}\phi _{1}\left( \tau ^{\ast }\right) +\ln \left\vert \eta
_{1}\right\vert  &>\kappa -\kappa \ln \left( -\kappa /\cos \beta \right)
+\kappa \ln \left\vert \frac{\tau _{+}}{\tau _{+}+\kappa }\right\vert  
\notag \\
&+\ln \left\vert \tau _{+}+1\right\vert +\left( \kappa -1\right) \ln \left(
\kappa -1\right)   \notag \\
&>\kappa -\kappa \ln \left( -\kappa /\cos \beta \right) +\ln \left\vert 
\bar{\tau}_{+}+1\right\vert +\left( \kappa -1\right) \ln \left( \kappa
-1\right) ,  \label{another inequality}\tag{SM158}
\end{align}%
where $\bar{\tau}_{+}=-\kappa +\left\vert \tau _{0}+\kappa \right\vert \exp
\left\{ \mathrm{i}\pi /2\right\} $ is the point where $\mathcal{C}_{2}$ and $%
\mathcal{C}_{1}$ meet. On the other hand,%
\begin{equation*}
\cos \beta \leq \cos \arg \bar{\tau}_{+}=-\kappa /\left\vert \bar{\tau}%
_{+}\right\vert 
\end{equation*}%
and thus%
\begin{eqnarray*}
\kappa -\kappa \ln \left( -\kappa /\cos \beta \right)  &>&\kappa -\kappa \ln
\left\vert \bar{\tau}_{+}\right\vert =\kappa +\kappa \ln \left( \left\vert 
\bar{\tau}_{+}+\kappa \right\vert /\left\vert \bar{\tau}_{+}\right\vert
\right) -\kappa \ln \left\vert \bar{\tau}_{+}+\kappa \right\vert  \\
&>&\kappa +\kappa \ln \left( 1/\sqrt{2}\right) -\kappa \ln \left\vert \bar{%
\tau}_{+}+\kappa \right\vert >-\kappa \ln \left\vert \bar{\tau}_{+}+\kappa
\right\vert .
\end{eqnarray*}%
Using this inequality and (\ref{another inequality}), we obtain%
\begin{equation*}
\operatorname{Re}\phi _{1}\left( \tau ^{\ast }\right) +\ln \left\vert \eta
_{1}\right\vert >-\kappa \ln \left\vert \bar{\tau}_{+}+\kappa \right\vert
+\ln \left\vert \bar{\tau}_{+}+1\right\vert +\left( \kappa -1\right) \ln
\left( \kappa -1\right) .
\end{equation*}%
Since $\left\vert \tau +\kappa \right\vert $ stays constant for $\tau \in 
\mathcal{C}_{1}$ whereas $\left\vert \tau +1\right\vert $ is strictly
decreasing as $\tau $ moves along $\mathcal{C}_{1}$ from $\tau _{0}$ to $%
\bar{\tau}_{+},$ there exists $\alpha _{2}>0$ which is independent of $\tau
_{+}$ and $\theta ,$ such that%
\begin{eqnarray*}
\operatorname{Re}\phi _{1}\left( \tau ^{\ast }\right) +\ln \left\vert \eta
_{1}\right\vert  &>&-\kappa \ln \left( \tau _{0}+\kappa \right) +\ln \left(
\tau _{0}+1\right) +\left( \kappa -1\right) \ln \left( \kappa -1\right)
+\alpha _{2} \\
&>&-\tau _{0}-\kappa \ln \left( \tau _{0}+\kappa \right) +\ln \left( \tau
_{0}+1\right) +\left( \kappa -1\right) \ln \left( \kappa -1\right) +\alpha
_{2} \\
&=&\operatorname{Re}\phi _{1}\left( \tau _{0}\right) +\ln \left\vert \eta
_{10}\right\vert +\alpha _{2},
\end{eqnarray*}%
where $\eta _{10}$ is the value of $\eta _{1}$ that corresponds to $z_{0}$.
Therefore, by (\ref{connections}), we have%
\begin{equation*}
-m\left( \operatorname{Re}\phi _{1}\left( \tau ^{\ast }\right) +\ln \left\vert \eta
_{1}\right\vert \right) <-m\left( \frac{p}{2m}f_{\mathrm{h}}(z_{0})+\alpha
_{2}-k\right) .
\end{equation*}%
Using this inequality together with (\ref{changed variables}) and (\ref{Stir}%
), we obtain that, for some $\tilde{\alpha}>0$ that does not depend on $\tau
_{+}$ and $\theta ,$%
\begin{equation*}
\left\vert F_{1}\right\vert \leq \tilde{\alpha}\sqrt{p}\exp \left\{ -\frac{p%
}{2}\left( f_{\mathrm{h}}(z_{0})+\frac{2m}{p}\alpha _{2}\right) \right\}
\int_{\mathcal{K}^{\ast }}\left\vert \chi _{1}\left( \tau \right) \mathrm{d}%
\tau \right\vert .
\end{equation*}%
Analysing the integral $\int_{\mathcal{K}^{\ast }}\left\vert \chi _{1}\left(
\tau \right) \mathrm{d}\tau \right\vert $ as above, we conclude that there
exist $\alpha ,\alpha _{1}>0$ that do not depend on $\tau _{+}$ and $\theta ,
$ such that (\ref{LemmasEquation}) holds. $\square $

Using Lemma \ref{BoundF1} and equation (\ref{L2}), we obtain the following
bound on $\left\vert L_{2}\left( \theta ;\Lambda \right) \right\vert $%
\begin{equation}
\left\vert L_{2}\left( \theta ;\Lambda \right) \right\vert \leq \alpha
_{1}p\exp \left\{ -\frac{p}{2}\alpha \right\} \left\vert g_{\mathrm{c}}g_{%
\mathrm{e}}(z_{0})\right\vert \int_{\mathcal{K}_{2}}\left\vert \eta
_{1}\prod_{j=1}^{p}\left( \frac{z_{0}-\lambda _{j}}{z-\lambda _{j}}\right)
^{1/2}\mathrm{d}z\right\vert .  \label{firstbound}\tag{SM159}
\end{equation}%
On the other hand, for any $\lambda _{j}$ from the support of $F_{\mathbf{c}%
},$ we have%
\begin{equation}
\left\vert \frac{z_{0}-\lambda _{j}}{z-\lambda _{j}}\right\vert <\left\vert 
\frac{z_{0}}{z}\right\vert =\left\vert \frac{\eta _{10}}{\eta _{1}}%
\right\vert =\left\vert \frac{\tau _{0}\left( \tau _{0}+1\right) \left( \tau
_{+}+\kappa \right) }{\left( \tau _{0}+\kappa \right) \tau _{+}\left( \tau
_{+}+1\right) }\right\vert   \label{multifloor}\tag{SM160}
\end{equation}%
and 
\begin{equation}
\mathrm{d}z=\frac{c_{1}\left( 1-c_{1}\right) }{\theta c_{2}}\mathrm{d}\eta
_{1}=\frac{c_{1}\left( 1-c_{1}\right) }{\theta c_{2}}\left( 1-\frac{\kappa
\left( \kappa -1\right) }{\left( \tau _{+}+\kappa \right) ^{2}}\right) 
\mathrm{d}\tau _{+}.  \label{dz}\tag{SM161}
\end{equation}%
Note that, for any $\tau _{+}\in \mathcal{C}_{2}$%
\begin{equation*}
\left\vert \frac{\kappa \left( \kappa -1\right) }{\left( \tau _{+}+\kappa
\right) ^{2}}\right\vert <\frac{\kappa \left( \kappa -1\right) }{\left( \tau
_{0}+\kappa \right) ^{2}}.
\end{equation*}%
A direct calculation based on the definitions 
\begin{eqnarray*}
\tau _{0} &=&\tfrac{1}{2}\left\{ \eta _{10}-1+\sqrt{\left( \eta
_{10}-1\right) ^{2}+4\kappa \eta _{10}}\right\} , \\
\eta _{10} &=&\frac{z_{0}\theta c_{2}}{c_{1}\left( 1-c_{1}\right) },\text{ }%
\kappa =\frac{c_{1}+c_{2}-c_{1}c_{2}}{c_{2}\left( 1-c_{1}\right) },\text{ and%
} \\
z_{0} &=&\frac{\left( 1+\theta \right) \left( \theta +c_{1}\right) }{\theta
\left( 1+\left( 1+\theta \right) c_{2}/c_{1}\right) }
\end{eqnarray*}%
yields%
\begin{eqnarray*}
\tau _{0} &=&\frac{\theta +c_{1}}{1-c_{1}}\text{, }\eta _{10}=\frac{%
c_{2}\left( \theta +1\right) \left( \theta +c_{1}\right) }{\left(
c_{1}+c_{2}+\theta c_{2}\right) \left( 1-c_{1}\right) },\text{ and} \\
\frac{\kappa \left( \kappa -1\right) }{\left( \tau _{0}+\kappa \right) ^{2}}
&=&c_{1}\frac{c_{1}+c_{2}-c_{1}c_{2}}{\left( c_{1}+c_{2}+\theta c_{2}\right)
^{2}}.
\end{eqnarray*}%
The latter two equalities together with (\ref{dz}) imply that there exists a
constant $\alpha _{2}>0$ that does not depend on $\theta \in \left( 0,\bar{%
\theta}-\varepsilon \right] $ such that (for sufficiently large $\mathbf{n},p
$ as $\mathbf{n},p\rightarrow _{\boldsymbol{\gamma }}\infty $) 
\begin{equation*}
\left\vert \eta _{10}\mathrm{d}z\right\vert <\frac{\alpha _{2}}{\theta }%
\left\vert \mathrm{d}\tau _{+}\right\vert .
\end{equation*}%
Using this and (\ref{multifloor}) in (\ref{firstbound}), we obtain%
\begin{equation*}
\left\vert L_{2}\left( \theta ;\Lambda \right) \right\vert \leq \alpha
_{1}\alpha _{2}p\exp \left\{ -\frac{p}{2}\alpha \right\} \left\vert \frac{g_{%
\mathrm{c}}}{\theta }g_{\mathrm{e}}(z_{0})\right\vert \int_{\mathcal{C}%
_{2}}\left\vert \frac{\tau _{0}\left( \tau _{0}+1\right) \left( \tau
_{+}+\kappa \right) }{\left( \tau _{0}+\kappa \right) \tau _{+}\left( \tau
_{+}+1\right) }\right\vert ^{p/2-1}\left\vert \mathrm{d}\tau _{+}\right\vert
.
\end{equation*}%
Note that for any $\tau _{+}\in \mathcal{C}_{2}$ we have $\left\vert \tau
_{+}+1\right\vert >\left\vert \tau _{+}+\kappa \right\vert .$ On the other
hand, $\tau _{0}+1<\tau _{0}+\kappa .$ Therefore,%
\begin{equation*}
\left\vert \frac{\left( \tau _{0}+1\right) \left( \tau _{+}+\kappa \right) }{%
\left( \tau _{0}+\kappa \right) \left( \tau _{+}+1\right) }\right\vert <1
\end{equation*}%
and%
\begin{equation*}
\left\vert L_{2}\left( \theta ;\Lambda \right) \right\vert \leq \alpha
_{1}\alpha _{2}p\exp \left\{ -\frac{p}{2}\alpha \right\} \left\vert \frac{g_{%
\mathrm{c}}}{\theta }g_{\mathrm{e}}(z_{0})\right\vert \int_{\mathcal{C}%
_{2}}\left\vert \frac{\tau _{0}}{\tau _{+}}\right\vert ^{p/2-1}\left\vert 
\mathrm{d}\tau _{+}\right\vert .
\end{equation*}%
Using parameterization (\ref{parametrization of tau+}), we obtain%
\begin{eqnarray*}
\left\vert L_{2}\left( \theta ;\Lambda \right) \right\vert  &\leq &\alpha
_{1}\alpha _{2}p\exp \left\{ -\frac{p}{2}\alpha \right\} \left\vert \frac{g_{%
\mathrm{c}}}{\theta }g_{\mathrm{e}}(z_{0})\right\vert \int_{0}^{\infty
}\left\vert \frac{\tau _{0}}{x+\kappa -\mathrm{i}\left\vert \tau _{0}+\kappa
\right\vert }\right\vert ^{p/2-1}\mathrm{d}x \\
&\leq &\alpha _{1}\alpha _{2}\alpha _{3}p\exp \left\{ -\frac{p}{2}\alpha
\right\} \left\vert \frac{g_{\mathrm{c}}}{\theta }g_{\mathrm{e}%
}(z_{0})\right\vert 
\end{eqnarray*}%
for some $\alpha _{3}>0$ that does not depend on $\theta .$ Finally, note
that $g_{\mathrm{c}}/\theta =O(1)$ and $g_{\mathrm{e}}(z_{0})=O_{\mathrm{P}%
}(1)$, so that the above display implies equation (\ref{restInequalities}). Since 
\begin{equation*}
L_{1}\left( \theta ;\Lambda \right) =\frac{g(z_{0})}{\sqrt{-\textrm{d}%
^{2}f(z_{0})/\textrm{d}z^{2}}}+O_{\mathrm{P}}\left( p^{-1}\right) ,
\end{equation*}%
we see that $L_{2}\left( \theta ;\Lambda \right) $ is asymptotically
dominated by $L_{1}\left( \theta ;\Lambda \right) $. \vspace{2 mm}

\textbf{\small{CCA.}}
Let $A$ be an arbitrarily large positive constant. Split the
contour $\mathcal{K}_{2}$ into $\mathcal{K}_{21}$ and $\mathcal{K}_{22},$
where 
\begin{equation*}
\mathcal{K}_{21}=\left\{ z:z\in \mathcal{K}_{2},\operatorname{Re}z>-A\right\} .
\end{equation*}%
Note that the approximation 
\begin{equation*}
F_{2}=C_{m}\psi _{2}\left( t_{2}\right) e^{-\mathrm{i}\omega
_{2}/2}\left\vert 2\pi m\varphi _{2}^{\prime \prime }\left( t_{2}\right)
\right\vert ^{-1/2}\exp \left\{ -m\varphi _{2}\left( t_{2}\right) \right\}
\left( 1+o(1)\right) 
\end{equation*}%
derived in Lemma JO\ref{1F1approximation} remains valid for $z\in \mathcal{K}_{21}.$
Therefore, the representation%
\begin{equation*}
L\left( \theta ;\Lambda \right) =\sqrt{\pi p}\frac{1}{2\pi \mathrm{i}}\int_{%
\mathcal{K}}\exp \left\{ -\frac{p}{2}f(z)\right\} g(z)\mathrm{d}z
\end{equation*}%
is valid for $z\in \mathcal{K}_{21}\cup \mathcal{K}_{1}.$ Hence, if we show
that $\mathcal{K}_{21}\cup \mathcal{K}_{1}$ is a contour of steep descent
for $-\operatorname{Re}f(z),$ then%
\begin{equation*}
L_{21}\left( \theta ;\Lambda \right) +L_{1}\left( \theta ;\Lambda \right) 
\end{equation*}%
must be asymptotically equivalent to $L_{1}\left( \theta ;\Lambda \right) $,
where 
\begin{equation*}
L_{21}\left( \theta ;\Lambda \right) =\sqrt{\pi p}\frac{1}{2\pi \mathrm{i}}%
\int_{\mathcal{K}_{21}\cup \mathcal{\bar{K}}_{21}}\exp \left\{ -\frac{p}{2}%
f(z)\right\} g(z)\mathrm{d}z,
\end{equation*}%
and thus, $L_{21}\left( \theta ;\Lambda \right) $ must be asymptotically
dominated by $L_{1}\left( \theta ;\Lambda \right) .$

Obviously, $-\operatorname{Re}f_{\mathrm{e}}(z)$ is decreasing as $z$ moves along $%
\mathcal{K}_{21}$ so that $\operatorname{Re}z$ becomes more and more negative. Let
us consider the behavior of 
\begin{equation}
-\operatorname{Re}f_{\mathrm{h}}(z)=\frac{1-c_{1}}{c_{1}}\left( -\varphi _{2}\left(
t_{2}\right) -\kappa \ln \kappa +\left( \kappa -1\right) \ln \left( \kappa
-1\right) \right) .  \label{identityagain}\tag{SM162}
\end{equation}%
Recall (\ref{newf2CCA}), that states 
\begin{equation*}
\operatorname{Re}\varphi _{2}\left( t_{2}\right) =-2\kappa \ln \left\vert
t_{2}\right\vert +\left( 2\kappa -1\right) \ln \left\vert t_{2}-1\right\vert
+\kappa \ln \frac{\kappa }{\kappa -1}.
\end{equation*}%
Parametrize $z\in \mathcal{K}_{21}$ as%
\begin{equation*}
z=z_{1}-x\left\vert z_{0}-z_{1}\right\vert +\left\vert
z_{0}-z_{1}\right\vert \mathrm{i,}\text{ }x\in \left[ 0,\left(
A+z_{1}\right) /\left\vert z_{0}-z_{1}\right\vert \right] 
\end{equation*}%
where 
\begin{equation*}
z_{1}=-\frac{c_{1}\left( 1-c_{1}\right) ^{2}l\left( \theta \right) }{4\theta
r^{2}}.
\end{equation*}%
For the corresponding $\eta _{2}=z\theta c_{2}^{2}/\left[ c_{1}^{2}l\left(
\theta \right) \right] $ we have 
\begin{equation*}
\eta _{2}=R_{0}-xR_{1}+R_{1}\mathrm{i,}\text{ }x\in \left[ 0,\left(
A+z_{1}\right) /\left\vert z_{0}-z_{1}\right\vert \right] ,
\end{equation*}%
where%
\begin{equation*}
R_{0}=-\frac{1}{4\kappa \left( \kappa -1\right) }\text{ and }%
R_{1}=\left\vert z_{0}-z_{1}\right\vert \frac{\theta c_{2}^{2}}{%
c_{1}^{2}l\left( \theta \right) }.
\end{equation*}

From the definition of $t_{2}$ we obtain%
\begin{equation*}
t_{2}=\frac{2\kappa }{1+\sqrt{1+4\kappa \left( \kappa -1\right) \left(
R_{0}-xR_{1}+R_{1}\mathrm{i}\right) }},
\end{equation*}%
which implies that%
\begin{equation}
t_{2}=\frac{2\kappa }{1+\rho \sqrt{-x+\mathrm{i}}},  \label{t2local}\tag{SM163}
\end{equation}%
where%
\begin{equation*}
\rho =\sqrt{4\kappa \left( \kappa -1\right) R_{1}}.
\end{equation*}

\begin{lemma}
\label{aboutt2}Let (\ref{t2local}) hold. Then $\frac{\mathrm{d}}{\mathrm{d}x%
}\left( -\operatorname{Re}\varphi _{2}\left( t_{2}\right) \right) <0$ for $x\geq 0$.
\end{lemma}

\textbf{\small{Proof:}} Since%
\begin{equation*}
\operatorname{Re}\sqrt{-x+\mathrm{i}}=\sqrt{\frac{\sqrt{x^{2}+1}-x}{2}}\text{ and }%
\operatorname{Im}\sqrt{-x+\mathrm{i}}=\sqrt{\frac{\sqrt{x^{2}+1}+x}{2}},
\end{equation*}%
we obtain%
\begin{eqnarray*}
\frac{\mathrm{d}}{\mathrm{d}x}\left( -\operatorname{Re}\varphi _{2}\left(
t_{2}\right) \right)  &=&-\frac{1}{2\sqrt{x^{2}+1}}\frac{\rho ^{2}x-\rho 
\operatorname{Re}\sqrt{-x+\mathrm{i}}}{\left\vert 1+\rho \sqrt{-x+\mathrm{i}}%
\right\vert ^{2}} \\
&&-\frac{2\kappa -1}{2\sqrt{x^{2}+1}}\frac{\rho ^{2}x+\left( 2\kappa
-1\right) \rho \operatorname{Re}\sqrt{-x+\mathrm{i}}}{\left\vert 2\kappa -1-\rho 
\sqrt{-x+\mathrm{i}}\right\vert ^{2}}.
\end{eqnarray*}%
For $x\geq 0$ this is no larger than%
\begin{equation*}
-\frac{\rho \operatorname{Re}\sqrt{-x+\mathrm{i}}}{2\sqrt{x^{2}+1}}\left( \frac{-1}{%
\left\vert 1+\rho \sqrt{-x+\mathrm{i}}\right\vert ^{2}}+\frac{\left( 2\kappa
-1\right) ^{2}}{\left\vert 2\kappa -1-\rho \sqrt{-x+\mathrm{i}}\right\vert
^{2}}\right) ,
\end{equation*}%
which is negative because%
\begin{equation*}
\left\vert 1+\rho \sqrt{-x+\mathrm{i}}\right\vert >\left\vert 1-\frac{\rho }{%
2\kappa -1}\sqrt{-x+\mathrm{i}}\right\vert . \square 
\end{equation*}

Lemma \ref{aboutt2} and identity (\ref{identityagain}) imply that $-\operatorname{Re}%
f_{\mathrm{e}}(z)$ is decreasing as $z$ moves along $\mathcal{K}_{21}.$
Hence $\mathcal{K}_{21}\cup \mathcal{K}_{1}$ is indeed a contour of steep
descent for $-\operatorname{Re}f(z),$ and therefore $L_{21}\left( \theta ;\Lambda
\right) $ is asymptotically dominated by $L_{1}\left( \theta ;\Lambda
\right) .$ It remains to be shown that $L_{22}\left( \theta ;\Lambda \right)
=L_{2}\left( \theta ;\Lambda \right) -L_{1}\left( \theta ;\Lambda \right) $
is asymptotically dominated by $L_{1}\left( \theta ;\Lambda \right) .$

For any $z\in \mathcal{K}_{22}$ and the corresponding $\eta _{2}=z\theta
c_{2}^{2}/\left[ c_{1}^{2}l\left( \theta \right) \right] $, consider the
integral representation%
\begin{equation}
F_{2}=\frac{C_{m}}{2\pi \mathrm{i}}\int_{0}^{(1+)}\exp \left\{ -m\varphi
_{2}\left( t\right) \right\} \psi _{2}\left( t\right) \mathrm{d}t,
\label{contourrep}\tag{SM164}
\end{equation}%
where%
\begin{equation*}
\varphi _{2}\left( t\right) =-\kappa \ln \left( t\right) +\left( \kappa
-1\right) \ln \left( t-1\right) +\kappa \ln \left( 1-\eta _{2}t\right) 
\end{equation*}%
\begin{equation*}
\psi _{2}\left( t\right) =\left( t-1\right) ^{-1}\left( 1-\eta _{2}t\right)
^{-1}.
\end{equation*}%
For a fixed contour $\mathcal{K}_{\ast }$ in (\ref{contourrep}), it is
clearly possible to make $\operatorname{Re}\varphi _{2}\left( t\right) $ arbitrarily
large and $\left\vert \psi _{2}\left( t\right) \right\vert $ arbitrarily
close to zero, uniformly with respect to $t\in \mathcal{K}_{\ast }$ by
choosing $A$ sufficiently large (so that $\left\vert \eta _{2}\right\vert $
is sufficiently large). Therefore, by choosing $A$ sufficiently large, we
shall have inequality%
\begin{equation*}
\left\vert F_{2}\right\vert \leq \tilde{\alpha}\sqrt{p}\exp \left\{ -\frac{p%
}{2}\left( \operatorname{Re}f_{\mathrm{h}}(z_{0})+\alpha \right) \right\} 
\end{equation*}%
for some $\tilde{\alpha},\alpha >0$ (that do not depend on $\theta $) and
any $z\in \mathcal{K}_{22}.$ Using this upper bound in (\ref{L2}), we obtain%
\begin{equation*}
L_{22}\left( \theta ;\Lambda \right) \leq \alpha _{1}p\exp \left\{ -\frac{p}{%
2}\alpha \right\} \left\vert g_{\mathrm{c}}g_{\mathrm{e}}(z_{0})\right\vert
\int_{\mathcal{K}_{22}}\left\vert \prod_{j=1}^{p}\left( \frac{z_{0}-\lambda
_{j}}{z-\lambda _{j}}\right) ^{1/2}\mathrm{d}z\right\vert 
\end{equation*}%
for some $\alpha _{1}>0$ that does not depend on $\theta $.

Clearly for any $z\in \mathcal{K}_{22}$ and any $\lambda _{j}$ from the
support of $F_{\mathbf{c}}$ we have%
\begin{equation*}
\left\vert \frac{z_{0}-\lambda _{j}}{z-\lambda _{j}}\right\vert \leq
\left\vert \frac{z_{0}}{z}\right\vert =\left\vert \frac{\eta _{20}}{\eta _{2}%
}\right\vert ,
\end{equation*}%
where $\eta _{20}$ is the value of $\eta _{2}$ that correspond to $z=z_{0}.$
Therefore, we have for some $\alpha _{2}>0$ that does not depend on $\theta $%
\begin{equation*}
L_{22}\left( \theta ;\Lambda \right) \leq \alpha _{2}p\exp \left\{ -\frac{p}{%
2}\alpha \right\} \left\vert \frac{g_{\mathrm{c}}}{\theta }g_{\mathrm{e}%
}(z_{0})\right\vert \int_{\mathcal{K}_{22}}\left\vert \left\vert \frac{\eta
_{20}}{\eta _{2}}\right\vert ^{p/2}\mathrm{d}\eta _{2}\right\vert ,
\end{equation*}%
and thus, for some $\alpha _{3}>0$ that does not depend on $\theta ,$%
\begin{equation*}
L_{22}\left( \theta ;\Lambda \right) \leq \alpha _{3}p\exp \left\{ -\frac{p}{%
2}\alpha \right\} \left\vert \frac{g_{\mathrm{c}}}{\theta }g_{\mathrm{e}%
}(z_{0})\right\vert .
\end{equation*}%
Finally, note that $g_{\mathrm{c}}/\theta =O(1)$ and $g_{\mathrm{e}%
}(z_{0})=O_{\mathrm{P}}(1)$, so that the above display implies (\ref{restInequalities})
with $L_{22}$ replacing $L_{2}$. Since 
\begin{equation*}
L_{1}\left( \theta ;\Lambda \right) =\frac{g(z_{0})}{\sqrt{-\mathrm{d}%
^{2}f(z_{0})/\mathrm{d}z^{2}}}+O_{\mathrm{P}}\left( p^{-1}\right) ,
\end{equation*}%
we see that $L_{22}\left( \theta ;\Lambda \right) $ is asymptotically
dominated by $L_{1}\left( \theta ;\Lambda \right) $.

\section{Asymptotics of LR}

\subsection{Derivations for Theorem JO\ref{LR asymptotics} (limiting LR)} \label{sec: derivations for theorem 10}

We record details to verify that 
\begin{equation*}
\frac{g(z_{0})}{\sqrt{-f^{\prime \prime }(z_{0})}}=\exp \left\{ -\tfrac{1}{2}%
\Delta _{p}(\theta )+\tfrac{1}{2}\log [1-\delta _{p}^{2}(\theta )]\right\}
(1+o(1)),
\end{equation*}%
where, perhaps surprisingly, our six cases reduce to the three values for $%
\delta _{p}(\theta )$ given in Theorem JO\ref{LR asymptotics}. Recall the decomposition $g=g_{%
\mathrm{c}}g_{\mathrm{e}}g_{\mathrm{h}}$ and note from the definitions
(JO\ref{Deltap}) that $g_{\mathrm{e}}(z_{0})=\exp \{-\tfrac{1}{2}\Delta _{p}(\theta )\}$%
. Consequently, from the definition of $D_{2}$ in Table JO\ref{Table 5}, the left side
of the previous display may be written as 
\begin{equation*}
\theta ^{-1}g_{\mathrm{c}}g_{\mathrm{h}}\sqrt{D_{2}}\exp \{-\tfrac{1}{2}%
\Delta _{p}(\theta )\},
\end{equation*}%
% we have
% \begin{equation*}
%   \frac{g(z_0)}{\sqrt{-2f''(z_0)}}
%     = \theta^{-1} \grc \cdot \grh \cdot \sqrt{D_2} \exp \{ - \hf
%     \Delta_p(\theta) \},
% \end{equation*}
so our task is to verify that 
\begin{equation}
P=\theta ^{-1}g_{\mathrm{c}}g_{\mathrm{h}}\sqrt{D_{2}}=(1-\delta
_{p}^{2}(\theta ))^{1/2}(1+o(1)).  \label{eq:prod}\tag{SM165}
\end{equation}%
To this end, Table \ref{tab:ThJO9} collects values for 
% $\theta^{-1} \grc$ (Table JO1),
% $\grh$ (Section JO4.1),
% and $\sqrt{D_2}$ (Table JO4).
$\theta ^{-1}\check{g}_{\mathrm{c}}$, $g_{\mathrm{h}}$, and $\sqrt{D_{2}}$
from Table JO\ref{Table 2a}, Section JO4.1 and Table JO\ref{Table 5} respectively. Cases SMD and PCA
require no further comment. For the remaining cases, we add remarks on the
evaluation of $g_{\mathrm{h}}(z_{0})$ and then the product \eqref{eq:prod}.

\begin{table}[h]
\centering 
\begin{tabular}{lccc}
& $\theta ^{-1}\check{g}_{c}$ & $g_{h}$ & $\sqrt{D_{2}}$ \\ \hline
&  &  &  \\ 
SMD \qquad  & $1$ & $1$ & $\sqrt{1-\theta ^{2}}$ \\[10pt]
PCA & $\dfrac{1}{c_{1}(1+\theta )}$ & $1$ & $c_{1}(1+\theta )\sqrt{%
h_{0}/c_{1}}$ \\[10pt]
SigD & $\dfrac{r\sqrt{c_{1}+c_{2}}}{c_{1}^{2}(1+\theta )}$ & $\dfrac{%
c_{1}l(\theta )}{r^{2}}$ & $\dfrac{r(1+\theta )\sqrt{h}}{l^{2}(\theta )}$ \\%
[10pt]
REG$_{0}$ & $\dfrac{1}{c_{1}\sqrt{1-c_{1}}}$ & $\dfrac{\sqrt{1-c_{1}}}{\sqrt{%
K_{0}}}$ & $c_{1}\sqrt{K_{0}h_{0}/c_{1}}$ \\[10pt]
REG & $\dfrac{r\sqrt{c_{1}+c_{2}}}{c_{1}^{2}\sqrt{1-c_{1}}}$ & $\dfrac{\sqrt{%
c_{1}(1-c_{1})}l(\theta )}{r\sqrt{K_{1}}}$ & $\dfrac{\sqrt{c_{1}K_{1}h}}{%
l^{2}(\theta )}$ \\[10pt]
CCA & $\dfrac{r^{2}(c_{1}+c_{2})}{c_{1}^{3}\sqrt{1-c_{1}}l(\theta )}$ & $%
\dfrac{c_{1}\sqrt{1-c_{1}}l^{3/2}(\theta )}{r^{2}\sqrt{K_{2}}}$ & $\dfrac{%
c_{1}\sqrt{K_{2}h}}{\sqrt{c_{1}+c_{2}}l^{3/2}(\theta )}$%
\end{tabular}
\caption{Components of the product $P=\protect\theta ^{-1}g_{\mathrm{c}}g_{%
\mathrm{h}}\protect\sqrt{D_{2}}$. The CCA entry for $g_{\mathrm{h}}$ is
shown for completeness -- it is derived, post facto, from the calculations
above.}
\label{tab:ThJO9}
\end{table}

\textbf{\small{SigD.}}
First observe that since $z_{0}\theta =(1+\theta )(c_{1}+\theta )/l(\theta )$%
, 
\begin{equation*}
g_{\mathrm{h}}(z_{0})=\left( 1-\frac{c_{2}z_{0}\theta }{c_{1}(1+\theta )}%
\right) ^{-1}=\frac{c_{1}l(\theta )}{r^{2}},
\end{equation*}%
and we get the claimed expression for $P$, 
\begin{equation}
P^{2}=\frac{(c_{1}+c_{2})h}{c_{1}^{2}l^{2}}=1-\frac{\theta ^{2}r^{2}}{%
c_{1}^{2}l^{2}},  \label{eq:prod-sigD}\tag{SM166}
\end{equation}%
after using the identity 
\begin{equation}
(c_{1}+c_{2})h=c_{1}^{2}l^{2}-\theta ^{2}r^{2}  \label{eq:id1}\tag{SM167}
\end{equation}
\vspace{3 mm}

\textbf{\small{REG$_{0}$.}}
From (JO\ref{g2}) and \eqref{eq:eta0-D0}, we have 
\begin{equation*}
g_{\mathrm{h}}(z_{0})\sim (1+4\eta _{0})^{-1/4}\sim \sqrt{1-c_{1}}/\sqrt{%
K_{0}}.
\end{equation*}
\vspace{3 mm}

\textbf{\small{REG.}}
We use (JO\ref{g2}) to evaluate $g_{\mathrm{h}}(z_{0})$. Using \eqref{eq:t1} to
evaluate $t_{1}(z_{0})$, we have 
\begin{align}
\varphi _{1}^{\prime \prime }(t_{1})& =\frac{\kappa }{t_{1}^{2}}-\frac{%
\kappa -1}{(t_{1}-1)^{2}}=\frac{c_{2}^{2}(1+\theta )^{2}}{c_{2}(1-c_{1})}%
\left[ \frac{r^{2}}{L^{2}(\theta )}-\frac{1}{c_{1}}\right]   \notag \\
& =-\frac{c_{2}^{2}(1+\theta )^{2}}{c_{1}^{2}(1-c_{1})}\frac{K_{1}(\theta )}{%
l^{2}(\theta )},  \label{eq:phi1pp}\tag{SM168}
\end{align}%
using the identity 
\begin{equation*}
L^{2}(\theta )-c_{1}r^{2}=c_{1}c_{2}K_{1}(\theta ).
\end{equation*}%
Since $t_{1}-1>0$ and $\varphi _{1}^{\prime \prime }(t_{1})<0$, we can take $%
\omega _{1}=0$. Together with $\psi
_{1}(t_{1})=(t_{1}-1)^{-1}=c_{2}(1+\theta )/c_{1}$, we obtain from (JO\ref{g2})
and \eqref{eq:t1} 
\begin{equation*}
g_{\mathrm{h}}(z_{0})\sim \sqrt{\frac{c_{1}}{r^{2}}}\,|\varphi _{1}^{\prime
\prime }(t_{1})|^{-1/2}\psi _{1}(t_{1})=\frac{\sqrt{c_{1}(1-c_{1})}l(\theta )%
}{r\sqrt{K_{1}(\theta )}}
\end{equation*}%
The product $P$ then reduces to the first expression in \eqref{eq:prod-sigD}.
\vspace{3 mm}

\textbf{\small{CCA.}}
We show that $P_{\mathrm{CCA}}=P_{\mathrm{REG}}(1+o(1))$. From Table \ref%
{tab:ThJO9} and (JO\ref{gammaj}), we have 
\begin{equation*}
\frac{\theta ^{-1}g_{\mathrm{c,C}}}{\theta ^{-1}g_{\mathrm{c,R}}}=\frac{r%
\sqrt{c_{1}+c_{2}}}{c_{1}l(\theta )},\qquad \sqrt{\frac{D_{2,C}}{D_{2,R}}}=%
\sqrt{\frac{c_{1}l(\theta )}{c_{1}+c_{2}}\frac{K_{2}}{K_{1}}},\qquad \frac{%
\psi _{2}(t_{2})}{\psi _{1}(t_{1})}=\frac{1}{1-\eta _{2}t_{2}}=\frac{%
c_{1}l(\theta )}{r^{2}}.
\end{equation*}%
% \begin{align*}
%   \theta^{-1} g_{\rm c, C} 
%     & = \theta^{-1} g_{\rm c, R} \frac{r\sqrt{c_1+c_2}}{c_1 l(\theta)}
%     \\
%   \frac{D_{2,C}}{D_{2,R}}
%     & = \frac{c_1 l(\theta)}{c_1+c_2} \frac{K_2}{K_1},  \\
%   \frac{\psi_2(t_2)}{\psi_1(t_1)} 
%     & = \frac{1}{1-\eta_2 t_2}  
%       = \frac{c_1 l(\theta)}{r^2}.
% \end{align*}
Multiplying these ratios and referring to (JO\ref{g2}), we obtain 
\begin{equation}
\frac{P_{CCA}}{P_{REG}}\sim \left\vert \frac{\varphi _{1}^{\prime \prime
}(t_{1})}{\varphi _{2}^{\prime \prime }(t_{2})}\right\vert ^{1/2}\left[ 
\frac{c_{1}l}{r^{2}}\frac{K_{2}(\theta )}{K_{1}(\theta )}\right]
^{1/2}e^{i\omega _{2}/2}.  \label{eq:P-rat}\tag{SM169}
\end{equation}%
We now compare $\varphi _{2}^{\prime \prime }(t_{2})$ to $\varphi
_{1}^{\prime \prime }(t_{1})$, recalling that $t_{2}=t_{1}$. First, from %
\eqref{eq:phi2p}, 
\begin{equation*}
\varphi _{2}^{\prime \prime }(t)=-\frac{\kappa \eta _{2}^{2}}{(1-\eta
_{2}t_{2})^{2}}+\varphi _{1}^{\prime \prime }(t).
\end{equation*}%
In particular, $\varphi _{2}^{\prime \prime }(t_{2})<0$ and, as with $\omega
_{1}$, also $\omega _{2}=0$. From \eqref{eq:t2eta2}, we evaluate 
\begin{equation*}
-\frac{\kappa \eta _{2}^{2}}{(1-\eta _{2}t_{2})^{2}}=-\frac{%
c_{2}^{3}(1+\theta )^{2}}{c_{1}^{2}(1-c_{1})}\frac{(c_{1}+\theta )^{2}}{%
l^{2}r^{2}},
\end{equation*}%
so that from \eqref{eq:phi1pp}, 
\begin{equation*}
\frac{\varphi _{2}^{\prime \prime }(t_{2})}{\varphi _{1}^{\prime \prime
}(t_{2})}=1+\frac{c_{2}(c_{1}+\theta )^{2}}{r^{2}K_{1}(\theta )}=\frac{%
c_{1}l(\theta )K_{2}(\theta )}{r^{2}K_{1}(\theta )},
\end{equation*}%
where the second identity follows after some algebra. The latter display and %
\eqref{eq:P-rat} show that $P_{\mathrm{CCA}}=P_{\mathrm{REG}}(1+o(1))$.

\subsection{Proof of Theorem JO\ref{AsymptoticNormality} (Gaussian process limit)} \label{sec: proof of theorem 11}

\textit{Some general considerations}

\textbf{\small{Almost sure continuity of $\ln L\left( \protect\theta ;\Lambda
\right) $.}}
Let $\varepsilon >0$ be a fixed small number. First, let us show that $\ln
L\left( \theta ;\Lambda \right) $ are continuous functions of $\theta \in
\lbrack 0,\bar{\theta}-\varepsilon ]$ for each of the six cases under study.
Recall equation (JO\ref{LRgeneral})%
\begin{equation}
L^{(Case)}\left( \theta ;\Lambda \right) =\alpha \left( \theta \right)
\left. _\mathsf{p}F_\mathsf{q}\left( a,b;\Psi ,\Lambda \right) \right. ,
\label{likelihood}\tag{SM170}
\end{equation}%
where $\Psi $ is a $p$-dimensional matrix $\operatorname{diag}\left\{ \Psi
_{11},0,...,0\right\} ,$ and the values of $\Psi _{11},$ $\alpha \left(
\theta \right) ,$ $\mathsf{p}$$,$ $\mathsf{q}$$,$ $a,$ and $b$ are as given
in Table \ref{Table 2}. Consider the series representation%
\begin{eqnarray*}
_\mathsf{p}F_\mathsf{q}\left( a,b;\Psi ,\Lambda \right) 
&=&\sum_{k=0}^{\infty }\frac{1}{k!}\sum_{\kappa \vdash k}\frac{\left(
a_{1}\right) _{\kappa }...\left( a_{\mathsf{p}}\right) _{\kappa }}{\left(
b_{1}\right) _{\kappa }...\left( b_{\mathsf{q}}\right) _{\kappa }}\frac{%
C_{\kappa }\left( \Psi \right) C_{\kappa }\left( \Lambda \right) }{C_{\kappa
}\left( I_{p}\right) } \\
&=&\sum_{k=0}^{\infty }\frac{1}{k!}\frac{\left( a_{1}\right) _{k}...\left(
a_{\mathsf{p}}\right) _{k}}{\left( b_{1}\right) _{k}...\left( b_{\mathsf{q}%
}\right) _{k}}\frac{\Psi _{11}^{k}C_{k}\left( \Lambda \right) }{C_{k}\left(
I_{p}\right) },
\end{eqnarray*}%
where the second equality follows from the fact that $C_{\kappa }\left( \Psi
\right) =0$ unless partition $\kappa \vdash k$ is trivial, that is $\kappa
=k,$ in which case $C_{\kappa }\left( \Psi \right) =\Psi _{11}^{k}$ (see
definition 7.2.1 iii in Muirhead (1982)). James (1968) shows that the
coefficients of zonal polynomials are positive. Therefore, for non-negative $%
\Psi _{11}$ and $\lambda _{j},$ $j=1,...,p,$ we have%
\begin{equation*}
0\leq \frac{\Psi _{11}^{k}C_{k}\left( \Lambda \right) }{C_{k}\left(
I_{p}\right) }\leq \left( \Psi _{11}\lambda _{1}\right) ^{k}.
\end{equation*}%
This implies that $_\mathsf{p}F_\mathsf{q}\left( a,b;\Psi ,\Lambda \right) $
is an analytic function of $\theta \in \lbrack 0,\bar{\theta}-\varepsilon ]$
and $_\mathsf{p}F_\mathsf{q}\left( a,b;\Psi ,\Lambda \right) \geq 1$ (the
first term in the expansion of $_\mathsf{p}F_\mathsf{q}\left( a,b;\Psi
,\Lambda \right) $ is 1) when $\mathsf{p}\leq \mathsf{q}$, that is for SMD,
PCA, REG$_{0}$, and REG cases. For SigD and CCA, $_\mathsf{p}F_\mathsf{q}%
\left( a,b;\Psi ,\Lambda \right) $ is an analytic function of $\theta $ in
the domain%
\begin{equation*}
\Psi _{11}\lambda _{1}<1.
\end{equation*}%
But for SigD and CCA $\lambda _{j}$ are solutions to 
\begin{equation*}
\det \left( H-\lambda \left( E+\frac{n_{1}}{n_{2}}H\right) \right) =0,
\end{equation*}%
and hence, with probability 1, $\lambda _{1}\leq n_{2}/n_{1}$ because $H$
and $E$ are positive definite. Therefore, for SigD we have 
\begin{equation*}
\Psi _{11}\lambda _{1}=\frac{\theta n_{1}}{n_{2}\left( 1+\theta \right) }%
\lambda _{1}\leq \frac{\theta }{1+\theta }<1
\end{equation*}%
for any $\theta \in \lbrack 0,\bar{\theta}-\varepsilon ],$ and for CCA we
have%
\begin{equation*}
\Psi _{11}\lambda _{1}=\frac{\theta n_{1}^{2}}{n_{2}^{2}+n_{2}n_{1}\left(
1+\theta \right) }\lambda _{1}\leq \frac{\theta n_{1}}{n_{2}+n_{1}(1+\theta )%
}<1
\end{equation*}%
for any $\theta \in \lbrack 0,\bar{\theta}-\varepsilon ].$ Thus, $%
_\mathsf{p}F_\mathsf{q}\left( a,b;\Psi ,\Lambda \right) $ is an analytic
function of $\theta \in \lbrack 0,\bar{\theta}-\varepsilon ]$ and $%
_\mathsf{p}F_\mathsf{q}\left( a,b;\Psi ,\Lambda \right) \geq 1$ for all six
cases that we consider. Using (\ref{likelihood}) we conclude that $\ln
L\left( \theta ;\Lambda \right) $ are continuous functions of $\theta \in
\lbrack 0,\bar{\theta}-\varepsilon ]$ with probability one. In particular
(see Bosq (2000) p. 22) $\ln L\left( \theta ;\Lambda \right) $ can be
interpreted as random element of the space $C_{\left[ 0,1-\varepsilon \right]
}$ of continuous functions on $\left[ 0,1-\varepsilon \right] $ equipped
with the supremum norm.

%TCIMACRO{\TeXButton{B}{\begin{table}[tbp] \centering}}%
%BeginExpansion
\begin{table}[tbp] \centering%
%EndExpansion
\begin{tabular}{llllll}
Case & $_\mathsf{p}F_\mathsf{q}$ & $\alpha \left( \theta \right) $ & $a$ & $b
$ & $\Psi _{11}$ \\ \hline
SMD & $_{0}F_{0}$ & $\exp \left( -p\theta ^{2}/4\right) $ & \_ & \_ & $%
\theta p/2\bigskip $ \\ 
PCA & $_{0}F_{0}$ & $\left( 1+\theta \right) ^{-n_{1}/2}$ & \_ & \_ & $%
\theta n_{1}/(2\left( 1+\theta \right) )\bigskip $ \\ 
SigD & $_{1}F_{0}$ & $\left( 1+\theta \right) ^{-n_{1}/2}$ & $n/2$ & \_ & $%
\theta n_{1}/\left( n_{2}\left( 1+\theta \right) \right) \bigskip $ \\ 
REG$_{0}$ & $_{0}F_{1}$ & $\exp \left( -n_{1}\theta /2\right) $ & \_ & $%
n_{1}/2$ & $\theta n_{1}^{2}/4\bigskip $ \\ 
REG & $_{1}F_{1}$ & $\exp \left( -n_{1}\theta /2\right) $ & $n/2$ & $n_{1}/2$
& $\theta n_{1}^{2}/\left( 2n_{2}\right) \bigskip $ \\ 
CCA & $_{2}F_{1}$ & $\left( 1+n_{1}\theta /n\right) ^{-n/2}$ & $\left(
n/2,n/2\right) $ & $n_{1}/2$ & $\theta n_{1}^{2}/\left(
n_{2}^{2}+n_{2}n_{1}\left( 1+\theta \right) \right) $ \\ \hline
\end{tabular}
\caption{Parameters of the JO's explicit expression (JO\ref{LRgeneral})  for the
likelihood ratios. Here $n\equiv n_1+n_2$.}\label{Table 2}%
%TCIMACRO{\TeXButton{E}{\end{table}}}%
%BeginExpansion
\end{table}%
%EndExpansion
\vspace{3 mm}

\textbf{\small{Reduction to a linear spectral statistic.}}
By Theorem JO\ref{LR asymptotics} we have%
\begin{equation}
\ln L\left( \theta ;\Lambda \right) =-\tfrac{1}{2}\Delta _{p}(\theta )+\frac{1%
}{2}\ln \left( 1-\left[ \delta _{p}\left( \theta \right) \right] ^{2}\right)
+o_{\mathrm{P}}(1),  \label{additive form}\tag{SM171}
\end{equation}%
where%
\begin{equation*}
\delta _{p}\left( \theta \right) =\left\{ 
\begin{array}{ll}
\theta  & \text{for SMD} \\ 
\theta /\sqrt{c_{1}} & \text{for PCA and REG}_{0} \\ 
\theta r/\left( c_{1}l\left( \theta \right) \right)  & \text{for SigD, REG,
and CCA}%
\end{array}%
\right. 
\end{equation*}%
and 
\begin{equation}
\Delta _{p}(\theta )=p\int \ln \left( z_{0}-\lambda \right) \mathrm{d}\left( 
\hat{F}\left( \lambda \right) -F_{\mathbf{c}}\left( \lambda \right) \right) 
\label{deltaptheta}\tag{SM172}
\end{equation}%
with%
\begin{equation}
z_{0}=\left\{ 
\begin{array}{ll}
\theta +1/\theta  & \text{for SMD} \\ 
\left( 1+\theta \right) \left( \theta +c_{1}\right) /\theta  & \text{for PCA
and REG}_{0} \\ 
\left( 1+\theta \right) \left( \theta +c_{1}\right) /\left[ \theta l\left(
\theta \right) \right]  & \text{for SigD, REG, and CCA}%
\end{array}%
\right.   \label{zcase}\tag{SM173}
\end{equation}%
and $F_{\mathbf{c}}$ equals the semicircle distribution for SMD, the
Marchenko-Pastur distribution for PCA and REG$_{0},$ and the scaled Wachter
distribution for SigD, REG, and CCA. As explained in JO, the statistic $%
\Delta _{p}(\theta )$ should be interpreted as zero whenever $z_{0}\leq
\lambda _{1}$.

Since both $\ln L\left( \theta ;\Lambda \right) $ and $\Delta _{p}(\theta )$
are random element of $C_{\left[ 0,1-\varepsilon \right] },$ $o_{\mathrm{P}%
}(1)$ is also a random element of $C_{\left[ 0,1-\varepsilon \right] },$ and 
$\left\Vert o_{\mathrm{P}}(1)\right\Vert \overset{\mathrm{P}}{\rightarrow }0$%
. Therefore by the standard argument, see for example Theorem 3.1 of
Billingsley (1999), p. 27, the weak limits of $\ln L\left( \theta ;\Lambda
\right) $ and of $-\tfrac{1}{2}\Delta _{p}(\theta )+\tfrac{1}{2}\ln \left( 1-%
\left[ \delta _{p}\left( \theta \right) \right] ^{2}\right) $ coincide. Note
that $\tfrac{1}{2}\ln \left( 1-\left[ \delta _{p}\left( \theta \right) \right]
^{2}\right) $ is converging in the space $C_{\left[ 0,1-\varepsilon \right]
} $ to%
\begin{equation*}
\delta \left( \theta \right) =\left\{ 
\begin{array}{ll}
\theta & \text{for SMD} \\ 
\theta /\sqrt{\gamma _{1}} & \text{for PCA and REG}_{0} \\ 
\theta \rho /\left( \gamma _{1}+\gamma _{2}+\theta \gamma _{2}\right) & 
\text{for SigD, REG, and CCA}%
\end{array}%
\right. .
\end{equation*}%
Therefore, we only need to establish the weak convergence of $\Delta
_{p}(\theta ).$ There are two facts to be established. First, the tightness
of $\Delta _{p}(\theta ),$ and second, the convergence of its finite
dimensional distributions. \vspace{3 mm}

\textbf{\small{Tightness of $\Delta _{p}(\protect\theta )$.}}
There are three cases to consider: $F_{\mathbf{c}}$ is the semicircle, the
Marchenko-Pastur, and the Wachter distribution. Whether the Marchenko-Pastur 
$F_{\mathbf{c}}$ corresponds to PCA or REG$_{0}$ cases is of no importance
because we consider the tightness under the null hypothesis so that $\hat{F}$
is the same for PCA and REG$_{0}.$ Similarly, the differences between SigD,
REG and CCA cases are of no importance here.

\paragraph{Tightness, Semi-circle $F_{\mathbf{c}}$}

The tightness of $\Delta _{p}(\theta )$ in this case is a direct
consequence of Theorem 1.1 of Bai and Yao (2005).

\paragraph{Tightness, Marchenko-Pastur $F_{\mathbf{c}}$}

Following Bai and Silverstein (2004), let us represent the linear spectral
statistic $\Delta _{p}(\theta )$ in the following form%
\begin{equation*}
\Delta _{p}(\theta )=-\frac{1}{2\pi \mathrm{i}}\oint_{\mathcal{R}}\ln
\left( z_{0}-z\right) p\left[ \hat{s}(z)-s_{\mathbf{c}}(z)\right] \mathrm{d}%
z,
\end{equation*}%
where $\mathcal{R}$ is contour that does not intersect the supports of $\hat{%
F}$ and $F_{\mathbf{c}}$ and does not encircle $z_{0}.$ Here 
\begin{equation*}
\hat{s}(z)=\int \left( \lambda -z\right) ^{-1}\mathrm{d}\hat{F}\left(
\lambda \right) \text{ and }s_{\mathbf{c}}\left( z\right) =\int \left(
\lambda -z\right) ^{-1}\mathrm{d}F_{\mathbf{c}}\left( \lambda \right) .
\end{equation*}%
With asymptotically negligible probability the above requirements for $%
\mathcal{R}$ are impossible to satisfy. We will therefore condition our
arguments on the high probability event that ensures the existence of
required $\mathcal{R}$.

Precisely, recall that the supports of $F_{\boldsymbol{\gamma }}$ and $F_{%
\mathbf{c}}$ are given by%
\begin{eqnarray*}
\left[ \beta _{-},\beta _{+}\right] &=&\left[ \left( 1-\sqrt{\gamma _{1}}%
\right) ^{2},\left( 1+\sqrt{\gamma _{1}}\right) ^{2}\right] \text{ and} \\
\left[ b_{-},b_{+}\right] &=&\left[ \left( 1-\sqrt{c_{1}}\right) ^{2},\left(
1+\sqrt{c_{1}}\right) ^{2}\right] ,
\end{eqnarray*}%
respectively, and the threshold $\bar{\theta}$ equals $\sqrt{\gamma _{1}}$.
Furthermore, $\mathbf{c}\rightarrow \boldsymbol{\gamma }$. Using these facts and
the definition of $z_{0}$, it is straightforward to verify that there exists 
$\eta >0$ that depends on $\varepsilon $ such that%
\begin{equation*}
\min_{\theta \in \left[ 0,\bar{\theta}-\varepsilon \right] }\left(
z_{0}-\beta _{+}-\eta \right) >0
\end{equation*}%
for all sufficiently large $n_{1},n_{2},p$ along the sequence $\mathbf{n}%
,p\rightarrow _{\boldsymbol{\gamma }}\infty $. Further, note that $\lambda _{1}%
\overset{a.s.}{\rightarrow }\beta _{+}$ and $\lambda _{p}\overset{a.s.}{%
\rightarrow }\beta _{-}$ when $\mathbf{n},p\rightarrow _{\boldsymbol{\gamma }%
}\infty $.

Consider the event%
\begin{equation}
Q_{p}=\left\{ \max \left\{ \lambda _{1},b_{+}\right\} \leq \beta _{+}+\eta
/2<z_{0}-\eta /2,\text{ }\min \left\{ \lambda _{p},b_{-}\right\} \geq \beta
_{-}-\eta /2\right\} .  \label{Qp}\tag{SM174}
\end{equation}%
The discussion above implies that 
\begin{equation}
\lim_{p\rightarrow \infty }\Pr \left\{ Q_{p}\right\} =1.
\label{frequent event}\tag{SM175}
\end{equation}%
Let $\mathcal{R}$ be the rectangular contour with the vertices at $\left(
\beta _{+}+\eta \right) \pm \mathrm{i}v$ and $\left( \beta _{-}-\eta \right)
\pm \mathrm{i}v$ for an arbitrary fixed positive $v.$ Conditional on the
event $Q_{p},$ $\mathcal{R}$ does not intersect the supports of $\hat{F}$
and $F_{\mathbf{c}}$ and does not encircle $z_{0}$ as required. Since $\Pr
\left\{ Q_{p}\right\} \rightarrow \infty ,$ it is sufficient to establish
the tightness of $\Delta _{p}(\theta )$ conditional on $Q_{p}$. Therefore,
in what follows we shall assume that $Q_{p}$ holds.

Let $\mathcal{C}$ be the part of $\mathcal{R}$ that lies in the upper half
complex plane. Then%
\begin{equation*}
\Delta _{p}(\theta )=-\frac{1}{\pi }\operatorname{Im}\int_{\mathcal{C}}\ln \left(
z_{0}-z\right) p\left[ \hat{s}(z)-s_{\mathbf{c}}(z)\right] \mathrm{d}z.
\end{equation*}%
Since the mapping%
\begin{equation*}
f(z)\mapsto g(\theta )=-\frac{1}{\pi }\operatorname{Im}\int_{\mathcal{C}}\ln \left(
z_{0}-z\right) f(z)\mathrm{d}z
\end{equation*}%
is a continuous mapping from the space $C_{\mathcal{C}}$ of the
complex-valued continuous functions on $\mathcal{C}$ (with the supremum
norm) to the space $C_{\left[ 0,1-\varepsilon \right] },$ the tightness of $%
\Delta _{p}(\theta )$ would follow from that of 
\begin{equation*}
M_{p}(z)\equiv p\left[ \hat{s}(z)-s_{\mathbf{c}}(z)\right] .
\end{equation*}

As in Bai and Silverstein (2004) p. 561, choose sequence $\left\{
\varepsilon _{p}\right\} $ such that $\varepsilon _{p}\rightarrow 0$ as $%
\mathbf{n},p\rightarrow _{\boldsymbol{\gamma }}\infty $ and 
\begin{equation*}
\varepsilon _{p}\geq p^{-\alpha }
\end{equation*}%
for some $\alpha \in \left( 0,1\right) $. Further, let%
\begin{eqnarray*}
\mathcal{C}_{u} &=&\left\{ x+\mathrm{i}v:x\in \left[ \beta _{-}-\eta ,\beta
_{+}+\eta \right] \right\} , \\
\mathcal{C}_{l} &=&\left\{ \left( \beta _{-}-\eta \right) +\mathrm{i}y:y\in %
\left[ p^{-1}\varepsilon _{p},v\right] \right\} , \\
\mathcal{C}_{r} &=&\left\{ \left( \beta _{+}+\eta \right) +\mathrm{i}y:y\in %
\left[ p^{-1}\varepsilon _{p},v\right] \right\} ,
\end{eqnarray*}%
and let $\mathcal{C}_{p}=\mathcal{C}_{l}\cup \mathcal{C}_{u}\cup \mathcal{C}%
_{r}$. Define the process $\hat{M}_{p}\left( z\right) $ on $\mathcal{C}$ as
follows%
\begin{equation*}
\hat{M}_{p}\left( z\right) =\left\{ 
\begin{array}{ll}
M_{p}(z) & \text{for }z\in \mathcal{C}_{p} \\ 
M_{p}(\beta _{+}+\eta +\mathrm{i}p^{-1}\varepsilon _{p}) & \text{for }%
z=\beta _{+}+\eta +\mathrm{i}y,\text{ }y\in \left[ 0,p^{-1}\varepsilon _{p}%
\right] \\ 
M_{p}(\beta _{-}-\eta +\mathrm{i}p^{-1}\varepsilon _{p}) & \text{for }%
z=\beta _{-}-\eta +\mathrm{i}y,\text{ }y\in \left[ 0,p^{-1}\varepsilon _{p}%
\right]%
\end{array}%
\right. .
\end{equation*}

Note that 
\begin{equation*}
\hat{M}_{p}\left( z\right) =p\left[ \hat{s}(z)-s_{\mathbf{c}}(z)\right] +o_{%
\mathrm{P}}(1),
\end{equation*}%
where $o_{\mathrm{P}}(1)$ is uniform over $z\in \mathcal{C}$. Indeed, for
any $z\in \mathcal{C}_{p}$ we have%
\begin{equation*}
\hat{M}_{p}\left( z\right) =p\left[ \hat{s}(z)-s_{\mathbf{c}}(z)\right] ,
\end{equation*}%
whereas by the definition of $\hat{s}(z)$ and (\ref{Qp})%
\begin{equation*}
\sup_{y\in \left[ 0,p^{-1}\varepsilon _{p}\right] }p\left\vert \hat{s}(\beta
_{\pm }\pm \eta +\mathrm{i}y)-\hat{s}(\beta _{\pm }\pm \eta +\mathrm{i}%
p^{-1}\varepsilon _{p})\right\vert \leq p\frac{p^{-1}\varepsilon _{p}}{%
\left( \eta /2\right) ^{2}}\rightarrow 0,
\end{equation*}%
and similarly%
\begin{equation*}
\sup_{y\in \left[ 0,p^{-1}\varepsilon _{p}\right] }p\left\vert s_{\mathbf{c}%
}(\beta _{\pm }\pm \eta +\mathrm{i}y)-s_{\mathbf{c}}(\beta _{\pm }\pm \eta +%
\mathrm{i}p^{-1}\varepsilon _{p})\right\vert \mathbf{1}\left\{ Q_{p}\right\}
\leq p\frac{p^{-1}\varepsilon _{p}}{\left( \eta /2\right) ^{2}}\rightarrow 0.
\end{equation*}%
Therefore, it is sufficient to prove the tightness of $\hat{M}_{p}\left(
\cdot \right) $ as a sequence of random elements of $C_{\mathcal{C}}$. Lemma
1.1 of Bai and Silverstein (2004) establishes this result along with the
weak convergence of $\hat{M}_{p}\left( \cdot \right) $ to a Gaussian process.

\paragraph{Tightness, Wachter $F_{\mathbf{c}}$}

We shall base our arguments on the results established in Zheng (2012). He
establishes a CLT for linear spectral statistics of multivariate $F$ and $%
Beta$ matrices via representing those statistics in the form of a contour
integral that involves a process related to $M_{p}\left( z\right) $ (see the
previous section). The CLT follows from his proving the convergence of the
process to a Gaussian process.

In contrast to JO, whose attention is focused on the eigenvalues of $H\left(
E+\frac{n_{1}}{n_{2}}H\right) ^{-1}$, Zheng's (2012) primary focus is on the
eigenvalues of $HE^{-1}.$ Let $\hat{F}$ and $\hat{G}$ be the empirical
distributions of the eigenvalues of $H\left( E+\frac{n_{1}}{n_{2}}H\right)
^{-1}$ and $HE^{-1},$ respectively. If $x$ is an eigenvalue of $HE^{-1},$
then $x\left( 1+c_{2}x/c_{1}\right) ^{-1}$ is an eigenvalue of $H\left( E+%
\frac{n_{1}}{n_{2}}H\right) ^{-1}$, and thus 
\begin{equation*}
\hat{G}\left( x\right) =\hat{F}\left( \frac{x}{1+c_{2}x/c_{1}}\right) .
\end{equation*}%
A similar equality holds for the corresponding limiting distributions $G_{%
\mathbf{c}}$ and $F_{\mathbf{c}}$. Therefore,%
\begin{eqnarray*}
\Delta _{p}(\theta ) &\equiv &p\int \ln \left( z_{0}-\lambda \right) \mathrm{%
d}\left( \hat{F}\left( \lambda \right) -F_{\mathbf{c}}\left( \lambda \right)
\right) \\
&=&p\int \ln \left( z_{0}-\frac{x}{1+c_{2}x/c_{1}}\right) \mathrm{d}\left( 
\hat{G}\mathcal{(}x)-G_{\mathbf{c}}\left( x\right) \right) .
\end{eqnarray*}%
Denote the Stieltjes transform of $\hat{G}$ as $\hat{m}(z)$ and that of $G_{%
\mathbf{c}}$ as $m_{\mathbf{c}}(z).$ Then, similarly to the Marchenko-Pastur
case, we have%
\begin{equation*}
\Delta _{p}(\theta )=-\frac{1}{2\pi \mathrm{i}}\oint_{\mathcal{R}}\ln
\left( z_{0}-\frac{z}{1+c_{2}z/c_{1}}\right) p\left[ \hat{m}(z)-m_{\mathbf{c}%
}(z)\right] \mathrm{d}z,
\end{equation*}%
where $\mathcal{R}$ is contour that does not intersect the supports of $\hat{%
G}$ and $G_{\mathbf{c}}$ and does not encircle $z_{0}/\left(
1-c_{2}z_{0}/c_{1}\right) .$ As above, the existence of such a contour
requires conditioning on a large probability event, which we shall assume.

Zheng (2012) pp. 467--470 sketches a proof of the weak convergence of $p\left[
\hat{m}(z)-m_{\mathbf{c}}(z)\right] .$ Such a weak convergence implies the
tightness, which in its turn implies the tightness of $\Delta _{p}(\theta ).$

For the reader's convenience, we provide here a brief description of the
main steps in Zheng's proof. The proof is based on the decomposition%
\begin{equation*}
p\left[ \hat{m}(z)-m_{\mathbf{c}}(z)\right] =p\left[ \hat{m}(z)-m_{\mathbf{c}%
}^{\left( E\right) }(z)\right] +p\left[ m_{\mathbf{c}}^{\left( E\right)
}(z)-m_{\mathbf{c}}(z)\right] ,
\end{equation*}%
where $m_{\mathbf{c}}^{\left( E\right) }(z)$ is the Stieltjes transform of $%
G_{\mathbf{c}}^{(E)},$ the limiting spectral distribution (as $\mathbf{n}%
,p\rightarrow _{\mathbf{c}}\infty $) of $HA_{p}^{-1}$ where the empirical
spectral distribution of symmetric positive definite matrix $A_{p}$
converges to that of $E$ as $\mathbf{n},p\rightarrow _{\mathbf{c}}\infty $.
First, Zheng establishes the weak convergence of $p\left[ \hat{m}(z)-m_{%
\mathbf{c}}^{\left( E\right) }(z)\right] $ conditional on $\left\{
E,p=1,2...\right\} $ by appealing to Lemma 1.1 of Bai and Silverstein
(2004). Since the limiting process does not depend on $\left\{
E,p=1,2...\right\} ,$ the unconditional convergence also follows. Next,
Zheng represents $p\left[ m_{\mathbf{c}}^{\left( E\right) }(z)-m_{\mathbf{c}%
}(z)\right] $ as a product of a continuous function of $z$ that converges in 
$C_{\mathcal{R}}$ and the term $p\left[ \hat{m}_{E}\left( -\underline{m}_{%
\mathbf{c}}(z)\right) -m_{c_{2}}\left( -\underline{m}_{\mathbf{c}}(z)\right) %
\right] ,$ where $\hat{m}_{E}$ is the Stieltjes transform of the empirical
spectral distribution of $E,$ $m_{c_{2}}$ is that of the corresponding
limiting distribution as $\mathbf{n},p\rightarrow _{\mathbf{c}}\infty ,$ and 
$\underline{m}_{\mathbf{c}}$ is defined via the Stieltjes transform $m_{%
\mathbf{c}}$ of $G_{\mathbf{c}}$ by%
\begin{equation*}
\underline{m}_{\mathbf{c}}(z)=-\frac{1-c_{1}}{z}+c_{1}m_{\mathbf{c}}\left(
z\right) .
\end{equation*}%
Then, she points out that $-\underline{m}_{\mathbf{c}}(z)$ converges to $-%
\underline{m}_{\boldsymbol{\gamma }}(z),$ which is defined analogously with $%
\mathbf{c}$ replaced by $\boldsymbol{\gamma }$. Function $-\underline{m}_{%
\boldsymbol{\gamma }}(z)$ transforms $\mathcal{R}$ to a contour encircling the
support of the limiting spectral distribution of $E.$ Zheng appeals to Lemma
1.1 of Bai and Silverstein (2004) to establish the weak convergence of $p%
\left[ \hat{m}_{E}\left( z\right) -m_{c_{2}}\left( z\right) \right] $ as a
random continuous function on such a contour. Zheng's proof omits some
details, probably for the sake of saving the space. For example, she does
not mention that to be able to view $p\left[ \hat{m}(z)-m_{\mathbf{c}%
}^{\left( E\right) }(z)\right] $ and $p\left[ m_{\mathbf{c}}^{\left(
E\right) }(z)-m_{\mathbf{c}}(z)\right] $ as continuous random functions on $%
\mathcal{R}$, a conditioning on some event of increasing probability in
needed. Having a detailed proof would be useful, but requires a separate
research effort. \vspace{3 mm}

\textbf{\small{Finite dimensional convergence.}}
The convergence of the finite dimensional distributions to Gaussian
distributions follow from Theorem 1.1 of Bai and Yao (2005) for the
semicircle $F_{\mathbf{c}}$, from Theorem 1.1 of Bai and Silverstein (2004)
for the Marchenko-Pastur $F_{\mathbf{c}},$ and from Theorem 4.1 of Zheng
(2012) for the Wachter $F_{\mathbf{c}}$. We now use the results in the above
mentioned papers to compute the means and covariance matrices of the
asymptotic finite-dimensional distributions of $\Delta _{p}(\theta )$.

\paragraph{Finite dimensional asymptotics, Semi-circle $F_{\mathbf{c}}$}

Recall that $z_{0}=\theta +1/\theta ,$ we obtain%
\begin{equation*}
\Delta _{p}\left( \theta \right) =p\int \ln \left( \theta ^{2}-\lambda
\theta +1\right) \mathrm{d}\left( \hat{F}\left( \lambda \right) -F_{\mathbf{c%
}}\left( \lambda \right) \right) .
\end{equation*}%
Theorem 1.1 Bai and Yao (2005) implies that the random vector $\left( \Delta
_{p}\left( \theta _{1}\right) ,...,\Delta _{p}\left( \theta _{k}\right)
\right) $ with $\theta _{i}\in \left[ 0,\bar{\theta}-\varepsilon \right] $
converges in distribution to a Gaussian vector $(\mathcal{D}\left( \theta
_{1}\right) ,...,\mathcal{D}\left( \theta _{k}\right) )$ with 
\begin{equation}
\mathbb{E}\mathcal{D}\left( \theta _{i}\right) =\frac{1}{4}\left[ \ln \left[
\left( 1-\theta _{i}\right) ^{2}\right] +\ln \left[ \left( 1+\theta
_{i}\right) ^{2}\right] \right] -\tfrac{1}{2}\tau _{0}\left( \theta \right)
\label{MeanSc}\tag{SM176}
\end{equation}%
and%
\begin{equation}
Cov\left( \mathcal{D}\left( \theta _{i}\right) ,\mathcal{D}\left( \theta
_{j}\right) \right) =2\sum_{l=1}^{\infty }l\tau _{l}\left( \theta
_{i}\right) \tau _{l}\left( \theta _{j}\right) ,  \label{CovSc}\tag{SM177}
\end{equation}%
where 
\begin{equation*}
\tau _{l}\left( \theta \right) =\frac{1}{2\pi }\int_{-\pi }^{\pi }\ln \left(
1+\theta ^{2}-2\theta \cos \varphi \right) \cos \left( l\varphi \right) 
\mathrm{d}\varphi .
\end{equation*}

\begin{lemma}
\label{tau integral}For any $\theta $, such that $\left\vert \theta
\right\vert <1$, and any integer $l>0,$ we have $\tau _{l}\left( \theta
\right) =-\theta ^{l}/l$ and $\tau _{0}\left( \theta \right) =0$.
\end{lemma}

\textbf{\small{Proof:}} Changing the variable of integration from $\varphi $ to $%
z=e^{\mathrm{i}\varphi },$ we obtain%
\begin{equation*}
\tau _{l}\left( \theta \right) =\frac{1}{2\pi \mathrm{i}}\oint \ln \left[
\left( 1-\theta z\right) \left( 1-\theta z^{-1}\right) \right] z^{l-1}%
\mathrm{d}z,
\end{equation*}%
where the contour integral is taken over the counter-clockwise oriented unit
circle. Representing the logarithm of a product as the sum of logarithms, we
obtain%
\begin{equation*}
\tau _{l}\left( \theta \right) =\frac{1}{2\pi \mathrm{i}}\oint \ln \left(
1-\theta z\right) z^{l-1}\mathrm{d}z+\frac{1}{2\pi i}\oint \ln \left[
1-\theta z^{-1}\right] z^{l-1}\mathrm{d}z.
\end{equation*}%
Since, for $\left\vert \theta \right\vert <1,$ $\ln \left( 1-\theta z\right) 
$ is analytic in the unit circle and equal to zero at $z=0,$ we have%
\begin{equation*}
\frac{1}{2\pi \mathrm{i}}\oint \ln \left( 1-\theta z\right) z^{l-1}\mathrm{d%
}z=0
\end{equation*}%
for any integer $l\geq 0$. Hence,%
\begin{equation*}
\tau _{l}\left( \theta \right) =\frac{1}{2\pi \mathrm{i}}\oint \ln \left[
1-\theta z^{-1}\right] z^{l-1}\mathrm{d}z.
\end{equation*}%
Changing the variable of integration from $z$ to $\zeta =z^{-1}$, and noting
that $\mathrm{d}z/z=-\mathrm{d}\zeta /\zeta $, we get%
\begin{equation*}
\tau _{l}\left( \theta \right) =\frac{1}{2\pi \mathrm{i}}\oint \ln \left[
1-\theta \zeta \right] \zeta ^{-l-1}\mathrm{d}\zeta .
\end{equation*}%
On the other hand, for $\left\vert \zeta \right\vert \leq 1$, we have the
following power series expansion%
\begin{equation*}
\ln \left[ 1-\theta \zeta \right] =-\sum_{j=1}^{\infty }\frac{\theta ^{j}}{j%
}\zeta ^{j}.
\end{equation*}%
Thus, by Cauchy's residue theorem, $\tau _{l}\left( \theta \right) =-\theta
^{l}/l$ for $l>0$ and $\tau _{0}\left( \theta \right) =0$. $\square $

Lemma \ref{tau integral} together with (\ref{MeanSc}) and (\ref{CovSc}) yield%
\begin{equation*}
\mathbb{E}\mathcal{D}\left( \theta _{i}\right) =\tfrac{1}{2}\ln \left(
1-\theta _{i}^{2}\right)
\end{equation*}%
and%
\begin{equation*}
Cov\left( \mathcal{D}\left( \theta _{i}\right) ,\mathcal{D}\left( \theta
_{j}\right) \right) =-2\ln \left( 1-\theta _{i}\theta _{j}\right) .
\end{equation*}

\paragraph{Finite dimensional asymptotics, Marchenko-Pastur $F_{\mathbf{c}}$}

For PCA and REG$_{0}$ the finite dimensional distributions of $\Delta
_{p}\left( \theta \right) $ are derived in Lemma 12 of Onatski et al (2013).
They show that the random vector $\left( \Delta _{p}\left( \theta
_{1}\right) ,...,\Delta _{p}\left( \theta _{k}\right) \right) $ with $\theta
_{i}\in \left[ 0,\bar{\theta}-\varepsilon \right] $ converges in
distribution to a Gaussian vector $(\mathcal{D}\left( \theta _{1}\right)
,...,\mathcal{D}\left( \theta _{k}\right) )$ with%
\begin{equation*}
\mathbb{E}\mathcal{D}\left( \theta _{i}\right) =\tfrac{1}{2}\ln \left(
1-\theta _{i}^{2}/\gamma _{1}\right)
\end{equation*}%
and%
\begin{equation*}
Cov\left( \mathcal{D}\left( \theta _{i}\right) ,\mathcal{D}\left( \theta
_{j}\right) \right) =-2\ln \left( 1-\theta _{i}\theta _{j}/\gamma
_{1}\right) .
\end{equation*}

\paragraph{Finite dimensional asymptotics, Wachter $F_{\mathbf{c}}$}

Let 
\begin{equation*}
\hat{G}\left( x\right) =\hat{F}\left( \frac{x}{1+c_{2}x/c_{1}}\right) \text{
and }G_{\mathbf{c}}\left( x\right) =F_{\mathbf{c}}\left( \frac{x}{%
1+c_{2}x/c_{1}}\right) .
\end{equation*}%
Then%
\begin{equation}
\Delta _{p}\left( \theta \right) =p\int \ln \left( z_{0}-\frac{x}{%
1+c_{2}x/c_{1}}\right) \mathrm{d}\left( \hat{G}\left( x\right) -G_{\mathbf{c}%
}\left( x\right) \right) .  \label{FDelta}\tag{SM178}
\end{equation}%
Recall that%
\begin{equation*}
z_{0}=\frac{\left( 1+\theta \right) \left( \theta +c_{1}\right) }{\theta
\left( 1+\left( 1+\theta \right) c_{2}/c_{1}\right) }.
\end{equation*}%
Let 
\begin{equation}
z_{\gamma 0}=\frac{\left( 1+\theta \right) \left( \theta +\gamma _{1}\right) 
}{\theta \left( 1+\left( 1+\theta \right) \gamma _{2}/\gamma _{1}\right) }.
\label{zgamma0}\tag{SM179}
\end{equation}%
Since $z_{0}\rightarrow z_{\gamma 0}$ and $\mathbf{c}\rightarrow \mathbf{%
\gamma }$ as $\mathbf{n},p\rightarrow _{\boldsymbol{\gamma }}\infty ,$ the
asymptotic distribution of the random vector $\left( \Delta _{p}\left(
\theta _{1}\right) ,...,\Delta _{p}\left( \theta _{k}\right) \right) $ must
be the same as that of $\left( \Delta _{\gamma p}\left( \theta _{1}\right)
,...,\Delta _{\gamma p}\left( \theta _{k}\right) \right) ,$ where%
\begin{equation}
\Delta _{\gamma p}\left( \theta \right) =p\int \ln \left( z_{\gamma 0}-\frac{%
x}{1+\gamma _{2}x/\gamma _{1}}\right) \mathrm{d}\left( \hat{G}\left(
x\right) -G_{\mathbf{c}}\left( x\right) \right) .  \label{Dgammap}\tag{SM180}
\end{equation}%
This can be formally shown by considering the representation%
\begin{equation*}
\Delta _{p}\left( \theta \right) -\Delta _{\gamma p}\left( \theta \right) =-%
\frac{1}{2\pi \mathrm{i}}\oint_{\mathcal{R}}\ln \left[ \frac{z_{0}-\frac{z}{%
1+c_{2}z/c_{1}}}{z_{\gamma 0}-\frac{z}{1+\gamma _{2}z/\gamma _{1}}}\right] p%
\left[ \hat{m}(z)-m_{\mathbf{c}}(z)\right] \mathrm{d}z
\end{equation*}%
(see subsection \textquotedblleft Tightness, Wachter $F_{\mathbf{c}}$%
\textquotedblright ) and using the convergence of $p\left[ \hat{m}(z)-m_{%
\mathbf{c}}(z)\right] $ established by Zheng (2012) to demonstrate that 
\begin{equation}
\Delta _{p}\left( \theta \right) -\Delta _{\gamma p}\left( \theta \right)
=o_{\mathrm{P}}(1).  \label{negligeability}\tag{SM181}
\end{equation}%
Theorem 3.1 of Zheng (2012) implies that the random vector $\left( \Delta
_{\gamma p}\left( \theta _{1}\right) ,...,\Delta _{\gamma p}\left( \theta
_{k}\right) \right) $ with $\theta _{i}\in \left[ 0,\bar{\theta}-\varepsilon %
\right] $ converges in distribution to a Gaussian vector $(\mathcal{D}\left(
\theta _{1}\right) ,...,\mathcal{D}\left( \theta _{k}\right) )$. Let us find
the asymptotic mean and covariances.

Equations (\ref{zgamma0}), (\ref{Dgammap}) and some elementary algebra yield%
\begin{equation}
\Delta _{\gamma p}\left( \theta \right) =\Delta _{\gamma p}^{(1)}\left(
\theta \right) -\Delta _{\gamma p}^{(2)}\left( \theta \right) ,
\label{decompDgp}\tag{SM182}
\end{equation}%
where%
\begin{equation*}
\Delta _{\gamma p}^{(1)}\left( \theta \right) =p\int \ln \left( \theta
+\gamma _{1}+\left( \gamma _{2}-\frac{\theta }{1+\theta }\right) x\right) 
\mathrm{d}\left( \hat{G}\left( x\right) -G_{\mathbf{c}}\left( x\right)
\right) ,
\end{equation*}%
and%
\begin{equation*}
\Delta _{\gamma p}^{(2)}\left( \theta \right) =p\int \ln \left( \gamma
_{1}/\gamma _{2}+x\right) \mathrm{d}\left( \hat{G}\left( x\right) -G_{%
\mathbf{c}}\left( x\right) \right) .
\end{equation*}%
Note that both $\Delta _{\gamma p}^{(1)}\left( \theta \right) $ and $\Delta
_{\gamma p}^{(2)}\left( \theta \right) $ have form%
\begin{equation*}
Y_{ab}=p\int \ln \left( a+bx\right) \mathrm{d}\left( \hat{G}\left( x\right)
-G_{\mathbf{c}}\left( x\right) \right) .
\end{equation*}%
For $a,b,a^{\prime },b^{\prime }>0,$ Zheng (2012), Example 4.1 proves that $%
\left( Y_{ab},Y_{a^{\prime }b^{\prime }}\right) $ converge to a Gaussian
vector $\left( X_{ab},X_{a^{\prime }b^{\prime }}\right) $ with%
\begin{equation}
\mathbb{E}X_{ab}=\tfrac{1}{2}\log \frac{\left( c^{2}-d^{2}\right) \rho ^{2}}{%
\left( c\rho -\gamma _{2}d\right) ^{2}}  \label{meanZ}\tag{SM183}
\end{equation}%
and%
\begin{equation}
Cov\left( X_{ab},X_{a^{\prime }b^{\prime }}\right) =2\log \frac{cc^{\prime }%
}{cc^{\prime }-dd^{\prime }},  \label{CovZ}\tag{SM184}
\end{equation}%
where $c>d>0$ satisfy%
\begin{equation}
c^{2}+d^{2}=a+b\frac{1+\rho ^{2}}{\left( 1-\gamma _{2}\right) ^{2}}\text{
and }cd=\frac{b\rho }{\left( 1-\gamma _{2}\right) ^{2}}  \label{cd}\tag{SM185}
\end{equation}%
and $c^{\prime }>d^{\prime }>0$ satisfy 
\begin{equation}
c^{\prime 2}+d^{\prime 2}=a^{\prime }+b^{\prime }\frac{1+\rho ^{2}}{\left(
1-\gamma _{2}\right) ^{2}}\text{ and }c^{\prime }d^{\prime }=\frac{b^{\prime
}\rho }{\left( 1-\gamma _{2}\right) ^{2}}.  \label{cdprime}\tag{SM186}
\end{equation}

A direct inspection reveals that Zheng's proof of (\ref{meanZ}) and (\ref%
{CovZ}) remains valid for any real $a,b,a^{\prime },$ and $b^{\prime }$ such
that $\log \left( a+bz\right) $ and $\log \left( a^{\prime }+b^{\prime
}x\right) $ are analytic in an open domain containing the support of $G_{%
\boldsymbol{\gamma }}$ as long as there exist real $c$ and $d$ satisfying (\ref%
{cd}) and real $c^{\prime }$ and $d^{\prime }$ satisfying (\ref{cdprime})
such that $\left\vert c\right\vert >\left\vert d\right\vert $ and $%
\left\vert c^{\prime }\right\vert >\left\vert d^{\prime }\right\vert .$ Such 
$c,d,c^{\prime }$ and $d^{\prime }$ do exist for $Y_{ab}=\Delta _{\gamma
p}^{(1)}\left( \theta \right) $ and $Y_{a^{\prime }b^{\prime }}=\Delta
_{\gamma p}^{(2)}\left( \theta \right) $. Indeed, the values of $a$ and $b$
for $Y_{ab}=\Delta _{\gamma p}^{(1)}\left( \theta \right) $ are 
\begin{equation*}
a=\theta +\gamma _{1}\text{ and }b=\gamma _{2}-\frac{\theta }{1+\theta }.
\end{equation*}%
The corresponding $c$ and $d$ that satisfy (\ref{cd}) are%
\begin{equation}
c=\frac{\rho }{\sqrt{\theta +1}\left( 1-\gamma _{2}\right) }\text{ and }d=%
\frac{\gamma _{2}-\theta \left( 1-\gamma _{2}\right) }{\sqrt{\theta +1}%
\left( 1-\gamma _{2}\right) }.  \label{cdexplicit}\tag{SM187}
\end{equation}%
Since $\gamma _{2}<\rho ,$ $\left\vert c\right\vert $ is clearly larger than 
$\left\vert d\right\vert $ for positive $d.$ For non-positive $d,$ $%
\left\vert c\right\vert >\left\vert d\right\vert $ if and only if%
\begin{equation*}
\theta \left( 1-\gamma _{2}\right) -\gamma _{2}<\rho ,
\end{equation*}%
But this inequality folds for any $\theta \in \left[ 0,\bar{\theta}%
-\varepsilon \right] $ because $\bar{\theta}=\left( \gamma _{2}+\rho \right)
/\left( 1-\gamma _{2}\right) $ (see Table JO\ref{Table 3}).

Further, the values of $a^{\prime }$ and $b^{\prime }$ for $Y_{a^{\prime
}b^{\prime }}=\Delta _{\gamma p}^{(2)}\left( \theta \right) $ are%
\begin{equation*}
a^{\prime }=\gamma _{1}/\gamma _{2}\text{ and }b^{\prime }=1.
\end{equation*}%
The corresponding $c^{\prime }$ and $d^{\prime }$ that satisfy (\ref{cdprime}%
) are%
\begin{equation}
c^{\prime }=\frac{\rho }{\left( 1-\gamma _{2}\right) \sqrt{\gamma _{2}}}%
\text{ and }d^{\prime }=\frac{\sqrt{\gamma _{2}}}{1-\gamma _{2}}.
\label{cdprimeexplicit}\tag{SM188}
\end{equation}%
Since $\gamma _{2}<\rho ,$ we have $c^{\prime }>d^{\prime }>0.$

Using (\ref{decompDgp}), (\ref{meanZ}), (\ref{cdexplicit}), and (\ref%
{cdprimeexplicit}), we find that 
\begin{eqnarray*}
\mathbb{E\mathcal{D}}\left( \theta \right)  &\mathbb{=}&\tfrac{1}{2}\log 
\frac{\left( c^{2}-d^{2}\right) \left( c^{\prime }\rho -\gamma _{2}d^{\prime
}\right) ^{2}}{\left( c\rho -\gamma _{2}d\right) ^{2}\left( c^{\prime
2}-d^{\prime 2}\right) } \\
&=&\tfrac{1}{2}\log \left( 1-\frac{\rho ^{2}\theta ^{2}}{\left( \gamma
_{1}+\gamma _{2}\left( 1+\theta \right) \right) ^{2}}\right) 
\end{eqnarray*}%
and%
\begin{eqnarray*}
Cov\left( \mathbb{\mathcal{D}}\left( \theta _{i}\right) ,\mathbb{\mathcal{D}}%
\left( \theta _{j}\right) \right)  &=&2\log \frac{\rho ^{2}}{\rho
^{2}-\left( \gamma _{2}-\theta _{i}\left( 1-\gamma _{2}\right) \right)
\left( \gamma _{2}-\theta _{j}\left( 1-\gamma _{2}\right) \right) } \\
&&-2\log \frac{\rho ^{2}}{\rho ^{2}-\left( \gamma _{2}-\theta _{j}\left(
1-\gamma _{2}\right) \right) \gamma _{2}} \\
&&-2\log \frac{\rho ^{2}}{\rho ^{2}-\left( \gamma _{2}-\theta _{i}\left(
1-\gamma _{2}\right) \right) \gamma _{2}} \\
&&+2\log \frac{\rho ^{2}}{\rho ^{2}-\gamma _{2}^{2}} \\
&=&-2\log \left( 1-\frac{\rho ^{2}\theta _{i}\theta _{j}}{\left( \gamma
_{1}+\gamma _{2}\left( 1+\theta _{i}\right) \right) \left( \gamma
_{1}+\gamma _{2}\left( 1+\theta _{j}\right) \right) }\right) .
\end{eqnarray*}

\section{Concluding remarks}

\subsection{Power of the LR test under multi-spike alternatives}
\label{power-lr-test}

Consider the likelihood ratio test that rejects the null hypothesis 
of no spikes when the supremum of $\ln L\left( \theta
;\Lambda \right) $ over $\theta \in \left[ 0,\bar{\theta}-\varepsilon \right]
$ is above an asymptotic critical value. In this section, we study the power of 
such a test in the situation where the rank-one assumption on the alternative 
is wrong and there are multiple spikes, the highest of which is at least as high 
as the spike under our rank-one setting. 

Intuitively, the power should increase under such a multi-spike alternative
because it is ``further away'' from the null than the one-spike alternative.
Below, we confirm this intuition for SMD and PCA cases.

First let us show that, in any of James' cases, the corresponding likelihood
ratio test has a monotone acceptance region. That is, the null is accepted
if and only if $g\left( \lambda _{1},...,\lambda _{p}\right) <\mathsf{const}
$ for a function $g$ which is non-decreasing in each argument. Recall that
the likelihood ratio has the following form 
\begin{equation}
L\left( \theta ;\Lambda \right) =\alpha (\theta )\left. _{\mathsf{p}}F_{%
\mathsf{q}}\right. \left( a,b;\Psi ,\Lambda \right) ,
\label{lkelihood ratio 1}\tag{SM189}
\end{equation}%
where $\Psi =\mathsf{diag}\left\{ \Psi _{11},0,...,0\right\} ,$ $\Lambda =%
\mathsf{diag}\left\{ \lambda _{1},...,\lambda _{p}\right\} ,$ and the
values of $\Psi _{11},\alpha (\theta ),a,b,\mathsf{p}$, and $\mathsf{q}$ for
the different cases are given in Table JO\ref{Table 2}. As explained in
Section \ref{sec: proof of theorem 11}, we have the following expansion%
\[
L\left( \theta ;\Lambda \right) =\alpha (\theta )\sum_{k=0}^{\infty }\frac{1%
}{k!}\frac{\left( a_{1}\right) _{k}...\left( a_{\mathsf{p}}\right) _{k}}{%
\left( b_{1}\right) _{k}...\left( b_{\mathsf{q}}\right) _{k}}\frac{\Psi
_{11}^{k}C_{k}\left( \Lambda \right) }{C_{k}\left( I_{p}\right) },
\]%
where $C_{k}$ are zonal polynomials. James (1968) shows that zonal
polynomials have positive coefficients. Therefore, $C_{k}\left( \Lambda
\right) $ and $L\left( \theta ;\Lambda \right) $ are nondecreasing in each $%
\lambda _{j}$ for any fixed $\theta \in \left[ 0,\bar{\theta}-\varepsilon %
\right] .$ As a consequence, the supremum of $\ln L\left( \theta ;\Lambda
\right) $ over $\theta \in \left[ 0,\bar{\theta}-\varepsilon \right] $ is a
non-decreasing function in each $\lambda _{j}.$

Next, recall that SMD refers to the problem of testing $H_{0}:\Phi =0$
against $H_{1}:\Phi =\theta _{1}\psi _{1}\psi _{1}^{\prime }$ using the
eigenvalues $\lambda _{j},$ $j=1,...,p$ of matrix $X=\Phi +Z/\sqrt{p},$
where $Z$ is a noise matrix from the Gaussian Orthogonal Ensemble. Now
suppose that the actual situation corresponds to the alternative 
\[
H_{\mathrm{mult}}:\Phi =\sum_{j=1}^{r}\theta _{j}\psi _{j}\psi _{j}^{\prime
},
\]%
where $\theta _{1}\geq ...\geq \theta _{r}>0$ and $\psi _{1},...,\psi _{r}$
is a set of orthonormal nuisance vectors. Since $\Phi $ under $H_{\mathsf{%
mult}}$ is no smaller than under $H_{1},$ the $j$-th largest eigenvalue of $X
$ under $H_{\mathsf{mult}}$ is no smaller than under $H_{1}.$ But as shown
above, the likelihood ratio test has a monotone acceptance region. Hence,
its power to reject $H_{0}$ in favour of $H_{\mathsf{mult}}$ is at least as
high as its power to reject $H_{0}$ in favour of $H_{1}.$

Similarly, recall that PCA refers to the problem of testing $H_{0}:\Omega
=I_{p}$ against $H_{1}:\Omega =I_{p}+\theta _{1}\psi _{1}\psi _{1}^{\prime }$
using the eigenvalues of $YY^{\prime }/n_{1},$ where $Y=\Omega
^{1/2}\varepsilon $ and $\varepsilon $ is a $p\times n_{1}$ matrix with
i.i.d. standard normal entries. Suppose that the actual situation
corresponds to the alternative 
\[
H_{\mathsf{mult}}:\Omega =I_{p}+\sum_{j=1}^{r}\theta _{j}\psi _{j}\psi
_{j}^{\prime }.
\]%
Note that the non-zero eigenvalues of $YY^{\prime }/n_{1}$ coincide with
those of $\varepsilon ^{\prime }\Omega \varepsilon /n_{1}.$ Since $\Omega $
under $H_{\mathsf{mult}}$ is no smaller than under $H_{1},$ the $j$-th
largest eigenvalue of $\varepsilon ^{\prime }\Omega \varepsilon /n_{1}$
under $H_{\mathsf{mult}}$ is no smaller than under $H_{1}.$ Therefore, using
the monotonicity of the acceptance region of the test, we conclude that the
power corresponding to $H_{\mathsf{mult}}$ is no smaller than that
corresponding to $H_{1}.$

Unfortunately, for the remaining cases, the above logic does not go through.
For example, for SigD, we test $H_{0}:\Omega =I_{p}$ against $H_{1}:\Omega
=I_{p}+\theta _{1}\psi _{1}\psi _{1}^{\prime }$ using eigenvalues of $\left(
XX^{\prime }/n_{2}\right) ^{-1}\left( YY^{\prime }/n_{1}\right) ,$ where $%
Y=\Omega ^{1/2}\varepsilon $ is as above, and $X$ is a $p\times n_{2}$
matrix with i.i.d. standard normal entries independent from $Y$. It is
conceivable that%
\[
\varepsilon ^{\prime }\Omega ^{1/2}\left( XX^{\prime }/n_{2}\right)
^{-1}\Omega ^{1/2}\varepsilon /n_{1},
\]%
as opposed to $\varepsilon ^{\prime }\Omega \varepsilon /n_{1}$ (cf. the PCA
case above), has some of its eigenvalues under $H_{\mathsf{mult}}$ smaller
than the corresponding eigenvalues under $H_{1}.$ Of course, on average over
the distribution of $X,$ the situation will be exactly the same as in the
PCA case. Therefore, although we cannot prove the increase in power, it
remains intuitively plausible.

Perlman and Olkin (1980) study the unbiasedness and power monotonicity of
tests with monotone acceptance regions in cases that correspond to our REG$%
_{0},$ REG, and CCA. Although they prove the unbiasedness of such tests, the
power monotonicity remains a \textquotedblleft strong
conjecture\textquotedblright\ (see p. 1329 of their paper). Their
Proposition 2.6 (ii) formulates conditions on the likelihood ratio
(corresponding to general alternatives) that guarantee the power
monotonicity. However, as shown in Richards (2004), these conditions do not
hold for likelihoods of form (\ref{lkelihood ratio 1}), in general. Of course,
this does not mean that Perlman and Olkin's conjecture is wrong, it just
cannot be established directly via Proposition 2.6.

\end{document}